\documentclass[reqno,oneside]{amsart} 
\pdfoutput=1
\usepackage{amsmath}
\usepackage{amsthm,
amssymb,amsmath,enumerate,stackrel,
verbatim,
}
\usepackage{xr-hyper}
\usepackage[all,cmtip]{xy}
\usepackage{enumitem}
\usepackage{soul}
\usepackage{mathrsfs}
\usepackage[T1]{fontenc} % for \k{e}
\usepackage[pagebackref,hyperindex,linktocpage=true]{hyperref}
\hypersetup{
    colorlinks,
    linkcolor={red!50!black},
    citecolor={red!50!black}, 
    urlcolor={red!50!black}, 
    filecolor={red!50!black} 
}

\makeatletter
\newcommand{\fakephantomsection}{%
  \Hy@MakeCurrentHref{\@currenvir.\the\Hy@linkcounter}
  \Hy@raisedlink{\hyper@anchorstart{\@currentHref}\hyper@anchorend}%
}
\makeatother

\usepackage{theoremref}

\usepackage{tikz}
\usetikzlibrary{matrix,arrows}

\usepackage{tikz-cd}

\usepackage{leftidx} 

\makeatletter
\@addtoreset{equation}{section}
\makeatother

\numberwithin{equation}{section}

\theoremstyle{definition}
\newtheorem{Defi}{Definition}[section] \newcommand{\defi}{\begin{Defi}} \newcommand{\xdefi}{\end{Defi}} \newtheorem{DefiLemm}[Defi]{Definition and Lemma} \newcommand{\defilemm}{\begin{DefiLemm}} \newcommand{\xdefilemm}{\end{DefiLemm}} 
\newtheorem{Bsp}[Defi]{Example} \newcommand{\exam}{\begin{Bsp}} \newcommand{\xexam}{\end{Bsp}} 
\newtheorem{Syno}[Defi]{Synopsis} \newcommand{\syno}{\begin{Syno}} \newcommand{\xsyno}{\end{Syno}} 
\newtheorem{Bem}[Defi]{Remark} \newcommand{\rema}{\begin{Bem}} \newcommand{\xrema}{\end{Bem}} 
\newtheorem{Notation}[Defi]{Notation} \newcommand{\nota}{\begin{Notation}} \newcommand{\xnota}{\end{Notation}} 
\newtheorem{Convention}[Defi]{Convention} \newcommand{\conv}{\begin{Convention}} \newcommand{\xconv}{\end{Convention}} 
\newtheorem{Warning}[Defi]{Warning} \newcommand{\warn}{\begin{Warning}} \newcommand{\xwarn}{\end{Warning}} 
\newtheorem{Situation}[Defi]{Situation} \newcommand{\situ}{\begin{Situation}} \newcommand{\xsitu}{\end{Situation}}
\newtheorem{Assumption}[Defi]{Assumption} \newcommand{\assu}{\begin{Assumption}} \newcommand{\xassu}{\end{Assumption}} 

\theoremstyle{plain}
\newtheorem{Theo}[Defi]{Theorem} \newcommand{\theo}{\begin{Theo}} \newcommand{\xtheo}{\end{Theo}} 
\newtheorem{Satz}[Defi]{Proposition} \newcommand{\prop}{\begin{Satz}} \newcommand{\xprop}{\end{Satz}} 
\newtheorem{Lemm}[Defi]{Lemma} \newcommand{\lemm}{\begin{Lemm}} \newcommand{\xlemm}{\end{Lemm}} 
\newtheorem{Coro}[Defi]{Corollary} \newcommand{\coro}{\begin{Coro}} \newcommand{\xcoro}{\end{Coro}}
\newtheorem{Ques}[Defi]{Question} \newcommand{\ques}{\begin{Ques}} \newcommand{\xques}{\end{Ques}}
\newtheorem{Conj}[Defi]{Conjecture} \newcommand{\conj}{\begin{Conj}} \newcommand{\xconj}{\end{Conj}}

\newcommand{\refsect}[1]{Section \ref{sect--#1}}
\newcommand{\refit}[1]{(\ref{item--#1})}
\newcommand{\refeq}[1]{\eqref{eqn--#1}}
\renewcommand{\eqref}[1]{(\ref{#1})}

\newcommand{\eqn}{\begin{equation}} \newcommand{\xeqn}{\end{equation}}
\newcommand{\eqnarr}{\begin{eqnarray*}} \newcommand{\xeqnarr}{\end{eqnarray*}}
\newcommand{\eqnarra}{\begin{eqnarray}} \newcommand{\xeqnarra}{\end{eqnarray}}

\newcommand{\pf}{\begin{proof}} \newcommand{\xpf}{\end{proof}}

\newcounter{heyheyCounter}[section]

\newif \ifHideSomeComments
\HideSomeCommentstrue

\ifHideSomeComments
  \newcommand{\heyheyX}[1]{}
\else
  \newcommand{\heyheyX}[1]{\heyhey{#1}}
\fi
\newcommand{\heyhey}[1]{\refstepcounter{heyheyCounter}$\bigstar^\text{\theheyheyCounter}$\marginpar{\footnotesize $\bigstar^\text{\theheyheyCounter}$ #1}}

\newif \ifDraft % whether it is a draft or otherwise meant to be sent to other people etc.
\Draftfalse

\ifDraft
  \usepackage[top=1.5cm, bottom=1.5cm, left=.5cm, right=7.5cm, heightrounded, marginparwidth=7cm, marginparsep=.2cm]{geometry}
\else
  \usepackage[top=1.5cm, bottom=1.5cm, left=2cm, right=2cm, heightrounded]{geometry}
  \usepackage{letltxmacro}
  \LetLtxMacro\Oldfootnote\footnote
  \renewcommand{\footnote}[2][]{\relax}
  \renewcommand{\heyhey}[1]{}
\fi

\newcommand{\nc}{\newcommand}
\nc{\StP}[1]{\cite[Tag~\href{http://stacks.math.columbia.edu/tag/#1}{#1}]{StacksProject}}
\nc{\StPd}[2]{\cite[Tags~\href{http://stacks.math.columbia.edu/tag/#1}{#1}, \href{http://stacks.math.columbia.edu/tag/#2}{#2}]{StacksProject}} 

\nc{\on}{\operatorname}
\nc{\aff}{{\on{aff}}}
\nc{\modi}{{\on{mod}}} 
\nc{\even}{{\on{even}}}
\nc{\odd}{{\on{odd}}}
\nc{\naive}{{\on{naive}}}
\nc{\hofib}{\on{hofib}}
\nc{\Bun}{\on{Bun}}
\nc{\ad}{{\on{ad}}}
\nc{\lft}{{\on{lft}}}
\nc{\modulo}{\on{mod}} 
\nc{\FinSet}{\on{FinSet}} 
\nc{\surj}{\on{surj}} 

\nc{\Weil}{{\on{Weil}}} 
\nc{\FWeil}{{\on{FWeil}}} 
\nc{\cons}{{\on{cons}}} 
\nc{\tot}{{\on{Tot}}} 

\nc{\str}{\on{-}}
\nc{\perf}{{\on{perf}}}
\nc{\Rel}{{\on{Pos}}}
\nc{\lan}{\langle}
\nc{\ran}{\rangle}

\nc{\PM}{{\on{PM}}}
\nc{\bbf}{{\mathbf f}}
\nc{\calH}{{\mathcal H}}
\nc{\calP}{{\mathcal P}} 
\nc{\calF}{{\mathcal F}} 
\nc{\calB}{{\mathcal B}} 
\nc{\calI}{{\mathcal I}}
\nc{\calJ}{{\mathcal J}} 
\nc{\calO}{{\mathcal O}}
\nc{\calM}{{\mathcal M}}
\nc{\co}{\colon}
\nc{\GL}{GL} 
\nc{\zero}{\underline{0}}
\nc{\ux}{\underline{x}} 
\nc{\uT}{\underline{T}} 
\nc{\uP}{\underline{P}}
\nc{\uS}{\underline{S}}

\newcommand{\category}[1]{\mathrm{#1}}

\newcommand\restr[2]{{ 
  \left.\kern-\nulldelimiterspace 
  #1 
  \vphantom{\big|} 
  \right|_{#2} 
  }}

\makeatletter 
\newcommand*{\doublerightarrow}[2]{\mathrel{
  \settowidth{\@tempdima}{$\scriptstyle#1$}
  \settowidth{\@tempdimb}{$\scriptstyle#2$}
  \ifdim\@tempdimb>\@tempdima \@tempdima=\@tempdimb\fi
  \mathop{\vcenter{
    \offinterlineskip\ialign{\hbox to\dimexpr\@tempdima+1em{##}\cr
    \rightarrowfill\cr\noalign{\kern.5ex}
    \rightarrowfill\cr}}}\limits^{\!#1}_{\!#2}}}
\newcommand*{\triplerightarrow}[1]{\mathrel{ 
  \settowidth{\@tempdima}{$\scriptstyle#1$}
  \mathop{\vcenter{
    \offinterlineskip\ialign{\hbox to\dimexpr\@tempdima+1em{##}\cr
    \rightarrowfill\cr\noalign{\kern.5ex}
    \rightarrowfill\cr\noalign{\kern.5ex}
    \rightarrowfill\cr}}}\limits^{\!#1}}}
\makeatother

\def\gr{\category {gr}} 
\newcommand{\Ho}{\category{Ho}} 
\newcommand{\Mod}{\category{Mod}} 
\newcommand{\Grp}{\category{Grp}} 
\newcommand{\Mon}{\category{Mon}} 
\newcommand{\Ani}{\category{Ani}} 
\newcommand{\coMod}{\category{coMod}}

\newcommand{\Perf}{\category{Perf}}

\renewcommand{\Bar}{\operatorname{Bar}} 

\newcommand{\OldLambda}{E} % used to be Lambda in the version before

\renewcommand{\L}{\category{L}}
\renewcommand{\Pr}{\category{Pr}}
\newcommand{\Stable}{\mathrm{St}}
\newcommand{\PrL}{\Pr^\L} 
\newcommand{\PrLomegaSt}{\PrL_{\omega,\Stable}}
\newcommand{\PrLRomegaSt}{\Pr^{\category {LR}}_{[\omega],\Stable}} 
\newcommand{\PrLSt}{\PrL_{\mathrm{St}}} 
\newcommand{\PrLst}{\PrLSt} 
\newcommand{\PrLgr}{\Pr^\category{L}_\gr}

\newcommand{\Cat}{\category{Cat}} 
\newcommand{\Fun}{\category{Fun}} 
 
\newcommand{\Corr}{\category{Corr}} 
\newcommand{\Sch}{\category{Sch}} 
\newcommand{\AffSch}{\category{AffSch}} 
 
\newcommand{\IndSch}{\category{IndSch}}

\newcommand{\Fin}{\category{Fin}} 
\newcommand{\dual}{\vee} 
\newcommand{\Zar}{\mathrm{Zar}} 
\newcommand{\Set}{\category{Set}}
\newcommand{\Sat}{\category{Sat}}

\newcommand{\Fl}{Fl}
\newcommand{\et}{\mathrm{et}} 
\newcommand{\opp}{\mathrm{op}}
\def\proper{\mathrm {proper}} 
\def\placid{\mathrm {pl}}

\newcommand{\free}{\mathrm{free}} % free (abelian groups) 
\newcommand{\grAb}{\Ab^\Z} % graded abelian groups 
\newcommand{\grD}{\D^\Z} % graded derived category
\newcommand{\grCh}{\Ch^\Z} % graded chain complexes
\newcommand{\grModZ}{\grAb}% {\gr\Mod_\Z}
\newcommand{\grModQ}{\Mod^\Z_\Q}
 
\def\Gm{\mathbf {G}_\mathrm m}
\def\GmC{\mathbf{G}_{\mathrm{m}, C}}
\def\GmCI{\mathbf{G}_{\mathrm{m}, C^I}}
\def\GmS{\mathbf{G}_{\mathrm{m}, S}}
\def\Gmmono{\Gm\mathrm{-mono}}

\newcommand{\GmX}[  1]{\mathbf {G}_{\mathrm {m}, #1}} 
\def\IC{\mathrm{IC}} 
\def\BD{\mathrm{BD}} 
\def\Supp{\operatorname{Supp}} %Support

\newcommand{\colim}{\operatornamewithlimits{colim}}

\def\id{{\rm id}}

\def\pr{{\rm pr}} 
\def\opp{{\rm op}} 
 
\def\To#1#2{\mathop{\count0=#1 \loop\ifnum\count0>0 \smash-\mkern-7mu \advance\count0 -1 \repeat \mathord\rightarrow}\limits^{#2}}

\def\Maps{\mathop{\rm Maps}\nolimits}

\def\cofib{\mathop{\rm cofib}\nolimits} 
 
\def\Hom{\mathop{\rm Hom}\nolimits} 
\def\rightind{{\rm r}}  % right-induced
\def\leftind{{\rm l}}  % left-induced
\def\Tens{\mathop{\rm Tens}\nolimits} 
\def\Nilp{\mathop{\rm Nilp}\nolimits} 
\def\NilpQ{\Nilp_\Q} 
\def\Ran{\mathop{\rm Ran}\nolimits} 
\def\Lan{\mathop{\rm Lan}\nolimits} 
 
\def\PGL{\mathrm {PGL}} 
\def\GL{\mathrm {GL}}

\def\Fil{\mathop{\rm Fil}\nolimits}
\def\Grad{\mathop{\rm Grad}\nolimits}

\def\Ab{\mathop{\rm Ab}\nolimits}
\def\Ind{\mathop{\category{Ind}}} 
 
\def\Gr{\mathop{\rm Gr}\nolimits} 
\def\Hck{\mathop{\rm Hck}\nolimits} 
 
\def\Fl{\mathop{\rm Fl}\nolimits}

\def\Sym{\mathop{\rm Sym}\nolimits} 
\def\CAlg{\mathop{\rm CAlg}\nolimits}

\def\IHom{\underline{\Hom}}

\def\Map{\Maps} 
\def\Rep{\category{Rep}} 
 
\def\dualizable{\mathrm{dualizable}}

\def\et{\mathrm{\acute et}} 
\def\Nis{\mathrm{Nis}}

\def\der{\mathrm{der}} 
\def\CT{\mathrm{CT}} 
 
\def\sph{\mathrm{sph}} 
 
\def\res{\mathrm{res}} 
\def\act{\mathrm{act}} 
\def\aug{\mathrm{aug}}

\def\Z{{\mathbf Z}} 
\def\Fp{{{\mathbf F}_p}} 
\def\Fq{{{\mathbf F}_q}}

\def\Q{{\mathbf Q}}

\def\C{{\mathbf {C}}} 
\def\A{{\bf A}} 
\renewcommand{\P}[1][1]{\mathbf P^{#1}} 
\def\Gm{\mathbf {G}_\mathrm m}

\def\qq{{\mathbf{q}}}

\def\Cc{{\mathcal{C}}}
\def\Dd{{\mathcal{D}}}

\def\Ff{{\mathcal{F}}}
\def\calH{{\mathcal{H}}}

\def\Uu{{\mathcal{U}}}

\def\H{{\rm H}} 
\def\pe{{{}^\mathrm{p}} \! } 
\def\pH{\pe\H}

\def\SH{\category{SH}}

\def\R{\mathrm{R}} 
\def\Betti{\mathrm{B}} 
\def\red{\mathrm{r}} 
\def\redx{(\red)} 
\def\DM{\category{DM}} 
\def\DMr{\DM_\red}

\def\DMrx{\DM_{\redx}}

\def\DATM{\category{DATM}} 
\def\DTM{\category{DTM}} 
\def\DTMr{\DTM_\red} 
\def\DTMrx{\DTM_{\redx}}

\def\anti{\mathrm{anti}} 
\def\xantic{(\anti),\comp} % possibly anti-effective and compact 
\def\antic{\anti,\comp} % anti-effective and compact 
\def\antilc{\anti,\locc} % anti-effective and locally compact 
\def\redu{\mathrm{red}} 

\def\Wak{\mathrm{Wak}} % Category generated by Wakimoto motives
\def\JWak{\mathbf{J}^B} % Wakimoto objects; can still decide whether to remove B from notation

\def\Perv{\category{Perv}} 
\def\MTM{\category{MTM}}

\def\Sat{\category{Sat}}

\def\ii{$\infty$}

\def\Spec{\mathop{\rm Spec}}

\newcommand{\M}{\mathrm{M}}

\newcommand{\comp}{\mathrm{c}} 
\newcommand{\locc}{\mathrm{lc}}  % locally compact 

\newcommand{\PreStk}{\category{PreStk}}
\newcommand{\adm}{\mathrm{Adm}}

\newcommand{\an}{\mathrm{an}} 
 
\newcommand{\D}{\category{D}} 
 
\newcommand{\Ch}{\category{Ch}}

\def\sm{{\rm sm}} 
\def\proper{{\rm pr}} 
\def\sbuildrel#1\over#2{\mathrel{\smash{\mathop{\kern0pt #2}\limits^{#1}}}}

\let\x\times
\let\ol\overline
\renewcommand{\t}{\otimes}

\renewcommand{\r}{\rightarrow}

\newcommand{\pot}[1]{ [\hspace{-0,5mm}[ {#1} ]\hspace{-0,5mm}] }
\newcommand{\rpot}[1]{ (\hspace{-0,7mm}( {#1} )\hspace{-0,7mm}) }

\mathchardef\mhyphen="2D

\def\matrix#1{\null\,\vcenter{\normalbaselines
    \ialign{\hfil$##$\hfil&&\quad\hfil$##$\hfil\crcr
      \mathstrut\crcr\noalign{\kern-\baselineskip}
      #1\crcr\mathstrut\crcr\noalign{\kern-\baselineskip}}}\,}

\newdimen\harrowsize
\def\mapright#1{\smash{\mathop{\hbox to\harrowsize{\rightarrowfill}}\limits^{#1}}}

{\catcode`@=11
\gdef\cal{\fam\tw@}

\global\let\over\@@over
\global\let\atop\@@atop
\global\let\above\@@above
\global\let\overwithdelims\@@overwithdelims
\global\let\atopwithdelims\@@atopwithdelims
\global\let\abovewithdelims\@@abovewithdelims
\gdef\eqalign#1{\null\,\vcenter{\openup\jot\m@th
\ialign{\strut\hfil$\displaystyle{##}$&$\displaystyle{{}##}$\hfil
      \crcr#1\crcr}}\,}
\newskip\xcentering \global\xcentering=0pt plus 1000pt minus 1000pt
\gdef\eqalignno#1{\displ@y \tabskip\xcentering
  \halign to\displaywidth{\hfil$\@lign\displaystyle{##}$\tabskip\z@skip
    &$\@lign\displaystyle{{}##}$\hfil\tabskip\xcentering
    &\llap{$\@lign##$}\tabskip\z@skip\crcr
    #1\crcr}}
\gdef\eqlabel#1{\refstepcounter{equation}\label{eqn--#1}\eqno\hbox{\@eqnnum}}
}

\pagecolor{white}

\begin{document}

\title{Central motives on parahoric flag varieties}
\author{Robert Cass, Thibaud van den Hove, Jakob Scholbach}

\begin{abstract}
	We construct a refinement of Gaitsgory's central functor for integral motivic sheaves, and show it preserves stratified Tate motives.
	Towards this end, we develop a reformulation of unipotent motivic nearby cycles, which also works over higher-dimensional bases.
	We moreover introduce Wakimoto motives and use them to show that our motivic central functor is t-exact.
	A decategorification of these functors yields a new approach to generic Hecke algebras for general parahorics.
\end{abstract}

\maketitle

\setcounter{tocdepth}{1}
\tableofcontents

\section{Introduction}
For a split reductive group $G$, the study of the spherical Hecke algebra and the Iwahori--Hecke algebra  are at the center of the Langlands program. These are defined as
\begin{align*}
\calH^{\sph} & := \Fun_c(G (\Fq \pot t) \setminus G (\Fq \rpot t) / G (\Fq \pot t), \C), \\
\calH^{\calI} & := \Fun_c(\mathcal{I} \setminus G (\Fq \rpot t) / \mathcal{I}, \C),
\end{align*}
where $\mathcal{I} \subset G (\Fq \pot t)$ denotes an Iwahori subgroup.

These algebras are categorified by various categories of sheaves on the affine Grassmannian $\Gr$ and the affine flag variety $\Fl$, respectively.
Moreover, certain properties of Hecke algebras can be upgraded to categorical statements.
For example, the commutativity of $\calH^\sph$ can be refined to the existence of a symmetric monoidal structure on the category $\Perv_{L^+G} (\Gr)$ of $L^+G$-equivariant perverse sheaves on \(\Gr\).
This property is one of the salient aspects of the geometric Satake equivalence, a cornerstone of the geometric Langlands program, which establishes an equivalence between that category and the representations of the Langlands dual group $\hat G$.

Beginning with Mirkovi\'{c}--Vilonen's pioneering paper \cite{MirkovicVilonen:Geometric}, geometric Satake  has attracted interest from many authors including \cite{Richarz:New,Zhu:Ramified, Zhu:Affine,Mautner:ExoticTiltingSheaves2018,FarguesScholze:Geometrization, XuZhu:Bessel}. These papers treat various contexts, such as Betti sheaves, $\ell$-adic sheaves, $D$-modules, as well as Witt vector and \(B_{\mathrm{dR}}^+\)-affine Grassmannians.
In \cite{CassvdHScholbach:MotivicSatake}, the present authors refined several of these equivalences to a statement for integral motivic sheaves, building upon earlier work for rational motives in \cite{RicharzScholbach:Motivic}.

Bernstein \cite[Theorem 2.13]{Bernstein:Centre} gave a description of the center of $\calH^{\calI}$ as 
$$Z(\calH^{\calI}) = \calH^\sph.$$
Following a suggestion of Beilinson, Gaitsgory \cite{Gaitsgory:Central} achieved a geometrization of this isomorphism, which was improved upon by Zhu \cite{Zhu:Coherence}.
Namely, there is a geometric object over \(\A^1\) whose fibers over \(x\neq 0\) are \(\Gr\), and whose fiber over \(0\) is \(\Fl\), i.e., a deformation 
$$\Gr \leadsto \Fl.\eqlabel{GrleadstoFl}$$
Taking nearby cycles for this family then yields the desired geometrization.

The purpose of this paper is to refine Gaitsgory's central functor to the level of integral motivic sheaves.
In order to achieve this goal, we use new insights into the geometry of such \(\A^1\)-deformations and related Beilinson--Drinfeld Grassmannians, and combine these with a much-simplified approach to motivic unipotent nearby cycles.

In the \(\ell\)-adic setting, Gaitsgory's central functor is also a key ingredient in Bezrukavnikov's tamely ramified local Langlands correspondence \cite{ArkhipovBezrukavnikov:Perverse, Bezrukavnikov:Two}.
As such, the constructions of the present paper, along with the motivic Satake equivalence from \cite{CassvdHScholbach:MotivicSatake}, should allow us to enhance \cite{ArkhipovBezrukavnikov:Perverse,Bezrukavnikov:Two} to the motivic setting. In particular, this will provide a mixed upgrade of these works as anticipated in \cite[\S 11.1]{Bezrukavnikov:Two}, and it will eliminate the dependence on $\ell$ inherent in the use of $\ell$-adic cohomology.

\subsection{Geometry of Beilinson--Drinfeld Grassmannians}
Let \(G/\Spec \Z\) be a split reductive group and let \(\mathbf{f}\) be a facet in the Bruhat-Tits building of \(G\) contained in the closure of a standard alcove \(\mathbf{a}_0\). Let $\Fl_{\mathbf{f}}$ be the partial affine flag variety associated to the corresponding parahoric $\Z[\![t]\!]$-model \(G_{\mathbf{f}}\) of \(G\) ($\Gr$ and $\Fl$ are special cases). All results below also hold over arbitrary bases; for this introduction we restrict to \(\Spec \Z\) for simplicity.

There is a smooth affine group scheme \(\mathcal{G}_{\mathbf{f}}/\A^1_\Z\) which generically agrees with $G$, and whose fiber over the completed local ring at \(0\) agrees with the parahoric \(G_{\mathbf{f}}\).
The associated quotient \(\Gr_{\mathcal{G}_{\mathbf{f}}} := L \mathcal G_{\mathbf{f}}/ L^+\mathcal G_{\mathbf{f}} \) of the loop group by the positive loop group for $\mathcal G_{\mathbf f}$, known as a \emph{Beilinson--Drinfeld Grassmannian}, provides a degeneration from $\Gr$ to $\Fl_{\mathbf{f}}$.
The family in \refeq{GrleadstoFl} is the special case when $\mathbf f$ is a standard alcove.

The $\Gm$-fiber $(\Gr_{\mathcal{G}_{\mathbf{f}}})_{\eta} \cong \Gr \times \Gm$ carries the stratification by \(L^+G\)-orbits, while the special fiber $(\Gr_{\mathcal{G}_{\mathbf{f}}})_{s} = \Fl_{\mathbf f}$ is stratified by \(L^+G_{\mathbf{f}}\)-orbits in \(\Fl_{\mathbf{f}}\).
These sets of strata parametrize a basis of the spherical Hecke algebra, respectively the Hecke algebra for $\mathbf f$.

Previous work shows that categories of stratified Tate motives, i.e., motives whose restriction to each stratum lies in the subcategory generated by the ``constant'' motive $\Z$ and its Tate twists, are a small, but highly interesting subcategory of the category of all motivic sheaves, e.g.~\cite{RicharzScholbach:Motivic, CassvdHScholbach:MotivicSatake}. The following theorem ensures this approach is again applicable to \(\Gr_{\mathcal{G}_{\mathbf{f}}}\), cf.~\thref{theo-WT-BD}.
Joint with \thref{Ups.intro}\refit{Ups.WT}, it ensures that unipotent nearby cycles preserve stratified Tate motives.

\theo
The above stratifications of \((\Gr_{\mathcal{G}_{\mathbf{f}}})_\eta\) and \((\Gr_{\mathcal{G}_\mathbf{f}})_{s}\) determine a universally anti-effective Whitney--Tate stratification of \(\Gr_{\mathcal{G}_{\mathbf{f}}}\), in the sense of \cite[Definition 2.6]{CassvdHScholbach:MotivicSatake}.
\xtheo

Anti-effectivity ensures that no positive Tate twists arise, which will be useful when considering generic Hecke algebras by decategorifying the motivic central functor.
The key to proving Whitney--Tateness is the compatibility of push-pull functors from $\Gm$ to the special fiber of $\A^1$ with hyperbolic localization.
By additionally using that certain families of constant term functors preserve and reflect Tateness, we can reduce to the case where \(G=T\) is a torus, in which case the degeneration is trivial.

\subsection{Revisiting unipotent nearby cycles}
As was mentioned above, Gaitsgory's work crucially hinges on the unipotent nearby cycles functor (in the context of $\ell$-adic sheaves).
Such a functor has been developed by Ayoub for motivic categories \cite{Ayoub:Six1,Ayoub:Motivic}.
Taking our cue from Campbell's work \cite{Campbell:Nearby}, we construct this functor independently, using an \ii-categorical and much simplified approach. 
This is also necessary in view of the need of at least the rudiments of nearby cycles over 2-dimensional bases.

In the definition below, $S$ is a connected scheme that is smooth of finite type over a Dedekind ring or a field. For a scheme of finite type $X/S$, we have the category $\DM(X)$ of integral motivic sheaves as constructed by Spitzweck \cite{Spitzweck:Commutative} and the full subcategory $\DTM(X)$ of Tate motives, i.e., the presentable stable subcategory generated by $\Z_X(k)$, for $k \in \Z$.

\defi
For a scheme $X \r \A^1_S$, we consider the generic and special fibers 
$$X_\eta := X \x_{\A^1} \Gm \stackrel j \r X \stackrel i \gets X_s := X \times_{\A^1, 0} S.$$
The unipotent nearby cycles functor is 
$$\Upsilon : \DM(X_\eta) \r \DTM(\GmX S) \t_{\DTM(S)} \DM(X_\eta) \stackrel {1^* \t i^* j_*} \r \NilpQ \DM(X_s).$$
\xdefi

Here $\t$ denotes Lurie's tensor product of presentable stable \ii-categories.
The first functor exhibits a natural action of the motivic cohomology of $\GmX S$ (relative to $S$) on any motivic sheaf on $X_\eta$, akin to the natural $R$-action on any quasi-coherent sheaf $\mathcal F$ on a scheme $Z / \Spec R$; see the beginning of \refsect{unipotent} for more details.
Objects in the category $\NilpQ \DM(X_s)$ are pairs consisting of some $M \in \DM(X_s)$ and a locally nilpotent map $\phi: M \t \Q \r M \t \Q(-1)$. This map $\phi$ can be thought of as the ``logarithm of monodromy''. We stress, however, that no logarithm (of any sort) appears in this construction; instead this operator arises as the action of the Koszul dual of the cohomology ring of $\Gm$.
The appearance of rational coefficients at this point is caused by the non-formality of the motive of $\Gm$: while the motive $\M(\Gm)$ is $\Z \oplus \Z(1)[1]$, its natural coalgebra structure is square-zero (on the summand $\Z(1)[1]$) only after passing to rational coefficients.
By contrast, for reduced motives as introduced in \cite{EberhardtScholbach:Integral}, our construction gives a similarly defined functor $\Upsilon$ which carries a ``logarithm of monodromy'' map already integrally.

\theo
\thlabel{Ups.intro} 
The unipotent nearby cycles functor $\Upsilon$ enjoys the following properties.
\begin{enumerate}
  \item The functor $\Upsilon$ is naturally functorial (in a highly structured manner) with respect to smooth pullbacks and proper pushforwards. In addition, it is lax compatible with exterior products (cf.~\thref{Upsilon.all.in.one}). 
  \item \label{item--Ups.WT} If $X$ carries a Whitney--Tate stratification (compatible with the strata $\Gm$ and $0$ in $\A^1$), then $\Upsilon$ preserves stratified Tate motives (cf.~\thref{Upsilon preserves Tate}). 
  \item \label{item--Ups.CT} 
  If $\Gm$ (regarded as a constant group scheme over $\A^1$) acts on $X$, with fixed point locus $X^0 / \A^1$, then $\Upsilon_{X / \A^1}$ (applied to $\Gm$-monodromic motivic sheaves) is compatible with $\Upsilon_{X^0 / \A^1}$ under hyperbolic localization (cf.~\thref{prop-Gm-Upsilon}).
\end{enumerate}
\xtheo

Statement \refit{Ups.WT} is crucial to be able to use $\Upsilon$ in the context of Beilinson--Drinfeld Grassmannians, while \refit{Ups.CT} is the key to actual computations, since it allows us to eventually replace $G$ by a maximal torus.
The \(\ell\)-adic version of \refit{Ups.CT} appeared in \cite{Richarz:Spaces} and has already seen many applications, such as \cite{HainesRicharz:TestFunctions}.
We refer to \refsect{unipotent} for further results, including the setup of (unipotent) nearby cycles for families over $\A^n$.

\subsection{The central functor} 
We construct the central functor in \thref{defi--nearby-hecke} as
$$\mathsf{Z}_{\mathbf{f}} \colon \DM(L^+G \backslash \Gr) \rightarrow \NilpQ \DM(L^+G_{\mathbf f} \backslash \Fl_{\mathbf f}), \quad \mathsf{Z}_{\mathbf{f}}(\mathcal{F}) = \Upsilon (\mathcal{F} \boxtimes \Z_{\Gm}),$$ 
where $\Upsilon$ denotes the (unipotent) nearby cycles functor for the Beilinson--Drinfeld Grassmannian $\Gr_{\mathcal G_{\mathbf f}}$ introduced above (or, more precisely its extension to $L^+G$-equivariant sheaves developed in \refsect{KanInd}).
Its main properties are as follows, cf.~\thref{CentralityIso} and \thref{theo-moniodal-functor}.
In the following theorem, $\DTM(L^+ G \setminus \Gr)$ denotes the full subcategory of $\DM(L^+G \setminus \Gr_G)$ spanned by objects whose underlying non-equivariant motive restricts to a Tate motive on each $L^+G$-orbit, cf.~\refsect{stratified Tate motives} and \thref{notation DTM Gr} (and the references there).
The category $\DTM(L^+ G_{\mathbf f} \setminus \Fl_{\mathbf f})$ is defined similarly.

\theo\thlabel{thm--central functor}
Let $\mathcal F, \mathcal F' \in \DTM(L^+ G \setminus \Gr)$ and $\mathcal M \in \DTM(L^+ G_{\mathbf f} \setminus \Fl_{\mathbf f})$.
\begin{enumerate}
  \item \label{item--Z Tate} The functor $\mathsf{Z}_{\mathbf{f}}$ preserves (anti-effective) stratified Tate motives, and hence restricts to a functor
  \[\mathsf{Z}_{\mathbf{f}}\colon \DTM(L^+ G \setminus \Gr_G)^{(\anti)} \to \DTM(L^+G_{\mathbf{f}} \backslash \Fl_{\mathbf{f}})^{(\anti)}.\]
  \item  The functor $\mathsf{Z}_{\mathbf{f}}$ takes values in central sheaves with respect to the convolution product $\star$ in the sense that there are isomorphisms 
  $$\mathsf Z_\bbf (\calF) \star \calM \cong \calM \star \mathsf Z_\bbf(\calF).$$
  \item \label{item--Z monoidal} There is an isomorphism
  $$\mathsf{Z}_{\bbf}(\calF) \star \mathsf{Z}_{\bbf}(\calF') \cong \mathsf{Z}_{\bbf}(\calF \star \calF').$$
  \end{enumerate}
At least on the level of the homotopy categories, these isomorphisms endow $\mathsf Z_{\mathbf f}$ with the structure of a monoidal functor to the Drinfeld center,
  $$\Ho(\DTM(L^+ G \setminus \Gr )^{(\anti)}) \r \mathcal{Z}(\Ho(\DTM(L^+ G_\bbf \setminus \Fl_{\mathbf f})^{(\anti)})).$$
\xtheo

The proof follows similar reasoning as in the work of Achar--Riche \cite{AcharRiche:Central}. To state further properties of $\mathsf{Z}_{\mathbf{f}}$, we recall the abelian category of mixed Tate motives $\MTM(L^+G_{\mathbf{f}} \setminus \Fl_{\mathbf{f}})$, which arises as the heart of a t-structure on $\DTM(L^+G_{\mathbf{f}} \setminus \Fl_{\mathbf{f}})$. At least with $\Q$-coefficients, a compact object lies in the heart if and only if its Betti realization is a perverse sheaf.
The main result of \cite{CassvdHScholbach:MotivicSatake} identifies $\MTM(L^+G \setminus \Gr_G)$ with a motivic category of integral representations of $\hat{G}$.
For the restriction of $\mathsf{Z}_{\mathbf{f}}$ to \(\MTM(L^+G \setminus \Gr_G)\), we upgrade \thref{thm--central functor} with a compatibility between the centrality and monoidality isomorphisms (\thref{theo-E2}), as well as a compatibility with the commutativity constraint (\thref{prop--fusion-v-nearby}) coming from fusion \cite[§5.3]{CassvdHScholbach:MotivicSatake}. The former compatibility uses nearby cycles over a 2-dimensional base, which was one of our motivations to revisit the nearby cycles functor. 
Upgrading the monoidality of $\mathsf Z_\bbf$ to an \ii-categorical statement does not seem to require new geometric ideas, but rather an \ii-categorical study of the compatibility of Koszul dualities for varying (co)algebras, which we defer to another paper. 

\subsection{Wakimoto motives}
One property of \(\mathsf{Z}_{\mathbf{f}}\) that is not yet addressed in \thref{thm--central functor} is t-exactness.
For $\ell$-adic and Betti sheaves, t-exactness for the perverse t-structure is a general property of nearby cycles. Such a general result is currently not available for motivic sheaves, but we are able to prove the t-exactness of  \(\mathsf{Z}_{\mathbf{f}}\) by other methods.
The t-exactness when \(\mathbf{f}=\mathbf{f}_0\) is hyperspecial follows from \thref{prop--satake-exact}. On the opposite end, when \(G_{\mathbf{a}_0} = \mathcal{I}\) is an Iwahori group scheme, we adapt arguments of \cite{ALWY:Gaitsgory}, by using a motivic refinement of Wakimoto sheaves as follows.

Classically, Wakimoto sheaves are perverse sheaves introduced by Mirkovi\'{c} to categorify a maximal commutative subalgebra of \(H^{\mathcal{I}}\).
Geometrically, this manifests by the fact that every perverse central sheaf on \(\Fl\) admits a filtration by Wakimoto sheaves.
Motivically, we then construct a Wakimoto functor
\(\JWak_\mu \colon \DTM( \mathcal{I} \setminus \Spec \Z) \rightarrow \DTM( \mathcal{I} \setminus \Fl)\), cf.~\cite{ArkhipovBezrukavnikov:Perverse}, where \(\mu\in X_*(T)\) and $T \subset B$ are a maximal torus and Borel, and prove it is t-exact.
We also use the t-exact constant term functor \(\CT_{B^-}[-\deg] := \oplus \CT_{B^-}^\mu[-\langle 2\rho,\mu\rangle]\) on $\DTM(L^+G \setminus \Gr)$  from \cite[Definition 5.3]{CassvdHScholbach:MotivicSatake}, where $B^-$ is the opposite Borel and $2\rho$ is the sum of the positive roots, which under the Satake equivalence corresponds to the restriction functor $\Rep_{\hat G} \rightarrow \Rep_{\hat T}$.

\theo
For any \(\Ff\in \DTM(L^+G \setminus \Gr)\), the central motive \(\mathsf{Z}_{\mathbf{a}_0}(\Ff)\) admits a (necessarily unique) \(X_*(T)\)-filtration with graded pieces given by \(\JWak_{\mu}(\CT_{B^-}^\mu(\Ff))(-\langle 2\rho,\mu\rangle)[-\langle 2\rho,\mu\rangle]\) (cf.~\thref{Grad-CT}).
In particular, \(\mathsf{Z}_{\mathbf{a}_0}\) is t-exact.
\xtheo

This t-exactness will be a crucial ingredient in a motivic refinement of \cite{ArkhipovBezrukavnikov:Perverse,Bezrukavnikov:Two}, as alluded to above. Along the way we must also give alternative proofs of two related t-exactness results in \thref{StandardExact} (cf.~\cite[Lemma 4.21]{EberhardtKelly:Mixed}) and \thref{lemm--leftright} (cf.~\cite[Proposition 4.6]{Achar:ModularPerverseSheavesII}), whose analogues in \cite{AcharRiche:Central, ALWY:Gaitsgory} are proved using Artin vanishing.
Finally, following a suggestion of Achar and discussions with Louren\c{c}o, we use the above Wakimoto filtration to show that the central functor is t-exact for general facets in \thref{t-exactness general facets}.

 \subsection{Generic Hecke algebras}
In \cite[Definition 6.34]{CassvdHScholbach:MotivicSatake}, the present authors defined the generic spherical Hecke algebra \(\mathcal{H}^{\sph}(\qq)\). 
This is a \(\Z[\qq]\)-algebra, where \(\qq\) is a formal variable, which specializes to the spherical Hecke algebra of functions on $G (\Fq \rpot t)$ under the map $\qq \mapsto q$. 
The motivic Satake equivalence gives an isomorphism between \(\mathcal{H}^{\sph}(\qq)\) and the representation ring of a certain Vinberg monoid for $\hat{G}$ appearing in Zhu's work  \cite{Zhu:Integral}, and hence with the Grothendieck ring of a suitable category of anti-effective motives on \(\Gr\). 
On the other hand, Vign\'{e}ras has defined the generic Iwahori--Hecke algebra \(\mathcal{H}^{\mathcal{I}}(\qq)\) \cite{Vigneras:Algebres}. 
This can also be realized as the Grothendieck ring of a suitable category of anti-effective motives on the full affine flag variety, cf.~\thref{comparison with Vigneras}.

These definitions of generic Hecke algebras (which also appear in \cite{PepinSchmidt:Generic} for $G  = \GL_2$) rely on combinatorics: Kazhdan--Lusztig theory for 
$\mathcal{H}^{\sph}$, and the Iwahori--Matsumoto presentation for \(\calH^{\mathcal{I}}\).
However, the description in terms of Grothendieck rings suggests that, for a general facet \(\mathbf{f}\), we can \emph{define} the generic parahoric Hecke algebra \(\calH^{\mathbf{f}}(\qq)\) at \(G_{\mathbf{f}}\)-level as the Grothendieck ring \(K_0(\DTM(L^+G_{\mathbf{f}} \backslash \Fl_{\mathbf{f}})^{\antilc})\), cf.~\thref{Defi generic parahoric} for details.
These are \(\Z[\qq]\)-algebras which specialize to the usual parahoric Hecke algebras by \thref{generic parahoric}.
Moreover, we can again decategorify \(\mathsf{Z}_{\mathbf{f}}\) to relate \(\calH^{\sph}(\qq)\) to the center of \(\calH^{\mathbf{f}}(\qq)\). The situation can be summarized as follows:
\[\begin{tikzcd}
	\DTM(L^+G\backslash \Gr)^{\antilc} \arrow[r, "\mathsf{Z}_{\mathbf{f}}"] \arrow[d, mapsto, "K_0"] & \DTM(L^+G_{\mathbf{f}}\backslash \Fl_{\mathbf{f}})^{\antilc} \arrow[d, mapsto,"K_0"]\\
	\calH^{\sph}(\qq) \arrow[d, mapsto, "-\otimes_{\Z[\qq]{,}\qq\mapsto q} \Z"] \arrow[r] & \calH^{\mathbf{f}}(\qq) \arrow[d, mapsto, "-\otimes_{\Z[\qq]{,}\qq\mapsto q} \Z"]\\
	\calH^{\sph} \arrow[r] & \calH^{\mathbf{f}},
\end{tikzcd}\]
where the two lower horizontal arrows have central image.

For the particular case \(\mathbf{f}=\mathbf{a}_0\), we get a \emph{generic Bernstein isomorphism}, generalizing Bernstein's isomorphism \(Z(\calH^{\mathcal{I}}) \cong \calH^{\sph}\), cf.~\thref{generic Bernstein iso}.

\theo If \(\mathbf{f}=\mathbf{a}_0\) is a standard alcove, then decategorifying the central functor $\mathsf{Z}_{\mathbf{a}_0}$ induces an isomorphism between \(\mathcal{H}^{\sph}(\qq)\) and the center of \(\mathcal{H}^{\mathcal{I}}(\qq)\).
\xtheo

Along the way to this theorem, we also obtain independence-of-$\ell$ for related results in  \cite[§5.3]{AcharRiche:Central} which used mixed $\ell$-adic sheaves.
We note that we also expect the theorem to hold without the assumption \(\mathbf{f}=\mathbf{a}_0\), and in that case the natural map \(\mathcal{H}^{\sph}(\qq)\to \mathcal{H}^{\mathbf{f}}(\qq)\) is still injective with central image.
However, to show it surjects onto the center, we use structural results of Vignéras \cite{Vigneras:Algebres} about the center of \(\mathcal{H}^{\mathcal{I}}(\qq)\), which are not available for general parahoric Hecke algebras.
In fact, it seems that when working with \(\Z\)-coefficients, a Bernstein isomorphism relating the centers of parahoric Hecke algebras with the spherical Hecke algebra is not yet known for general parahorics.

\subsection*{Acknowledgements} We thank Pramod Achar, Johannes Anschütz, Joseph Ayoub, Katia Bogdanova, Rızacan Çiloğlu, Sean Cotner, Ishai Dan-Cohen, Dennis Gaitsgory, Martin Gallauer, Adeel Khan, Achim Krause, Jo\~{a}o Louren\c{c}o, Fabien Morel, Guglielmo Nocera, Cédric Pépin, Simon Pepin Lehalleur, Timo Richarz and Alberto Vezzani for helpful discussions around the subject matter and / or valuable feedback on earlier drafts of this paper. We are also grateful to the referees for their helpful comments and corrections. 

J.S. was supported by the EuropeanUnion – Project 20222B24AY (subject area: PE - Physical Sciences and Engineering) ``The arithmetic of motives and L-functions'', and also by the Max-Planck-Institute for Mathematics in Bonn.
T.v.d.H. was supported by the European Research Council (ERC) under the European Union’s Horizon 2020 research and innovation programme (grant agreement 101002592), and by the Deutsche Forschungsgemeinschaft (DFG), through the TRR 326 GAUS (project number 444845124). R.C. was supported by the National Science Foundation under Award No. 2103200 and No. 1840234.

\section{Preliminaries}
\label{sect--prelim}
\subsection{Higher algebra}

In this section, we recall and establish some basic notions around tensor products in the context of \ii-categories and Koszul duality.
\subsubsection{Monoidality and adjoints}
\label{sect--monoidal aspects}
We will work with the following symmetric monoidal \ii-categories \cite[§3.1]{BlumbergGepnerTabuada:Universal}
$$\Cat^\perf \stackrel[\cong]\Ind \r \PrLomegaSt \r  
\PrLSt.\eqlabel{PrLSt etc}$$
From right to left, we have the category of stable presentable \ii-categories, with colimit-preserving (equivalently, left adjoint) functors; 
then its (non-full) subcategory of compactly generated stable categories, with functors preserving compact objects \cite[Notation~5.3.2.8]{Lurie:HA}, and then the \ii-category of small stable idempotent complete \ii-categories with exact functors.

For two stable and compactly generated \ii-categories $C$ and $D$ a continuous functor $f: C \r D$ preserves compact objects iff its right adjoint $f^R$ is continuous, i.e., if $f^R$ is also a left adjoint.
Therefore, the operation of passing to right adjoints yields an equivalence
\begin{equation}
\PrLomegaSt \stackrel \cong \r (\PrLRomegaSt)^\opp, C \mapsto C, f \mapsto f^R, \label{PrLR equivalence}
\end{equation}
where at the right $\PrLRomegaSt$ denotes the \ii-category of presentable compactly generated stable \ii-categories with functors being \emph{both left and right adjoints} (but not necessarily preserving compact objects).

Recall Lurie's tensor product on $\PrL$ and on its subcategory $\PrLSt$ \cite[§4.8.1, Proposition~4.8.2.18]{Lurie:HA}.

\lemm
\thlabel{PrLomega.equivalence}
The categories $\PrLomegaSt \subset \PrLSt$ and $\PrLRomegaSt \subset \PrLSt$ are symmetric monoidal subcategories.
Furthermore, with respect to these symmetric monoidal structures, the equivalence in \eqref{PrLR equivalence} is symmetric monoidal.
\xlemm

\pf
The condition of being compactly generated is preserved under Lurie's tensor product \cite[Lemma~5.3.2.11]{Lurie:HA}.

Let $f_i : C_i \r C'_i$ be finitely many functors in $\PrLomegaSt$, i.e., $f_i$ preserve colimits and compact objects (and $C_i$ and $C'_i$ are stable compactly generated categories).
Since their right adjoints $f_i^R$ are also left adjoints, we can consider the functor $\bigotimes_i (f_i^R)$.
Let $f := \bigotimes f_i$.
There is a natural map $\bigotimes f_i^R \r f^R$.
To see it is an isomorphism let $c_i \in C_i$ be arbitrary objects. The elementary tensors $\boxtimes c_i$ generate the category $\bigotimes C_i$. 
Thus it suffices to note that
$$\Maps_{\bigotimes C_i}(\boxtimes c_i, \boxtimes f_i^R(c'_i)) = \Maps_{\bigotimes C_i}(\boxtimes c_i, f^R(\boxtimes c'_i)).$$
Indeed, by \cite[Proposition~I.7.4.2]{GaitsgoryRozenblyum:StudyI}, both sides are canonically identified with $\bigotimes_i \Maps_{C_i}(f_i (c_i), c'_i)$.

Let $\Cat_\Stable^\L$ (resp.~$\Cat_\Stable^\R)$ be the (non-full) subcategory of $\Cat_\Stable$ (stable \ii-categories and exact functors) where only left (resp.~right adjoints) are allowed.
The cartesian symmetric monoidal structure on $\Cat_\Stable (\subset \Cat_\infty)$ \cite[Theorem~1.1.4.4]{Lurie:HA} respects these subcategories. We indicate these by $\Cat_\Stable^{\L / \R, \x}$ (even though on these subcategories it is no longer cartesian).
We have an equivalence $\Cat_\infty^\L \cong (\Cat_\infty^\R)^\opp$ given by the unicity of adjoints \cite[Proposition~5.2.6.2]{Lurie:Higher}. According to the construction of this equivalence, this equivalence is symmetric monoidal.
This respects the (symmetric monoidal) ``$\Stable$''-subcategories, giving rise to the commutative triangle at the right below:
$$\xymatrix{
\Pr_{\omega,\Stable}^{\L, \t} \ar[d] \ar[rr] & & \Cat_\Stable^{\L, \x} \ar[d]^\cong \ar[r] & \Fin_* \\
(\Pr_{\omega,\Stable}^{\L\R, \t})^\opp \ar[rr] & & (\Cat_\Stable^{\R, \x})^\opp \ar[ur].
}$$
The left horizontal maps are the canonical inclusions \cite[Notation~4.8.1.2, Proposition~4.8.1.15]{Lurie:HA}.
The left vertical arrow exists, and is an equivalence in each fiber over $\langle n \rangle \in \Fin_*$, by the above (applied to each $f_i$ individually). It is symmetric monoidal by the isomorphisms $(\bigotimes f_i)^R = \bigotimes (f_i^R)$.
\xpf

Recall that a $\t$-category $D$ is \emph{rigid} if it is stable, and compactly generated by a symmetric monoidal subcategory $D_0 \subset D$ in which every object is dualizable.
I.e., $D = \Ind(D_0)$, and $(-)^\dual : D_0 \stackrel \cong \r D_0^\opp$ is an equivalence.
The Ind-completion of this equivalence gives an equivalence $D^\dual = D$, where $\dual$ refers to the dualization with respect to the tensor product in $\PrLSt$ \cite[Ch. I, Proposition~7.3.5]{GaitsgoryRozenblyum:StudyI}.

\lemm
\thlabel{ModD.equivalence}
For a rigid $\t$-category $D$, there are equivalences
$$\Mod_D(\PrLSt) = \Mod_{D^\dual}(\PrLSt) = \coMod_D(\PrLSt),$$
which are given by the identity on the level of the underlying objects in $\PrLSt$.
\xlemm

\pf
The first equivalence arises from the isomorphism (of commutative algebra objects) $D \cong D^\dual$.
The second is a generality about (co)modules over dualizable objects: if $d$ is a dualizable object in a symmetric monoidal \ii-category $C$, the forgetful functor $\coMod_d C \r C$ admits a \emph{left} adjoint, given by $d^\dual \t - $.
It is therefore monadic, and thus admits a functor to $\Mod_{d^\dual} C$, which is an equivalence.
\xpf

\rema
Let $C$ be an object in $\Mod_D(\PrLSt)$, with action map $a : D\t C \r C$. Under the above equivalence, it corresponds to the $D$-comodule whose coaction is given by $a^R$, the right adjoint of $a$ \cite[Lemma~I.9.3.2]{GaitsgoryRozenblyum:StudyI}.
Furthermore, we have $a = \mu \t_D \id$, for the multiplication $\mu: D \t D \r D$.
Using again that passing to right adjoints is compatible with $\t$, we obtain the following description of the $D$-coaction on $C$: 
$$a^R = \mu^R \t_D \id.\eqlabel{coaction}$$
\xrema

\lemm
For a rigid category $D$, the right adjoint of the multiplication map, $\mu^R : D \r D \t D$, is a lax symmetric monoidal functor. In particular, there is a natural transformation like so:
$$\xymatrix{
D \t D \ar[d]^\mu \ar[r]^{\mu^R \t \mu^R} & D^2 \t D^2 \ar[d]^{\mu_{D^2}} \ar@{=>}[dl] \\
D \ar[r]^{\mu^R} & D^2.
}$$
\xlemm

\pf
The multiplication $\mu$ is symmetric monoidal, so its right adjoint is lax symmetric monoidal.
\xpf

\coro
\thlabel{unit map lax monoidal}
Let $D$ be rigid and let
$$U : \Mod_D (\PrLSt) \r \PrLSt$$
be the forgetful functor.
According to the above, its left adjoint and its right adjoint agree up to equivalence, and are denoted by $L$.
The unit map (for it being the right adjoint)
$$\alpha: \id \r LU$$
is then a natural transformation of lax symmetric monoidal endofunctors of $\Mod_D(\PrLSt)$.
\xcoro

\pf
The map $\alpha$ is given by applying $\id \t_D - $ to $\mu^R : D \r D \t D$. 
\xpf

\lemm
\thlabel{generalities.tensor} \cite[Proposition~4.1]{BenZviFrancisNadler:Integral}
For a symmetric monoidal presentable \ii-category $C$, a commutative algebra object $A \in \CAlg(C)$ and a $C$-module $M$ (in $\PrLSt$), there is an equivalence 
$$\Mod_A(C) \t_C M \stackrel \cong \r \Mod_A(M).\eqlabel{Mod.A.C}$$
\xlemm

\subsubsection{Koszul duality}
\label{sect--Koszul duality}

The following discussion of an element of Koszul duality will later be used in order to establish the monodromy map on unipotent nearby cycles.

Let $\grAb$ denote the (ordinary) category of $\Z$-graded abelian groups, and $\grAb_\free \subset \grAb$ its full subcategory of degreewise free abelian groups.
We also use $\grCh$, the category of unbounded chain complexes of $\Z$-graded abelian groups.
These categories are symmetric monoidal for the usual (underived) tensor product of graded objects, given by Day convolution along $\Z$.
Finally, let $\grD := \D(\grAb)$ be the derived \ii-category of $\Z$-graded abelian groups \cite[Definition~1.3.5.8]{Lurie:HA}.
This is a symmetric monoidal \ii-category when endowed with the derived tensor product \cite[Remark~7.1.2.12]{Lurie:HA}, and again Day convolution along $\Z$.
The composite $\grAb_\free \subset \grAb \subset \grD$ is symmetric monoidal with respect to these tensor structures.

Let $\Sigma := \bigoplus_{n \ge 0} \Z(n) \in \grAb_\free$. It is the free commutative algebra generated by $\Z(1)$ (i.e. $\Z$ placed in graded degree $-1$), in the abelian category $\grAb$.
Endowed with the counit being the natural augmentation map $\Sigma \r \Z$ and the comultiplication 
\begin{align}
\nabla_\Sigma : \Sigma \r & \Sigma \t \Sigma, \nonumber \\
\Z(1) \ni t \mapsto & t \t 1 + 1 \t t,\label{comultiplication Sigma}
\end{align}
$\Sigma$ is a bicommutative (commutative and cocommutative) bialgebra object in $\grAb_\free$. Using the symmetric monoidal functor to $\grD$, we will also regard $\Sigma$ as a bicommutative bialgebra object in $\grD$. 

Note that $\Sigma$ is the free associative algebra in $\grD$ generated by $\Z(1)$. The map $\Sigma \r \Sym_{\grD} \Z(1) = \bigoplus_{n \ge 0} \Z(n) / S_n$ to the free commutative algebra on $\Z(1)$, taken in $\grD$, is only an isomorphism after taking 0-th cohomology, or after tensoring with $\Q$.
Indeed, the \emph{derived} coinvariants of the trivial action of $S_n$ on $\Z(n)$ have torsion cohomologies in positive cohomological degrees. 

Since $\Sigma$ is a commutative algebra in $\grD$, the \ii-category of left $\Sigma$-modules $\Mod_\Sigma(\grD)$ carries the derived tensor product $- \t_\Sigma -$.
In the sequel, however, we need a symmetric monoidal structure whose tensor product is the derived tensor product over $\Z$ instead.
To formally construct that symmetric monoidal structure, we will present $\Mod_\Sigma(\grD)$ as the \ii-category underlying a symmetric monoidal model category whose construction is as follows: $\grCh$ carries the projective model structure, e.g. \cite[Proposition~7.1.2.8]{Lurie:HA}.
Being degreewise free, $\Sigma$ is cofibrant in $\grCh$, so the model structure on $\Mod_\Sigma(\grCh)$ right-induced from the forgetful functor $U : \Mod_\Sigma(\grCh) \r \grCh$ exists \cite[Remark~4.2]{SchwedeShipley:Algebras}. We indicate it by $\Mod_\Sigma(\grCh)^\rightind$. 
Its weak equivalences and fibrations are created by $U$; in particular weak equivalences are ($\Sigma$-module) maps which are quasi-isomorphisms in each graded degree.
Its generating (acyclic) cofibrations are of the form $\Sigma (k) \t c$, where $c$ is a generating (acyclic) cofibration in $\Ch(\Ab)$ (for those, see, e.g.~\cite[Definition~2.3.3]{Hovey:Model}). 

On $\Mod_\Sigma(\grCh)$, we now also consider the model structure \emph{left-induced} along the forgetful functor $U$ to $C$, i.e., $U$ creates \emph{co}fibrations and weak equivalences. We denote it by $\Mod_\Sigma(\grCh)^{\leftind}$.
To ensure its existence we use \cite[Corollary~2.7]{GarnerKedziorekRiehl:Lifting}: if $f \in \Mod_\Sigma(\grCh)$ is a map having the right lifting property with respect to all maps $g \in \Mod_\Sigma(\grCh)$ such that $U(g)$ is a cofibration in $\grCh$, we need to show $f$ is a weak equivalence. 
Indeed, for any cofibration $h$ in $\grCh$, the map $g := \Sigma \t_\Z h$ is such that $U(g)$ is a cofibration, since $\Sigma$ is cofibrant in $\grCh$.
Such maps $g$ generate the cofibrations in $\Mod_\Sigma(\grCh)^\rightind$. Maps having the right lifting property against such $g$ are acyclic fibrations in $\Mod_\Sigma(\grCh)^\rightind$, in particular these are degreewise quasi-isomorphisms. This confirms the existence of the left-induced model structure.

Given that the \ii-category underlying $\grCh$ is $\grD$, and given \cite[Theorem 4.3.3.17]{Lurie:HA}, the \ii-category $\Mod_\Sigma(\grD)$ is equivalent to the \ii-category underlying $\Mod_\Sigma(\grCh)^\leftind$ or the Quillen equivalent $\Mod_\Sigma(\grCh)^\rightind$.

We regard $\Mod_\Sigma(\grCh)$ as a symmetric monoidal category whose tensor product is computed on the underlying objects in $\grCh$, endowed with the $\Sigma$-action 
$$\Sigma \t (M \t N) \stackrel {\nabla_\Sigma} \r \Sigma \t \Sigma \t M \t N \stackrel{\act_M \t \act_N} \r M \t N.$$
Equivalently, expressing the $\Sigma$-action on $M$ as a map $e_M : M(1) \r M$, the map $(M \t N)(1) \r M \t N$ is the sum $e_M \t \id_N + \id_M \t e_N$.
With this symmetric monoidal structure, $\Mod_\Sigma(\grCh)^{\leftind}$ is a symmetric monoidal model category since both (acyclic) cofibrations and the tensor product are created by the forgetful functor $U$.

\defi
We denote by $(\Mod_\Sigma(\grD), \t_\Z)$ the underlying symmetric monoidal \ii-category \cite[Definition~4.1.7.6]{Lurie:HA} of $\Mod_\Sigma(\grCh)^\leftind$.
\xdefi

\rema
By the setup, the forgetful functor $(\Mod_\Sigma(\grD), \t_\Z) \r \grD$ is symmetric monoidal, i.e., the (derived) tensor product of $\Sigma$-modules is computed on the underlying objects in $\grD$. In particular, $\Z$ (not $\Sigma$) is the monoidal unit.
We expect that the symmetric monoidal \ii-category structure constructed above is compatible with the one constructed by Beardsley \cite[Theorem~3.18]{Beardsley:Bialgebras} in much greater generality, but we have not investigated this in detail.
\xrema

We let 
\begin{equation}
\OldLambda := \Z \oplus \Z(-1)[-1] \in \grCh.
\label{notation OldLambda}
\end{equation}
Endowed with the multiplication such that $\Z(-1)[-1]$ is square zero, it is a cdga, i.e., a commutative algebra object in $\grCh$.
It is cofibrant in $\grCh$.
In parallel to $\Mod_\Sigma(\grCh)^\rightind$, we have the right-induced model structure $\Mod_\OldLambda(\grCh)^\rightind$, which is a symmetric monoidal model category with respect to the (usual) tensor product $- \t_\OldLambda -$.
Its underlying symmetric monoidal \ii-category is equivalent to $\Mod_\OldLambda(\grD)$ (with its usual derived tensor product over $\OldLambda$) by applying \cite[Theorem~7.1.2.13]{Lurie:HA} to $\OldLambda$.

The augmentation maps turn $\Z$ into a $\Sigma$-$\OldLambda$-bimodule (in $\grCh$). 
Moreover, the resulting functor
$$\Z \t_\OldLambda - : \Mod_\OldLambda(\grCh)^\rightind \r \Mod_\Sigma(\grCh)^\leftind\eqlabel{adjunction model categories}$$
is a symmetric monoidal left Quillen functor, i.e., preserves the monoidal unit and tensor products up to isomorphism, and preserves (acyclic) cofibrations. Its right adjoint is the (underived) functor $\IHom_{\Sigma}(\Z, -)$.

\defilemm
\thlabel{Nilp}
The derived tensor product functor
$$\Z \t_\OldLambda - : (\Mod_{\OldLambda}(\grD), \t_\OldLambda)  \r (\Mod_\Sigma(\grD), \t_\Z),\eqlabel{adjunction infty categories}$$
is fully faithful and symmetric monoidal.
Therefore, if we write $\Nilp$ for the essential image of this functor, we obtain an equivalence (of presentably symmetric monoidal \ii-categories)
$$(\Mod_\OldLambda(\grD), \t_\OldLambda) \stackrel \cong \r \Nilp \subset (\Mod_\Sigma(\grD), \t_\Z). \eqlabel{Nilp diagram}$$
\xdefilemm

Henceforth, we abbreviate $\Mod_\Sigma := \Mod_\Sigma(\grD)$ and $\Mod_\OldLambda := \Mod_\OldLambda(\grD)$, with the symmetric monoidal \ii-category structures constructed above.

\pf
We obtain the functor \refeq{adjunction infty categories} from the left Quillen functor \refeq{adjunction model categories} by passing to the underlying \ii-categories.
The latter, and therefore also the former, functor is symmetric monoidal by the above discussion.

The usual short exact sequence in $\Mod_\Sigma(\grAb)$ 
$$0 \r \Sigma(1) \r \Sigma \r \Z \r 0,\eqlabel{resolution Z Sigma}$$
is a resolution by cofibrant $\Sigma$-modules. It shows that $\Z$ is compact in $\Mod_\Sigma(\grD)$, so that the derived right adjoint $\IHom_{\Sigma}(\Z, -)$ (the derived inner Hom in $\Mod_{\Sigma}(\grD)$) preserves colimits. It also implies that the unit map $\OldLambda \r \IHom_\Sigma(\Z, \Z)$ for the adjunction given by \refeq{adjunction infty categories} is a quasi-isomorphism.
This implies that this unit map is a quasi-isomorphism for all objects in $\Mod_{\OldLambda}(\grD)$.
\xpf

\rema
\thlabel{remarks Nilp}
For $\Q$-linear (but otherwise arbitrary) categories, a similar symmetric monoidal Koszul-type equivalence appears in \cite[Proposition~2.32]{BindaGallauerVezzani:Motivic}, which in turn refines non-monoidal results due to Raskin and Barkan (see op.~cit.~for references).

The objects of $\Mod_{\Sigma}(\grD)$ can be described as pairs $(c \in \D(\grAb), \phi : c(1) \r c)$. 
Being equivalent to $\Mod_\OldLambda(\grD)$, the full subcategory $\Nilp \subset \Mod_\Sigma(\grD)$ is compactly generated by $\Z(k)$, $k \in \Z$, regarded as a trivial $\Sigma$-module (i.e., $\Sigma$ acts via its augmentation).
Thus, for a compact object $(c, \phi) \in \Nilp$  (which is a retract of a finite colimit of such generators), $\phi$ is nilpotent.
An arbitrary object $(c, \phi) \in \Nilp$ is a filtered colimit of compact objects, so one may think of $\phi$ as being locally nilpotent.

If we replace $\D(\grAb)$ by $\D(\Ab)$ (i.e., suppress all gradings), then $\Nilp$ can be characterized geometrically as the full subcategory of the derived category of quasi-coherent sheaves on $\A^1_\Z = \Spec \Z[t]$ spanned by objects supported at $0$, i.e., annihilated by inverting $t$, cf.~\cite[Theorem~2.1]{DwyerGreenlees:Complete}.
From this perspective, the symmetric monoidal structure on $\Nilp$ is the convolution along the addition on $\A^1$.
\xrema

\subsection{Motives}

\nota
\thlabel{notation S}
Throughout we fix a connected base scheme $S$ that is smooth of finite type over a Dedekind ring or a field.
Unless specifically mentioned otherwise, all schemes are supposed to be of finite type over $S$.
We often abbreviate $\GmX S$ by $\Gm$ etc.
\xnota
For a scheme $X / S$, we denote the stable presentable (symmetric) monoidal \ii-category of motives (or motivic sheaves) with integral coefficients by $\DM(X)$.
The corresponding category with rational coefficients is denoted by $\DM(X, \Q)$.

In the generality of \thref{notation S}, its construction is due to Spitzweck \cite{Spitzweck:Commutative} as the category of modules over a certain $E_\infty$-ring spectrum in the stable $\A^1$-homotopy category $\SH(X)$.
If $S$ is smooth over a field, that spectrum is (isomorphic to) Voevodsky's spectrum representing motivic cohomology (with integral coefficients).

\subsubsection{Six functors}
The six functors on motives can be packaged into a lax (symmetric) monoidal functor
$$\DM := \DM^*_! : \Corr \r \PrLSt \eqlabel{DM.blah}$$
Its domain is the \ii-category of correspondences. Its objects are (finite type $S$-)schemes. Informally, its morphisms are correspondences of schemes over $S$, i.e., diagrams of the following form:
$$Y \stackrel g \gets Z \stackrel f \r X.\eqlabel{correspondences}$$
The above functor sends $X \mapsto \DM(X)$, while a correspondence as above is mapped to $g_! f^*$. 
For $X, X' / S$, the map $(X, X') \r X \x_S X'$ in the symmetric monoidal \ii-category $\Corr$ is sent to the exterior product functor $\DM(X) \t \DM(X') \stackrel{\boxtimes} \r \DM(X \x_S X')$.
Part of the lax monoidality of the functor is the assertion that $\boxtimes$ is compatible with *-pullback and !-pushforward.
We refer to \cite{GaitsgoryRozenblyum:StudyI,LiuZheng:Enhanced,Mann:p-adic} for in-depth discussions of the category of correspondences, to \cite{Khan:Motivic,Hoyois:Six} for its application in the context of motivic sheaves, and to \cite[§A]{RicharzScholbach:Motivic} for an exposition also including the lax monoidality of $\DM$.

Any $S$-scheme $X$ is a coalgebra object in $\Corr$ (the comultiplication being the diagonal $X \r X \x_S X$).
The lax symmetric monoidality of $\DM$ in \refeq{DM.blah} then endows $\DM(X)$ with the structure of a symmetric monoidal \ii-category.
In addition, for any scheme $T / S$ (e.g., $T = \GmX S$), any $T$-scheme $f : X \r T$ is naturally a coalgebra object over $T$ (by means of the graph $\Gamma_f : X \r X \x_S T$). The lax symmetric monoidal functor above leads naturally to a lax symmetric monoidal functor 
$$\DM^*_! : \Corr_T \r \Mod_{\DM(T)} \PrLSt.\eqlabel{DM Mod}$$
Here $\Mod$ denotes the indicated category of modules with respect to the Lurie tensor product in $\PrLSt$, cf.~\refsect{monoidal aspects}.

We will also use the \ii-category $\DM(X)$ of motives on an ind-scheme $X$ (such as the affine flag variety), as defined and discussed in \cite[§2.4]{RicharzScholbach:Intersection}.

\subsubsection{Stratified Tate motives}
\label{sect--stratified Tate motives}
Recall from \cite{RicharzScholbach:Intersection, CassvdHScholbach:MotivicSatake} the following basic notation related to stratified Tate motives:
The presentable stable subcategory $\DTM(X) \subset \DM(X)$ generated by $\Z(n)$, $n \in \Z$ is called the category of \emph{Tate motives}. If we only allow $n \le 0$, we speak of \emph{anti-effective Tate motives}; the category is denoted by $\DTM(X)^\anti$.
In the sequel, we will write ``$(\anti)$'' to denote either the presence or the absence of the condition of anti-effectivity.

For a stratified (ind-)scheme $\iota: X^\dagger = \bigsqcup_{w \in W} X_w \r X$, we denote by $\DTM(X, X^\dagger)^{(\anti)} := \{M \in \DM(X), \iota^* M \in \DTM(X^\dagger)^{(\anti)}\}$ the category of (anti-effective) stratified Tate motives.
We will only use this category if the stratification is \emph{(anti-effective) Whitney--Tate}, which means that $\iota^* \iota_*$ preserves $\DTM(X^\dagger)^{(\anti)}$.
If $\iota$ is clear from the context, we also denote this category by $\DTM(X)^{(\anti)}$.
We call the objects therein \emph{(anti-effective) stratified Tate motives}. If furthermore the base change map \cite[Eqn.~(2.4)]{CassvdHScholbach:MotivicSatake} is an isomorphism then we say the stratification is \emph{universally} Whitney--Tate. 

\subsubsection{Reduced motives}
\label{sect--reduced motives}
For $X / S$, the category $\DM(X)$ is a module over $\DTM(S)$, i.e., is endowed with an action of the motivic cohomology of $S$.
The action of the higher motivic cohomology of $S$ often has no representation-theoretic significance. This led to the introduction of reduced motives in \cite{EberhardtScholbach:Integral}.
By definition, the category of reduced Tate motives over $S$ is $\DTMr(S) = \D(\grAb)$, the derived \ii-category of $\Z$-graded abelian groups. 
There is a canonical reduction functor $\red : \DTM(S) \r \DTMr(S)$ considered in op.~cit.
This functor can be thought of as modding out the higher motivic cohomology of $S$.
It admits a section $i : \DTMr(S) \r \DTM(S)$, sending $\Z(k)$ ($\Z$ in graded degree $-k$) to the $k$-fold Tate twist $\Z(k)$.

For $X / S$, the \ii-category of reduced Tate motives is defined as $\DTMr(X) := \DTM(X) \t_{\DTM(S), \red} \DTMr(S)$ and similarly $\DMr(X) := \DM(X) \t_{\DTM(S)} \DTMr(S)$.
The formalism of reduced motives enjoys the same six functor properties as $\DM$.
We will write $\DMrx$ or $\DTMrx$ to indicate that a statement holds for reduced and non-reduced (Tate) motives.

We let $\PrLgr := \Mod_{\D(\grAb)}(\PrLSt)$ be the \ii-category of stable presentable $\Z$-linear $\Z$-graded categories (with functors preserving that structure). 
Using the natural restriction functor along $i$, we will sometimes consider $\DM^*_!$ as the lax symmetric monoidal functor
$$\DM^*_! : \Corr \r \Mod_{\DM(S)}(\PrLSt) \stackrel \res \r \Mod_{\DTM(S)}(\PrLSt) \stackrel{\res_i} \r \Mod_{\D(\grAb)}(\PrLSt).\eqlabel{DM Prgr}$$

\subsubsection{Mixed Tate motives}
\label{sect--mixed Tate motives}
As in \cite[§2.3]{CassvdHScholbach:MotivicSatake}, we will make use of motivic t-structures. 
We consider the t-structure on $\DTM(S)$ whose aisle $\DTM(S)^{\le 0}$ is generated under colimits by the objects $\Z(k)$, $k \in \Z$.
Recall that if $S$ satisfies the Beilinson--Soulé conjecture (e.g., $S = \Spec \Z$), these generators lie in the heart of the t-structure. \emph{Whenever we consider mixed Tate motives, we will assume this condition to hold.}

If $X$ is an ind-scheme with an admissible Whitney--Tate stratification (for the definition of admissibility see \cite[Definition~2.9]{CassvdHScholbach:MotivicSatake}, for example this holds if the strata are isomorphic to some $\Gm^n \x_S \A^m$), then $\DTM(X, X^\dagger)$ carries a t-structure whose aisle $\DTM(X, X^\dagger)^{\le 0}$ is generated by the objects $\iota_{w!} \Z_{X_w}(k)[\dim_S X_w]$, for all strata $\iota_w : X_w \r X$ and $k \in \Z$.
The heart of this t-structure is denoted by $\MTM(X) := \MTM(X, X^\dagger)$. It is generated under extensions by intersection motives, as detailed in \cite[Lemma~2.15]{CassvdHScholbach:MotivicSatake}.

The same notions exist for categories of (stratified) reduced Tate motives, except that the Beilinson--Soulé condition on $S$ is not needed.

\subsubsection{Equivariant motives}
In order to consider equivariant motives on affine flag varieties, we recall the general formalism of equivariant motives, following \cite{RicharzScholbach:Intersection, CassvdHScholbach:MotivicSatake}.
Namely, via Kan extensions we can extend \(\DM\) to a functor
\[\DM\colon (\PreStk_S)^{\mathrm{op}} \to \PrLSt,\]
cf.~\cite[§2.1.5]{CassvdHScholbach:MotivicSatake}.
Here \(\PreStk_S\) is the \ii-category of prestacks, i.e., the (\(\infty\)-categorical) free cocompletion of the category of affine (but not necessarily finite type) \(S\)-schemes, cf.~\cite[§5.1]{Lurie:Higher}.

For a group prestack \(G\) acting on an ind-scheme \(X\), the category of equivariant motives on \(X\) is by definition \(\DM(G\backslash X)\), where \(G\backslash X\) is the prestack quotient, i.e., the colimit of the cosimplicial diagram 
$
\dots
\mathrel{
	\substack{
		\xrightarrow{} \\[-0.4ex]
		\xrightarrow{} \\[-0.4ex]
		\xrightarrow{}
	}
}
G \x X
\mathrel{
	\substack{
		\xrightarrow{} \\[-0.4ex]
		\xrightarrow{}
	}
}
X
$
.

If \(X\) is equipped with a Whitney--Tate stratification, the category of \(G\)-equivariant stratified Tate motives on \(X\) is defined as
\[\DTM_G(X):= \DM(G\backslash X) \times_{\DM(X)} \DTM(X).\]
In other words, they are those equivariant motives which are stratified Tate after forgetting the equivariance.
Sometimes, we also write \(\DTM(G\backslash X)=\DTM_G(X)\), but we emphasize that this really depends on the presentation, rather than just the prestack \(G\backslash X\).
Under certain assumptions on the \(G\)-action on \(X\), cf.~\cite[Lemma 2.26]{CassvdHScholbach:MotivicSatake}, there is moreover a (necessarily unique) t-structure on \(\DTM_G(X)\) such that the forgetful functor \(u^! : \DTM_G(X)\to \DTM(X)\) is t-exact.
(As in \refsect{mixed Tate motives}, we need to assume \(S\) satisfies the Beilinson-Soulé vanishing conjecture, if we consider non-reduced motives.)
We denote its heart by \(\MTM_{G}(X)=\MTM(G\backslash X)\).

\subsubsection{Grothendieck groups}\label{sect--Grothendieck}

Eyeing applications to Iwahori--Hecke algebras in \refsect{generic IH algebra}, we provide some general computations of Grothendieck groups of categories of Tate motives, including the equivariant and stratified  case.
To this end, we use the slice filtration \cite{HuberKahn:Slice}: if $\DTM(X)^{(\ge n)}$ denotes the full stable presentable subcategory generated by $\Z_X(m)$, $m \ge n$, the inclusion into $\DTM(X)$ admits a right adjoint $\nu^{\ge n}$.
The kernel of this functor is denoted by $\DTM(X)^{(\le n-1)}$.
We write $\DTM(X)^{(=0)} := \DTM(X)^{(\le 0)} \cap \DTM(X)^{(\ge 0)}$ for the ``heart'' of the slice structure.

We begin with unstratified (anti-effective) Tate motives, where the superscript $\comp$ stands for compact objects.
\lemm
\thlabel{K0 DTM}
Let $X / S$ be a smooth connected scheme.
Then the following maps given by $\mathbf q \mapsto \Z(-1)$ are ring isomorphisms: 
\begin{align*}
\Z[\mathbf q] \r & K_0(\DTM(X)^{\antic}), \\
\Z[\mathbf q^{\pm 1}] \r & K_0(\DTM(X)^{\comp}).
\end{align*} 
\xlemm

\pf
According to \cite[Theorem~3]{Spitzweck:Mixed}, the group $\Hom_{\DTM(X)}(\Z, \Z(r)[n])$ vanishes for $r < 0$ and $n \in \Z$ and also for $r = 0$, $n \ne 0$.
For $r = n = 0$, it is isomorphic to $\Z$, since $X$ is connected.
Since $\DTM(X)$ is compactly generated by $\Z(m)$, and 
$$\nu^{\ge n} \Z(m) = \left \{ \begin{tabular}{cc} 0 & $m < n$ \\ $\Z(m)$ & $m \ge n,$ \end{tabular} \right .$$
we obtain that $\nu^{\ge n}$ preserves compact objects.
In addition, the cofiber of the natural map $\nu^{\ge n} M \r M$ lies in the stable presentable subcategory generated by $\Z(s)$, $s < n$, so that $\DTM(X)^{(\le n-1)}$ is the stable presentable subcategory generated by $\Z(s)$, $s < n$.
Finally, any compact object $F \in \DTM(X)^\comp$ lies in some $\DTM(X)^{(\ge -n)}$, and is such that $\nu^{\ge n} F = 0$ for $n \gg 0$.
Thus, $F$ is obtained in finitely many steps, using cofiber sequences, from Tate twists of compact objects in $\langle \Z \rangle \subset \DTM(X)$.
The subcategory spanned by such compact objects is equivalent to $\Perf_\Z$, cf.~\cite[Proposition~4.5]{HuberKahn:Slice}, whose $K_0$-group is isomorphic to $\Z$.
\xpf

\thref{K0 DTM} and the localization fiber sequence $j_! j^* \r \id \r i_! i^*$ (for complementary open and closed embeddings $j$ and $i$, respectively) imply the following result.

\coro
\thlabel{coro_K0 DTM}
If $X^\dagger = \bigsqcup_{w \in W} X_w \r X$ is an (anti-effective) Whitney--Tate stratified ind-scheme whose strata $X_w$ are smooth over $S$ and connected, then 
\begin{align*}
K_0(\DTM(X, X^\dagger)^{\antic}) & = \bigoplus_{w \in W} \Z[\mathbf q], \\
K_0(\DTM(X, X^\dagger)^\comp) & = \bigoplus_{w \in W} \Z[\mathbf q^{\pm 1}].
\end{align*}
\xcoro

Recall the following notion from \cite[Definition 5.52]{CassvdHScholbach:MotivicSatake}: if \(H\) is a smooth group scheme acting on an ind-scheme \(X\), we say that a motive $\calF \in \DM_H(X)$ is \emph{locally compact} if its underlying motive $u^! \calF \in \DM(X)$ is compact. We denote the full subcategory of locally compact equivariant Tate motives by $\DTM_H(X)^{\locc}$.

\lemm
\thlabel{K0 forget group action}
Suppose $H$ is a smooth, fiberwise connected group scheme over $S$.
If $H$ acts on a smooth $S$-scheme $X$, the forgetful functor induces an isomorphism for the Grothendieck groups of equivariant vs.~non-equivariant (unstratified) Tate motives:
$$K_0(\DTM_H(X)^{(\anti), \locc}) \stackrel \cong \r K_0(\DTM(X)^{\xantic}).$$
\xlemm

\pf
The smoothness of $H$ allows us to invoke relative purity (i.e., the isomorphism $f^! = f^*(d)[2d]$, for a smooth map of relative dimension $d$). We obtain an equivalence
\begin{equation}
\DTM_H(X) = \lim (\DTM(X) \stackrel[p^*]{a^*} \rightrightarrows \DTM(X \x_S H) \substack{\rightarrow\\[-.9em] \rightarrow \\[-.9em] \rightarrow} \dots).\label{DTM*}
\end{equation}
It therefore makes sense to consider anti-effective equivariant Tate motives to begin with.
The functors $a^*$ and $p^*$ are ``slice-exact'', i.e., preserve the subcategories $\DTM(-)^{(\ge n)}$ and $\DTM(-)^{(\le n-1)}$ and thus also the heart $\DTM(-)^{(=0)}$. Both $a$ and $p$ have a section, so that the functors are conservative; and therefore also detect the property of being in either subcategory.
Therefore, \eqref{DTM*} holds mutatis mutandis for these three subcategories.

By \StPd{0385}{0378}, $\pi_0(X) = \pi_0(X \x_S H^{\x n})$ for each $n \ge 0$.
This gives equivalences $\DTM(X \x H^n)^{\comp, (=0)} = \DTM(X)^{\comp, (=0)}$. Passing to the limit gives an equivalence $\DTM_H(X)^{\locc, (=0)} = \DTM(X)^{\comp, (=0)}$. This implies our claim by functoriality of the slice filtration.
\xpf

\coro\thlabel{Grothendieck-stratified-equivariant}
Suppose an ind-scheme $X$ carries  an (anti-effective) Whitney--Tate stratification $X^\dagger = \bigsqcup_{w \in W} X_w$ whose strata $X_w$ are smooth.
Let a pro-algebraic group $H$ act on $X$ in such a way that each stratum $X_w$ is preserved by the $H$-action, and such that the $H$-action induces an isomorphism $X_w = (H/H_w)_\Zar$, where $H_w \subset H$ is an extension of a split pro-unipotent group by a connected split reductive group.
Then the forgetful functors again induce an isomorphism on Grothendieck groups:
$$K_0(\DTM_H(X, X^\dagger)^{(\anti), \locc}) \stackrel \cong \r K_0(\DTM(X, X^\dagger)^{\xantic}).$$
\xcoro

\pf
Localization reduces the case of $X/H$ to the case $X_w / H$.
By \cite[Proposition~3.1.23]{RicharzScholbach:Intersection}, \cite{RicharzScholbach:IntersectionCorrigendum} there is an equivalence $\DTM_H(X_w) = \DTM_{H_w}(S)$, so that we can apply \thref{K0 forget group action}.
\xpf

\rema\thlabel{Rationalization and K0}
In any of the situations above, the rationalization functor induces an isomorphism on Grothendieck rings.
Indeed, if \(X/S\) is smooth connected, we still have an isomorphism \(\Z[\qq]\cong K_0(\DTM(X,\Q)^{\anti,\comp})\) as in \thref{K0 DTM}, and this carries over to equivariant and/or stratified Tate motives.
\xrema

\subsubsection{Tate motives on $\Gm$}

Our approach to (unipotent) nearby cycles rests on applying some higher categorical algebra to various categories of motivic sheaves.
This is facilitated by a module-theoretic interpretation of the category of Tate motives on $\Gm$ that we study in this section. It is worth noting that the multiplication on $\Gm$ plays no role at this point.

Let $p : \GmX S \r S$ be the structural map.
A basic, but key object, is the ``cohomology'' of $\Gm$,
$$\Lambda := p_* p^* \Z.\eqlabel{PI}$$
It is a commutative algebra object in $\DM(S)$ whose underlying object in $\DM(S)$ is $\Z \oplus \Z(-1)[-1]$. 
In particular, it lies in the subcategory $\DTM(S)^\anti$.
The counit map 
$$p^* \Lambda = p^* p_* p^* \Z \r p^* \Z = \Z_{\Gm}$$
is a map of commutative algebras in $\DTM(\Gm)^\anti$.
The point $S \stackrel 1 \r \GmX S$ equips $\Lambda$ with an augmentation, $\aug: \Lambda \r \Z$.
The following result appears independently in \cite[Proposition~2.56]{BindaGallauerVezzani:Motivic}.

\lemm
\thlabel{DTM.Gm}
The restriction of $p_*$ to $\DTM(\Gm)$ is monadic and therefore gives rise to an equivalence
$$\DTM(\Gm) \cong \Mod_\Lambda(\DTM(S))\eqlabel{DTM Gm}$$
under which $p_*$ corresponds to the functor forgetting the $\Lambda$-action.
The same statement also holds for reduced Tate motives, for anti-effective Tate motives, and for reduced anti-effective Tate motives.
\xlemm

\pf
The functor $p_{*}$, restricted to $\DTMrx(\Gm) \r \DTMrx(S)$, is conservative since $p^* \DTMrx(S)$ generates $\DTMrx(\Gm)$.
By general theory, $p_{*}$ preserves colimits so that $p_{*}$ is monadic.
For anti-effective motives, the same proof works, since $p^*$ and $p_*$ preserve anti-effective motives.
\xpf

Recall that for any object $M$ in a stable presentable symmetric monoidal \ii-category $C$, there is the split square zero extension, which is an object in $\CAlg(C)$ whose underlying object in $C$ is isomorphic to $1 \oplus M$ \cite[§7.3.4]{Lurie:HA}.

\lemm
\thlabel{Lambda square zero}
There is an isomorphism in $\CAlg(\DTM(S, \Q))$
$$\Lambda \t \Q = \Q \oplus \Q(-1)[-1].$$
Similarly, if we let $\Lambda_\red = p_* p^* \Z \in \DTMr(S) (=\D(\grAb))$ be the dual of the reduced motive of $\Gm$, there is an isomorphism in $\CAlg(\DTMr(S))$:
$$\Lambda_\red = \Z \oplus \Z(-1)[-1].\eqlabel{Lambda red}$$
(Thus, under the equivalence $\DTMr(S) = \D(\grAb)$, $\Lambda_\red$ corresponds to the object $\OldLambda$ considered in \eqref{notation OldLambda}.)
\xlemm

\pf
The statement with rational coefficients holds by \cite[Proposition~2.63]{BindaGallauerVezzani:Motivic} which constructs isomorphisms between the free symmetric algebra $\Sym(\Q(-1)[-1])$ and both terms. We also refer to \cite[Theorem~7.1.1]{AnconaEnrightWardHuber:Motive} for the identification of the Hopf algebra structure on $\M(\Gm)$ (again with rational coefficients).

For the statement for reduced motives (with integral coefficients), note that both objects lie in the full subcategory spanned by $\Z(-n)[-n]$ for $n \ge 0$, which is an \emph{ordinary} symmetric monoidal category.
Let $I = \Z(-1)[-1]$  be the augmentation ideal of $\Lambda_\red$.
The space of maps $\Map_{\DTMr(S)}(I^{\t n}, I^{\t m}) = \Map_{\D(\grAb)}(\Z(-n)[-n], \Z(-m)[-m])$ is trivial for any $n \ne m$ and is isomorphic to $\Z$ for $n = m$.
Therefore, the multiplication maps are all null-homotopic.
The switching maps (for the symmetric monoidal structure) $I \t I \stackrel \cong \r I \t I$ are the identity both for $\Lambda_\red$ (since this is true for the unreduced object $\Lambda$), and also for the square zero extension.
\xpf

\rema
\thlabel{Gm integral}
The analogous statement for $\DM$ and integral or $\Fp$-coefficients fails.
In this sense the passage to rational coefficients in the above formality result \cite[Proposition~2.63]{BindaGallauerVezzani:Motivic} is optimal.
To see this note that the multiplication $\Lambda \t \Lambda \r \Lambda$ is dual to the comultiplication on $\M(\Gm) = \Lambda^\dual$, given by the diagonal $\M(\Gm) \stackrel{\Delta_{\Gm}} \r \M(\Gm \x \Gm) = \M(\Gm) \t \M(\Gm)$. The restriction of this map to the direct summands $\Z(1)[1] \r \Z(2)[2]$ is nonzero since applying $\Hom_{\DM(\Spec \Q)}(\Z, -)$ gives the non-zero map $K_1^M(\Q) = \Q^\x \r K_2^M(\Q)$, $[a] \mapsto [a] \cdot [a]$. (Instead, the Steinberg relation implies $[a] \cdot [-a] = 0 \in K_2^M(\Q)$.) Inverting 2 in the coefficients, we do have $\Map_{\DM(\Spec \Z)}(\Z[\frac 12](1)[1], \Z[\frac 12](2)[2]) = 0$, but higher order multiplications still form an obstruction to $\M(\Gm) \t \Z[\frac 12]$ being the square zero coalgebra (or $\Lambda [\frac 12]$ the square zero algebra).
In fact, the Betti realization functor $\rho_\Betti$ is symmetric monoidal and therefore preserves square zero extensions.
Now $\rho_\Betti(\Lambda \t \Fp) = C^*(\Gm^{\an}, \Fp) = C^*(S^1, \Fp) \in \CAlg(\D(\Mod_{\Fp}))$. We claim it is not isomorphic, as a commutative algebra object, to the square zero extension $\Fp \oplus \Fp(-1)[-1]$.
Indeed, any such isomorphism would be compatible with the Steenrod operations, but $P^0 = \id$ on $C^*(S^1, \Fp)$ \cite[Proposition~8.1]{May:General}, while $P^0 = 0$ on the (strictly commutative) cdga $\Fp \oplus \Fp[-1]$.
Thus, there is no isomorphism, in $\CAlg(\DM(\Spec \Z, \Fp))$, between $\Lambda \t \Fp$ and the split square zero extension $\Fp \oplus \Fp(-1)[-1]$.
\xrema

\coro
\thlabel{DTM Gm Q}
There are equivalences (of commutative algebra objects in $\Mod_{\D(\grAb)}(\PrLSt)$)
\begin{align*}
\DTM(\Gm, \Q) & = \DTMr(\Gm) \t_{\DTMr(S)} \DTM(S, \Q). \\
\DTMr(\Gm) & = \Nilp.
\end{align*}
\xcoro

\pf
The statement about $\DTMr(\Gm)$ follows from the construction of $\Nilp$ in \thref{Nilp} and from \refeq{Lambda red}.
To deduce the first statement, we consider the symmetric monoidal functor $i : \DTMr(S) = \D(\grAb) \r \DTM(S, \Q)$ given by $\Z(k) \mapsto \Q(k)$.
It maps the split square zero extension $\Lambda_\red = \Z \oplus \Z(-1)[-1]$ to the split square zero extension $\Q \oplus \Q(-1)[-1]$, which by \thref{Lambda square zero} is isomorphic to $\Lambda \t \Q$.
Thus using \thref{DTM.Gm} and \thref{generalities.tensor}, we get a chain of equivalences of \ii-categories 
\begin{align*}
\DTM(\Gm, \Q) & = \Mod_{\Lambda \t \Q}(\DTM(S, \Q)) \\
& = \Mod_{i (\Lambda_\red)} (\DTM(S, \Q)) \\
& = \Mod_{\Lambda_\red}(\DTMr(S)) \t_{\DTMr(S)} \DTM(S, \Q).
\end{align*}
\xpf

\rema
\thlabel{1^*}
In other words, Tate motives $M$ on $\Gm$ with rational coefficients (and, likewise, reduced motives with integral coefficients) are equivalent to the datum of their restriction $1^* M$ together with a locally nilpotent map $(1^* M) (1) \r 1^* M$.
Concerning (unreduced) motives with integral coefficients, we still have a natural map $\Z(1) \stackrel \alpha \r \IHom_{\Lambda}(\Z, \Z)$, and thus a map of $E_1$-algebras $\Sigma \r \IHom_{\Lambda}(\Z, \Z)$.
Note that $1^*M = \Z \t_\Lambda M$ is naturally a module over $\IHom_\Lambda(\Z, \Z)$ and therefore by restriction a $\Sigma$-module.
In other words, $1^*M$ still comes naturally with the map $\alpha$, but we have little control over it (cf.~\thref{lemm-monodromy-triangle}) and therefore do not use it in this paper.
\xrema

We end this section by introducing a category that will eventually allow to keep track of some motivic sheaf $M \in \DM(X)$ together with a map $M \t \Q(1) \r M \t \Q$.

\defi
\thlabel{NilpQ}
Let $\NilpQ := \D(\grModZ) \x_{\D(\grModZ) \t \Q} (\Nilp \t \Q),$
where we write $C \t \Q := C \t_{\D(\Ab)} \D(\Mod_\Q)$ for the rationalization of any $\Z$-linear stable \ii-category, and the pullback is formed using the forgetful functor $\Nilp \subset \Mod_{\Sigma}(\D(\grModZ)) \r \D(\grModZ)$.

Both $\Nilp$ and $\NilpQ$ are modules over $\D(\grAb)$. For some $C \in \Mod_{\D(\grAb)}(\PrLst)$ (such as some category of motives, cf.~\refeq{DM Prgr}), we will abbreviate $\NilpQ C := \NilpQ \t_{\D(\grAb)} C$ and $\Nilp C := \Nilp \t_{\D(\grAb)} C$.
\xdefi

\rema
\thlabel{NilpQ explained}
Thus, $\NilpQ$ is a full subcategory of $\D(\grModZ) \x_{\D(\grModQ)} \Mod_{\Sigma \t \Q}(\D(\grModQ))$;
objects in the latter category can be described as $(A \in \D(\grAb), \phi : A \t \Q(1) \r A \t \Q)$.
The category $\NilpQ$ is compactly generated by $(\Z(k), \Q(k+1) \stackrel 0 \r \Q(k))$, $k \in \Z$.
Thus, if $C$ is compactly generated, then $\NilpQ C = C \x_{C \t \Q} (C \t \Nilp \t \Q)$ and this category is compactly generated as well. If an object $(c \in C, \phi : c \t \Q(1) \r c \t \Q)$ in this category is compact, then the map $\phi$ is nilpotent.
Similar statements hold for $\Nilp C$.
The natural functor $\Nilp \r \NilpQ$ induces an equivalence after rationalization, i.e., $\Nilp \t \Q \stackrel \cong \r \NilpQ \t \Q$.
\xrema

In the setup of the unipotent nearby cycles functor we will use the following functor.

\defi
\thlabel{DTM Gm Nilp}
The symmetric monoidal functors
\begin{align*}
1^* : \DTM(\Gm) \r & \DTM(S), \\
\DTM(\Gm) \r & \DTM(\Gm, \Q) \stackrel{\text{\ref{DTM Gm Q}}} = \Nilp \t_{\D(\grAb)} \DTM(S, \Q)
\end{align*}
induce a symmetric monoidal functor
$$\DTM(\Gm) \r \NilpQ \DTM(S).$$
(After passing to rationalizations, this functor recovers the equivalence in \thref{DTM Gm Q}.)
\xdefi

\subsubsection{Convolution structures} 
\label{sect--convolution structures}
For a coherent construction of Wakimoto functors (\thref{Definition Wakimoto}), we refine the construction of the convolution monoidal structure on categories such as $\DM(\mathcal{I} \setminus LG / \mathcal{I})$ from \cite[§3]{RicharzScholbach:Motivic} to the level of \ii-categories (as opposed to their homotopy categories). Related statements appear, e.g., in \cite[Ch.~5, §5]{GaitsgoryRozenblyum:StudyI} and \cite[§4.2]{ALWY:Gaitsgory}, so we will keep this brief.

First, the restriction of $\DM^!$ to the category $\Sch_S^{\placid}$ of placid schemes (and placid morphisms) is lax symmetric monoidal \cite[§A.2]{RicharzScholbach:Motivic}.
Recall that this encodes the existence of $\boxtimes$ and its compatibility with !-pullback along smooth (and then placid) maps.
Second, the Kan extension to placid ind-schemes $\IndSch_S^\placid$ \cite[§C.1]{Gaitsgory:Local}, is therefore again lax symmetric monoidal.
Third, the Kan extension to the free completion under sifted colimits is lax symmetric monoidal, see, e.g., the discussion in \cite[§3.2]{Robalo:Theorie}. Denote this by $\IndSch_S^\placid(\mathrm{sift})$.
Finally, consider the category $\Corr(\IndSch_S^\placid(\mathrm{sift}))$ of correspondences, where as ``horizontal'' maps we only allow maps $f$ such that for any pullback $f'$ of $f$, the functor $f'^!$ admits a left adjoint, denoted $f'_!$, which is moreover adjointable (i.e., commutes with all !-pullbacks).
Then $\DM^!$ extends naturally to a lax symmetric monoidal functor, see, e.g., \cite[Proposition~A.5.10]{Mann:p-adic}.

These facts can be applied as follows:
let $L$ be a group object in the category of placid ind-schemes over $S$, and fix a subgroup object $L^+ \subset L$ such that the quotient $(L/L^+)_\Nis$ is representable by an ind-proper ind-scheme.
Then the convolution groupoid $X := L^+ \setminus L / L^+$ is an object in the sifted completion.
The multiplication map $m : X \x_{S / L^+} X \r X$ is a pullback of the projection $L / L^+ \r S$. By assumption, up to sheafification, it is an ind-proper map, so that the same is true for any pullback of it. For ind-proper maps $f$, the left adjoint $f_!$ of $f^!$ exists and commutes with all !-pullbacks \cite[Proposition~2.3.3]{RicharzScholbach:Intersection}.
Thus, $X$ is a (non-commutative) algebra object in $\Corr(\IndSch_S^\placid(\mathrm{sift}))$ \cite[Part I, Ch.~9, Corollary~4.4.5]{GaitsgoryRozenblyum:StudyI}, so that the (symmetric) monoidal functor $\DM^!$ maps it to an algebra object $\DM(X)$ in $\PrLSt$ or $\Mod_{\DM(S)}(\PrLSt)$.

\section{Unipotent nearby cycles}
\label{sect--unipotent}

The following definition of unipotent nearby cycles is a motivic adaptation of the approach in \cite[§3]{Campbell:Nearby}, see also \cite[§2.1]{Chen:Nearby}.
Throughout, we consider the following geometric situation, where $X$ is a scheme and the two top squares are cartesian:
$$\xymatrix{
X_\eta \ar[r]^j \ar[d]^{f_\eta} & X \ar[d]^f & X_s \ar[d]_{f_s} \ar[l]_i \\
\GmX S \ar[r] \ar[dr]_{p} & \A^1_S \ar[d]^{q} & S \ar[dl]^\id \ar[l]_0 \\ & S.}\eqlabel{diagram}$$
We write $\Gamma_{f_\eta} : X_\eta \r X_\eta \x_S \Gm$ for the graph of $f_\eta$ and consider the functor
$$\DTM(\Gm) \t_{\DTM(S)} \DM(X_\eta) \stackrel \boxtimes \r \DM(X_\eta \x \Gm) \stackrel{\Gamma_{f_\eta}^*} \r \DM(X_\eta).$$
This composite, which is a map in $\PrLSt$, has a right adjoint, which we denote by $\Gamma_{f_\eta?}$.

\defi
\thlabel{unipotent nearby cycles}
The \emph{unipotent nearby cycles functor} is the composite
$$\Upsilon_X : \DM(X_\eta) \stackrel{\Gamma_{f_\eta ?}} \r \DTM(\Gm) \t_{\DTM(S)} \DM (X_\eta) \stackrel {1^* \t i^* j_*} \r \DM(X_s).$$
\xdefi

We now make this definition more explicit, using the equivalence $\DTM(\Gm) \t_{\DTM(S)} \DM(X_\eta) = \Mod_\Lambda (\DM(X_\eta))$, the category of $\Lambda$-modules (where $\Lambda$ is pulled back from $S$ to $X_\eta$ along $p f_\eta$) in $\DM(X_\eta)$.
This equivalence results from \thref{generalities.tensor} and \thref{DTM.Gm}.
Here and throughout, all tensor products of categories of (Tate) motives are carried out in $\Mod_{\DTM(S)}(\PrLSt)$. I.e., unless otherwise mentioned, $C \t D$ is a shorthand for $C \t_{\DTM(S)} D$.

Every object $\calF \in \DM(X_\eta)$ becomes naturally a $\Lambda$-module: its $\Lambda$-action is just the map
$$f_\eta^* p^* \Lambda \t \calF = f_\eta^* p^* p_* \Z \t \calF \stackrel{adj.} \r f_\eta^* \Z \t \calF = \calF.$$
The functor $\Gamma_{f_\eta?}$ agrees with the resulting natural functor $\DM(X_\eta) \r \Mod_\Lambda(\DM(X_\eta))$, since their left adjoints, given by $\Gamma_{f_\eta}^*$ and $- \t_\Lambda \Z$ agree. Here the tensor product is formed using the counit map $\Lambda (:= p^* \Lambda)  = p^* p_* \Z_{\Gm} \r \Z_{\Gm}$.
The standard functor $i^* j_*$ extends naturally to a functor $\Mod_\Lambda(\DM(X_\eta)) \r \Mod_\Lambda(\DM(X_s))$ (in a way that is compatible with forgetting the $\Lambda$-action).
This can be identified with $\id_{\DTM(\Gm)} \t i^* j_*$. In other words, that functor is just keeping track of the extra $\Lambda$-action, but otherwise is the standard $i^* j_*$.
Finally, the functor $1^* \t \id_{\DM(X_s)}$ identifies with $- \t_\Lambda \Z : \Mod_\Lambda(\DM(X_s)) \r \DM(X_s)$, where the tensor product is formed using the augmentation map, which geometrically corresponds to pullback along $1 \in \Gm$.

The fiber sequence
$$I := \Z(-1)[-1] \stackrel \iota \r \Lambda \stackrel \aug \r \Z\eqlabel{cofiber I}$$
gives rise to fiber sequences
$$I^{\t n+1} \stackrel{\iota \t \id} \r \Lambda \t I^{\t n} = \Lambda(-n)[-n] \stackrel{\aug \t \id} \r I^{\t n}.$$
This gives a ``resolution'' of $\Z$ by free $\Lambda$-modules:
$$\colim \left ( \dots \Lambda(-2)[-2] \r \Lambda(-1)[-1] \r \Lambda \right ) \stackrel \cong \r \Z.\eqlabel{resolution.Z.P}$$
Using this resolution, we obtain the following explicit formula for $\Upsilon$:
$$\Upsilon(M) = \colim \left ( \dots \r i^* j_* M (-2)[-2] \r i^* j_* M(-1)[-1] \r i^* j_* M \right).\eqlabel{Upsilon explicit}$$
By the construction of \refeq{resolution.Z.P}, the maps in this diagram are obtained from the map $\Z(-1)[-1] \t M \r \Lambda \t M \r M$ induced by the $\Lambda$-action on $M$. 

\subsection{Functoriality, Künneth formula and monodromy}

Our goal in this section is to enhance this definition in the following ways: 
we establish the compatibility of $\Upsilon$ with proper pushforward, smooth pullback and its lax compatibility with exterior products. We do this in a highly structured way, having as a payoff an effortless approach to, say, equivariant nearby cycles, cf.~\thref{Upsilon equivariant}.
We also enhance $\Upsilon$ by the additional datum of a monodromy map, which works integrally for reduced motives, but requires rational coefficients for (unreduced) motives.

We will cast $\Upsilon$ as a natural transformation between two functors, the domain of which is the following lax symmetric monoidal functor
$$\DM_\eta : \Corr_{\A^1} \stackrel{j^*} \r \Corr_{\Gm} \stackrel {\DM^*_!} \r \Mod_{\DTM(\Gm)} \PrLSt \stackrel U \r \Mod_{\DTM(S)}(\PrLSt).$$
It sends $X / \A^1$ to $\DM(X_\eta)$. Its lax symmetric monoidality expresses the maps
$$\boxtimes_{\Gm} : \DM(X_{1\eta}) \t \DM(X_{2\eta}) \r \DM(X_{1\eta} \x_{\Gm} X_{2\eta})$$
and their functoriality with respect to !-pushforward and *-pullback.

As for the codomain of $\Upsilon$, recall from \refsect{Koszul duality} and \thref{NilpQ} the category $\NilpQ$. It is a commutative algebra object in $\PrLgr := \Mod_{\D(\grAb)}(\PrLSt)$. Objects in the category $\NilpQ \DM(Y) := \NilpQ \t_{\D(\grAb)} \DM(Y)$ are pairs $(M \in \DM(Y), \phi : M(1) \t \Q \r M \t \Q)$, where $\phi$ is a locally nilpotent map of the rationalized motives, as is explained in \thref{remarks Nilp} and \thref{NilpQ explained}.
The codomain of $\Upsilon$ is now the lax symmetric monoidal functor
$$\NilpQ \DM_s : \Corr_{\A^1} \stackrel {i^*} \r \Corr \stackrel {\DM^*_!} \r \PrLgr \stackrel {\NilpQ \t_{\D(\grAb)} -} \r \PrLgr.$$
It maps $X / \A^1$ to $\NilpQ \DM(X_s)$.
The lax symmetric monoidality of this records the functors
$$\boxtimes: \NilpQ \DM(X_{1s}) \t \NilpQ \DM(X_{2s}) \r \NilpQ \DM(X_{1s} \x X_{2s}),$$
sending $(M_1, \phi_1), (M_2, \phi_2)$ to $(M_1 \boxtimes M_2, \phi_1 \boxtimes \id + \id \boxtimes \phi_2)$.

We consider the wide subcategory $\Corr_{\A^1}^{\proper, \sm} \subset \Corr_{\A^1}$ in which we require the correspondences in \refeq{correspondences} to consist of a smooth map $f$ and a proper map $g$. Thus, the restriction $\DM^*_!|_{\Corr_{\A^1}^{\proper, \sm}}$ only records *-pullback along smooth maps and !-pushforward along proper maps.

The following theorem will be proven in \refsect{proof of Upsilon theorem}.
\theo
\thlabel{Upsilon.all.in.one}
There is a natural transformation of lax symmetric monoidal functors $\Corr_{\A^1}^{\proper, \sm} \r \PrLgr$:
$$\Upsilon : \DM_\eta \r \NilpQ \DM_s,$$
whose evaluation at a scheme $X / \A^1$, composed with the forgetful functor $\NilpQ \DM(X_s) \r \DM(X_s)$, is the functor $\Upsilon$ in \thref{unipotent nearby cycles}.
\xtheo

\rema
\thlabel{remarks Upsilon}
We explain the content of \thref{Upsilon.all.in.one} by probing $\Upsilon$ in various ways.
\begin{enumerate}
  \item 
  Evaluating $\Upsilon$ on some scheme $X / \A^1$, we have a colimit-preserving functor
  $$\Upsilon_X : \DM(X_\eta) \r \NilpQ \DM(X_s).$$
  That is, for $M \in \DM(X_\eta)$, $\Upsilon(M) \in \DM(X_s)$ is naturally equipped with a locally nilpotent morphism $\Upsilon(M)(1) \t \Q \r \Upsilon(M) \t \Q$, that we refer to as the \emph{monodromy}.
  In \thref{lemm-monodromy-triangle} we will compute its cofiber to be $i^* j_* M(1)[1] \t \Q$.

  \item Evaluating $\Upsilon$ at a correspondence $Y \stackrel g \gets Z \stackrel f \r X$ (of schemes over $\A^1$), with $f$ smooth and $g$ proper, we have a natural isomorphism
  $$g_{s!} f_s^* \Upsilon_X \stackrel \cong \r \Upsilon_Y g_{\eta!} f_\eta^*\eqlabel{Upsilon lax functoriality}$$
  or, more diagramatically, a commutative diagram
  $$\xymatrix{
  \DM(X_\eta) \ar[r]^-{\Upsilon_X} \ar[d]^{g_{\eta!} f_\eta^*} & \NilpQ \DM(X_s) \ar[d]^{g_{s!} f_s^*} \\
  \DM(Y_\eta) \ar[r]^-{\Upsilon_Y} & \NilpQ \DM(Y_s).
  }$$
  Here $g_\eta$ is the pullback of $g$ along $\Gm \subset \A^1$ etc.
  In other words, the formation of (unipotent) nearby cycles, including the monodromy, is compatible with proper pushforward and smooth pullback.  

  For $g$ non-proper or $f$ non-smooth, we still have a map in \refeq{Upsilon lax functoriality}, but it need not be an isomorphism, cf.~\thref{Upsilon lax functor}.

  \item
  For a pair $X_1, X_2$ of schemes over $\A^1$, there is a morphism
  $$(X_1, X_2) \r X_1 \x_{\A^1} X_2$$
  in (the symmetric monoidal \ii-category) $\Corr_{\A^1}$.
  The lax symmetric monoidality of $\Upsilon$ then gives a natural transformation in the diagram
  $$\xymatrix{
  \DM(X_{1\eta}) \t \DM(X_{2\eta}) \ar[rr]^(.45){\Upsilon_{X_1} \t \Upsilon_{X_2}} \ar[d]^{\boxtimes_{\Gm}} & & \NilpQ \DM(X_{1s}) \t \NilpQ \DM(X_{2s}) \ar@{=>}[dll] \ar[d]^{\boxtimes}\\
  \DM(X_{1\eta} \x_{\Gm} X_{2\eta}) \ar[rr]_{\Upsilon_{X_1 \x_{\A^1} X_2}} & & \NilpQ \DM(X_{1s} \x X_{2s}).}$$
  Thus, for $M_i \in \DM(X_{i\eta})$, there is a natural map
  $$\kappa_{(X_1, X_2)} : \Upsilon_{X_1} (M_1) \boxtimes \Upsilon_{X_2} (M_2) \r \Upsilon_{X_1 \x_{\A^1} X_2} (M_1 \boxtimes_{\Gm} M_2).$$
  We call this the \emph{Künneth map} for unipotent nearby cycles.
  (In general it is not an isomorphism, but see \thref{theo-WT-BD}\refit{theo-WT-BD3} for a case where it is.)
  This map is functorial and compatible with arbitrary colimits in both $M_i$. It is compatible with the monodromy maps in the sense that it is a map of (locally nilpotent) $\Sigma$-modules, where $\Sigma$ acts on the source of $\kappa$ via the comultiplication of the bialgebra $\Sigma$, cf.~\eqref{comultiplication Sigma}.
Further, $\kappa$ is compatible with functoriality in $\Corr_{\A^1}^{\proper, \sm}$.
  For example, for a pair of smooth maps $f_k : Z_k \r X_k$, there is a commutative square in $\Corr_{\A^1}$
  $$\xymatrix{
  (Z_1, Z_2) \ar[r] \ar[d]_{(f_1, f_2)} & Z_1 \x_{\A^1} Z_2 \ar[d]^{f_1 \x_{\A^1} f_2} \\
  (X_1, X_2) \ar[r] & X_1 \x_{\A^1} X_2.
  }$$ 
  Then
  $$(f_{1s} \x f_{2s})^* \kappa_{(X_1, X_2)} \left (M_1, M_2 \right ) = \kappa_{(Z_1, Z_2)} \left ((f_{1\eta})^* M_1, (f_{2\eta})^* M_2 \right ).$$ 
  A similar formula holds for !-pushforwards along proper maps.

  \item
  \label{item--Upsilon reduced}
  For reduced motives, we have a similarly defined functor
  $$\Upsilon : \DM_{\red,\eta} \r \Nilp \DM_{\red, s},$$
  i.e., the monodromy map $\phi$ exists already integrally. Indeed, the point of using $\NilpQ$ above is the usage of \thref{DTM Gm Nilp} (which goes back to \thref{Lambda square zero}), which already holds integrally for reduced motives.
  All properties of $\Upsilon$ proved in this section (for regular motives) then hold verbatim for reduced motives.

\end{enumerate}
\xrema

\coro
\thlabel{Upsilon.tensor}
For fixed $X / \A^1$, $\Upsilon : \DM(X_\eta) \r \NilpQ \DM(X_s)$ is a lax symmetric monoidal functor.
\xcoro

\pf
The symmetric monoidal structure on individual categories $\DM(Y)$ arises by appending $\Delta^*$ (for $Y \stackrel \Delta \r Y \x Y$) to $\boxtimes$.
The symmetric monoidal structure is then given by the composite:
\begin{align*}
\Upsilon (-) \t \Upsilon (-) & = \Delta_{X_s}^* (\Upsilon (-) \boxtimes \Upsilon (-)) \\
& \stackrel * \r \Upsilon (\Delta_{X_\eta}^*  (- \boxtimes -)) & \\
& = \Upsilon (- \t -).
\end{align*}
The map labelled ``$*$'' exists by the lax compatibility of $\Upsilon$ with respect to (non-smooth) pullback, cf.~\thref{Upsilon lax functor}.
Alternatively, using Zariski descent to reduce to the case of separated $X$, one can use the compatibility of $\Upsilon$ with pushforward along the proper map $\Delta$.
\xpf

\subsection{Further remarks}

The above definition is related to various works: a precise comparison with Ayoub's unipotent nearby cycles appears below in \refsect{comparison Ayoub}.
The approach in \thref{unipotent nearby cycles} can be regarded as a motivic refinement of Campbell's approach to Beilinson's unipotent nearby cycles \cite{Campbell:Nearby}. It also bears some similarity with the one in \cite[Scholie~1.3.26]{Ayoub:Motifs}, but is simpler in that it avoids the usage of rigid analytic motives.
It also shares a kinship with the recent work of Binda--Gallauer--Vezzani \cite{BindaGallauerVezzani:Motivic}, who have independently carved out the close relation between unipotent nearby cycles (in the context of rigid-analytic motives) and Koszul duality.

As was highlighted to us by Vezzani, the existence of an integral nearby cycles functor (without monodromy), and the existence of a monodromy map after passing to rational coefficients matches the outcome of the approach in Ayoub's work \cite[§3.6]{Ayoub:Six2} as well as the situation for rigid analytic motives in the work of Binda--Gallauer--Vezzani \cite{BindaGallauerVezzani:Motivic}, with the salient point being \cite[Theorem~3.3.3]{AyoubGallauerVezzani:Six}. An integral monodromy map for \emph{reduced} motives is discussed in \thref{remarks Upsilon}\refit{Upsilon reduced}.

If we consider $\Upsilon^! := (1^* \t i^! j_!) \circ \Gamma_{f_\eta?}$ instead of $\Upsilon$, we get $\Upsilon^! = \Upsilon[1]$.

An interesting variant is to consider other categories in place of $\DTM(\Gm)$.
For example, one can consider the functor
$$\DM(X_\eta) \r \DATM(\Gm) \t_{\DTM(S)} \DM(X_\eta) \stackrel{1^* \t i^* j_*} \r \DM(X_s),$$
where $\DATM(\Gm)$ denotes the subcategory of $\DM(\Gm)$ spanned by $f_{n*} \Z(n)$, for $f_n : \Gm \r \Gm$ being the $n$-th power map.
Over a field of characteristic 0, does this functor agree with Ayoub's full nearby cycles functor?
Other interesting variants might be to allow $f_* \Z(n)$ for all finite étale $f$, or even replace $\DTM(\Gm)$ by the subcategory $\DM(\Gm)^\dualizable \subset \DM(\Gm)$ of dualizable objects.

\subsection{Proof of \thref{Upsilon.all.in.one}}
\label{sect--proof of Upsilon theorem}
The following diagram gives an overview of the functors and natural transformations appearing in the construction of $\Upsilon$:
$$\xymatrix{
\Corr^{\sm, \proper}_{\A^1} \ar[r]^{j^*} \ar[dr]^{i^*} &
\Corr^{\sm, \proper}_{\Gm} \ar[r]^(.4){\DM} &
\Mod_{\DTM(\Gm)}(\PrLSt) \ar[d]_U \ar[r]^F \ar@{=>}[dl]_{i^* j_*} &
\Mod_{\NilpQ \DTM(S)}(\PrLSt) \ar[d]_{U'} \\
&
\Corr^{\sm, \proper} \ar[r]^(.4){\DM} &
\Mod_{\DTM(S)}(\PrLSt) \ar@{=}[r] \ar@{=>}@/_1pc/[u]_{\alpha} \ar@/_2pc/[u]_L & 
\Mod_{\DTM(S)}(\PrLSt) \ar@/_2pc/[u]_{L'}.}$$
Here $i^*$ and $j^*$ at the left are the pullbacks along the respective inclusions; $\DM$ is the functor in \refeq{DM Mod}.
We have abbreviated $F := \NilpQ \DTM(S) \t_{\DTM(\Gm)} -$, where the tensor product is formed using the functor in \thref{DTM Gm Nilp}. The functors $U$ and $U'$ are the forgetful functors, whose left adjoints are given by $L := \DTM(\Gm) \t_{\DTM(S)} -$ and $L' := \NilpQ \DTM(S)\t_{\DTM(S)} -$.
In this diagram, $i^*$, $j^*$, $F$, $L$ and $L'$ are symmetric monoidal and all other functors are lax symmetric monoidal.
As before, we abbreviate $\DM_\eta := \DM \circ j^*$, $\DM_s = \DM \circ i^*$.

We now discuss the pertinent natural transformations between these functors, including their lax monoidality.
In the above diagram, the two arrows $\Rightarrow$ stand for natural transformations 
\begin{eqnarray*}
i^* j_* & : & U \DM_\eta \r \DM_s \text{ and } \\
\alpha & : & \id \r LU
\end{eqnarray*}
which are explained shortly.

First of all, there is a natural transformation of lax symmetric monoidal endofunctors on $\Mod_{\DTM(S)}(\PrLSt)$, $UL \r U'L'$. It arises from $L' = F L$, so if $G$ denotes the right adjoint of $F$ (which is again a forgetful functor), the unit map $\id \r GF$ gives rise to
$$UL \r U'L' = UGFL.$$
For some $C \in \Mod_{\DTM(S)}(\PrLSt)$ (such as $C = \DM(X_\eta)$), the evaluation on $C$ is the functor $\DTM(\Gm) \t_{\DTM(S)} C \r \NilpQ \DTM(S) \t_{\DTM(S)} C = \NilpQ \t_{\D(\grAb)} C$ that on the first factor sends $M \in \DTM(\Gm)$ to the pair consisting of $1^* M \in \DTM(S)$ and $\phi : M \t \Q(1) \r M \t \Q$, and on the second factor is $\id_C$.

\lemm
\thlabel{right.adjoint.action}
The adjunction $(L, U)$ is ambidextrous, and the unit map 
$$\alpha : \id  \r  LU$$
is a natural transformation of lax symmetric monoidal endofunctors of $\Mod_{\DTM(\Gm)}(\PrLSt)$.
The composition with $\DM_\eta$ thus gives a natural transformation
$$\alpha : \DM_\eta \r LU \DM_\eta$$
of lax symmetric monoidal functors $\Corr_{\A^1} \r \Mod_{\DTM(\Gm)}(\PrLSt)$.
\xlemm

\pf
The categories $\DTM(\Gm)$ and $\DTM(S)$ are compactly generated by the objects $\Z(n)$, and are therefore rigid.
By \thref{ModD.equivalence}, the right and left adjoint functor (denoted $L$) of $U: \Mod_{\DTM(\Gm)} \PrLSt \r \Mod_{\DTM(S)}(\PrLSt)$ agree up to equivalence. 
The unit map $\id \r LU$ is a natural map of lax symmetric monoidal functors, and therefore so is its composition with $\DM_\eta$.
\xpf

\rema The functor $\Gamma_{f_\eta?}$ has the following properties:
\begin{itemize}
  \item 
For any scheme $X / \A^1$, the evaluation $\alpha(X)$ is the right adjoint
$$\Gamma_{f_\eta?} : \DM(X_\eta) \r \DTM(\Gm) \t_{\DTM(S)} \DM(X_\eta)$$
of the functor $\DTM(\Gm) \t \DM(X_\eta) \stackrel \boxtimes \r \DM(\Gm \x_S X_\eta) \stackrel{\Gamma_{f_\eta}^*} \r \DM(X_\eta)$.
Indeed, this latter functor encodes the action of $\DTM(\Gm)$ on $\DM(X_\eta)$.
\item
We can disentangle the rôle of $\Gm$ and $X_\eta$ using the formula \refeq{coaction}:
$$\Gamma_{f_\eta?} = \Delta_? \t \DM(X_\eta),\eqlabel{Gamma? Delta?}$$ 
where $\Delta_?$ is the right adjoint of the $\DTM(\Gm)$-action on itself, which is more concretely given by
$$\DTM(\Gm) \t_{\DTM(S)} \DTM(\Gm) \stackrel[\cong]\boxtimes \r \DTM(\Gm^2) \stackrel{\Delta^*} \r \DTM(\Gm).$$
The first functor $\boxtimes$ is an equivalence by \thref{DTM.Gm} and \thref{generalities.tensor}.
Therefore, $\Delta_?$ is essentially the right adjoint of $\Delta^*$.
The functor $\Delta_?$ is \emph{distinct} from the right adjoint $\Delta_*$ of $\Delta^* \colon \DM(\Gm \x \Gm) \r \DM(\Gm)$. Instead, the inclusion $\DTM(\Gm \x \Gm) \subset \DM(\Gm \x \Gm)$ admits a right adjoint $R$, and $\Delta_? = R \circ \Delta_*$. The functor $R$ was studied by Totaro \cite{Totaro:Adjoint}.

More generally, if $X_\eta$ carries a universal Whitney--Tate stratification (in the sense of \cite[Definition~2.6]{CassvdHScholbach:MotivicSatake}), then one can show an equivalence
$$\DTM(\Gm) \t_{\DTM(S)} \DTM(X_\eta) \stackrel[\cong]{\boxtimes} \r \DTM(\Gm \x_S X_\eta).$$
Thus, in this case $\Gamma_{f_\eta?}$ can be computed on $\DTM(X_\eta)$ as $\Gamma_{f_\eta?} = R \circ \Gamma_{f_\eta*}$, where $R$ is again a right adjoint to the inclusion $\DTM(X_\eta) \subset \DM(X_\eta)$ of stratified Tate motives.
\item
Among other things, \thref{right.adjoint.action} asserts that $\Gamma_{f_\eta?}$ is compatible with *-pullbacks and !-pushforwards (along the $X$-direction), as well as with exterior products in the following sense. For two schemes $f_k : X_k \r \A^1$, the following diagram commutes, where we abbreviate $D := \DTMr(\Gm)$:
$$\xymatrix{
\DM(X_{1\eta}) \t \DM(X_{2\eta}) 
\ar[d]^{\Gamma_{f_{1\eta?}} \t \Gamma_{f_{2\eta?} }} 
\ar[r]^{\boxtimes_{\Gm}} & \DM(X_{1\eta} \x_{\Gm} X_{2\eta}) \ar[d]^{\Gamma_{(f_1 \x f_2)_{\eta ?}}} \\
(D \t \DM(X_{1\eta})) \t_D (D \t \DM(X_{2\eta}))  \ar[r]^(.6){\boxtimes_{\Gm}} & 
D \t \DM(X_{1\eta} \x_{\Gm} X_{2\eta}).}$$
\end{itemize}
\xrema

\lemm
There is a natural transformation
$$i^* j_* : U \DM_\eta \r \DM_s$$
of lax symmetric monoidal functors $\Corr_{\A^1}^{\sm, \proper} \r \Mod_{\DTM(S)}(\PrLSt)$.
\xlemm

\pf
The functors $i^* : \DM \r \DM_s$ and $j^* : \DM \r \DM_\eta$ are clearly transformations of lax symmetric monoidal functors.

By construction, the restriction of $\DM^*_!$ to $\Corr_{\A^1}^{\sm, \proper}$ is taking values in the non-full $\t$-subcategory $\PrLRomegaSt \subset \PrLSt$ (cf.~\thref{PrLomega.equivalence}), since the categories $\DM(-)$ are compactly generated, and the functors $f^*$ (for $f$ smooth) and $g_!$ (for $g$ proper) admit left adjoints, namely $f_\sharp$ and $g^*$, respectively.
Thus, $\DM^*_!$ is a lax $\t$-functor $\Corr_{\A^1}^{\sm, \proper} \r \PrLRomegaSt$.
Appending the $\t$-equivalence $\PrLRomegaSt \cong (\PrLomegaSt)^\opp$ obtained by passing to left adjoints (\thref{PrLomega.equivalence}), we obtain a symmetric monoidal functor denoted by
$$\DM_\sharp^* : \Corr_{\A^1}^{\sm, \proper} \r (\PrLomegaSt)^\opp.$$
It encodes $\sharp$-pushforwards along smooth maps $f$ and *-pullbacks along proper maps $g$.

Precomposition with $j_\sharp j^* \r \id$, which is a transformation of symmetric monoidal endofunctors on $\Corr_{\A^1}$ (its evaluation on $X/ \A^1$ is the map $j : X_\eta \r X$) induces the natural transformation of lax $\t$-functors $j^* : \DM^*_\sharp \r \DM^*_{\eta,\sharp}$.
Passing back to right adjoints, we obtain the requested transformation, of lax $\t$-functors, $j_* : \DM^*_{\eta!} \r \DM^*_!$.
\xpf

This accomplishes the construction of unipotent nearby cycles in the form stated in \thref{Upsilon.all.in.one}.

\rema
\thlabel{Upsilon lax functor}
The restriction to $\Corr_{\A^1}^{\sm, \proper} \subset \Corr_{\A^1}$ is only necessary to ensure that $j_* : \DM_\eta \r \DM$ is functorial. Without that restriction, $j_*$ and therefore also $\Upsilon$ is a \emph{lax} natural transformation of functors $\Fun(\Corr_{\A^1}, \PrLSt)$ obtained from the (non-lax) natural transformation $j^*$ by passing to adjoints, cf.~e.g.~\cite{Haugseng:Lax}.
\xrema

\subsection{Monodromy triangle}

The following fiber sequence has been achieved for motives with rational coefficients by Ayoub in \cite[Theorem~4.28]{Ayoub:Motivic}.
The approach taken here using Koszul duality highlights the fact that the relevance of rational coefficients is related to the non-formality of the motive of $\Gm$, cf.~\thref{Gm integral}, \thref{1^*}.
Recall the object $\Sigma = \bigoplus_{n \ge 0} \Z(n)$ from \refsect{Koszul duality}.

\prop \thlabel{lemm-monodromy-triangle}
Let $X / \A^1$, $M \in \DM(X_\eta, \Q)$ or $M \in \DMr(X_\eta)$.
There is a functorial isomorphism
$$\Upsilon(M) \t_\Sigma \Z = i^* j_* M(1)[1].\eqlabel{Upsilon tensor Z}$$
In other words, there is a cofiber sequence (in $\DM(X_s, \Q)$, resp.~$\DMr(X_s)$),
$$i^* j_* M \r \Upsilon(M) \r \Upsilon(M)(-1).\eqlabel{cofiber sequence Upsilon}$$
\xprop

\pf
Indeed, the cofiber sequence \refeq{resolution Z Sigma} shows that $\Z \t_\Sigma \Z$ is the split square zero extension $\Z \oplus \Z(1)[1]$. For reduced motives (and likewise after passing to rational coefficients for regular motives), this identifies with $\Lambda(1)[1]$, by \thref{Lambda square zero}. This yields
\begin{align*}
\Upsilon(M) \t_\Sigma \Z & = (i^* j_* M \t_\Lambda \Z) \t_{\Sigma} \Z \\
& = i^* j_* M \t_\Lambda \Lambda(1)[1].
\end{align*}
The cofiber sequence is then an immediate consequence.
\xpf

\rema
\thref{lemm-monodromy-triangle} shows that $i^* j_*$ can be recovered from $\Upsilon$ (while the converse holds by definition of $\Upsilon$).
This implies that Beilinson's glueing theorem for perverse sheaves (e.g. \cite[Theorem~8.1]{Morel:Beilinson}) can be achieved motivically as soon as one has a motivic t-structure, e.g., for mixed stratified Tate motives over an appropriately stratified scheme, cf.~\cite[Lemma~2.15]{CassvdHScholbach:MotivicSatake}. 
\xrema

\rema
\thlabel{trivial monodromy}
Both maps in \refeq{cofiber sequence Upsilon} exist already integrally, but we only claim it is a cofiber sequence for rational motives (or reduced motives).
For the latter map, see \thref{1^*}. The map $\alpha : i^* j_* M \r \Upsilon(M)$ arises from applying $i^* j_* M \t_{\Lambda} -$ to the unit map $\Lambda = p_* p^* \Z \r p_* 1_* 1^* p^* \Z = \Z$, where $S \stackrel 1 \r \GmX S \stackrel p \r S$. We say that $\Upsilon(M)$ has \emph{trivial monodromy} if $\alpha$ splits. By functoriality, this condition on $M$ is preserved under proper pushforward and smooth pullback.
\xrema

\subsection{Computation for trivial families}
\lemm \thlabel{lemm-A1-Upsilon}
For the trivial family $f = \id_{\A^1}$ we have $\Upsilon_\id(\Z) = \Z$, and the monodromy is trivial (in the sense of \thref{trivial monodromy}).
\xlemm

\pf
We use the following diagram, in which the functors $\Delta^*$, $\id \t j^*$ and $\id \t q^*$ are the left adjoints of the adjacent functors in the opposite direction.
All tensor products are over $\DTM(S)$, and $\DTM(\A^1)$ denotes stratified Tate motives with respect to the stratification given by $\Gm \sqcup 0$.
Note that $q_*$ maps these motives to $\DTM(S)$, so $\id \t q_*$ below takes values in $\DTM(\Gm) = \DTM(\Gm) \t \DTM(S)$.
$$\xymatrix{
\DTM(\Gm) \t \DTM(\Gm) \ar@<-.5ex>[d]_{\Delta^*} \ar@<-.5ex>[r]_{\id \t j_*} & \DTM(\Gm) \t \DTM(\A^1) \ar@<-.5ex>[l]_{\id \t j^*} \ar@<.5ex>[d]^{\id \t q_*} \ar[r]^(.6){\id \t i^*} & \DTM(\Gm) \ar[d]^{1^*}  \\
\DTM(\Gm) \ar@{=}[r] \ar@<-.5ex>[u]_{\Delta_?} 
& \DTM(\Gm) \ar@{=}[ur] \ar[r]^{1^*} \ar@<.5ex>[u]^{\id \t q^*} &
\DTM(S).
}\eqlabel{diagram yeah}$$
At the left half, we have $(\id \t q_*) (\id \t j_*) \Delta_? = \id_{\DTM(\Gm)}$ since it is right adjoint to $\Delta^* (\id \t j^*) (\id \t q^*) = (q \circ (\id \x j) \circ \Delta)^* = \id_{\Gm}^*$.

The unit map 
$$q_* \r q_* i_* i^* = i^*\eqlabel{q* i* iso}$$ 
is an isomorphism when evaluated on objects in $\DTM(\A^1)$: it is enough to check this for $i_* \Z$, which is clear, and for $\Z$, where it follows from homotopy invariance, i.e., $q_* q^* = \id$.

The triviality of the monodromy also follows from this: under the equivalence $\DTM(\Gm) = \Mod_\Lambda$, $i^* j_* \Z$ corresponds to $\Lambda$, and the unit map in \thref{trivial monodromy} is the augmentation map $\Lambda \r \Z$, which is split by the unit map for $\Lambda$.
\xpf

\rema
The importance of the isomorphism $q_* \Z = i^* \Z$ is hard to overstate.
First, it is also the basis of the same result in Ayoub's approach \cite[Proposition~4.9]{Ayoub:Motivic} (or \cite[Proposition~3.4.9]{Ayoub:Six2}).
It is also the core computation in the proof of hyperbolic localization, cf.~\cite[Corollary~2.9]{Richarz:Spaces} or \cite[Proposition~3.2.2]{DrinfeldGaitsgory:Theorem}.
The full force of $\A^1$-invariance is actually not needed, instead it suffices to use that the pullback along the stack quotients $\overline q : \A^1 / \Gm \r B\Gm (:= S / \Gm)$ is fully faithful.
This should prove useful in applying the above paradigm to prismatic cohomology, cf.~\cite[Lemma~5.3.2]{KubrakPrikhodko:p-adic}, and also to the motivic spaces considered in \cite{AnnalaHoyoisIwasa:Algebraic}.
\xrema

\coro
\thlabel{nearby trivial family}
Suppose $f : X = Y \x \A^1 \r \A^1$ is a trivial family.
Let $M \in \DM(Y)$ be such that the natural map
$$M \boxtimes j_* \Z \r (\id_Y \x j)_* (M \boxtimes \Z) \eqlabel{universal WT condition}$$
is an isomorphism.
Then there is a natural isomorphism (and the monodromy is trivial)
$$\Upsilon(M \boxtimes \Z_{\Gm}) = M.$$
\xcoro

\pf
We have the following diagram,
where the bottom right horizontal functor is ${\id \t (i \x \id)^* (j \x \id)_*}$. (Recall the convention that $\t$ is meant to be over $\DTM(S)$.)
$$\xymatrix{
\DTM(\Gm) \t \DM(Y) \ar[d]^\boxtimes \ar[r]^{\Delta_? \t \id} & 
\DTM(\Gm)^{\t 2} \t \DM(Y) \ar[r]^{\id \t i^* j_* \t \id} \ar[d]^{\id \t \boxtimes} & 
\DTM(\Gm) \t \DM(Y) \ar@{=}[d] \ar@{=>}[dl]
\\
\DM(\Gm \x Y) \ar[r]^(.4){(\Delta \x \id_Y)_?} &
\DTM(\Gm) \t \DM(\Gm \x Y) \ar[r] &
\DTM(\Gm) \t \DM(Y).
}$$
The left square commutes since the left vertical $\boxtimes$-functor is a map of $\DTM(\Gm)$-modules, and thus also a map of $\DTM(\Gm)$-comodules (\thref{ModD.equivalence}).
The two left horizontal arrows are just these coaction maps.

The indicated natural transformation at the right stems from the commutativity of $i^*$ with $\boxtimes$ (which holds for any motive in $\DM(\A^1) \t \DM(Y)$), and the map \refeq{universal WT condition} (which in general need not be an isomorphism). 

If we append $1^* \t \id$ to the top row, we get $\Upsilon_{\id} \t \id_{\DM(Y)}$.
Using $\Upsilon(\Z_{\Gm}) = \Z$ (\thref{lemm-A1-Upsilon}) our claim follows.
\xpf

\rema
\thlabel{KuennethYinJang}
If the base scheme $S$ is a field of exponential characteristic $n$, then \refeq{universal WT condition} (and also \refeq{j* box} below) is an isomorphism for arbitrary motives with $\Z[1/n]$-coefficients \cite[Theorem~2.4.6]{JinYang:Kuenneth}.

The motivation for the tailored statement above is that we will perform the construction of central motives over $S = \Spec \Z$, where that property does not hold in general.

We also note that Jin--Yang \cite{JinYang:Some} take a similar approach by singling out a ULA-ness condition that ensures compatibility of tame nearby cycles with exterior products.
\xrema

\subsection{Preservation of Tate motives}

\prop
\thlabel{Upsilon preserves Tate}
Suppose $X / \A^1$ carries an (anti-effective) Whitney--Tate stratification that is compatible with the stratification of $\A^1_S$ by $\GmX S \sqcup S$.
Then $\Upsilon$ preserves (anti-effective) stratified Tate motives, i.e., it restricts to a functor
$$\Upsilon: \DTM(X_\eta, X_\eta^+)^{(\anti)} \r \DTM(X_s, X_s^+)^{(\anti)}.$$
\xprop

\pf
By \refeq{Gamma? Delta?}, we can compute $\Gamma_{f_\eta?}$ as 
$$\Delta_? \t \id : \DTM(\Gm) \t_{\DTM(\Gm)} \DM(X_\eta) \r 
\DTM(\Gm^2) \t_{\DTM(\Gm)} \DM(X_\eta),$$
where $\Delta_?$ is the right adjoint of $\Delta^*$ (restricted to Tate motives!).
In particular, $\Gamma_{f_\eta?}$ maps $\DTM(X_\eta)^{(\anti)}$ to the full subcategory $\DTM(\Gm) \t \DTM(X_\eta)^{(\anti)}$.
Since $i^* j_*$ preserves (anti-effective) Tate motives by assumption, and $1^* : \DTM(\Gm) \r \DTM(S)$ also does, this shows that $\Upsilon$ preserves (anti-effective) Tate motives.
\xpf

\subsection{Compatibility with hyperbolic localization}
We address the compatibility of unipotent nearby cycles with hyperbolic localization. Given a scheme $X / \A^1$ endowed with a Nisnevich-locally linearizable $\Gm$-action (where $\A^1$ is acted upon trivially), we consider the hyperbolic localization diagram (e.g., \cite{DrinfeldGaitsgory:Theorem, Richarz:Spaces})
$$\xymatrix{
X^0 \ar[dr]_{f^0} & X^\pm \ar[d]^{f^\pm} \ar[l]_{p^\pm} \ar[r]^{q^\pm} & X \ar[dl]^f \\ & \A^1.}\eqlabel{hyp.loc}$$
Here $X^+ = \Map^{\Gm}(\A^1, X)$ are the attractors of the $\Gm$-action, $X^- = \Map^{\Gm}(\A^1_\anti, X)$ are the repellers and $X^0 = \Map^{\Gm}(S, X)$ are the fixed points of the $\Gm$-action.
We will write subscripts $s$ and $\eta$ for pullbacks along $0 \subset \A^1$ and $\Gm \subset \A^1$ and superscripts $0$ and $+$ for the corresponding objects on the level of fixed points and attractors.
Note that $(X_\eta)^+ = X^+ \x_{\A^1} \Gm$ etc.

The subcategory $\DM(X_\eta)^{\Gmmono} \subset \DM(X_\eta)$ of \emph{$\Gm$-monodromic motives} is defined to be the presentable subcategory generated by $\Gm$-equivariant motives, i.e., by the image of the forgetful functor $\DM(X_\eta / \Gm) \r \DM(X_\eta)$ (see \eqref{DM equivariant} for the category of equivariant motives).
Braden's hyperbolic localization asserts that the natural map
$$q^-_* p^{-!} \r q^+_! p^{+*} \ (\in \Fun(\DM(X), \DM(X^0))),\eqlabel{hyperbolic loc map}$$
is an isomorphism when restricted to $\DM(X)^{\Gmmono}$ (cf.~\cite[Proposition~2.5]{CassvdHScholbach:MotivicSatake} for this motivic statement and the references there for precursor statements for $\ell$-adic sheaves and D-modules).

\prop \thlabel{prop-Gm-Upsilon}
Unipotent nearby cycles commute with hyperbolic localization in the sense that the following diagram commutes:
$$\xymatrix{
\DM(X_\eta)^{\Gmmono} \ar[d]^{p^+_{\eta !} q_\eta^{+*}} \ar[r]^{\Upsilon_f}  &
\NilpQ \DM(X_s) \ar[d]^{p^+_{s!} q_s^{+*}} \\
\DM(X^0_\eta) \ar[r]^{\Upsilon_{f^0}} &
\NilpQ \DM(X^0_s).
}$$
\xprop

\pf
Similarly to the proof of \thref{Upsilon preserves Tate}, the functor $\Gamma_{f_\eta?} = \Delta_? \t \id$ maps $\DM(X_\eta)^{\Gmmono}$ to $\DTM(\Gm) \t \DM(X_\eta)^{\Gmmono}$.
On this category, Braden's theorem ensures that $\id_{\DTM(\Gm)} \t i^* j_*$ commutes with the hyperbolic localization functors in \refeq{hyperbolic loc map}.
And so does, for trivial reasons, $1^*$.
\xpf

\subsection{Compatibility with realization}

Suppose given a \emph{realization functor}, i.e., a natural transformation 
$$\rho: \DM \r D$$
(of lax symmetric monoidal functors $\Corr \r \PrLSt$; note this forces $\rho$ to be compatible with *-pullbacks, !-pushforwards, and exterior products).
Then the construction of $\Upsilon$ can be repeated verbatim for $D$ if we replace $\DTM(X)$ by the subcategory $DT(X) \subset D(X)$ generated by the objects $\rho(\Z(n))$.

\lemm
Suppose $\rho$ is right adjointable, i.e., that it commutes with *-pushforwards and !-pullbacks.
Then, for $X / \A^1$, there is a natural commutative diagram
$$\xymatrix{
\DM(X_\eta) \ar[d]^{\rho} \ar[r]^{\Upsilon} & \DM(X_s) \ar[d]^{\rho} \\
D(X_\eta) \ar[r]^\Upsilon & D(X_s).
}$$
\xlemm

\pf
We will prove that each functor in the composite
$$\Upsilon_X : \DM(X_\eta) \stackrel{\Gamma_{f_\eta ?}} \r \DTM(\Gm) \t_{\DTM(S)} \DM (X_\eta) \stackrel {1^* \t i^* j_*} \r \DM(X_s)$$
is compatible with $\rho$.
This holds by the naturality and right adjointability of $\rho$ for $1^* \t i^* j_*$.
It remains to check that 
$$\xymatrix{
\DTM(\Gm^2) \ar[r]^{\Delta^*} \ar[d]^\rho &
\DTM(\Gm) \ar[d]^\rho \\
DT(\Gm^2) \ar[r]^{\Delta^*} & DT(\Gm)
}$$
is right adjointable, i.e., that $\Delta_?$ commutes with $\rho$. This implies our claim since $\Gamma_{f_\eta?} = \id \t \Delta_?$.
Using \thref{DTM.Gm} (and the identical statement for $DT$), this means that
$$\xymatrix{
\Mod_{\Lambda^{\t 2}} (\DTM(S)) \ar[d]^{\rho} & \Mod_\Lambda (\DTM(S)) \ar[d]^\rho \ar[l]^{\res} \\
\Mod_{\rho(\Lambda^{\t 2})} (DT(S)) & \Mod_{\rho(\Lambda)} (DT(S)) \ar[l]^\res
}$$
should commute. This is the case as can be seen after applying the conservative forgetful functor to $DT(S)$.
\xpf 

\exam
According to \cite[Théorème 3.7]{Ayoub:Note}, \cite[§6]{Ayoub:Realisation}, and \cite[Proposition~3.9]{EberhardtScholbach:Integral}, respectively, examples of such realization functors are the Betti and $\ell$-adic realization functor and the reduction functor $\DM \r \DMr$.

Beilinson's approach to unipotent nearby cycles (see, e.g., \cite[Corollary~4.3]{Morel:Beilinson} for a recent exposition) is
$$\Upsilon^B (F) := \colim_n i^* j_* (L_n \t F),$$
where $L_n = \Z \oplus \dots \oplus \Z(-n)$, where $\Z(1)$ is the local system on $\Gm^\an$ associated to $\pi_1(\Gm^\an)$.
We have a natural isomorphism
$$L_n = \colim (\Z(-n)[-n] \r \dots \r \Z(-1)[-1] \r \Z)$$
so that the definition of $\Upsilon$ (for sheaves in the analytic topology), cf.~\refeq{Upsilon explicit} agrees with $\Upsilon^B$.

In the context of analytic sheaves, it is known that $\Upsilon^B$ is a direct summand of the full (i.e., not necessarily unipotent) nearby cycles functor, and is therefore t-exact and preserves constructible sheaves.
Without using a notion of full motivic nearby cycles, we will prove such properties for Tate motives in the special situation of the deformation of, say, $\Gr_G \x \Gr_G \x \Gm$ to $\Fl$ in \thref{prop--satake-exact} and \thref{theo-WT-BD}. 
\xexam

\subsection{Comparison with Ayoub's approach}
\label{sect--comparison Ayoub}
We now show that our definition of $\Upsilon$ agrees with the one due to Ayoub \cite{Ayoub:Six2}, \cite[§10]{Ayoub:Realisation}. Ayoub's construction is based on the geometric bar construction associated to the diagram
$$\xymatrix{& Y := \Gm \ar[d]^{\id \x 1} \\ 
X := \Gm \ar[r]^\Delta & B := \Gm^2.}\eqlabel{XYB}$$
Here $\Delta$ denotes the diagonal map; for clarity we will denote by $\mathbf{\Delta}$ the category of finite linearly ordered sets below.
We denote this geometric bar construction, which is a cosimplicial $\Gm$-scheme, by $A: \mathbf{\Delta} \r \Sch_{\Gm}$.
In low degrees, it is given by
$$\xymatrix{\Gm \ar@<.5ex>[rr]^\Delta \ar@<-.5ex>[rr]_{\id \x 1} & &
\Gm^2 
\ar@<1ex>[rr]^{\Delta \x \id} 
\ar[rr]|{\id \x \Delta} 
\ar@<-1ex>[rr]_{\id \x \id \x 1}
&  & 
\Gm^3 
\ar@<1ex>[rr] \ar@<.33ex>[rr] \ar@<-.33ex>[rr] \ar@<-1ex>[rr] & & \dots}$$ 
We refer to, say, \cite[Lemme~3.4.1]{Ayoub:Six2} for its full definition.
Below, we will use that the terms in this cosimplicial scheme are $\Gm$-isomorphic to schemes of the form $\Gm \x \Gm^n$.
We also have the functor 
$$\Gamma : \Sch_{\Gm}^\opp \r \DM(\Gm), (X \stackrel f \r \Gm) \mapsto f_* f^* \Z.\eqlabel{Gamma}$$

\defi Ayoub's unipotent nearby cycles functor is defined as
$$\Upsilon^A : \DM(X_\eta) \r \DM(X_s), M \mapsto \Upsilon^A (M) := i^* j_* (\mathscr U \t M),$$
where
$$\mathscr U := \colim (\mathbf{\Delta}^\opp \stackrel A \r \Sch_{\Gm}^\opp \stackrel \Gamma \r \DM(\Gm)).$$
\xdefi

\rema
The above presentation of $\mathscr U$, which is somewhat easier to digest than Ayoub's original approach using categories of motives over diagrams of schemes, can be found in \cite[§2.2]{JinYang:Some}. 

Previously to this definition, Ayoub had proposed another, mildly different definition \cite[Definition~4.3]{Ayoub:Motivic},
$$\Upsilon^{A'}(M) := i^* j_* \IHom(\mathscr U', M),$$ where $\mathscr U' := \lim (\mathbf{\Delta} \stackrel A \r \Sch_{\Gm} \stackrel \Gamma \r \DM(\Gm) \stackrel{(-)^\dual} \r \DM(\Gm)^\opp)$.
That is, instead of tensoring termwise with $f_{n*} f_n^* \Z$ (for $f_n : A_n \r \Gm$ the structural map in the $n$-th term in the cosimplicial scheme given by the bar construction) one can also take $\IHom(\lim f_{n!} f_n^! \Z, -)$.
Given that $\mathscr U'$ is the dual of $\mathscr U$ we have a natural map $\Upsilon^A \r \Upsilon^{A'}$.
However, since $\mathscr U$ is not dualizable (unlike the individual terms in that diagram), it is not an isomorphism.
From a structural point of view, the former functor is preferable since it preserves colimits.
\xrema

\prop
Ayoub's unipotent nearby cycle is isomorphic to the one defined in \thref{unipotent nearby cycles}:
$$\Upsilon^A = \Upsilon.$$
\xprop

\pf
Applying $- \t_{\DTM(\Gm)} \DM(X_\eta)$ to the functors
$$\DTM(\Gm) \stackrel {\Delta_?} \r \DTM(\Gm^2) \stackrel {(\id, 1)^*} \r \DTM(\Gm),\eqlabel{Delta ? etc}$$
gives us, according to \refeq{Gamma? Delta?},
$$\DM(X_\eta) \stackrel{\Gamma_{f_\eta}?} \r \DTM(\Gm) \t_{\DTM(S)} \DM(X_\eta) \stackrel {1^* \t \id} \r \DM(X_\eta).$$
By definition, $\Upsilon$ is the composition of this with $i^* j_*$. Thus, it suffices to prove that the functor $(\id, 1)^* \Delta_?$ agrees with $- \t \mathscr U$.

Let us write $C^\t := (\Sch_{\Gm}^\opp)^{\t}$ for the ordinary symmetric monoidal category of $\Gm$-schemes (with their fiber product over $\Gm$). We regard any $\Gm$-scheme as an algebra-object (in the opposite category) via the diagonal, and $X$ as a $B$-module via $\Gamma_\Delta : X \r X \x_{\Gm} B$ etc.
In the notation of \cite[§4.4.1, 4.4.2]{Lurie:HA}, the diagram \refeq{XYB} then amounts to a map $\Tens_{[2]}^\t \r C^\t$.
By inspection of the definitions, the composite $\mathbf{\Delta}^\opp \stackrel U \r \Tens_{[2]}^\t \r C^\t$, with $U$ as in \cite[Notation~4.4.2.4]{Lurie:HA}, is such that its pushforward along the maps $[m] \stackrel \alpha \r [n] \r [1]$ in $\Fin_*$ agrees with the $m$-th, and $n$-th terms of the geometric bar construction, including the simplicial structure mentioned above.

The functor $\Gamma$ in \refeq{Gamma} is lax symmetric monoidal.
When restricted to the full subcategory $C' \subset \Sch_{\Gm}$ consisting of $\Gm$-schemes that are isomorphic to $\Gm \x \Gm^n$, the functor is symmetric monoidal and produces Tate motives, i.e., we have a symmetric monoidal functor
$$\Gamma : (C'^\opp)^\t \r \DTM(\Gm) \stackrel{\text{\ref{DTM.Gm}}} = \Mod_\Lambda(\DTM(S)).$$
Applying $\Gamma$ to \refeq{XYB}, gives 
$$\xymatrix{
& \Lambda \\
\Lambda & \Lambda \t \Lambda \ar[l]_\mu \ar[u]_{\id \t \aug},
}$$
i.e., the multiplication map of the algebra object $\Lambda$, and the augmentation map. 
The colimit of the bar construction $\Bar(\Lambda, \Lambda \t \Lambda, \Lambda)$ (formed in $\Mod_\Lambda(\DTM(S))$, using these maps) agrees with $\Lambda \t_{\mu, \Lambda \t \Lambda, \id \t \aug} \Lambda$.
This is precisely the image of $\Lambda$ under the functor \refeq{Delta ? etc}.
\xpf

\subsection{Higher-dimensional unipotent nearby cycles}
Let $n \ge 0$ and consider the following ``$n$-dimensional'' variant of $\Nilp$ and $\NilpQ$ obtained by redoing the construction of \thref{Nilp} and \thref{NilpQ} for $p^n_* p^{n*} \Z$, where $p^n : \Gm^n \r S$ is the structural map.
By the Künneth formula (for Tate motives), this is isomorphic to $\Lambda^{\t n} = (p_* p^* \Z)^{\t n}$.
If $\Lambda_\red = \Z \oplus \Z(-1)[-1]$ denotes the square zero extension in $\D(\grAb)$, we again have a fully faithful, symmetric monoidal functor
$$\Mod_{\Lambda_\red^{\t n}} \D(\grAb) \stackrel {\Z \t_{\Lambda_\red^{\t n}} -} \r \left (\Mod_{(\bigoplus_{k \ge 0} \Z(k))^{\t n}}(\D(\grAb)), \t_\Z \right)$$
and we define $\Nilp^n$ to be the essential image of this functor.
As in \thref{NilpQ}, we let $\NilpQ^n := \D(\grAb) \x_{(\D(\grAb) \t \Q)} (\Nilp^n \t \Q)$.
The category $\NilpQ^n$ is linear over $\D(\grAb)$, and for some $\D(\grAb)$-module $C$, we write $\NilpQ^n C := \NilpQ^n \t_{\D(\grAb)} C$ as before.
We have a natural symmetric monoidal functor
$$\DTM(\Gm^n) \r \NilpQ^n \DTM(S),$$
which is an equivalence after passing to rationalizations (cf.~\thref{DTM Gm Nilp}).

The following functor serves as a replacement and far-reaching extension of the constructions in \cite[§9.4]{AcharRiche:Central}, which in turn goes back to \cite[§4]{Gaitsgory:Braiding}.

\defi
The \emph{$n$-dimensional unipotent nearby cycles functor} 
$$\Upsilon_n := \Upsilon_{X / \A^n} : \DM(X \x_{\A^n} \Gm^n) \r \Nilp^n_\Q \DM(X \x_{\A^n, 0} S)$$
is obtained by replacing $\Gm$ by $\Gm^n$ in \thref{unipotent nearby cycles}.
\xdefi

This functor has essentially the same good properties. In more detail, $\Upsilon_n$ enjoys the same functoriality as in \thref{Upsilon.all.in.one}.
The formula \refeq{Upsilon tensor Z} continues to hold as stated, but the cofiber sequence in \refeq{cofiber sequence Upsilon} becomes slightly more involved (since a resolution of $\Z$ by free $\Sigma^n$-modules is longer).
For the computation of $\Upsilon_n(\Z) = \Z$ for the trivial family $\id_{\A^n}$ (\thref{lemm-A1-Upsilon}) note that for any $\Gm$-equivariant sheaf on $\A^n$, the map in \refeq{q* i* iso} is still an isomorphism when applied to $j_* \Z$, where $j : \Gm^I \x 0^{J} \r \A^{I \sqcup J}$, as can be seen by applying \refeq{q* i* iso} iteratively.
The preservation of Tate motives and compatibility with hyperbolic localization can also be proven the exact same way.
Finally, as in \thref{remarks Upsilon}\refit{Upsilon reduced}, the same definition leads to a functor $\Upsilon_n : \DMr(X \x_{\A^n} \Gm^n) \r \Nilp^n \DMr(X \x_{\A^n, 0} S)$ between categories of \emph{reduced} motives with an integrally-defined monodromy map.

In this section, we collect a few basic compatibilities of this construction as $n$ varies. For concreteness, we highlight only the relation of $\Upsilon_1$ vs.~$\Upsilon_2$, which is going to be used in~\thref{lemm-conv-nearby}.

We consider the diagonal $\A^1 \stackrel \Delta \r \A^2$.
Taking the pullback to $\Gm^2$, and taking cohomology (as in \refeq{PI}) gives the map $\Lambda_\Q^2 = (\Sym \Q(-1)[-1])^{\t 2} \r \Lambda_\Q = \Sym (\Q(-1)[-1]).$
We will use the adjunction (pullback, and forgetful functor, respectively):
$$\Delta^* : \Corr_{\A^2} \rightleftarrows \Corr_{\A^1} : \Delta_*.$$
In the following, we consider the functor (for $X / \A^1$)
$$\res_{\Sigma^2}^\Sigma \Upsilon_2 \Delta_* : \DM(X_\eta) \r \Nilp_\Q \DM(X_s),$$ which regards $X$ as a scheme over $\A^2$ via the diagonal, then applies 2-dimensional nearby cycles, and then appends the restriction 
$$\res_{\Sigma^2}^\Sigma: \Mod_{\Sigma_\Q^{\t 2}}(\DTM(S)) \r \Mod_{\Sigma_\Q} (\DTM(S))\eqlabel{restriction Mod Sigma}$$
along the comultiplication in \eqref{comultiplication Sigma} (which maps $\NilpQ^2$ to $\NilpQ$).

\prop
\thlabel{restriction Upsilon}
For schemes over $\A^1$, the two-dimensional nearby cycles functors are related to the (one-dimensional) ones by means of a natural transformation
$$\res_{\Sigma^2}^{\Sigma} \Upsilon_2 \Delta_* \r \Upsilon$$
of functors $\Corr^{\pr, \sm}_{\A^1} \r \Fun(\Delta^1, \PrLSt)$.
Thus for any $X / \A^1$, and writing $X' := \Delta_* (X)$ for the same scheme, but now regarded as lying over (the diagonal in) $\A^2$, there is a natural transformation (which is compatible with smooth pullbacks and proper pushforwards):
$$\xymatrix{
\DM(X_\eta) \ar@{=}[d] \ar[r]^(.4){\Upsilon} & \NilpQ(\DM(X_s)) \\
\DM(X'_\eta) \ar[r]^(.4){\Upsilon_2} \ar@{=>}[ur] & \NilpQ^2(\DM(X'_s)) \ar[u]_{\res}
}$$
\xprop

\pf
Note that $X_\eta = X'_\eta := X \x_{\A^2} \Gm^2$, $X_s = X'_s := X \x_{\A^2, 0} S$.
We write $f_\eta^2 : X_\eta \r \Gm^2$ etc. 
Let again $\Delta_?$ be the right adjoint of $\Delta^* : \DTM(\Gm^2) \r \DTM(\Gm)$.
By passing to left adjoints, one sees that there is a diagram, where all tensor products are over $\DTM(S)$:
$$\xymatrix{
\DM(X_\eta) \ar@{=}[d] \ar[r]^(.3){(\Gamma_{f_\eta})_?} & 
\DTM(\Gm) \t \DM(X_\eta) \ar[r]^(.5){\id \t i^* j_*} \ar[d]^{\Delta_? \t \id} & 
\DTM(\Gm) \t \DM(X_s) \ar[d]^{\Delta_? \t \id} \ar[r]^(.6){1^* \t \id} & 
\NilpQ \DM(X_s) \\
\DM(X'_\eta) \ar[r]^(.3){(\Gamma_{f_\eta^2})_?} & \DTM(\Gm^2) \t \DM(X'_\eta) \ar[r]^{\id \t i^* j_*} & 
\DTM(\Gm^2) \t \DM(X'_s) \ar[r]^(.6){(1,1)^* \t \id} \ar@{=>}[ur] & 
\NilpQ^2 \DM(X'_s) \ar[u]^\res. 
}$$
The left square commutes, since their left adjoints commute: $\Gamma_{f_\eta}^* \circ (\Delta^* \t \id) = \Gamma_{f_\eta^2}^*$.
The middle square commutes trivially.
The indicated natural transformation at the right arises again by passing to adjoints, using the identity $1^* \Delta^* = (1,1)^*$.
\xpf

\rema
The natural transformation at the right is not invertible. Suppressing the $\DM(X_s)$-factor from the notation and forgetting the $\NilpQ$-action, the square identifies with
$$\xymatrix{
\Mod_{\Lambda}(\DTM(S)) \ar[r]^(.6){- \t_{\Lambda} \Z} \ar[d]^{\res} &
\DTM(S) \ar@{=}[d] \\
\Mod_{\Lambda^2} (\DTM(S)) \ar[r]^(.6){- \t_{\Lambda^2} \Z}  &
\DTM(S).}$$
Here $\res$ is the restriction functor along the map $\Lambda^2 = \Gamma(\Gm^2) \r \Gamma(\Gm) = \Lambda$ obtained by functoriality of $\Gamma$, cf.~\refeq{Gamma}.
\xrema

The following statement is a replacement and enhancement of \cite[Lemma~9.4.12]{AcharRiche:Central}.

\coro \thlabel{coro--nearby-pullback}
For $X / \A^2$, let us write $\Delta : X' := X \x_{\A^2} \A^1 \r X$ for the pullback along the diagonal etc. 
For $M \in \DM(X)$, there is a natural map in $\NilpQ \DM(X_s)$
$$\res_{\Sigma^2}^\Sigma \Upsilon_2 (M) \r \Upsilon (\Delta^* M).$$
It is functorial in $X$ (with respect to smooth-proper correspondences over $\A^2$).
\xcoro

\pf
As before, we write subscripts $\eta$ and $s$ for the pullbacks along $\Gm^2 \subset \A^2$ and $(0,0) \subset \A ^2$.
The unit map $M \r \Delta_{\eta*} \Delta_\eta^* M$ yields a map 
$$\Upsilon_{X / \A^2} (M) \r \Upsilon_{X / \A^2}(\Delta_{\eta*} \Delta_\eta^* M) = \Delta_{s*} \Upsilon_{X' / \A^2} (\Delta_\eta^* M).$$
The second isomorphism is the compatibility with proper pushforward. The map $\Delta_s : X' \x_{\A^2, (0,0)} S \r X \x_{\A^2, (0,0)} S$ is an isomorphism.
We then append the map supplied by \thref{restriction Upsilon}, $\res_{\Sigma^2}^\Sigma \Upsilon_{X'/\A^2} (\Delta_\eta^* M) \r \Upsilon_{X' / \A^1} (\Delta_\eta^* M)$.
\xpf

\subsubsection{Compatibility with exterior products}
In order to state the following result, we use the symmetric monoidal functor
$$\NilpQ \t_{\D(\grAb)} \NilpQ \r \NilpQ^2 := \D(\grAb) \x_{\D(\grAb) \t \Q} (\Nilp^2 \t \Q)$$
whose first component is induced by the equivalence $\D(\grAb) \t_{\D(\grAb)} \D(\grAb) \stackrel \cong \r \D(\grAb)$ and whose second component is induced by the (rationalization of the) equivalence $\Nilp \t_{\D(\grAb)} \Nilp = \Nilp^2$ which follows from \thref{Nilp} and \thref{generalities.tensor}.
The description of compact generators of $\NilpQ$ (which carries over to $\NilpQ^2$) in \thref{NilpQ explained} implies this functor is an equivalence.

\prop
\thlabel{Upsilon.boxtimes}
Let $f_k : X_k \r \A^1$, $k = 1, 2$ be two schemes, and $f : X := X_1 \x X_2 \r \A^2$ their product.
For $M_k \in \DM(X_{k\eta})$, there is a natural map (in $\NilpQ^2 \DM(X_{1s} \x X_{2s})$)
$$\Upsilon_{X_1 / \A^1} (M_1) \boxtimes \Upsilon_{X_2 / \A^1} (M_2) \r \Upsilon_{X / \A^2} (M_1 \boxtimes M_2)\eqlabel{Upsilon box map}$$
which is an isomorphism if the natural map
$$(j_{1*} M_1) \boxtimes (j_{2*} M_2) \r (j_1 \x j_2)_* (M_1 \boxtimes M_2)\eqlabel{j* box}$$
is an isomorphism (concerning this condition, cf.~\thref{KuennethYinJang}).
\xprop

\pf
We consider the following diagram, where all tensor products are over $\DTM(S)$: 
$$\xymatrix{
\bigotimes_{k=1}^2 \DM(X_{k\eta}) \ar[d]^\boxtimes \ar[r]^(.4){\Gamma_{f_k?}} &
\bigotimes_{k=1}^2 (\DTM(\Gm) \t \DM(X_{k\eta})) \ar[d]^{\boxtimes} \ar@{=>}[dl] \ar[r]^(.6){1^* \t i_k^* j_{k*}} &
\bigotimes_{k=1}^2 \NilpQ \DM(X_{ks}) \ar[d]^\boxtimes \ar@{=>}[dl] \\
\DM(X_{1\eta} \x X_{2\eta}) \ar[r]^(.4){\Gamma_{f?}} & 
\DTM(\Gm^{2}) \t \DM(X_{1\eta} \x X_{2\eta}) \ar[r]^(.6){(1,1)^* \t i^* j_*} &
\NilpQ^2 \DM(X_{1s} \x X_{2s})\\
}$$
The top horizontal composite is $\Upsilon \t \Upsilon$.
The left horizontal functors are the right adjoints of the action of $\DTM(\Gm)$ on $\DM(X_{k\eta})$ (resp.~likewise for $X_\eta \r \Gm^{2}$).
(Therefore, the bottom horizontal composite relates closely to the functor considered in \cite[9.4.1 Definition]{AcharRiche:Central}.)
These actions are clearly compatible with the indicated $\boxtimes$-functors, which gives rise to the natural transformation at the left.

The left hand vertical $\boxtimes$ is a map of $\DTM(\Gm^2)$-modules, and thus also of comodules (\thref{ModD.equivalence}). The left hand horizontal functors are just the coaction maps.
Thus, the left hand part commutes.
The natural transformation in the right hand square arises since $*$-pullbacks are compatible with $\boxtimes$, while *-pushforwards are only lax compatible as in \refeq{j* box}.
If that latter map is an isomorphism, then so is \refeq{Upsilon box map}. Indeed, $\Gamma_{f_k?}(M_k)$ is a colimit of objects of the form $\Z_{\Gm}(-i)[-i] \boxtimes M_k$, and $j_{k*}$ and $\boxtimes$ preserve colimits.
\xpf

\rema
Along the lines of the proof above (or the ones in \thref{nearby trivial family}), one sees that for $Y / S$, $X / \A^1$, $N \in \DM(Y)$, $M \in \DM(X_\eta)$ there is a natural map (as usual, compatible with monodromy, smooth pullback and proper pushforward):
$$N \boxtimes \Upsilon_{X / \A^1}(M) \r \Upsilon_{Y \x X / \A^1}(N \boxtimes M).\eqlabel{Upsilon vs outer box}$$
\xrema

\subsubsection{Decomposing higher-dimensional nearby cycles}
\lemm \thlabel{lemm-composition}
Let $X / \A^2$, write $X_\eta := X \x_{\A^2} \Gm^2$, and let $X' := X \x_{\A^2} \A^1 \x 0 / \A^1$ via the first coordinate in $\A^2$, so that $X'_\eta = X \x_{\A^2} (\Gm \x 0)$, and $X'_s = X \x_{\A^2} (0, 0) =: X_s$. Denote by 
\begin{align*}
\Upsilon_2 : & \DM(X_\eta) \r \NilpQ \DM(X'_\eta) \\
\Upsilon_1 : & \DM(X'_\eta) \r \NilpQ \DM (X'_s)
\end{align*}
the (one-dimensional) nearby cycles along the second and first coordinate, respectively.
Finally, write
$$\Upsilon : \DM(X_\eta) \r \NilpQ^2 \DM(X_s)$$
for the two-dimensional nearby cycle functor.
There is a natural transformation (of functors $\Corr_{\A^2}^{\proper, \sm} \r \Fun(\Delta^1, \PrLgr)$)
$$\Upsilon \r \Upsilon_1 \circ \Upsilon_2.$$
\xlemm

\pf
For $k = 1, 2$, let $\Gamma_k$ be the graph of $X_\eta \r \Gm^2 \stackrel{p_k} \r \Gm$, and $\Gamma$ the graph of $X_\eta \r \Gm^2$.
Then 
$$\DM(X_\eta) \stackrel{\Gamma_{1?}} \r \DTM(\Gm) \t \DM(X_\eta) \stackrel{\id \t \Gamma_{2?}} \r \DTM(\Gm)^2 \t \DM(X_\eta)$$
agrees with $\Gamma_?$.
Indeed their left adjoints agree: $\Gamma^* = \Gamma_1^* (\Gamma_2 \x \id)^*$.
From here, the proof of \cite[Lemma~9.4.4(2)]{AcharRiche:Central} carries over.
\xpf

\subsection{Kan extension to ind-schemes and placid prestacks} \label{sect--KanInd}

Taking advantage of the \ii-categorical setup of $\Upsilon$, we can now quickly extend its reach to more general objects of relevance in geometric representation theory.
As a preliminary observation note from \thref{NilpQ explained} that $\NilpQ$ is compactly generated, and therefore a dualizable object in $\PrLgr$ \cite[Propositions~I.7.3.2 and 9.4.4]{GaitsgoryRozenblyum:StudyI}.
Thus, tensoring with it preserves not only colimits but also limits: $\NilpQ (\lim C_i) = \lim \NilpQ (C_i)$.

We can apply a Kan extension to obtain a unipotent nearby cycles functor for motives on placid prestacks in two steps:
$$\xymatrix{
(\AffSch_{\A^1}^{\sm})^\opp \ar[d] \ar[r] & \Corr_{\A^1}^{\sm, \proper} \ar[r]^{\Upsilon} & \Fun(\Delta^1, \PrLgr) \\
(\AffSch_{\A^1}^{\kappa, \placid})^\opp \ar[d] \ar@{.>}[urr]_{\Lan \Upsilon} \\
(\PreStk_{\A^1}^{\placid})^\opp \ar@{.>}[uurr]_{\Ran (\Lan \Upsilon)}}.$$
We first form the left Kan extension (denoted by $\Lan$) to the category of $\kappa$-small (for some fixed regular cardinal $\kappa$) placid affine $\A^1_S$-schemes (i.e., $X = \lim X_i$, with smooth affine transition maps $X_i \r X_j$), so that $\DM(X) = \colim \DM(X_i)$, with transition functors given by *-pullback.
We then form the right Kan extension to the category of placid prestacks, which is the free cocompletion of $\AffSch_{\A^1}^{\kappa, \placid}$.
If $X$ is a finite type $\A^1_S$-scheme and $G / \A^1$ is a pro-smooth algebraic group, then the quotient prestack $(X/G) / \A^1$ is placid and we have the following description of equivariant motives (cf.~around \cite[Lemma~2.2.7]{RicharzScholbach:Intersection}):
\begin{align}
\DM(X/G) & = \lim \left (\DM(X) \stackrel[a^!]{p^!}\rightrightarrows \DM(G \x X)  \substack{\rightarrow\\[-.9em] \rightarrow \\[-.9em] \rightarrow} \dots \right ) \nonumber \\ 
 & = \lim \left (\DM(X) \stackrel[a^*]{p^*}\rightrightarrows \DM(G \x X) \substack{\rightarrow\\[-.9em] \rightarrow \\[-.9em] \rightarrow} \dots \right ). \label{DM equivariant}
\end{align}

Evaluating the functor $\Upsilon := \Ran (\Lan \Upsilon)$ at such quotient prestacks gives the following.

\coro
\thlabel{Upsilon equivariant}
Let $X$ be a finite type $\A^1_S$-scheme, acted upon by a pro-smooth algebraic group $G / \A^1$.
Then there is a natural \emph{equivariant unipotent nearby cycles functor}, also denoted $\Upsilon$, such that the following diagram commutes:
$$\xymatrix{
\DM(X_\eta / G_\eta) \ar[d]^{u^!} \ar[r]^-{\Upsilon} & \NilpQ \DM(X_s / G_s) \ar[d]^{u^!} \\
\DM(X_\eta) \ar[r]^-{\Upsilon} & \NilpQ \DM(X_s).
}$$
\xcoro

In a similar manner, it is possible to do a Kan extension along proper maps (see also \cite[§2.4]{RicharzScholbach:Intersection}):
For example, if $X = \colim_{i \in I} X_i$ is an ind-scheme over $\A^1$, using that $\Upsilon$ commutes with pushforwards along the closed immersions $X_i \r X_j$, we get a natural functor
$$\DM(X_\eta) = \colim \DM(X_{i\eta}) \stackrel{\colim \Upsilon_{X_i}} \r \colim \NilpQ \DM(X_{is}) = \NilpQ \DM(X_s).$$ 
This functor is compatible with the one given by \thref{Upsilon.all.in.one} under the canonical insertion maps $\DM(X_{i\eta}) \r \DM(X_\eta)$.
More formally, the above functor is the evaluation (at $X$) of the left Kan extension of $\Sch^{\proper} \r \Corr^{\sm, \proper}_{\A^1} \stackrel \Upsilon \r \Fun(\Delta^1, \PrLgr)$ along the inclusion $\Sch_{\A^1}^{\proper} \r \IndSch^{\proper}$ into the category of ind-schemes.

\section{Geometry of Beilinson--Drinfeld Grassmannians}

Recall that we fixed a connected base scheme $S$ which is smooth of finite type over a Dedekind ring or a field. 
For an $S$-scheme of finite type $X$, when working with t-structures on $\DTM(X)$ we must further assume that $S$ satisfies the Beilinson--Soul\'e vanishing condition (cf.~\refsect{MTM}). 
For reduced motives this last assumption is not necessary. We note that the fact that Tate motives are preserved under convolution (\thref{thm-convDTM}) also implicitly uses of the existence of a t-structure, cf.~\cite[Theorem 3.17]{RicharzScholbach:Motivic}.

\subsection{Affine flag varieties}
\subsubsection{Group-theoretic notation}
Let $G$ be a split Chevalley group scheme over $\Z$. Fix a maximal torus and Borel subgroup $T \subset B \subset G$. When working over a general base $S$ we will usually write $G$ instead of $G \times_{\Z} S$. Let $X_*(T)$ and $X^*(T)$ be the groups of characters and cocharacters. Let $X_*(T)^+ \subset X_*(T)$ be the submonoid of dominant cocharacters and let $R_+$ be the positive roots determined by $B$. 
This determines a standard apartment $\mathscr{A} = X_*(T) \otimes \mathbb{R}$, and a special facet $\mathbf{f}_0$ containing the origin.
We moreover fix a standard alcove $\mathbf{a}_0$, containing \(\mathbf{f}_0\) in its closure (see \cite[\S 4.1]{RicharzScholbach:Intersection} for more details). 

\rema\thlabel{rmk-choice of borel}
Although the Borel \(B\) containing \(T\) yields a preferred choice of alcove \(\mathbf{a}_0\), we do not make this specific choice.
Indeed, in the proof of \thref{t-exactness general facets}, it will be useful to consider different Borels, while keeping the alcove \(\mathbf{a}_0\) fixed. For later use, we let $B_{\mathbf{a}_0} \subset T$ be the Borel whose Weyl chamber contains $\mathbf{a}_0$.
\xrema

The finite Weyl group and Iwahori--Weyl group are, respectively, $W_0 = N_G(T)/T$ and $W = X_*(T) \rtimes W_0.$ 
To avoid potential confusion, we denote the image of a cocharacter \(\mu\) under the inclusion \(X_*(T)\subset W\) by \(t(\mu)\).
There is an isomorphism $W \cong N_G(T)(\Z (\!(t)\!))/T(\Z[ \![t]\!])$ which extends the isomorphism $X_*(T) \cong T(\Z(\!(t)\!))/T(\Z[\![t]\!])$, $\lambda \mapsto \lambda(t)$.  
Let  $W_{\aff}$ be the subgroup of $W$ generated by the reflections in the walls of $\mathbf{a}_0$. This is a Coxeter group. 
If $\Omega \subset W$ is the stabilizer of $\mathbf{a}_0$ we have $W = W_{\aff} \rtimes \Omega$. 
Declaring elements of $\Omega$ to have length $0$ makes $W$ a quasi-Coxeter group. 
We denote the extension of the partial Bruhat order on $W_{\aff}$ to $W$ by $\leq$ and the length function by $l$.

\subsubsection{Loop groups}
\label{sect--loop-groups} Suppose that $S$ is affine, so that we have the power series scheme $S[\![t]\!] = \Spec \mathcal{O}_S[\![t]\!]$. We note that when $S$ is not affine, the loop groups and and affine flag varieties defined below can still be obtained by base change from $\Spec \Z$, cf.~\cite[\S 4.4]{RicharzScholbach:Intersection}, so this is not a serious assumption.
For a group scheme $H$ over $S[\![t]\!]$ we have the loop group (resp. positive loop group) functor $LH \colon \AffSch_S^{\text{op}} \rightarrow \Set$, $LH(R) = H(R (\!(t)\!))$ (resp. $L^+H(R) = H(R[ \![t]\!])$). 
Applying this to $G$ viewed as an $S[\![t]\!]$-scheme by base change, the functor $LG$ is an ind-affine $S$-scheme and $L^+G$ is an $S$-scheme (typically not of finite type). 
For each facet $\mathbf{f} \subset \mathscr{A}$ there is a parahoric Bruhat--Tits group scheme $G_{\mathbf{f}}$ constructed in \cite[\S 4.2.2]{PappasZhu:Kottwitz} using $\Z[\![t]\!]$ as the base ring.  For general affine $S$ this is obtained by base change from $\Z[\![t]\!]$. The scheme $G_{\mathbf{f}}$ is smooth, affine, and has geometrically connected fibers.
Then $L^+G_{\mathbf{f}} \subset LG$ is a pro-algebraic $\Z$-group scheme with geometrically connected fibers \cite[Lemma 4.2.4]{RicharzScholbach:Intersection}.
We have $L^+G_{\mathbf{f}_0} = L^+G$, and $L^+G_{\mathbf{a}_0}$ is the preimage of $B_{\mathbf{a}_0}$ under the projection $L^+G \rightarrow G$, $t \mapsto 0$. 
We denote $\mathcal{I} = L^+G_{\mathbf{a}_0}$ and call it the Iwahori group. 

\subsubsection{The partial affine flag variety}

\defi \thlabel{defi--partial}
The \emph{partial affine flag variety} $\Fl_{\mathbf{f}}$ associated to a facet $\mathbf{f}$ is the \'etale-sheafification $(LG/L^+G_{\mathbf{f}})_{\et}$ of the presheaf quotient of $LG$ by $L^+G_{\mathbf{f}}$.
\xdefi

It is well-known that $\Fl_{\mathbf{f}}$ is represented by an ind-projective $S$-scheme (see \cite[\S 4.2]{RicharzScholbach:Intersection}). Since $G$ is split, it follows from work of Faltings \cite[Def. 5 ff.]{Faltings:Loops} that the quotient map $LG \rightarrow \Fl_{\mathbf{f}}$ admits sections Zariski-locally. In particular, $\Fl_{\mathbf{f}} = (LG/L^+G_{\mathbf{f}})_{\Zar}$. This fact is necessary to define convolution (see \S \ref{sect--convolution}) on $\DM(L^+G_{\mathbf{f}} \backslash \Fl_{\mathbf{f}})$ since \(\DM\) satisfies Nisnevich but not \'etale descent. It will be convenient to work with \'etale-sheafifications so that we have access to the general theory of \cite[\S 3]{HainesRicharz:TestFunctionsWeil}.

We denote $\Fl = \Fl_{\mathbf{a}_0}$ and $\Gr = \Fl_{\mathbf{f}_0}$, and we call these the (full) affine flag variety and the affine Grassmannian. Sometimes we use a subscript to emphasize the group, e.g. $\Gr_G$. Let $G_{\text{der}}$ be the derived subgroup of $G$ and let $G_{\text{sc}}$ be the simply connected cover of $G_{\text{der}}$.
Since $S$ is connected, the connected components of $\Fl_{\mathbf{f}}$ are indexed by $\pi_1(G) = X_*(T)/X_*(T_{\text{sc}})$ where $T_{\text{sc}}$ is the preimage of $T \cap G_{\text{der}}$ in $G_{\text{sc}}$. 

\subsubsection{Affine Schubert schemes} 
Let $W_{\mathbf{f}} := (N_G(T)(\Z (\!(t)\!)) \cap L^+G_{\mathbf{f}}(\Z))/T(\Z[ \![t]\!]) \subseteq W$ for a fixed facet $\mathbf{f}$.
For another facet $\mathbf{f}' \subset \mathscr{A}$, there is a natural length function and partial order on the double cosets  $W_{\mathbf{f}'} \backslash W / W_{\mathbf{f}}$ also denoted by $l$ and $\leq$, cf. \cite[Lemma 1.6 ff.]{Richarz:Schubert}. 
The Bruhat decomposition implies that left $L^+\mathcal{G}_{\mathbf{f}'}$-orbits in $\Fl_{\mathbf{f}}$ are indexed by 
 $W_{\mathbf{f}'} \backslash W / W_{\mathbf{f}}$. 
 For example, $W_{\mathbf{f}_0} = W_0$ and $W_{\mathbf{a}_0}$ is the trivial group. 
 Thus, $\mathcal{I}$-orbits in $\Fl$ are indexed by $W$ and $L^+G$-orbits in $\Gr$ are indexed by $X_*(T)/W_0 \cong X_*(T)^+$.
 These observations motivate the following definition over any base scheme $S$ (see also \cite[Definition 4.4.1]{RicharzScholbach:Intersection}).

\defi
Let $w \in W_{\mathbf{f}'} \backslash W / W_{\mathbf{f}}$. The \emph{affine Schubert scheme} $\Fl_w(\mathbf{f}', \mathbf{f})$ is the scheme-theoretic image of the map
$$L^+G_{\mathbf{f}'} \rightarrow \Fl_{\mathbf{f}}, \quad g \mapsto g \cdot \dot{w} \cdot e$$ where $\dot{w} \in LG(S)$ is any representative of $w$ under the map $LG(\Z) \rightarrow LG(S)$ and $e \in \Fl_{\mathbf{f}}(S)$ is the basepoint. The $L^+G_{\mathbf{f}'}$-orbit of $w$ in $\Fl_{\mathbf{f}}$, denoted $\Fl_w^\circ(\mathbf{f}', \mathbf{f})$, is the \'etale sheaf-theoretic image of the above map.
\xdefi
The following proposition summarizes the properties of affine Schubert schemes.

\prop \thlabel{prop--summary}
The affine Schubert schemes satisfy the following properties.
\begin{enumerate}
\item The affine Schubert scheme $\Fl_w(\mathbf{f}', \mathbf{f})$ is a reduced, projective $S$-scheme which is moreover $L^+G_{\mathbf{f}'}$-stable.
\item The $L^+G_{\mathbf{f}'}$-orbit $\Fl_w^\circ(\mathbf{f}', \mathbf{f}) \subset \Fl_w(\mathbf{f}', \mathbf{f})$ is an open $S$-subscheme which is smooth, and fiberwise geometrically connected and dense.
\item \label{item--summary3} We have $\Fl_{w}^\circ(\mathbf{a}_0, \mathbf{f}) \cong \A^{l(w)}_S$ where $l$ is the length function on $W/W_{\mathbf{f}}$.
\item If $v \leq w$ with respect to the partial order on $W_{\mathbf{f}'} \backslash W / W_{\mathbf{f}}$ then there is a closed immersion $\Fl_v(\mathbf{f}', \mathbf{f}) \rightarrow \Fl_w(\mathbf{f}', \mathbf{f})$. As sets there is a decomposition
$$\Fl_w(\mathbf{f}', \mathbf{f}) = \bigsqcup_{v \leq w} \Fl_v^\circ(\mathbf{f}', \mathbf{f}).$$
Moreover, the reduced locus of $\Fl_\mathbf{f}$ is $(\Fl_\mathbf{f})_{\redu} = \colim_{w} \Fl_w(\mathbf{f}', \mathbf{f})$. 
\item If $S=k$ is a field, $\Fl_w(\mathbf{f}', \mathbf{f})$ agrees with the usual affine Schubert varieties defined in the literature, e.g. \cite[Def.~8.3]{PappasRapoport:Twisted}, \cite[Def.~2.5]{Richarz:Schubert}.
\item The formation of $\Fl_w(\mathbf{f}', \mathbf{f})$ commutes with base change along a map $S' \rightarrow S$ up to a possible nilpotent thickening.
\end{enumerate}
\xprop

\pf
All of these facts can be found in \cite[\S 4.3, \S 4.4]{RicharzScholbach:Intersection}. The proofs combine standard techniques from the case where $S$ is a field as in \cite{Richarz:Schubert, Richarz:AffGrass} with the Bruhat--Tits group schemes $G_{\mathbf{f}}$ from \cite{PappasZhu:Kottwitz}.
\xpf

We introduce the following standard notation for the most important affine Schubert schemes in this paper.
\begin{align*}
\Gr^{\leq \mu} &= \Fl_{\mu}(\mathbf{f}_0, \mathbf{f}_0), \quad  \Gr^{\mu} = \Fl_{\mu}^\circ (\mathbf{f}_0, \mathbf{f}_0), \quad \mu \in X_*(T)^+, \\
\Fl^{\leq w} &= \Fl_{w}(\mathbf{a}_0, \mathbf{a}_0), \quad \Fl^{w} = \Fl_{w}^\circ(\mathbf{a}_0, \mathbf{a}_0),  \quad w \in W.
\end{align*}

\subsubsection{The admissible locus} 
Inside $\Fl_{\mathbf{f}}$ there is another important family of subschemes. 
\defi \thlabel{defi-admissible}
Let $\mu \in X_*(T)^+$ and let $W_0(\mu) \subset X_*(T)$ be the $W_0$-orbit of $\mu$. The $\mu$\emph{-admissible set} is $$\adm_{\{\mu\}} = \{w \in W \: : \: w \leq t(\lambda) \text{ for some } \lambda \in W_0(\mu)\}.$$ For a facet $\mathbf{f}$ the $\mu$\emph{-admissible set relative to }$\mathbf{f}$ \cite[Eq.~(6.17)]{HainesRicharz:TestFunctions} is 
$$\adm_{\{\mu\}}^{\mathbf{f}} = W_{\mathbf{f}} \backslash \adm_{\{\mu\}} / W_{\mathbf{f}} \subset W_{\mathbf{f}} \backslash W / W_{\mathbf{f}},$$ i.e., $\adm_{\{\mu\}}^{\mathbf{f}}$ is the image of $\adm_{\{\mu\}}$ under the quotient $W \r W_{\mathbf{f}} \backslash W / W_{\mathbf{f}}$.
The $\mu$\emph{-admissible locus }$\mathcal{A}^{\mathbf{f}}(\mu)$ \emph{relative to} $\mathbf{f}$
is the scheme-theoretic image of $\bigsqcup_{\lambda \in W_0(\mu)} \Fl_{t(\lambda)}(\mathbf{f}, \mathbf{f}) \rightarrow \Fl_{\mathbf{f}}.$ As a set, $\mathcal{A}^{\mathbf{f}}(\mu)$ is the union of the $\Fl_w(\mathbf{f}, \mathbf{f})$ for those $w \in \adm_{\{\mu\}}^{\mathbf{f}}$. 
\xdefi
The admissible loci arise as reduced special fibers of global Schubert schemes in  the Beilinson--Drinfeld Grassmannian, which we define next.

\subsection{Beilinson--Drinfeld Grassmannians}\label{sect--BD-Grassmannians}

\subsubsection{Bruhat--Tits group schemes} \label{sect--BT}
Throughout this paper we let $C = \A^1_S$. For a prestack $Y \rightarrow C$, let $Y_s \rightarrow S$ be the fiber over the origin and let $Y_{\eta} \rightarrow \GmS$ be the fiber over the open complement.
We define the constant group scheme $G_C = G \times_{\Z} C$ over $C$.
Regarding $T \subset B_{\mathbf{a}_0}$ as closed subschemes of $G_C$ supported over the origin, we may form the dilatation $\mathcal{G} = \text{Bl}_{B_{\mathbf{a}_0}}^G G_C$ of $G_C$ in $B_{\mathbf{a}_0}$ along $G = (G_C)_s$ in the sense of \cite[Definition 2.1]{MayeuxRicharzRomagny:Neron}. 
This is a smooth affine $C$-group scheme with geometrically connected fibers \cite[Theorem 3.2]{MayeuxRicharzRomagny:Neron}.

Suppose $S$ is affine, and let $\mathcal{O}_s$ be the formal neighborhood of the origin in $C$ (so $\mathcal{O}_{s} \cong S[\![t]\!]$). 
Then
$$\restr{\mathcal{G}}{\mathcal{O}_{s}} = G_{\mathbf{a}_0}, \quad  \mathcal{G}_\eta  = G \times_\Z C_\eta, \quad L^+(\restr{\mathcal{G}}{\mathcal{O}_{s}}) = \mathcal{I}.$$
Thus, if $k$ is any field, $\mathcal{G}(k[\![t]\!]) \subset G(k[\![t]\!])$ consists of those elements whose reduction modulo $t$ lies in $B_{\mathbf{a}_0}(k)$ \cite[Example 3.3]{MayeuxRicharzRomagny:Neron}. 
The fiber of $\mathcal{G}$ along $\Z[t] \rightarrow k[t]$ is precisely the group scheme used in the construction of central sheaves in \cite{AcharRiche:Central, Zhu:Coherence} (in the split case), cf. \cite[Remark 2.2.9]{AcharRiche:Central}. 
The fiber along $\Z[t] \rightarrow \Z_p$, $t \mapsto p$ for $p$ a prime is a parahoric group scheme in mixed characteristic, but we will not need this fact. 

More generally, for each facet $\mathbf{f}$ there is a smooth affine group scheme $\mathcal{G}_{\mathbf{f}}$ over $C$ such that $(\mathcal{G}_{\mathbf{f}})_\eta = G \times_{\Z} C_\eta$ and $\restr{\mathcal{G}_{\mathbf{f}}}{\mathcal{O}_{s}}=G_{\mathbf{f}}$. 
The group $\mathcal{G}_{\mathbf{f}}$ is characterized uniquely as a particular family of parahorics \cite[Theorem 1.3]{Lorenco:Entiers}. 
The primary objective in \emph{loc.~cit.} is to overcome obstacles related to wild ramification, and $\mathcal{G}_{\mathbf{f}}$ can also be constructed by the methods in \cite{BruhatTits:Groups2} and \cite{PappasZhu:Kottwitz} (these works restrict to local fields but in fact they work over $\Z$ when the group is split).
According to our conventions, $\mathcal{G}_{\mathbf{f}_0} = G_C$ and $\mathcal{G}_{\mathbf{a}_0} = \mathcal{G}$.

\subsubsection{Basic properties of Beilinson--Drinfeld Grassmannians} Suppose that $S$ is affine; again, the Beilinson--Drinfeld loop groups and Grassmannians below can be defined for general $S$ by base change from $\Spec \Z$, so this is not a serious assumption.
Let $I$ be a finite nonempty set and let $X = C^I \times_S C \cong C^{I \sqcup \{*\}}$. We have the divisor $D \subset X$ corresponding to the ideal generated by $\prod_{i \in I}(t-x_i)$, where $t$ is the coordinate function on $C$ (the right factor) and the $x_i$ are the coordinates on $C^I$. For $R \in \AffSch_{C^I}^{\text{op}}$ let $X_R$ and $D_R$ be the base changes to $R$. 
We denote by $R[\![D]\!] \in \AffSch_S^{\text{op}}$ the ring of functions on the formal neighborhood $\hat D_R$ of $D_R$ in $X_R$. Let $\hat{D}_R^\circ$ be the complement of $D_R$ in $\Spec R[\![D]\!]$. This is an affine $S$-scheme; let $R(\!(D)\!)$ be its ring of functions (see \cite[\S 3.1.1]{HainesRicharz:TestFunctionsWeil} for more details). In the following definition, both $\hat{D}_R$ and $\hat{D}_R^\circ$ are viewed as $X$-schemes.

\defi
\thlabel{BD.Gr.I}
Let $I$ be a finite nonempty set. Let $H$ be a smooth affine $C$-group scheme, which we view as an $X$-group scheme via pullback along the projection $X = C^I \times_S C \rightarrow C$. The \emph{loop group relative to} $I$ (resp. \emph{positive loop group relative to} $I$) is the functor $\AffSch_{C^I}^{\text{op}} \rightarrow \Set$,
$$L_IH(R) = H(R(\!(D)\!)), \quad (\text{resp. } L^+_IH(R) = H(R[\![D]\!])).$$
The \emph{Beilinson--Drinfeld Grassmannian of} $H$ \emph{relative to} $I$ is the \'etale quotient $$\Gr_{H, I} = (L_IH/L_I^+H)_{\et}.$$
\xdefi

By \cite[Lemma 3.2]{HainesRicharz:TestFunctionsWeil}, $L_IH$ is an ind-affine $C^I$-group-scheme and $L^+_IH$ is an affine, faithfully flat, pro-smooth $C^I$-group-scheme. When working with equivariant motives it is convenient to work with groups of finite type. Toward this end, for $m \in \mathbb{Z}_{\geq 0}$ we let $R[\![D^{(m)}]\!]$ be the ring of functions on the $m$-th infinitesimal neighborhood of $D_R$ in $X_R$. Let
$$L_{I}^{+, (m)}H(R) = H(R[\![D^{(m)}]\!]).$$
Then $L^+_IH = \lim_m L_{I}^{+, (m)}H$ and each $ L_{I}^{+, (m)}H$ is a smooth affine $C$-group scheme. The kernel of $L^+_IH \rightarrow  L_{I}^{+, (0)}H$ is split pro-unipotent by \cite[Proposition A.9]{RicharzScholbach:Intersection}.
 
If $\mathcal{G}_{\mathbf{f}}$ is the Bruhat--Tits group scheme over $C$ associated to the facet $\mathbf{f}$ as in \refsect{BT} we have the Beilinson--Drinfeld Grassmannian $\Gr_{\mathcal{G}_{\mathbf{f}}, I} = (L_I\mathcal{G}_{\mathbf{f}}/ L^+_I \mathcal{G}_{\mathbf{f}})_{\et}$. If $I = \{*\}$ is a singleton we omit it from the notation, so $\Gr_{\mathcal{G}_{\mathbf{f}}} = \Gr_{\mathcal{G}_{\mathbf{f}}, \{*\}}$.

It follows from the definitions that $LG_C = LG \times_S C$ and  $L^+G_C = L^+G \times_S C$. More generally, $\Gr_{G_C, I}$ is the Beilinson--Drinfeld Grassmannian denoted by $\Gr_{G,I}$ in \cite{CassvdHScholbach:MotivicSatake}, so it is ind-projective. The proof of ind-projectivity for general $\mathcal{G}_{\mathbf{f}}$ follows from work of Pappas--Zhu \cite{PappasZhu:Kottwitz}, Haines--Richarz \cite{HainesRicharz:TestFunctionsWeil}, and Achar--Riche \cite{AcharRiche:Central} (see also the work of Lourenço \cite[Theorem 5.1.3]{Lorenco:Entiers} which proves ind-projectivity when $I = \{*\}$);
we recall the main points in \thref{pro-ind-proj}.

We now describe the fibers of $\Gr_{\mathcal{G}_{\mathbf{f}}, I}$. Let $\phi \colon I \twoheadrightarrow J$ be a surjection of non-empty finite sets and let $j \in J \sqcup \{*\}$.
We set $$C^{\phi, j} = \{(x_i) \in C^I \: : \: x_i = x_{i'} \text{ iff } \phi(i) = \phi(i'), \text{ and } x_i = 0 \text{ iff } \phi(i) = j\}.$$ Here exactly $|\phi^{-1}(j)|$ coordinates in $C^{\phi, j}$ are zero if $j \in J$, and every coordinate in $C^{\phi, *}$ is nonzero.
By letting $\phi$ and $j$ vary, we have a stratification of $C^I$ by the locally closed subschemes $C^{\phi, j}$.
Then
$$
\Gr_{\mathcal{G}_{\mathbf{f}}, I} \times_{C^I} C^{\phi, j} = \begin{cases}
  \Gr^{|J|-1} \times \Fl_{\mathbf{f}} \times C^{\phi, j}  & j \in J \\
  \Gr^{|J|} \times C^{\phi, j} & j = *.
\end{cases}
$$
This follows from the fact that $L^+_I \mathcal{G}_{\mathbf{f}} \times_{C^I} C^{\phi, j} = (L^+G)^{|J|-1} \times L^+G_{\mathbf{f}} \times C^{\phi, j}$ if $j \in J$ and  $L^+_I \mathcal{G}_{\mathbf{f}} \times_{C^I} C^{\phi, j} = (L^+G)^{|J|} \times C^{\phi, j}$ if $j= *$, and a similar description of the fibers of $L_I \mathcal{G}_{\mathbf{f}}$. In particular,
$$(\Gr_{\mathcal{G}_{\mathbf{f}}})_\eta = \Gr \times C_\eta, \quad (\Gr_{\mathcal{G}_{\mathbf{f}}})_s = \Fl_{\mathbf{f}} \times C_s \cong \Fl_{\mathbf{f}}.$$

In order to prove properties of the central functor we will need to consider the restriction of $\Gr_{\mathcal{G}_{\mathbf{f}}, I}$ to certain copies of $C$. The notation in the following definition is borrowed from \cite{Zhu:Coherence, AcharRiche:Central}. 

\defi \thlabel{defi--BDloop}
Let $H$ be a smooth affine $C$-group scheme, and let $C \cong C \times \{0\} \subset C^2$ be the closed subscheme obtained by setting $x_2=0$. We define the following three functors $\AffSch_C^{\text{op}} \rightarrow \Set$,
$$LH^{\BD} = \restr{(L_{\{1,2\}}H)}{C \times \{0\}},  \quad L^+H^{\BD} = \restr{(L_{\{1,2\}}^+H)}{C \times \{0\}},$$ and
$$\Gr_{H}^\BD = (LH^{\BD}/L^+H^{\BD})_{\et} = \restr{\Gr_{H, \{1,2\}}}{C \times \{0\}}.$$
\xdefi

All representability results which will be proved in \thref{pro-ind-proj} apply to $\Gr_{H}^\BD$ by base change.
If $H = \mathcal{G}_{\mathbf{f}}$ for a facet  $\mathbf{f}$ then by the above there are canonical isomorphisms
\begin{align*}
&(L^+\mathcal{G}_{\mathbf{f}}^{\BD})_{\eta} = L^+G \times L^+G_{\mathbf{f}} \times C_\eta  & (\Gr_{\mathcal{G}_{\mathbf{f}}}^{\BD})_{\eta} = \Gr \times \Fl_{\mathbf{f}} \times C_\eta,  \\
&(L^+\mathcal{G}_{\mathbf{f}}^{\BD})_{s} = L^+G_{\mathbf{f}} & (\Gr_{\mathcal{G}_{\mathbf{f}}}^{\BD})_{s} = \Fl_{\mathbf{f}}.
\end{align*}
Restriction from the disc $\hat{D}_{R}$ to the disc defined by the ideal $(t-x_1)$ (resp. $(t-x_2)$) defines a map $$L^+\mathcal{G}_{\mathbf{f}}^{\BD} \rightarrow L^+\mathcal{G}_{\mathbf{f}}, \quad  \text{(resp. } L^+\mathcal{G}_{\mathbf{f}}^{\BD} \rightarrow L^+G_{\mathbf{f}} \times C). $$

\subsubsection{The Beilinson--Drinfeld convolution Grassmannians}

We now recall the moduli interpretation of $\Gr_{\mathcal{G}_{\mathbf{f}}, I}$.
By \cite[Lemma 3.4]{HainesRicharz:TestFunctionsWeil}, $\Gr_{\mathcal{G}_{\mathbf{f}}}(R)$ is the set of isomorphism classes of pairs $(\mathcal{E}, \alpha)$ where $\mathcal{E}$ is a $\mathcal{G}_{\mathbf{f}}$-torsor on $\hat{D}_R$ and $\alpha \colon \restr{\mathcal{E}}{\hat{D}^\circ_R} \xrightarrow{\sim} \restr{\mathcal{E}^0}{\hat{D}^\circ_R}$ is a trivialization. Here $\mathcal{E}^0$ is a fixed trivial $\mathcal{G}_{\mathbf{f}}$-torsor. Similarly, $L_I \mathcal{G}_{\mathbf{f}}(R)$ is the set of isomorphism classes of triples $(\mathcal{E}, \alpha, \beta)$ where $(\mathcal{E}, \alpha)$ is as before and $\beta \colon \restr{\mathcal{E}^0}{\hat{D}_R} \xrightarrow{\sim} \restr{\mathcal{E}}{\hat{D}_R}$.

The moduli interpretation of $\Gr_{\mathcal{G}_{\mathbf{f}}, I}$ allows us to construct a convolution Grassmannian as follows, generalizing \cite[\S 3.5.2]{AcharRiche:Central} (where $|I| = 2$).  If $|I| > 1$, choose a distinguished element $i_0 \in I$. For $R \in \AffSch_{C^I}^{\text{op}}$, let $\hat{D}^\circ_{R, I-\{i_0\}}$ be the complement in $\hat{D}_R$ of the divisor defined by the ideal generated by $\prod_{i \in I \setminus \{i_0\}}(t-x_i)$. Thus, $\hat{D}^\circ_{R}$ is obtained from $\hat{D}^\circ_{R, I-\{i_0\}}$ by further removing the divisor defined by the ideal $(t-x_{i_0})$. Similarly, let $\hat{D}^\circ_{R, \{i_0\}}$ be the complement in $\hat{D}_R$ of the divisor defined by the ideal $(t-x_{i_0})$.

\defi \thlabel{def--i_0BDconv}
The \emph{Beilinson--Drinfeld convolution Grassmannian relative to} $i_0 \in I$ \emph{and} $\mathbf{f}$ is the  functor $\AffSch_{C^I}^{\text{op}} \rightarrow \Set$ defined by $\widetilde{\Gr}_{\mathcal{G}_{\mathbf{f}}, I}(R) =\{(\mathcal{E}^1, \mathcal{E}^2, \alpha_1, \alpha_2)\}$ where the $\mathcal{E}^i$ are $\mathcal{G}_{\mathbf{f}}$-torsors on $\hat{D}_R$ and 
$$\alpha_1 \colon \restr{\mathcal{E}^1}{\hat{D}^\circ_{R, I-\{i_0\}}} \xrightarrow{\sim} \restr{\mathcal{E}^0}{\hat{D}^\circ_{R, I-\{i_0\}}}, \quad \alpha_2 \colon \restr{\mathcal{E}^2}{\hat{D}^\circ_{R, \{i_0\}}} \xrightarrow{\sim} \restr{\mathcal{E}^1}{\hat{D}^\circ_{R, \{i_0\}}}.$$ The \emph{convolution map} $m_{i_0} \colon \widetilde{\Gr}_{\mathcal{G}_{\mathbf{f}}, I} \rightarrow \Gr_{\mathcal{G}_{\mathbf{f}}, I}$ is defined on points by $$\{(\mathcal{E}^1, \mathcal{E}^2, \alpha_1, \alpha_2)\} \mapsto \{(\mathcal{E}^2, \restr{\alpha_1}{\hat{D}^\circ_{R}} \circ\restr{\alpha_2}{\hat{D}^\circ_{R}})\}.$$
\xdefi

The convolution Grassmannian $\widetilde{\Gr}_{\mathcal{G}_{\mathbf{f}}, I}$ will be used to prove the ind-projectivity of $\Gr_{\mathcal{G}_{\mathbf{f}}, I}$. Toward this end, we note that as a functor on $\AffSch_{C^I}^{\text{op}}$ we have $(\Gr_{\mathcal{G}_{\mathbf{f}}, I-\{i_0\}} \times C)(R) = \{(\mathcal{E}, \alpha)\}$ where $\mathcal{E}$ is a $\mathcal{G}_{\mathbf{f}}$-torsor on $\hat{D}_R$ and $\alpha \colon \restr{\mathcal{E}}{\hat{D}^\circ_{R, I-\{i_0\}}} \xrightarrow{\sim} \restr{\mathcal{E}^0}{\hat{D}^\circ_{R, I-\{i_0\}}}$ is a trivialization.  Similarly, $(\Gr_{\mathcal{G}_{\mathbf{f}},\{i_0\}} \times C^{I - \{i_0\}})(R) = \{(\mathcal{E}, \alpha)\}$ where now $\alpha \colon \restr{\mathcal{E}}{\hat{D}^\circ_{R, \{i_0\}}} \xrightarrow{\sim} \restr{\mathcal{E}^0}{\hat{D}^\circ_{R, \{i_0\}}}$. There is a right $L^+_I \mathcal{G}_{\mathbf{f}}$-torsor $(\Gr_{\mathcal{G}_{\mathbf{f}}, I-\{i_0\}} \times C)^{(\infty)} \rightarrow \Gr_{\mathcal{G}_{\mathbf{f}}, I-\{i_0\}} \times C$ whose $R$-points are given by $\{(\mathcal{E}, \alpha, \beta)\}$ where $\{(\mathcal{E}, \alpha)\} \in (\Gr_{\mathcal{G}_{\mathbf{f}}, I-\{i_0\}} \times C)(R)$ and $\beta \colon \restr{\mathcal{E}}{\hat{D}_R} \xrightarrow{\sim} \restr{\mathcal{E}^0}{\hat{D}_R}$. The group $L^+_I \mathcal{G}_{\mathbf{f}}$ also acts on $\Gr_{\mathcal{G}_{\mathbf{f}},\{i_0\}} \times C^{I - \{i_0\}}$ on the left by changing the trivialization $\alpha \colon \restr{\mathcal{E}}{\hat{D}^\circ_{R, \{i_0\}}} \xrightarrow{\sim} \restr{\mathcal{E}^0}{\hat{D}^\circ_{R, \{i_0\}}}$. These observations lead to the following alternative description of $\widetilde{\Gr}_{\mathcal{G}_{\mathbf{f}}, I}$.

\lemm \thlabel{lemm--bundle-prod}
The convolution Grassmannian $\widetilde{\Gr}_{\mathcal{G}_{\mathbf{f}}, I}$ is isomorphic to the following \'etale quotient with respect to the diagonal action of $L^+_I \mathcal{G}_{\mathbf{f}}$:
$$\widetilde{\Gr}_{\mathcal{G}_{\mathbf{f}}, I} \cong (\Gr_{\mathcal{G}_{\mathbf{f}}, I-\{i_0\}} \times C)^{(\infty)} \times_{C^I}^{L^+_I \mathcal{G}_{\mathbf{f}}}  (\Gr_{\mathcal{G}_{\mathbf{f}},\{i_0\}} \times C^{I - \{i_0\}}).$$ 
\xlemm

\pf
This follows from \cite[Lemma 3.4]{HainesRicharz:TestFunctionsWeil} which says that every $\mathcal{G}_{\mathbf{f}}$-torsor on $\hat{D}_R$ is trivializable after some \'etale cover $R \rightarrow R'$, cf.~ also the proof of \cite[Proposition 2.3.11]{AcharRiche:Central}.
\xpf

\subsubsection{Ind-projectivity} We are now ready to prove that $\Gr_{\mathcal{G}_{\mathbf{f}}, I}$ is ind-projective.

\prop \thlabel{pro-ind-proj}
Let $H$ be a $C$-group scheme arising from the base change along $C = \A^1_S \rightarrow \A^1_\Z$ of a smooth affine group scheme over $\A^1_\Z$ with geometrically connected fibers. Then the Beilinson--Drinfeld Grassmannian $\Gr_{H, I}$ is represented by a separated ind-finite type $C$-scheme. For any facet $\mathbf{f}$,  $\Gr_{\mathcal{G}_{\mathbf{f}}, I}$ is furthermore ind-projective over $C^I$.
\xprop

\pf
Since sheafification commutes with base change we may assume $S = \Spec \Z$. If $G=\GL_n$ and $\mathbf{f} = \mathbf{f}_0$ then $\mathcal{G}_{\mathbf{f}} = \GL_{n,C}$, the general linear group over $C$. 
Ind-projectivity of $\Gr_{\GL_{n,C}, I}$ is proved using Quot schemes in \cite[Lemma 3.8]{HainesRicharz:TestFunctionsWeil}. 
For general $H$, since $\A^1_\Z$ has dimension $2$ then by \cite[Corollary 11.7]{PappasZhu:Kottwitz} we may choose a closed subgroup embedding $H \rightarrow \GL_{n, C}$ for some $n$ such that $(\GL_{n, C}/H)_{\text{fppf}}$ is quasi-affine. 
The induced map $\Gr_{H} \rightarrow \Gr_{\GL_{n,C}}$ is a quasi-compact immersion by \cite[Proposition 3.10]{HainesRicharz:TestFunctionsWeil}, so $\Gr_{H, I}$ is represented by a separated ind-finite-type $C$-scheme \cite[Corollary 3.11]{HainesRicharz:TestFunctionsWeil}. 

For ind-projectivity of $\Gr_{\mathcal{G}_{\mathbf{f}}, I}$ it then suffices to show ind-properness. We start with the case $I=\{*\}$. The argument for ind-properness in \cite[Proposition 6.5]{PappasZhu:Kottwitz} (where the base is a complete discrete valuation ring with perfect residue field) applies over $C$ as well; we sketch the argument for the convenience of the reader.  First, the geometric fibers of $\Gr_{\mathcal{G}_{\mathbf{f}}} \rightarrow C$ are ind-proper because $\Gr$ and $\Fl_{\mathbf{f}}$ are ind-proper \cite{Faltings:Loops, PappasRapoport:Twisted}. 
 For $\mu \in X_*(T)^+$ let $\Gr^{\leq \mu}_{\mathcal{G}_{\mathbf{f}}}$ be the scheme-theoretic image of the map $\Gr^{\leq \mu} \times_S C_\eta \rightarrow \Gr_{\mathcal{G}_{\mathbf{f}}}$. 
 Using that the fibers of $\Gr_{\mathcal{G}_{\mathbf{f}}}$ are ind-proper and that $\Gr_{\mathcal{G}_{\mathbf{f}}} \rightarrow \Gr_{\GL_{n,C}}$ is an immersion with ind-proper target, it can be shown that $\Gr^{\leq \mu}_{\mathcal{G}_{\mathbf{f}}}  \rightarrow C$ is proper, cf. the proof of \cite[Proposition 6.5]{PappasZhu:Kottwitz} or \cite[Theorem 2.19]{Richarz:AffGrass}.
 To conclude that $\Gr_{\mathcal{G}_{\mathbf{f}}}$ is ind-proper it remains to show that the fibers $(\Gr^{\leq \mu}_{\mathcal{G}_{\mathbf{f}}})_{s}$ cover $(\Gr_{\mathcal{G}_{\mathbf{f}}})_{s}$. 
 This can be checked on reduced loci of the geometric fibers. 
 Then the claim follows from the fact that when $S=\Spec k$ is the spectrum of an algebraically closed field, the reduced locus of $(\Gr^{\leq \mu}_{\mathcal{G}})_{s}$ contains the $\mu$-admissible locus $\mathcal{A}^{\mathbf{f}}(\mu)$ 
 \cite[Lemma 3.12]{Richarz:AffGrass} as in \thref{defi-admissible}.
 We then conclude by noting that $\Fl_{\mathbf{f}}$ is covered by the
 closed Schubert schemes $\Fl_{t(\nu)}(\mathbf{f}, \mathbf{f})$ associated to translation elements  $\nu \in X_*(T) \subset W$.

 For $|I| >1$ we proceed by induction on $|I|$, using the technique in \cite[\S 3.5.2]{AcharRiche:Central} (where $|I| = 2$) to write $\Gr_{\mathcal{G}_{\mathbf{f}}, I}$ as the image of an ind-proper scheme. First, as in \cite[Lemma 2.3.10]{AcharRiche:Central} the ind-projective scheme $\Gr_{\mathcal{G}_{\mathbf{f}}, \{*\}}$ admits an $L^+_{\{*\}} \mathcal{G}_{\mathbf{f}}$-stable presentation by projective $C$-schemes each admitting an $L^+_{\{*\}} \mathcal{G}_{\mathbf{f}}$-equivariant relatively ample line bundle. (The line bundle is pulled back from $\Gr_{\GL_{n,C}}$, using that the latter is constructed from Quot schemes which embed into ordinary Grassmannians.) By \cite[Remark 2.1.10]{AcharRiche:Central}, \thref{lemm--bundle-prod} shows that  $\widetilde{\Gr}_{\mathcal{G}_{\mathbf{f}}, I}$ is represented by an ind-proper $C^I$-scheme, cf.~ also \cite[Proposition 2.3.10]{AcharRiche:Central}. Since the convolution map $m_{i_0}$ in \thref{def--i_0BDconv} is surjective, this implies that $\Gr_{\mathcal{G}_{\mathbf{f}}, I}$ is ind-proper.
 \xpf

\subsubsection{Iterated Grassmannians} \label{sect--iterated}
To define the iterated Grassmannians we follow the setup in \cite[\S 2.3.4]{AcharRiche:Central}. In this subsection we fix $I=\{1,2\}$ when discussing divisors. Let $R \in \AffSch_C^{\text{op}}$. We may view $R$ as an affine scheme over $C^2$ by composing with the inclusion $C \cong C \times \{0\} \subset C^2$, and under this identification we have the disc $\hat{D}_R$ associated to the ideal generated by $(t-x_1)t$. Let $V((t-x_1)t)$ be the divisor associated to this ideal. For $P$ equal to one of the four symbols
$\emptyset$, $\zero$, $\ux$, $\ux \cup \zero$, we define
$$
\hat{D}_P^\circ = \begin{cases}
  \hat{D}_R  & P = \emptyset \\
  \hat{D}_R \setminus V(t) & P= \zero \\
  \hat{D}_R \setminus V(t-x_1) & P = \ux \\
  \hat{D}_R \setminus V((t-x_1)t) = \hat{D}_R^\circ & P =  \ux \cup \zero.
\end{cases}
$$

\defi \thlabel{defi-iterated} Let $\uP = (P_1, \ldots, P_n)$ where $P_i \in \{\emptyset, \zero, \ux, \ux \cup \zero \}$ for all $i$. The iterated Grassmannian associated to $\uP$ is the functor $\Gr_{\mathcal{G}_{\mathbf{f}}}(\uP)  \colon \AffSch_C^{\text{op}} \rightarrow \Set$ such that $$\Gr_{\mathcal{G}_{\mathbf{f}}}(\uP)(R) = \{(\mathcal{E}^i, \alpha_i)\}_{1 \leq i \leq n}$$ where the $\mathcal{E}^i$ are $\mathcal{G}_\bbf$-bundles on $\hat{D}_R$ and the $\alpha_i$ are isomorphisms $\alpha_i \colon \restr{\mathcal{E}^i}{\hat{D}_{P_i}^\circ} \xrightarrow{\sim} \restr{\mathcal{E}^{i-1}}{\hat{D}_{P_i}^\circ}$ (so $\alpha_1$ is a trivialization).
\xdefi

The fibers of the iterated Grassmannians are as follows. 

\lemm \thlabel{lemm--iteratedSingleton} \thlabel{lemm-conv-fibers}
Let $\uP = (P_1, \ldots, P_n)$.
We have a canonical isomorphism
$$
\Gr_{\mathcal{G}_{\mathbf{f}}}(\uP)_{\eta} = \underbrace{LG \x^{L^+G} \cdots \x^{L^+G} \Gr}_{i \textnormal{ factors}} \: \x \: \underbrace{LG \x^{L^+G_{\mathbf{f}}} \cdots \x^{L^+G_{\mathbf{f}}} \Fl_{\mathbf{f}}}_{j \textnormal{ factors}} \times C_\eta
$$
where $i$ is the number of symbols belonging to $\{\ux, \ux \cup \zero\}$ and $j$ is the number of symbols belonging to $\{\zero, \ux \cup \zero\}$. We also have a canonical isomorphism
$$
\Gr_{\mathcal{G}_{\mathbf{f}}}(\uP)_{s} = \underbrace{LG \x^{L^+G_{\mathbf{f}}} \cdots \x^{L^+G_{\mathbf{f}}} \Fl_{\mathbf{f}}}_{k \textnormal{ factors}}
$$
where $k$ is the number of symbols $P_i \neq \emptyset$. Furthermore, if $\uP$ is a singleton we have canonical isomorphisms
$$ \Gr_{\mathcal{G}_{\mathbf{f}}}(\emptyset) = C, \quad \Gr_{\mathcal{G}_{\mathbf{f}}}(\zero) = \Fl_\mathbf{f} \times C, \quad \Gr_{\mathcal{G}_{\mathbf{f}}}(\ux) = \Gr_{\mathcal{G}_{\mathbf{f}}}, \quad  \Gr_{\mathcal{G}_{\mathbf{f}}}(\ux \cup \zero) =  \Gr_{\mathcal{G}_{\mathbf{f}}}^{\BD}.
$$
\xlemm

\pf This is immediate from the moduli description of these functors; see \cite[Lemma 2.3.6, Table 2.3.1 ff.]{AcharRiche:Central}.
\xpf

\thref{lemm--iteratedSingleton} implies that $\Gr_{\mathcal{G}_{\mathbf{f}}}(\uP)$ is representable if $\uP$ is a singleton. For general $\uP$ we consider the functor $\Gr_{\mathcal{G}_{\mathbf{f}}}^{(\infty)}(\uP) \colon \AffSch_C^{\text{op}} \rightarrow \Set$ such that $\Gr_{\mathcal{G}_{\mathbf{f}}}^{(\infty)}(\uP)(R)$ is the set of isomorphism classes $\{(\mathcal{E}^i, \alpha_i), \beta\}_{1 \leq i \leq n}$ where $\{(\mathcal{E}^i, \alpha_i)\}_{1 \leq i \leq n}$ is as in \thref{defi-iterated} and $\beta \colon \mathcal{E}^n \xrightarrow{\sim} \mathcal{E}^0$ is a trivialization of $\mathcal{G}$-bundles on $\hat{D}_R$. 
There is a natural right action of $L^+\mathcal{G}_{\mathbf{f}}^{\BD}$ on $\Gr_{\mathcal{G}_{\mathbf{f}}}^{(\infty)}(\uP)$ which changes $\beta$, cf. \cite[Eqn.~(2.3.9)]{AcharRiche:Central}. 
The map $\Gr_{\mathcal{G}_{\mathbf{f}}}^{(\infty)}(\uP) \rightarrow \Gr_{\mathcal{G}_{\mathbf{f}}}(\uP)$ which forgets $\beta$ is an \'etale-locally trivial $L^+\mathcal{G}_{\mathbf{f}}^{\BD}$-torsor by \cite[Lemma 2.3.9]{AcharRiche:Central}. 
There is also a left action of $L^+\mathcal{G}_{\mathbf{f}}^{\BD}$ on $\Gr_{\mathcal{G}_{\mathbf{f}}}^{(\infty)}(\uP)$ and $\Gr_{\mathcal{G}_{\mathbf{f}}}(\uP)$ which changes $\alpha_1$, cf. \cite[Eqn.~(2.3.5)]{AcharRiche:Central}.

\prop \thlabel{BD-rep}
The iterated Grassmannian $\Gr_{\mathcal{G}_{\mathbf{f}}}(\uP)$ is represented by an ind-proper $C$-scheme. Moreover, it admits a presentation as a colimit of proper $L^+\mathcal{G}_{\mathbf{f}}^{\BD}$-stable subschemes such that the action of $L^+\mathcal{G}_{\mathbf{f}}^{\BD}$ on each subscheme factors through some $L^{+,(m)}\mathcal{G}_{\mathbf{f}}^{\BD}$.
\xprop

\pf
This is proved in \cite[Proposition 2.3.11]{AcharRiche:Central} when $S=\Spec \C$ but the proof applies to the general case since \cite[\S 3]{HainesRicharz:TestFunctionsWeil} works in this generality.
\xpf

\rema \thlabel{rema--jet}
If $Z \subset \Gr_{\mathcal{G}_{\mathbf{f}}}(\uP)$ is a finite-dimensional closed subscheme such that the action of $L^+\mathcal{G}_{\mathbf{f}}^{\BD}$ factors through $L^{+,(m)}\mathcal{G}_{\mathbf{f}}^{\BD}$, then there is a canonical equivalence $\DM(L^+\mathcal{G}_{\mathbf{f}}^{\BD} \backslash Z ) = \DM(L^{+,(m)}\mathcal{G}_{\mathbf{f}}^{\BD} \backslash Z)$. 
This follows from the argument in \cite[Proposition 2.2.11]{RicharzScholbach:Intersection}, using that the kernel of $L^+\mathcal{G}_{\mathbf{f}}^{\BD} \rightarrow  L^{+,(m)}\mathcal{G}_{\mathbf{f}}^{\BD}$ is split pro-unipotent by \cite[Proposition A.9]{RicharzScholbach:Intersection}.  
When $\uP = P$ is a singleton, this equivalence also restricts to $\DTM^{(\anti)}$ by the arguments in \cite[Proposition 3.1.27]{RicharzScholbach:Intersection} once we show that the latter category exists in \thref{theo-WT-BD}\refit{theo-WT-BD1}.
\xrema

\subsubsection{Nearby cycles on the Hecke stack}

\defi
\thlabel{defi--nearby-hecke}
Let $\uP$ be as in \thref{defi-iterated} and let $\mathbf{f}$ be a facet. The \emph{Hecke stack} associated to $\uP$ and $\mathbf{f}$ is the prestack quotient $\Hck_{\mathcal{G}_{\mathbf{f}}}(\uP) = \mathcal{L}^+\mathcal{G}_{\mathbf{f}}^{\BD} \backslash \Gr_{\mathcal{G}_{\mathbf{f}}}(\uP)$.

We define the functor 
$$\Upsilon_{\uP} \colon \DM(\Hck_{\mathcal{G}_{\mathbf{f}}}(\uP)_{\eta}) \rightarrow \NilpQ \DM(\Hck_{\mathcal{G}_{\mathbf{f}}}(\uP)_{s})$$
to be the unipotent nearby cycles functor associated to the structure map of $\Hck_{\mathcal{G}_{\mathbf{f}}}(\uP)$.
We also define the functor $$\mathsf{Z}_{\mathbf{f}} \colon \DM(L^+G \backslash \Gr) \rightarrow \NilpQ \DM(L^+G_{\mathbf{f}} \backslash \Fl_{\mathbf{f}}), \quad \mathsf{Z}_{\mathbf{f}}(\mathcal{F}) = \Upsilon_{\ux} (\mathcal{F} \boxtimes_{S} \Z_{C_\eta}),$$ where we use the identification $(\Gr_{\mathcal{G}_{\mathbf{f}}})_\eta = \Gr \times_S C_\eta$.
\xdefi 

\rema
In the expression $\mathcal F \boxtimes_S \Z_{C_\eta}$, implicitly we apply a restriction functor, i.e., !-pullback along the map of prestack quotients $(L^+\mathcal{G}_{\mathbf{f}}^{\BD} \backslash \Gr_{\mathcal{G}_{\mathbf{f}}})_\eta \rightarrow (L^+\mathcal{G}_{\mathbf{f}} \backslash \Gr_{\mathcal{G}_{\mathbf{f}}})_\eta$.

The functor $\Upsilon$ on $\Hck_{\mathcal G_{\mathbf f}}$ exists for formal reasons, cf.~\refsect{KanInd}:
such a functor exists for the motives on the ind-scheme $\Gr_{\mathcal G_{\mathbf f}}(\uP)$, and then on the prestack quotient by the action of the pro-smooth $\A^1$-group scheme $L^+ \mathcal G_{\mathbf f}$.

For later use, e.g. in \refsect{convolution}, also recall this: let $X = \colim X_i$ be an ind-scheme with an action of a pro-smooth group scheme $G= \lim G_j$ such that the action of $G$ on each $X_i$ factors through a smooth quotient $G_j$. If the kernel of $G \rightarrow G_j$ is split pro-unipotent, by \cite[Proposition 2.4.4]{RicharzScholbach:Intersection} we have $\DM(X_i/G) = \DM(X_i/G_j)$. We furthermore have an external product $\boxtimes$ on $\DM$ of such prestacks if they are placid in the sense of \cite[§A]{RicharzScholbach:Motivic}. By \thref{BD-rep} and the discussion following \thref{BD.Gr.I}, all of these assumptions are satisfied for $\Hck_{\mathcal{G}_{\mathbf{f}}}(\uP)$.
\xrema

\subsubsection{Convolution} \label{sect--convolution}
In order to refine the considerations in \cite[\S 4.2]{CassvdHScholbach:MotivicSatake} to a fully \ii-categorical treatment, we apply the constructions in \refsect{convolution structures} to $L^+ G_{\bbf} \subset LG$. Indeed, $LG$ is a placid ind-scheme \cite[§C.3]{Gaitsgory:Local} and $L^+G_{\bbf}$ is a pro-smooth group scheme, and the Zariski (and thus Nisnevich) sheafification of the presheaf quotient $LG/L^+G_{\bbf}$ is representable by the ind-proper ind-scheme $\Fl_{\bbf}$.
By construction, the multiplication map of the (non-symmetric) monoidal structure on $\DM(L^+G_{\mathbf{f}} \backslash LG / L^+G_{\mathbf{f}})$ is the usual convolution functor
$$\star : \DM(L^+G_{\mathbf{f}} \backslash LG / L^+G_{\mathbf{f}}) \x \DM(L^+G_{\mathbf{f}} \backslash LG / L^+G_{\mathbf{f}}) \r \DM(L^+G_{\mathbf{f}} \backslash LG / L^+G_{\mathbf{f}})$$
$$\calF_1 \star \calF_2 := m_! \pi^! (\calF_1 \boxtimes \calF_2).$$ 
Here $\pi$ and $m$ are the natural quotient and multiplication maps of prestacks:
$$\xymatrix{L^+G_{\mathbf{f}} \backslash LG / L^+G_{\mathbf{f}} \x L^+G_{\mathbf{f}} \backslash LG / L^+G_{\mathbf{f}} & \ar[l]_-{\pi} L^+G_{\mathbf{f}} \backslash LG  \x^{L^+G_{\mathbf{f}}} LG / L^+G_{\mathbf{f}}  \ar[r]^-m & L^+G_{\mathbf{f}} \backslash LG / L^+G_{\mathbf{f}}. }$$
The underlying motive of  $\pi^! (\calF_1 \boxtimes \calF_2)$ on the Zariski quotient $(LG  \x^{L^+G_{\mathbf{f}}} \Fl_{\mathbf{f}})_{\Zar}$ is sometimes also denoted $\calF_1 \widetilde \boxtimes \calF_2$. Under $\ell$-adic realization one recovers the usual notion of convolution, cf.~\cite[Proposition 3.14]{RicharzScholbach:Motivic}.

\prop \thlabel{thm-convDTM}
The convolution product turns $\DTM_{L^+G_{\mathbf{f}}}(LG / L^+G_{\mathbf{f}})^{(\anti)}$ into a monoidal \ii-category.
\xprop

\pf
By \cite[Definition and Lemma 4.11]{CassvdHScholbach:MotivicSatake} the functor $\star$ preserves anti-effective (resp.~all) stratified Tate motives, so the monoidal structure on $\DM(L^+G_{\bbf} \backslash LG/ L^+G_{\bbf})$ established above restricts to this full subcategory.
\xpf

The notion of local compactness (as defined in Section \ref{sect--Grothendieck}) is also preserved by the convolution product.

\lemm
\thlabel{monoidal structure locc}
The convolution product preserves the subcategory \(\DTM_{L^+G_{\mathbf{f}}}(\Fl_{\bbf})^{\locc}\subseteq \DTM_{L^+G_{\mathbf{f}}}(\Fl_{\bbf})\) and therefore turns it into a monoidal \ii-category.
\xlemm
\pf
Let \(\Ff_1,\Ff_2\in \DTM_{L^+G_{\mathbf{f}}}(\Fl_{\mathbf{f}})^{\locc}\).
Then \(u^!(\Ff_1\star \Ff_2)\) agrees with the pushforward of the twisted exterior product \(\Ff_1\widetilde{\boxtimes} \Ff_2\) (cf.~\cite[Lemma 4.16]{CassvdHScholbach:MotivicSatake}) along the multiplication map \[\Fl_{\mathbf{f}} \widetilde{\times} \Fl_{\mathbf{f}} := LG \overset{L^+G_{\mathbf{f}}}{\times}\Fl_{\mathbf{f}} \to \Fl_{\mathbf{f}}.\]
Since this multiplication map is ind-proper, pushforward along it preserves compactness, so it suffices to show \(\Ff_1\widetilde{\boxtimes} \Ff_2\) is compact.

This twisted exterior product is defined such that its pullback along the \(L^+G_{\mathbf{f}}\)-torsor \(LG \times \Fl_{\mathbf{f}} \to \Fl_{\mathbf{f}} \widetilde{\times} \Fl_{\mathbf{f}}\) agrees with the pullback of \(u^!(\Ff_1) \boxtimes u^!(\Ff_2)\) along the \(L^+G_{\mathbf{f}}\)-torsor \(LG \times \Fl_{\mathbf{f}} \to \Fl_{\mathbf{f}} \times \Fl_{\mathbf{f}}\), or similarly for suitable truncations of \(L^+G_{\mathbf{f}}\).
The lemma then follows from the fact that pullback along such (Zariski-locally trivial) torsors preserves and reflects compact objects.
\xpf

We have the following global variant of convolution as in \cite[Eqn.~(2.3.8)]{AcharRiche:Central}. 
There is a natural union operation $\cup$ on the symbols $\emptyset, \zero, \ux, \ux \cup \zero$.
For $1 \leq i \leq j \leq n$ let
$$m_{i,j} \colon \Gr_{\mathcal{G}_{\mathbf{f}}}(P_1, \ldots, P_n) \rightarrow \Gr_{\mathcal{G}_{\mathbf{f}}}(P_1, \ldots, P_{i-1}, P_i \cup \cdots \cup P_j, P_{j+1}, \ldots, P_n)$$
be the map which sends $\{(\mathcal{E}^i, \alpha_i)\}_{1 \leq i \leq n}$ to
$$\{(\mathcal{E}^1, \ldots, \mathcal{E}^{i-1}, \mathcal{E}^j, \mathcal{E}^{j+1}, \ldots, \mathcal{E}^n, \alpha_1, \ldots, \alpha_{i-1}, \alpha_{j}', \alpha_{j+1}, \ldots, \alpha_n)\}.$$ 
Here $\alpha_j' = \restr{\alpha_j}{\hat{D}_{P'}^\circ} \circ \cdots \circ \restr{\alpha_i}{\hat{D}_{P'}^\circ}$ and $P' = P_i \cup \cdots \cup P_j$. 
The convolution maps $m_{i,j}$ are $\mathcal{L}^+\mathcal{G}_{\mathbf{f}}^{\BD}$-equivariant. 
Over $C_\eta$ and $C_s$ the $m_{i,j}$ restrict to local convolution maps for $\Gr$ and $\Fl_{\mathbf{f}}$. These convolution maps are used throughout the proofs in \refsect{centralfunctor}. 

We construct a convolution product over $C_\eta$ as follows. Consider the diagram
$$\xymatrix{
\Gr_{\mathcal{G}_{\mathbf{f}}}(P_1) \x_{C} \Gr_{\mathcal{G}_{\mathbf{f}}}(P_2) & \Gr^{(\infty)}_{\mathcal{G}_{\mathbf{f}}}(P_1) \times_C \Gr_{\mathcal{G}_{\mathbf{f}}}(P_2) \ar[l]_p \ar[r]^-q & 
  \Gr_{\mathcal{G}_{\mathbf{f}}}(P_1, P_2)}. \eqlabel{TypeII}$$
The map $p$ is the base change of the \'etale-locally trivial $\mathcal{L}^+\mathcal{G}_{\mathbf{f}}^{\BD}$-torsor $\Gr_{\mathcal{G}_{\mathbf{f}}}^{(\infty)}(P_1) \rightarrow \Gr_{\mathcal{G}_{\mathbf{f}}}(P_1)$ from \refsect{iterated}. 
The map $q$ is an \'{e}tale torsor for the diagonal action of $\mathcal{L}^+\mathcal{G}_{\mathbf{f}}^{\BD}$ as in \cite[Corollary 2.3.14]{AcharRiche:Central}. We also have the $\mathcal{L}^+\mathcal{G}_{\mathbf{f}}^{\BD}$-equivariant convolution map $m_{1,2} \colon \Gr_{\mathcal{G}_{\mathbf{f}}}(P_1, P_2) \rightarrow  \Gr_{\mathcal{G}_{\mathbf{f}}}(P_1 \cup P_2)$, so we get a map $\ol m_{1,2} \colon \Hck_{\mathcal{G}_{\mathbf{f}}}(P_1, P_2) \rightarrow \Hck_{\mathcal{G}_{\mathbf{f}}}(P_1 \cup P_2)$.
The map $p$ is equivariant for the action of $(\mathcal{L}^+\mathcal{G}_{\mathbf{f}}^{\BD})^2$ by left translation on both factors. After taking fibers over $\eta$ and the Zariski-quotient by this action on the target, $p_\eta$ descends to a projection map
$$\xymatrix{
\Hck_{\mathcal{G}_{\mathbf{f}}}(P_1)_{\Zar, \eta} \x_{C} \Hck_{\mathcal{G}_{\mathbf{f}}}(P_2)_{\Zar, \eta} & \Hck_{\mathcal{G}_{\mathbf{f}}}(P_1, P_2)_{\eta}  \ar[l]_-{\ol p_\eta}. }\eqlabel{TypeII.1}
$$
Zariski-sheafification is necessary, but \'{e}tale sheafification is not, since $p_\eta$ and $q_\eta$ are Zariski-locally trivial. A similar projection exists over $C_s$.

\rema
A priori the maps $p_\eta$ and $q_\eta$ are only \'etale quotients, but since $LG \rightarrow \Fl_{\mathbf{f}}$ admits sections Zariski-locally by \thref{defi--partial}~ff., then the Zariski-local triviality of $p_\eta$ and $q_\eta$ follows from \thref{lemm-conv-fibers}, cf.~also \cite[Lemma 4.16 ff.]{CassvdHScholbach:MotivicSatake}. One can avoid taking Zariski-sheafifications in \refeq{TypeII.1} at the cost of working with products of pre-stack quotients of local loop groups in both the source and target, e.g. $LG/L^+G_{\mathbf{f}}$ instead of $\Fl_{\mathbf{f}}$. This loop group description does not extend globally to $\Gr_{\mathcal{G}_{\mathbf{f}}}(P_1, P_2)$, and we do not assert that $q$ is Zariski-locally trivial.
\xrema

We define $$\star_{C_\eta} : \DM(\Hck_{\mathcal{G}_{\mathbf{f}}}(P_1)_\eta) \x \DM(\Hck_{\mathcal{G}_{\mathbf{f}}}(P_2)_\eta ) \r \DM(\Hck_{\mathcal{G}_{\mathbf{f}}}(P_1 \cup P_2)_\eta )$$
$$\calF_1 \star_{C_\eta} \calF_2 := (\ol m_{1,2})_ {\eta !} \circ \ol p_{\eta}^! (\calF_1 \boxtimes_{C_\eta} \calF_2).$$
The construction of this functor is analogous to that in \cite[\S 4.2.4]{CassvdHScholbach:MotivicSatake}. We are implicitly using the fact that $\DM$ satisfies Nisnevich descent so that we can replace Zariski quotients with pre-stack quotients, e.g. $\DM(\Hck_{\mathcal{G}_{\mathbf{f}}}(P)_{\Zar, \eta}) = \DM(\Hck_{\mathcal{G}_{\mathbf{f}}}(P)_{\eta})$.
Working over $C_s$ gives another convolution product, which we denote by $\star$ since it recovers the convolution product on $\DM(L^+G_{\mathbf{f}} \backslash LG / L^+G_{\mathbf{f}})$. There are natural associativity constraints on $ \star_{C_\eta}$ and $\star$ at the level of homotopy categories constructed as in  \cite[Eqn.~(2.4.5)]{AcharRiche:Central}.

\lemm \thlabel{lemm-conv-nearby}
Let $\uP = (P_1, P_2)$.
For  $\mathcal{F}_1 \in \DM(\Hck_{\mathcal{G}_{\mathbf{f}}}(P_1)_\eta)$ and $\mathcal{F}_2 \in \DM(\Hck_{\mathcal{G}_{\mathbf{f}}}(P_2)_\eta)$, there are canonical morphisms in $\NilpQ \DM(\Hck_{\mathcal G_{\mathbf f}}(P_1, P_2)_s)$ and 
$\NilpQ \DM(\Hck_{\mathcal G_{\mathbf f}}(P_1 \cup P_2)_s)$, functorial in $\calF_1$ and $\calF_2$:
$$\ol{p}_s^! (\Upsilon_{P_1}(\calF_1) \boxtimes_{C_s} \Upsilon_{P_2}(\calF_2)) \rightarrow  \Upsilon_{\uP} (\ol{p}_\eta^! (\calF_1 \boxtimes_{C_\eta} \calF_2))\eqlabel{Upsilon-box-map}.$$
$$\Upsilon_{P_1}(\calF_1) \star \Upsilon_{P_2}(\calF_2) \rightarrow \Upsilon_{P_1 \cup P_2}(\calF_1 \star_{C_\eta} \calF_2)\eqlabel{Upsilon-convolution-map}.$$
\xlemm 

\pf By continuity we may assume that $\calF_1$ and $\calF_2$ are bounded in the sense of \cite[Example 2.4]{CassvdHScholbach:MotivicSatake}, i.e., the underlying non-equivariant motives are supported on finite type $S$-subschemes of their respective ind-schemes.  Applying $(\ol m_{1,2})_ {s !}$ to both sides of \refeq{Upsilon-box-map} and using compatibility with proper pushforward along $\ol m_{1,2}$ gives \refeq{Upsilon-convolution-map}.
To construct \refeq{Upsilon-box-map}  we first replace the maps $p$ and $q$ in \refeq{TypeII} by $\mathcal{L}^{+,(n)}\mathcal{G}_{\mathbf{f}}^{\BD}$-torsors for some integer $n$. The functor $ \Upsilon_{\uP}$ is an equivariant nearby cycles functor along  $\Gr_{\mathcal{G}_{\mathbf{f}}}(P_1, P_2)$, so we may view the motives in \refeq{Upsilon-box-map} as  $L^+G_{\mathbf{f}}$-equivariant motives on $\Gr_{\mathcal{G}_{\mathbf{f}}}(P_1, P_2)_s$.

Since $q_s$ admits sections Zariski-locally, then to give a map as in \refeq{Upsilon-box-map}, it is equivalent to apply $q_s^!$ and then give an  $(L^+G_{\mathbf{f}})^2$-equivariant map of motives on the source of $q_s$. By compatibility of $\Upsilon$ with smooth pullback, we need an $(L^+G_{\mathbf{f}})^2$-equivariant map
$$
p_s^! (\Upsilon_{P_1}(\calF_1) \boxtimes_{C_s} \Upsilon_{P_2}(\calF_2)) \rightarrow \Upsilon_{ \Gr^{(\infty)}_{\mathcal{G}_{\mathbf{f}}}(P_1) \times_C \Gr_{\mathcal{G}_{\mathbf{f}}}(P_2)} (p_\eta^! (\mathcal{F}_1 \boxtimes_{C_\eta} \mathcal{F}_2)). \eqlabel{eqn--conv-nearby}
$$
The underlying map of motives in \refeq{eqn--conv-nearby} is obtained from applying $p_s^!$ to a K\"unneth map for nearby cycles along $\Gr_{\mathcal{G}_{\mathbf{f}}}(P_1) \x_{S} \Gr_{\mathcal{G}_{\mathbf{f}}}(P_2) \rightarrow C^2$. K\"unneth maps commute with smooth pullback, so this is equivalently a K\"unneth map for nearby cycles along  $\Gr^{(\infty)}_{\mathcal{G}_{\mathbf{f}}}(P_1) \times_C \Gr_{\mathcal{G}_{\mathbf{f}}}(P_2)$. The particular K\"unneth map we need is obtained from \thref{Upsilon.boxtimes} followed by \thref{coro--nearby-pullback}. By construction of these natural transformations, the resulting map \refeq{eqn--conv-nearby} is $\Sigma$-linear (where $\Sigma$ acts on the domain by restriction along the comultiplication \eqref{comultiplication Sigma}, cf.~\refeq{restriction Mod Sigma}). 

It remains to see that \refeq{eqn--conv-nearby} is $(L^+G_{\mathbf{f}})^2$-equivariant. In fact, $p_\eta^! (\mathcal{F}_1 \boxtimes_{C_\eta} \mathcal{F}_2)$ is $(\mathcal{L}^+\mathcal{G}_{\mathbf{f}}^{\BD})_\eta^3$-equivariant, so \refeq{eqn--conv-nearby} can be refined to a K\"unneth map for $(\mathcal{L}^+\mathcal{G}_{\mathbf{f}}^{\BD})^3$-equivariant nearby cycles along $\Gr^{(\infty)}_{\mathcal{G}_{\mathbf{f}}}(P_1) \times_C \Gr_{\mathcal{G}_{\mathbf{f}}}(P_2)$. The desired $(L^+G_{\mathbf{f}})^2$-equivariance comes from restriction along the map $(L^+G_{\mathbf{f}})^2 \rightarrow (L^+G_{\mathbf{f}})^3$, $(g_1, g_2) \mapsto (g_1, g_2, g_2)$.
\xpf

\rema \thlabel{rema--multiple-factors}
There are diagrams similar to \refeq{TypeII.1} and \refeq{TypeII} for $n > 1$ factors, in particular
$$\Hck_{\mathcal{G}_{\mathbf{f}}}(P_1)_{\Zar, \eta} \x_{C} \cdots \x_{C}  \Hck_{\mathcal{G}_{\mathbf{f}}}(P_n)_{\Zar, \eta} \xleftarrow{\; \; \ol p \; \;} \Hck_{\mathcal{G}_{\mathbf{f}}}(P_1, \ldots P_n )_\eta$$ where $\uP = (P_1, \ldots, P_{i-1}, P_i \cup \cdots \cup P_j, P_{j+1}, \ldots, P_n)$. Considerations similar to those in \thref{lemm-conv-nearby} lead to a natural map in $\NilpQ \DM$:
$$\ol{p}_s^! (\Upsilon_{P_1}(\calF_1) \boxtimes_{C_s} \cdots \boxtimes_{C_s} \Upsilon_{P_n}(\calF_n)) \rightarrow  \Upsilon_{(P_1, \ldots, P_n)} (\ol{p}_\eta^! (\calF_1 \boxtimes_{C_\eta} \cdots \boxtimes_{C_\eta} \calF_n))\eqlabel{Upsilon-box-map-multiple}.$$ We will show in \thref{theo-WT-BD}\refit{theo-WT-BD3} that \refeq{Upsilon-box-map-multiple} is an isomorphism for Tate motives. Applying pushforward along $\ol m_{i,j} \colon \Hck_{\mathcal{G}_{\mathbf{f}}}(P_1, \ldots, P_n )_\eta \rightarrow \Hck_{\mathcal{G}_{\mathbf{f}}}(\uP)$ then gives an isomorphism of convolution products which will be used throughout \refsect{centralfunctor}.
\xrema

\subsection{Multiplicative group actions}
\label{sect--GmActions}
In order to construct Whitney--Tate stratifications of Beilinson--Drinfeld Grassmannians we must first take a detour and consider certain multiplicative group actions.
The general setup is as follows. 
Fix a cocharacter $\lambda \in X_*(T)$. 
There is an inclusion $\GmCI \rightarrow L_I^+(\GmC)$ coming from the identifications $\GmCI(R) = R^\times \subset R[\![D]\!]^\times = {L}^+_{I}(\GmC)(R)$. 
Let $\lambda_C \colon \GmC \rightarrow T_C$ be the base change of $\lambda$ by $C$.
By composing the previous inclusion with $L^+_I(\lambda_C) \colon {L}^+_{I}(\GmC) \rightarrow {L}^+_{I}(T_C)$ and the natural map $L^+_I(T_C) \rightarrow L^+_I\mathcal{G}_{\mathbf{f}}$, we get a $\GmCI$-action on $\Gr_{\mathcal{G}_{\mathbf{f}}, I}$. 
This also induces a $\GmC$-action on $\Gr_{\mathcal{G}_{\mathbf{f}}}^{\BD}$ and on the $\Gr_{\mathcal{G}_{\mathbf{f}}}(\uP)$. 
The following lemma allows us to apply the hyperbolic localization isomorphism \cite[Theorem B]{Richarz:Spaces}.

\lemm \thlabel{lemm-Zar-lin}
The $\GmCI$-action on $\Gr_{\mathcal{G}_{\mathbf{f}}, I}$ associated to a cocharacter $\lambda \in X_*(T)$ is Zariski-locally linearizable, as is the $\GmC$-action on $\Gr_{\mathcal{G}_{\mathbf{f}}}(\uP)$.
\xlemm

\pf
By \cite[Corollary 11.7]{PappasZhu:Kottwitz} and \cite[Proposition 3.10]{HainesRicharz:TestFunctionsWeil} there is a closed subgroup embedding $\mathcal{G}_{\mathbf{f}} \rightarrow   \GL_{n, C}$ for some integer $n$ which induces a $\GmCI$-equivariant immersion $\Gr_{\mathcal{G}_{\mathbf{f}}, I} \rightarrow \Gr_{ \GL_{n, C}, I}$. This immersion is in fact closed because $\Gr_{\mathcal{G}_{\mathbf{f}}, I}$ is ind-proper. 
Hence the $\GmCI$-action on $\Gr_{\mathcal{G}_{\mathbf{f}}, I}$ is Zariski-locally linearizable since the proof of \cite[Lemma 3.16]{HainesRicharz:TestFunctionsWeil} shows that the same is true for the $\GmCI$-action on $\Gr_{ \GL_{n, C}, I}$. By \thref{lemm--iteratedSingleton} this also treats the case where $\uP$ is a singleton.

The general case can be deduced as follows. Let $m_i \colon \Gr_{\mathcal{G}_{\mathbf{f}}}(\uP) \rightarrow \Gr_{\mathcal{G}_{\mathbf{f}}}(P_1 \cup \cdots \cup P_i)$ be the partial convolution map which sends  $\{(\mathcal{E}^i, \alpha_i)\}_{1 \leq i \leq n}$ to $\{(\mathcal{E}^i, \alpha_1 \circ \dots \restr{\alpha_i}{\hat{D}^{\circ}_{P_1 \cup \cdots \cup P_i}})\}$ for $1 \leq i \leq n$. The product of the $m_i$ over $C$ induces a map
$$f \colon \Gr_{\mathcal{G}_{\mathbf{f}}}(\uP) \rightarrow \Gr_{\mathcal{G}_{\mathbf{f}}}(P_1) \times_C \Gr_{\mathcal{G}_{\mathbf{f}}}(P_1 \cup P_2) \times_C \cdots \times_C \Gr_{\mathcal{G}_{\mathbf{f}}}(P_1 \cup \cdots \cup P_n),$$ which is $\GmC$-equivariant for the diagonal action of $\GmC$ on the target. Hence, it suffices to show that $f$ is a closed embedding. Since both ind-schemes are ind-proper, it suffices to show that $f$ is a monomorphism. 
Thus, suppose that  $f$ sends two $R$-points, $\{(\mathcal{E}^i, \alpha_i)\}_{1 \leq i \leq n}$ and $\{(\mathcal{F}^i, \beta_i)\}_{1 \leq i \leq n}$, to the same element. These elements are equal in $\Gr_{\mathcal{G}_{\mathbf{f}}}(\uP)$ if and only if there exist isomorphisms $ \mathcal{E}^i \xrightarrow{\sim} \mathcal{F}^i$ over $\hat{D}_R$ for $1 \leq i \leq n$ making an appropriate diagram commute. 
Equality may be checked after passing to an \'etale cover of $R$, so by \cite[Lemma 3.4]{HainesRicharz:TestFunctionsWeil} we may assume that the $\mathcal{E}^i$ and $\mathcal{F}^i$ are trivializable. 
We can then choose equivalent representatives in $\Gr_{\mathcal{G}_{\mathbf{f}}}(\uP)$ such that $\mathcal{E}^i = \mathcal{F}^i = \mathcal{E}^0$ for all $i$. 
Note that for all $P \in \{\emptyset, \zero, \ux, \ux \cup \zero\}$, $\Gr_{\mathcal{G}_{\mathbf{f}}}(P)(\hat{D}_P^\circ)$ is a subset of $\Gr_{\mathcal{G}_{\mathbf{f}}}(\hat{D}_R^\circ) := \Gr_{\mathcal{G}_{\mathbf{f}}}(R (\!( D) \!))$. 
Assuming that $\mathcal{E}^i = \mathcal{F}^i = \mathcal{E}^0$ for all $i$, it is straightforward to check that $\{(\mathcal{E}^i, \alpha_i)\}_{1 \leq i \leq n} = \{(\mathcal{F}^i, \beta_i)\}_{1 \leq i \leq n}$ in $\Gr_{\mathcal{G}_{\mathbf{f}}}(\uP)$ if and only if $$(\beta_1 \circ \cdots \circ \beta_i)^{-1} (\alpha_1 \circ \cdots \circ \alpha_i) \in \mathcal{G}_{\mathbf{f}}(\hat{D}_R), \quad \text{for all } 1 \leq i \leq n,$$ i.e, if and only if the above automorphisms of the trivial bundle $\restr{\mathcal{E}^0}{\hat{D}^\circ_R}$, which may be associated with elements of $\mathcal{G}_{\mathbf{f}}(\hat{D}_R^\circ)$, extend over $\hat{D}_R$.
On the other hand, $\{(\mathcal{E}^0, \alpha)\} = \{(\mathcal{E}^0, \beta)\}$ in $\Gr_{\mathcal{G}_{\mathbf{f}}}(P_1 \cup \cdots \cup P_i)$ if and only if $\beta^{-1} \alpha \in  \mathcal{G}_{\mathbf{f}}(\hat{D}_R)$. Hence $\{(\mathcal{E}^i, \alpha_i)\}_{1 \leq i \leq n} = \{(\mathcal{F}^i, \beta_i)\}_{1 \leq i \leq n}$, and $f$ is a closed embedding.
\xpf

We use the notation $(-)^0$, $(-)^+$, and $(-)^-$ to denote the fixed points, attractors, and repellers for an ind-scheme with $\Gm$-action in the sense of \cite{Richarz:Spaces, HainesRicharz:TestFunctionsWeil}. 
The cocharacter $\lambda$ induces a $\Gm$-action on $G$ by conjugation, $t \cdot g = \lambda(t) g \lambda(t)^{-1}$. 
Then $P^+ = G^+$ and $P^-= G^-$ are opposite parabolic subgroups of $G$ with Levi factor $M = P^+ \cap P^- = G^0$. 

\lemm
Consider the $\GmC$-action on $\mathcal{G}_{\mathbf{f}}$ induced by $\lambda \in X_*(T)$. 
\begin{enumerate}
\item The functors $\mathcal{P}_{\mathbf{f}}^+ := \mathcal{G}_{\mathbf{f}}^+$, $\mathcal{P}_{\mathbf{f}}^-:=\mathcal{G}_{\mathbf{f}}^-$ and $\mathcal{M}_{\mathbf{f}} := \mathcal{G}_{\mathbf{f}}^0$ are smooth, affine $C$-groups with geometrically connected fibers.
\item The Beilinson--Drinfeld Grassmannian  $\Gr_{\mathcal{P}_{\mathbf{f}}^\pm, I}$ is represented by a separated ind-finite type $C$-scheme for all $I$.
\item Furthermore, $\mathcal{M}_{\mathbf{f}}$ is the Bruhat--Tits group scheme for $M$ associated to the facet $\mathbf{f}$ as in \refsect{BT}.
\end{enumerate}
\xlemm

\pf We may assume that $S = \Spec \Z$ by \cite[Corollary 1.16]{Richarz:Spaces}. The geometric properties in (1) follow from \cite[Proposition 2.1.8]{CGP:Pseudo}. 
By \thref{pro-ind-proj} this also implies (2).  For (3), we note that the apartments $\mathscr{A}(G,T)$ and $\mathscr{A}(M,T)$ are naturally identified. The facet \(\mathbf{f}\) lies in a unique facet for \(M\), which we denote in the same way for simplicity.
Then for each prime $p$ the fiber $\mathcal{M}_{\mathbf{f}} \times_C \Z_p[t]$ is a parahoric $\Z_p[t]$-group scheme for $M$ by \cite[Lemma 5.16]{HainesRicharz:TestFunctions}. 
This uses the unique characterization of these parahoric group schemes in \cite[Theorem 4.1, Corollary 4.2]{PappasZhu:Kottwitz}. 
There is a similar fiberwise characterization of parahoric $\Z[t]$-group schemes \cite[Theorem 1.3]{Lorenco:Entiers}, so this implies that $\mathcal{M}_{\mathbf{f}}$ is a parahoric $\Z[t]$-group scheme.
\xpf

We will often abbreviate $P^+ = P$ and $\mathcal{P}_{\mathbf{f}}^+ = \mathcal{P}_{\mathbf{f}}$. The following lemma describes the fibers of $(\Gr_{\mathcal{G}_{\mathbf{f}}})^0$ and $(\Gr_{\mathcal{G}_{\mathbf{f}}})^{\pm}$.

\lemm \thlabel{lemm-Gm-on-Flf}
Consider the following maps associated to the $\Gm$-action on $\Fl_{\mathbf{f}}$.
$$\xymatrix{
(\Fl_{\mathbf{f}})^0 \ar@/^1pc/[rr]^{\iota} & (\Fl_{\mathbf{f}})^{\pm} \ar[l]^{q^{\pm}} \ar[r]_{p^{\pm}} & \Fl_{\mathbf{f}}
}$$
\begin{enumerate}
\item \label{item--lemm-Gm-on-Flf1} The map $q^{\pm}$ is ind-affine, has geometrically connected fibers, and induces a bijection on connected components $\pi_0((\Fl_{\mathbf{f}})^{\pm}) \cong \pi_0((\Fl_{\mathbf{f}})^0)$.
\item \label{item--lemm-Gm-on-Flf2} The map $p^{\pm}$ is schematic and restricts to an immersion on each connected component of $(\Fl_\mathbf{f})^{\pm}$. The map  $p^{\pm}$ is also bijective on $k$-points for every field $k$.
\item \label{item--lemm-Gm-on-Flf3} The map $\iota$ is a closed immersion.
\item \label{item--lemm-Gm-on-Flf4} If $S$ is an algebraically closed field then there is a canonical isomorphism $\pi_0((\Fl_{\mathbf{f}})^{0}) = W_{M, \aff} \backslash W / W_{\mathbf{f}}$ where $W_{M, \aff}$ is the affine Weyl group of $M$ with respect to $\mathbf{a}_0$.
\item \label{item--lemm-Gm-on-Flf5} Suppose that $\mathbf{f}$ is in the closure of $\mathbf{a}_0$, and that $\lambda$ is a regular cocharacter. Then $\pi_0((\Fl_{\mathbf{f}})^{0}) = W / W_{\mathbf{f}}$ and there is a canonical isomorphism of reduced loci $$(\Fl_{\mathbf{f}})^{0}_{\text{\normalfont red}} = \bigsqcup_{W / W_{\mathbf{f}}} S.$$
\end{enumerate}
\xlemm

\pf
Parts (1) and (3) are generalities about \'etale-locally linearizable $\Gm$-actions \cite[Theorem 2.1]{HainesRicharz:TestFunctions}. 
For part (2), the map $p^\pm$ is schematic again by \cite[Theorem 2.1]{HainesRicharz:TestFunctions}. 
To see that $p^\pm$ restricts to an immersion on each connected component of $(\Fl_\mathbf{f})^{\pm}$ we use the fact that the $\Gm$-action is Zariski-locally linearizable, cf.~also the proof of \cite[Lemma~3.5]{CassvdHScholbach:MotivicSatake}. 
The claim then follows since on an affine $\GmC$-scheme $Y$ the map $Y^+ \rightarrow Y$ is an closed immersion \cite[Lemma 1.9]{Richarz:Spaces}, and the map $q^\pm$ commutes with taking $\Gm$-equivariant affine Zariski covers \cite[Lemma 1.11]{Richarz:Spaces}. Finally, $p^\pm$ is a monomorphism by \cite[Remark 1.19]{Richarz:Spaces} and it is surjective on $k$-points by \StP{0BXZ}.

For part (4), let $M_{\text{sc}}$ be the simply connected cover of $M_{\text{der}}$ and let $P_{\text{sc}}^{\pm} = M_{\text{sc}} \ltimes U^{\pm}$ where $P^\pm = M \ltimes U^{\pm}$. By the argument in \cite[Theorem 5.2]{AGLR:Local} (see also \cite[\S 4]{HainesRicharz:TestFunctions}), when $S$ is an algebraically closed field the reduced locus of $(\Fl_{\mathbf{f}})^{\pm}$ agrees with the reduced locus of the disjoint union of the connected semi-infinite orbits $LP_{\text{sc}} \cdot \dot{w}$ where $w \in W_{M, \aff} \backslash W / W_{\mathbf{f}}$ and $\dot{w}$ is any lift to $LG(S)$ (and the reduced locus of $(\Fl_{\mathbf{f}})^{0}$ agrees with the $LM_{\text{sc}} \cdot \dot{w}$).

For part (5), we have $M = T$ so that $W_{M, \aff}$ is trivial in part (4). 
By \cite[Corollary 1.16]{Richarz:Spaces} we may assume that $S = \Spec \Z$. 
Then there is a map $f\ \colon \sqcup_{W / W_{\mathbf{f}}} \Spec \Z \rightarrow (\Fl_{\mathbf{f}})^{0}_{\text{\normalfont red}}$ obtained by taking representatives of elements of $W$ in $LG(\Z)$.  
The map $f$ is a closed immersion by part (3), so by the proof of part (4) it is an isomorphism on reduced loci over fields. 
Next, we claim that $(\Fl_{\mathbf{f}})^{0}_{\text{\normalfont red}}$ is ind-smooth over $\Z$. If $\mathbf{f} = \mathbf{f}_0$ then ind-smoothness follows since $(\Gr_G)^0 = \Gr_T$ by \cite[Theorem 3.17]{HainesRicharz:TestFunctionsWeil}. 
If $\mathbf{c}$ is any facet in the closure of another facet $\mathbf{c}'$ then the argument in \cite[Lemma 4.9]{HainesRicharz:TestFunctions} shows that the natural map $(\Fl_{\mathbf{c'}})^0 \rightarrow (\Fl_{\mathbf{c}})^0$ is smooth and surjective. 
This uses that $\Fl_{\mathbf{c'}} \rightarrow \Fl_{\mathbf{c}}$ is smooth (it is a fibration for a homogeneous space for $G$) and a generality about $\Gm$-actions over a general base \cite[Lemma 2.2]{HainesRicharz:TestFunctions}. 
Taking $\mathbf{c}' = \mathbf{a}_0$ and $\mathbf{c} = \mathbf{f}_0$ shows that $(\Fl_{\mathbf{a_0}})^0_{\text{red}}$ is ind-smooth, and then taking $\mathbf{c} = \mathbf{f}$ shows that $(\Fl_{\mathbf{f}})^0_{\text{red}}$ is ind-smooth by \StP{02K5}. 
Hence the $\Q$-fiber of $(\Fl_{\mathbf{f}})^0_{\text{red}}$ is dense. Moreover, after restricting to any closed Schubert cell in $\Fl_{\mathbf{f}}$, the reduced fixed point locus is a disjoint union of finitely many smooth (in particular, normal) integral subschemes, each admitting a birational map from $\Spec \Z$ that is bijective on points. 
Thus, $f$ is an isomorphism by Zariski's main theorem \StP{05K0}.
\xpf

The following theorem due to Haines--Richarz describes the geometry of the $\GmCI$-action on $\Gr_{\mathcal{G}_{\mathbf{f}}, I}$.

\theo 
\thlabel{HypLocCentral}
We have a commutative diagram with the following properties.
$$
\xymatrix{
\Gr_{\mathcal{M}_{\mathbf{f}}, I} \ar[d]^{\iota^0} & \Gr_{\mathcal{P}_{\mathbf{f}}^\pm, I} \ar[d]^{\iota^{\pm}} \ar[l]_{q_{\mathbf{f}}^\pm} \ar[r]^{p_{\mathbf{f}}^\pm} & \Gr_{\mathcal{G}_{\mathbf{f}}, I} \ar[d]^{\id} \\
(\Gr_{\mathcal{G}_{\mathbf{f}}, I})^0 & (\Gr_{\mathcal{G}_{\mathbf{f}}, I})^{\pm} \ar[l]_{q^{\pm}} \ar[r]^{p^{\pm}} & \Gr_{\mathcal{G}_{\mathbf{f}}, I}
}
$$
\begin{enumerate}
\item \label{item--HypLocCentral1} The maps $q_{\mathbf{f}}^\pm$ and $p_{\mathbf{f}}^\pm$ are induced by the natural maps $\mathcal{M}_{\mathbf{f}} \leftarrow \mathcal{P}_{\mathbf{f}}^{\pm} \rightarrow \mathcal{G}_{\mathbf{f}}$ associated to the $\GmC$-action on $\mathcal{G}_{\mathbf{f}}$.
\item \label{item--HypLocCentral2} The maps $q^{\pm}$ and $p^{\pm}$ are the natural maps associated to the $\GmCI$-action on $\Gr_{\mathcal{G}_{\mathbf{f}}, I}$.
\item \label{item--HypLocCentral3} The map $q^{\pm}$ is ind-affine, has geometrically connected fibers, and induces a bijection on connected components $\pi_0((\Gr_{\mathcal{G}_{\mathbf{f}}, I})^{\pm}) \cong \pi_0((\Gr_{\mathcal{G}_{\mathbf{f}}, I})^0)$.
\item \label{item--HypLocCentral4} The map $p^{\pm}$ is schematic and restricts to an immersion on each connected component of $(\Gr_{\mathcal{G}_{\mathbf{f}}, I})^{\pm}$.
\item \label{item--HypLocCentral5} The restrictions of $\iota^0$ and $\iota^\pm$ to $(C_\eta)^I$ are isomorphisms. 
\item \label{item--HypLocCentral6} Furthermore, the natural map $(\Gr_{\mathcal{G}_{\mathbf{f}}, I})^0 \rightarrow \Gr_{\mathcal{G}_{\mathbf{f}}, I}$ is a closed immersion.
\end{enumerate}
\xtheo

\pf
Since the action of $\GmC$ on $\Gr_{\mathcal{M}_{\mathbf{f}}}$ is trivial the natural map $\Gr_{\mathcal{M}_{\mathbf{f}}} \rightarrow \Gr_{\mathcal{G}_{\mathbf{f}}}$ factors through $(\Gr_{\mathcal{G}_{\mathbf{f}}})^0$. 
This defines $\iota^0$. 
The Rees construction explained in \cite[\S 3.3.1]{HainesRicharz:TestFunctionsWeil}, which uses the description of $\Gr_{\mathcal{P}_{\mathbf{f}}^\pm}$ as a moduli space of $\mathcal{P}_{\mathbf{f}}^\pm$-torsors, gives a canonical isomorphism $\Gr_{\mathcal{P}_{\mathbf{f}}^\pm} \cong (\Gr_{\mathcal{P}_{\mathbf{f}}^\pm})^{\pm}$ (it is assumed in loc.~cit.~that $\mathcal{G}_{\mathbf{f}}$ is reductive but this assumption is only necessary for showing that $\iota^0$ and $\iota^\pm$ are isomorphisms).
Then $\iota^\pm$ is defined by composing with the natural map $(\Gr_{\mathcal{P}_{\mathbf{f}}^\pm})^{\pm} \rightarrow (\Gr_{\mathcal{G}_{\mathbf{f}}^\pm})^{\pm}$. 
This constructs the commutative diagram with the maps in (1) and (2). 
Parts (3), (4), and (6) are generalities about $\Gm$-actions and can be proved as in \thref{lemm-Gm-on-Flf}. Part (5) is a special case of \cite[Theorem 3.17]{HainesRicharz:TestFunctionsWeil}.
\xpf

\subsection{Whitney--Tateness of Beilinson--Drinfeld Grassmannians}
\subsubsection{Stratifications} From now on we assume that all facets are contained in the closure of $\mathbf{a}_0$. Then we may consider the stratification of $\Gr$ by the $\Gr^{\mu}$ for $\mu \in X_*(T)^+$, and the stratification of $\Fl_{\mathbf{f}}$ by the $\Fl^\circ_w(\mathbf{f},\mathbf{f})$ for $w\in W_\mathbf{f} \backslash W/W_\mathbf{f}$.  We build a stratification of $\Gr_{\mathcal{G}_{\mathbf{f}}}(P)$ using these strata as follows.

\defilemm \thlabel{defi-strat-Gr}
Let $P \in \{\emptyset, \zero, \ux, \ux \cup \zero \}$. We define a stratification of $\Gr_{\mathcal{G}_{\mathbf{f}}}(P)$ by taking the union of the following strata over $C_\eta$ and $C_s$ with respect to the isomorphisms in \thref{lemm-conv-fibers}, where $\mu \in X_*(T)^+$ and $w\in W_\mathbf{f} \backslash W/W_\mathbf{f}$.
\begin{center}
\begin{tabular}{ |c|c|c| } 
 \hline
 $P$ & $\Gr_{\mathcal{G}_{\mathbf{f}}}(P)_\eta$ & $\Gr_{\mathcal{G}_{\mathbf{f}}}(P)_s$ \\
 \hline
 $\emptyset$ & $C_\eta$ & $C_s$ \\ 
 $\zero$ & $\Fl^\circ_w(\mathbf{f},\mathbf{f}) \times C_\eta$  & $\Fl^\circ_w(\mathbf{f},\mathbf{f}) \times C_s $ \\
  $\ux$ & $\Gr^\mu \times C_\eta$ & $\Fl^\circ_w(\mathbf{f},\mathbf{f}) \times C_s$ \\ 
   $\ux \cup \zero$ & $\Gr^\mu \times \Fl^\circ_w(\mathbf{f},\mathbf{f}) \times C_\eta$ & $\Fl^\circ_w(\mathbf{f},\mathbf{f}) \times C_s$ \\ 
 \hline
\end{tabular}
\end{center}
\xdefilemm

\pf
We must show that the closure of a stratum is a union of strata. In general, the strata are the orbits for the actions of $(L^+\mathcal{G}_{\mathbf{f}}^{\BD})_{\eta}$ and $(L^+\mathcal{G}_{\mathbf{f}}^{\BD})_{s}$. The closures of the strata over $C_\eta$ are thus stable under the action of $L^+\mathcal{G}_{\mathbf{f}}^{\BD}$, since this group is pro-smooth. Hence the fibers over $C_s$ of these closures are also stable under the action of $(L^+\mathcal{G}_{\mathbf{f}}^{\BD})_{s}$. This implies the statement if $S$ is a field, but in general, we need to show that the strata appearing in the fibers over $C_s$ are independent of the base $S$. Using basic facts about scheme-theoretic images, one can reduce to the cases $S = \Spec \Z$ or $S = \Spec k$, for $k$ an arbitrary field, as in \cite[Proposition 4.4.3]{RicharzScholbach:Intersection}. 

If $P \in \{\emptyset, \zero\}$, the statement is clear since the family is constant. If $P = \ux$, we claim that the reduced locus of the fiber over $C_s$ of the closure of $\Gr{^\mu} \times C_\eta$ is the $\mu$-admissible locus $\mathcal{A}^{\mathbf{f}}(\mu)$. Indeed, if $S = \Spec k$ this follows from \cite[Theorem 6.12]{HainesRicharz:TestFunctions}, and if $S = \Spec \Z$ this follows from \cite[Theorem 5.2.2]{Lorenco:Entiers}. While \cite[Theorem 5.2.2]{Lorenco:Entiers} assumes that $G_\der$ is simply connected, this is only needed for other geometric properties of the closure; here we may pass to a $z$-extension $G' \r G$ with $G'_{\der}$ simply connected. There is a canonical facet $\mathbf{f}'$ for $G'$ associated to $\mathbf{f}$, giving rise to a universal homeomorphism on connected components of affine flag varieties for $G'$ and $G$. Furthermore, the induced map on closures of strata over $C_\eta$ for $G'$ and $G$ is birational, proper, and surjective, cf.~the proof of \cite[Corollary 2.3]{HainesRicharz:Cohen}, so the fact that the $\mu$-admissible locus is the reduced $C_s$-fiber for $G'$ (for any lift of $\mu$ to $G'$) implies the same for $G$. 

Finally, if $P = \ux \cup \zero$, we consider the proper convolution map $m \colon \Gr_{\mathcal{G}_{\mathbf{f}}}(\ux, \zero) \rightarrow \Gr_{\mathcal{G}_{\mathbf{f}}}(\ux \cup \zero) $, which is an isomorphism over $C_\eta$. We may first compute the closure of $\Gr^\mu \times \Fl^\circ_w(\mathbf{f},\mathbf{f}) \times C_\eta$ inside $\Gr_{\mathcal{G}_{\mathbf{f}}}(\ux, \zero)$, and then take its image under $m$. Using the previous cases, this implies that the reduced locus over $C_s$ of the closure inside $\Gr_{\mathcal{G}_{\mathbf{f}}}(\ux \cup \zero)$ is the image of the convolution map $\mathcal{A}^{\mathbf{f}}(\mu) \overset{L^+G_{\mathbf{f}}}{\times} \Fl_w(\mathbf{f}, \mathbf{f}) \r \Fl_{\mathbf{f}}$. The Schubert cells appearing are thus uniquely determined by the combinatorics of $W$, independently of $S$, cf.~\cite[\S 3.2]{RicharzScholbach:Motivic}.
\xpf

\rema
\thlabel{notation DTM Gr}
By \thref{lemm-conv-fibers} and \cite[Proposition 3.7]{CassvdHScholbach:MotivicSatake} the stratifications of $\Gr_{\mathcal{G}_{\mathbf{f}}}(P)_{\eta}$ and $\Gr_{\mathcal{G}_{\mathbf{f}}}(P)_{s}$ in \thref{defi-strat-Gr} are anti-effective universally Whitney--Tate. In particular, we have the categories $\DTM(\Gr_{\mathcal{G}_{\mathbf{f}}}(P)_{\eta})^{\antic}$ and $\DTM(\Gr_{\mathcal{G}_{\mathbf{f}}}(P)_{s})^{\antic}$ of compact, anti-effective stratified Tate motives.
\xrema
 
To show that categories of Tate motives glue over $C$, we will use hyperbolic localization in the next section. \thref{example--WT} below is illustrative of the general case; first we need the following result.

\lemm \thlabel{lemm--DMorbit}
Suppose the action of $L^+G_{\mathbf{f}'}$ on $\Fl_{w}^\circ(\mathbf{f}', \mathbf{f})$ factors through some jet group $H : = L^{+, (m)}G_{\mathbf{f}'}$. Let $H_w \subset H$ be the stabilizer of the basepoint of $\Fl_{w}^\circ(\mathbf{f}', \mathbf{f})$ and let $e \colon S / H_w \rightarrow  H \backslash \Fl_{w}^\circ(\mathbf{f}', \mathbf{f})$ be the induced map of pre-stacks.

Then there are canonical equivalences
$$\DM(L^+G_{\mathbf{f}'} \backslash \Fl_{w}^\circ(\mathbf{f}', \mathbf{f}) ) = \DM(H \backslash \Fl_{w}^\circ(\mathbf{f}', \mathbf{f}) ) \xrightarrow[\sim]{e^!} \DM(S/H_w).$$
These equivalences restrict to $\DTM$.
Furthermore, some $H$-equivariant Tate motive in $\DTM(\Fl_{w}^\circ(\mathbf{f}', \mathbf{f}))$ is anti-effective iff its *-pullback to the basepoint $S \subset \Fl_w^\circ(\mathbf f', \mathbf f)$ is anti-effective.
Finally, if $S$ satisfies the Beilinson--Soulé vanishing condition then forgetting the equivariance induces an equivalence $\MTM(S/H_w) = \MTM(S)$.
\xlemm

\pf
The first equivalence follows as in \thref{rema--jet}. For the second we recall that $\Fl_{w}^\circ(\mathbf{f}', \mathbf{f}) = (H/H_w)_{\Zar}$ by \cite[Proposition 3.2]{CassvdHScholbach:MotivicSatake}, so that $e$ is an equivalence of stacks after Zariski-sheafification. Since $\DM$ is invariant under Zariski-sheafification we conclude that $e^!$ is an equivalence, cf.~\cite[Lemma 2.2.21]{RicharzScholbach:Intersection}. The same arguments as in \cite[Propositions 3.1.23, 3.1.27]{RicharzScholbach:Intersection}, \cite{RicharzScholbach:IntersectionCorrigendum} show that these equivalences respect $\DTM^{(\anti)}$ (the only difference is that we use Zariski instead of \'{e}tale quotients here).  Note that $*$-pullback necessarily appears when checking anti-effectivity, since $!$-pullback introduces a twist. Finally, the equivalence $\MTM(S/H_w) = \MTM(S)$ follows as in \cite[Proposition 3.2.23]{RicharzScholbach:Intersection} using that $H_w$ is cellular and fiberwise connected by \cite[Lemma 4.3.7]{RicharzScholbach:Intersection}, cf.~also \cite[Proposition 3.4]{CassvdHScholbach:MotivicSatake}.
\xpf

\exam \thlabel{example--WT}
Let us recall the simplest non-trivial example of such degenerations, which already appeared in \cite[§1.2.3]{Gaitsgory:Central}.
Namely, let \(G=\PGL_2\) and let $\mu$ be the unique minuscule dominant cocharacter. Then $\Gr^{\mu} \cong \P$, and the reduced special fiber of $\overline{\Gr^{\mu} \times \Gm} \subset \Gr_{\mathcal{G}_{\mathbf{a}_0}}$ is isomorphic to two copies of \(\P\) which intersect transversally in a single point $0$, which we identify with zero in both copies of \(\P\). Let us denote this special fiber by $X$, and the points at infinity by $\infty_1$ and $\infty_2$.

Let $j \colon \Gr^{\mu} \times \Gm \rightarrow  \overline{\Gr^{\mu} \times \Gm} \leftarrow X \colon i$ be the inclusions. To show that $i^*j_* \Z$ is Tate, we consider the semi-infinite orbits on $\Fl$ with respect to a regular dominant cocharacter. The semi-infinite orbits divide $X$ into three subschemes: a copy of $\A^1$ with origin $0$, a copy of $\A^1$ with origin identified with $\infty_1$ (up to relabelling), and a copy of $S$ identified with $\infty_2$. Switching to an anti-dominant regular cocharacter swaps  $\infty_1$ and $\infty_2$.

Note that hyperbolic localization preserves Tate motives on both $\Gr$ and $\Fl$, and it commutes with $i^*j_*$. Additionally, the hyperbolic localization of $i^*j_* \Z$ along the semi-infinite orbit identified with $\infty_2$ is just the stalk  $(i^*j_*\Z)_{\infty_2}$, which is therefore Tate. By using an anti-dominant regular cocharacter, $(i^*j_*\Z)_{\infty_1}$ is also Tate. Now by Iwahori-equivariance and \thref{lemm--DMorbit}, $(i^*j_* \Z)_{X \setminus \{0\}}$ is Tate. It remains to show that $(i^*j_*\Z)_{0}$ is Tate. Viewing $(i^*j_*\Z)_{0}$ as a motive supported on $X$, this stalk is trivially identified with the hyperbolic localization of $(i^*j_*\Z)_{0}$ along the copy of $\A^1$ with origin $0$. The advantage is that $(i^*j_*\Z)_{0}$ is the cofiber of $(i^*j_* \Z)_{X \setminus \{0\}} \rightarrow i^*j_* \Z$, and both of these terms have Tate hyperbolic localizations. Thus, $(i^*j_*\Z)_{0}$ is Tate.
\xexam

\subsubsection{Hyperbolic localization on affine flag varieties}
Let us fix a regular cocharacter \(\lambda\in X_*(T)\). As in \refsect{GmActions} this induces an action of $\GmC$ on $\Gr_{\mathcal{G}_{\mathbf{f}}}$ and $\Gr_{\mathcal{G}_{\mathbf{f}}}^{\BD}$.

\lemm 
\thlabel{lemm-closure-of-fixed}
The reduced loci of the fixed points $(\Gr_{\mathcal{G}_{\mathbf{f}}})^0$ and $(\Gr_{\mathcal{G}_{\mathbf{f}}}^{\BD})^0$ can be described as follows.
\begin{enumerate}
\item \label{item--lemm-closure-of-fixed1} We have canonical isomorphisms $$(\Gr_{\mathcal{G}_{\mathbf{f}}})^0_{\eta, \redu}  \cong  \coprod_{X_*(T)} C_\eta, \quad  (\Gr_{\mathcal{G}_{\mathbf{f}}}^{\BD})^0_{\eta, \redu} \cong \coprod_{X_*(T) \times W/W_{\mathbf{f}}} C_\eta.$$
\item \label{item--lemm-closure-of-fixed2} The reduced loci of both $(\Gr_{\mathcal{G}_{\mathbf{f}}})^0_{s}$ and $(\Gr_{\mathcal{G}_{\mathbf{f}}}^{\BD})^0_{s}$ are canonically isomorphic to $\coprod_{W/W_{\mathbf{f}}} S$.
\item \label{item--lemm-closure-of-fixed3} The reduced closure of any connected component of $(\Gr_{\mathcal{G}_{\mathbf{f}}})^0_{\eta}$ (resp. $(\Gr_{\mathcal{G}_{\mathbf{f}}}^{\BD})^0_{\eta}$) is canonically isomorphic to $C$.
\end{enumerate}
\xlemm

\pf
Parts (1) and (2) follow from \thref{lemm-Gm-on-Flf}\refit{lemm-Gm-on-Flf5}.
For part (3), note that $(\Gr_{T_C})_{\redu} = \coprod_{X_*(T)} C$. Hence the natural map $(\Gr_{T_C})_{\redu} \rightarrow \Gr_{\mathcal{G}_{\mathbf{f}}}$ admits a section over each connected component of the source, and thus it is a closed immersion by ind-properness. Then by \thref{HypLocCentral} we have a closed embedding $(\Gr_{T_C})_{\redu} \rightarrow (\Gr_{\mathcal{G}_{\mathbf{f}}})^0$. 
This implies that (3) holds for $\Gr_{\mathcal{G}_{\mathbf{f}}}$. To prove that (3) holds for $\Gr_{\mathcal{G}_{\mathbf{f}}}^{\BD}$, we consider the convolution map $m \colon \Gr_{\mathcal{G}_{\mathbf{f}}}(\ux, \zero) \rightarrow \Gr_{\mathcal{G}_{\mathbf{f}}}^{\BD}$, which is an isomorphism over $C_\eta$. 
The reduced closure of a connected component of $(\Gr_{\mathcal{G}_{\mathbf{f}}}^{\BD})^0_{\eta}$ inside $\Gr_{\mathcal{G}_{\mathbf{f}}}(\ux, \zero)$ is isomorphic to $C$, as can be seen using the isomorphism $\Gr_{\mathcal{G}_{\mathbf{f}}}(\ux, \zero) \cong \Gr_{\mathcal{G}^{(\infty)}_{\mathbf{f}}}(\ux) \x^{\mathcal{L}^+\mathcal{G}_{\mathbf{f}}^{\BD}} \Gr_{\mathcal{G}_{\mathbf{f}}}(\zero)$ and the case of $\Gr_{\mathcal{G}_{\mathbf{f}}}(\ux) = \Gr_{\mathcal{G}_{\mathbf{f}}}$. 
Hence the reduced closure inside $\Gr_{\mathcal{G}_{\mathbf{f}}}^{\BD}$ is an integral $C$-scheme $C'$ whose structure map $C' \rightarrow C$ is birational. The map $C' \rightarrow C$ is also bijective on points by considering the case where $S$ is an algebraically closed field. Since $C$ is normal we have $C' = C$ by Zariski's main theorem.
\xpf

\rema
By \thref{lemm-Gm-on-Flf}\refit{lemm-Gm-on-Flf5} and \thref{lemm-closure-of-fixed}, when $\uP$ is a singleton there is a canonical stratification of $\Gr_{\mathcal{G}_{\mathbf{f}}}(\uP)^0$ by copies of $C_\eta$ and  $S$ (the latter strata live over $C_s$). This stratification is universally anti-effective Whitney--Tate, since the same is true for the stratification $C = \A_S^1 =\GmS \sqcup S$. Indeed, for the inclusions $j \colon \GmS \rightarrow \A_S^1 \leftarrow S \colon  i$ we have $i^* j_* \Z = \Z \oplus \Z(-1)[-1]$ by relative purity. Thus, when $\uP$ is a singleton we may consider the category $\DTM(\Gr_{\mathcal{G}_{\mathbf{f}}}(\uP)^0)^{\anti}$.
\xrema

We now prove that hyperbolic localization preserves Tate motives on $\Fl_{\mathbf{f}}$. In the proof we will use the Demazure resolutions, which we now recall. Let $w \in W_{\aff}$ and let $\dot{w} = s_1 \cdots s_n$ be a reduced decomposition by simple reflections in the walls of $\mathbf{a}_0$. 
Let $\mathcal{P}_i \subset LG$ be the minimal parahoric subgroup containing $\mathcal{I}$ and a representative for $s_i$ (see \cite[\S 8.a]{PappasRapoport:Twisted} for more details). The Demazure scheme associated to $\dot{w}$ is $D(\dot{w}) = \mathcal{P}_1 \x^{\mathcal{I}} \cdots \x^{\mathcal{I}} \mathcal{P}_n/ \mathcal{I}$.
Multiplication in $LG$ induces an $\mathcal{I}$-equivariant map $D(\dot{w})  \rightarrow \Fl_w(\mathbf{a}_0, \mathbf{a}_0)$ which is an isomorphism over  $\Fl_w^\circ(\mathbf{a}_0, \mathbf{a}_0)$, cf. \cite[Definition 5 ff.]{Faltings:Loops}. 
The $\Gm$-action on $D(\dot{w})$ induced by the regular cocharacter $\lambda \in X_*(T)$ has the following properties.

\lemm \thlabel{lemm-Dem-Gm}
Let $w \in W_{\aff}$ and let $\dot{w} = s_1 \cdots s_n$ be a reduced decomposition.
\begin{enumerate}
\item The reduced locus of the fixed points $D(\dot{w})^0 \subset D(\dot{w})$ is a disjoint union of $2^n$ copies of $S$.
\item The natural map $\tilde p^+ \colon D(\dot{w})^+ \rightarrow D(\dot{w})$ is bijective on points and restricts to an immersion on connected components.
\item The restriction of the natural map $\tilde q^+\colon D(\dot{w})^+ \rightarrow D(\dot{w})^0$ to reduced loci is a disjoint union of relative affine spaces.
\item There exists a filtrable decomposition of $D(\dot{w})$ by the connected components of $D(\dot{w})^+$.
\end{enumerate} 
\xlemm

\pf
The product of the partial convolution maps $(p_1, \ldots, p_n) \mapsto p_1 \cdots p_i$ induces a $\Gm$-equivariant closed embedding $D(\dot{w}) \rightarrow \Fl^n$. Hence the fact that $\tilde p^+$ restricts to an immersion on connected components follows from the case of $\Fl$ as in \thref{lemm-Gm-on-Flf}\refit{lemm-Gm-on-Flf2}. 
The map $\tilde p^+$ is a bijection on points because $D(\dot{w})$ is projective, so in particular all maps $\Gm \rightarrow D(\dot{w})$ defined over a field extend to $\A^1$ (and it is a monomorphism because $D(\dot{w})$ is separated). This proves (2). 
When $S = k$ is an algebraically closed field the remaining parts (except for the precise number of fixed points) are true more generally for smooth projective $\Gm$-schemes with isolated fixed points by work of Bia{\l}ynicki-Birula \cite{BB:Theorems, BB:Properties}. 
Rather than generalizing these results to an arbitrary base we will give a direct proof for $D(\dot{w})$ following \cite[Theorem 2.17]{CassXu:Geometrization}.

First, we note that $D(\dot{s}_i) \cong \P$. 
Since $\lambda$ is regular then $D(\dot{s}_i)_{\redu}^0 \cong S \sqcup S$, corresponding to the two fixed points $\{e\}, \{s_i\}$ where $e$ is the basepoint. Furthermore, $D(\dot{s}_i)_{\redu}^+ \cong \A^1_S \sqcup S$. 
The connected component of $D(\dot{s}_i)_{\redu}^+$ which is isomorphic to $\A_S^1$ attracts toward $\{s_i\}$ (resp. $\{e\}$) if $\lambda$ attracts (resp. repels) the affine root group corresponding to $s_i$. 
This proves the remaining parts in the case $l(w) =1$. 

For the general case we induct on $l(w)$ using the isomorphism $D(\dot{w}) = D(\dot{t}) \x^{\mathcal{I}} D(\dot{s}_n)$, where $\dot{t} = s_1 \cdots s_{n-1}$. 
The projection $\pi \colon D(\dot{w}) \rightarrow D(\dot{t})$ is a $\Gm$-equivariant, Zariski-locally trivial $\P$-bundle. By applying the case $l(w) = 1$ to $D(\dot{s}_n)$, it follows that for each connected component $Y_i \subset D(\dot{t})^+_{\redu}$, the fiber $\pi^{-1}(Y_i)$ is set-theoretically the union of two connected components $X_{i}^1$, $X_i^2$ of $D(\dot{w})^+_{\redu}$. 
One of these connected components, say $X_i^1$, maps isomorphically on reduced loci onto $Y_i$ by the map $\pi$. The other connected component $X_i^2$ has the structure of a line bundle over $Y_i$ on reduced loci. Now (1) follows by induction. 
It also follows by induction that $\tilde q^+$ is a disjoint union of vector bundles. To see that these bundles are trivial we can reduce by base change to the case $S = \Z$, where it follows from the Quillen--Suslin theorem. This proves (3). 
Finally, to prove (4) we may suppose by induction that we have a filtrable decomposition by closed subschemes $Z_1 \subset \cdots \subset Z_{2^{n-1}} =  D(\dot{t})$ with relative complements $Z_i \setminus Z_{i-1} = Y_i$. 
Then it is easy to see that we have a filtrable decomposition of $D(\dot{w})$ with relative complements $X_1^1, X_1^2, \ldots, X_{2^{n-1}}^1, X_{2^{n-1}}^2$.
\xpf

\prop
\thlabel{HypLocTate}
In the notation of \thref{HypLocCentral}, hyperbolic localization induces a well-defined functor
$$q^+_!(p^+)^*\colon \DTM((\Gr_{\mathcal{G}_{\mathbf{f}}})_\eta)^{\antic} \to \DTM((\Gr_{\mathcal{G}_{\mathbf{f}}})^0_{\eta})^{\antic}.$$ 
The same is true more generally for  $\Gr_{\mathcal{G}_{\mathbf{f}}}(\uP)_\eta$ when $\uP$ is a singleton, and also for the respective fibers over $C_s$.
\xprop
\pf
Since $(\Gr_{\mathcal{G}_{\mathbf{f}}}^{\BD})_{\eta} = \Gr \times \Fl_{\mathbf{f}} \times C_\eta$, we may deduce the case of $\Gr_{\mathcal{G}_{\mathbf{f}}}^{\BD}$ from the case of $\Gr_{\mathcal{G}_{\mathbf{f}}}$ using the K\"unneth formula \cite[Theorem 2.4.6]{JinYang:Kuenneth} for $!$-pushforward. By the description of the fibers over $C_\eta$ and $C_s$ in \thref{lemm-conv-fibers}, we are reduced to proving a similar statement about hyperbolic localization on $\Fl_{\mathbf{f}}$. 
Recall that $(\Fl_{\mathbf{f}})^0_{\redu}$ is a disjoint union of copies of $S$ by \thref{lemm-Gm-on-Flf}\refit{lemm-Gm-on-Flf5}. 
As the stratification of $\Fl_{\mathbf{f}}$ by the $\Fl_w^\circ(\mathbf{a}_0, \mathbf{f})$ is finer than the stratification by the $\Fl_w^\circ(\mathbf{f}, \mathbf{f})$, and taking fixed points and attractors is compatible with passing to equivariant closed subschemes \cite[Lemma 1.4.9]{Drinfeld:Gm}, we can restrict ourselves to motives supported on \(\Fl_w(\mathbf{a}_0,\mathbf{f})\), for some \(w\in W/W_\mathbf{f}\). 

Let $\tilde w \in W$ be the unique lift of $w$ such that $\Fl_{\tilde w}^\circ(\mathbf{a}_0, \mathbf{a}_0) \cong \Fl_{w}^\circ(\mathbf{a}_0, \mathbf{f})$ via the canonical projection. After applying a $\Gm$-equivariant right translation by an element of $\Omega$ (which we may do since $\Omega$ normalizes $\mathcal{I}$) we may assume that $\tilde w \in W_{\aff}$.
To ease the notation we simply write $w$ instead of $\tilde w$ for the rest of the proof.
Chose a reduced decomposition $\dot{w} = s_1 \cdots s_n$.
Then we have a \(\Gm\)-equivariant Demazure resolution \(m\colon D(\dot{w})\to \Fl_w(\mathbf{a}_0,\mathbf{f})\), which is an isomorphism over $\Fl_{w}^\circ(\mathbf{a}_0, \mathbf{f})$. 

We claim that $m_! \Z \in \DTM(\Fl_w(\mathbf{a}_0,\mathbf{f}))^{\antic}$. 
To prove this, first suppose that $\mathbf{f} = \mathbf{a}_0$. Then as in the proof of \cite[Lemma 9]{Faltings:Loops}, $m$ factors as a composition of $\mathcal{I}$-equivariant maps with fibers trivial or isomorphic to $\P$.
The structure map $f \colon \P_S \rightarrow S$ satisfies $f_!(\Z) = \Z \oplus \Z(-1)[-2]$, which is compact and anti-effective Tate. 
By $\mathcal{I}$-equivariance and \thref{lemm--DMorbit}, the properties of Tateness, compactness, and anti-effectiveness may be checked fiberwise. Thus, these three properties follow from the cohomology of $\P$. 
For general $\mathbf{f}$ it remains to show that pushforward along $\Fl_{w}(\mathbf{a}_0, \mathbf{a}_0) \rightarrow \Fl_{w}(\mathbf{a}_0, \mathbf{f})$ preserves compact anti-effective stratified Tate motives. This follows by similar reasoning using that the reduced fibers of this map are unions of Schubert varieties in classical flag varieties, which in particular are stratified by affine spaces (the Bruhat decomposition), cf. the proof of \cite[Lemma 2.8]{CassXu:Geometrization}.

We now proceed similarly to \cite[Proposition 3.11]{CassvdHScholbach:MotivicSatake},  using an induction on \(l(w)\).
Let \(T\subseteq \Fl_w(\mathbf{a}_0,\mathbf{f})^+\) be a connected component, and denote by \(f\colon \tilde{T}:=m^{-1}(T)\to T\) the restriction of \(m\).
Then $\tilde T$ is set-theoretically a union of connected components of $D(\dot{w})^+$.
Consider the open immersion \(j\colon T\cap \Fl_w^\circ(\mathbf{a}_0,\mathbf{f})\subseteq T\) with closed complement \(i\colon Z\subseteq T\). 
By induction on $l(w)$ we need only show that $q^+_!(j_! \Z) \in \DTM(\Fl_w(\mathbf{a}_0,\mathbf{f})^0)^{\antic}$.
As \(m\) is an isomorphism over \(\Fl_w^\circ(\mathbf{a}_0,\mathbf{f})\), applying \(q^+_!\) to the localization sequence associated to \(f_!(\Z_{\tilde{T}})\) gives an exact triangle \(q^+_!(j_! \Z ) \to q^+_!f_!(\Z_{\tilde{T}}) \to q^+_!i_!i^*f_!(\Z_{\tilde{T}})\). 
\thref{lemm-Dem-Gm} shows that $\tilde T$ has a filtrable decomposition with cells isomorphic to affine spaces. Since the structure map $g \colon \A^1_S \rightarrow S$ satisfies $g_! \Z = \Z(-1)[-2]$, then \(q^+_!f_!(\Z_{\tilde{T}})\in \DTM(\Fl_w(\mathbf{a}_0,\mathbf{f})^0)^{\antic}\). It remains to show that \(q^+_!i_!i^*f_!(\Z_{\tilde{T}})\in \DTM(\Fl_w(\mathbf{a}_0,\mathbf{f})^0)^{\antic}\).
By base change, $i_!i^*f_!(\Z_{\tilde{T}})$ is a $*$-restriction of $m_! \Z$. Since $m_! \Z \in \DTM(\Fl_w(\mathbf{a}_0,\mathbf{f}))^{\antic}$ then we conclude by induction, as $Z$ is contained in a union of smaller affine Schubert cells.
\xpf

\subsubsection{Whitney--Tateness}\label{subsubsec-WT}
By \thref{lemm-Gm-on-Flf}, $\pi_0((\Fl_{\mathbf{f}})^+) = W/W_{\mathbf{f}}$. Let $\mathcal{S}_w$ be the reduced locus of the corresponding connected component, which we view as a locally closed sub ind-scheme $\mathcal{S}_w \subset \Fl_{\mathbf{f}}$. 

The following lemma is the equal characteristic analogue of \cite[Lemma 5.3]{AGLR:Local}. 
It can be proven in a similar way, but we sketch the proof for convenience of the reader.

\lemm
\thlabel{TrivialIntersection}
For any \(w\in W/W_\mathbf{f}\), there exists a regular cocharacter \(\Gm\to T\subseteq G\) such that for the induced \(\Gm\)-action, \(\mathcal{S}_w\cap \Fl_{w}(\mathbf{a}_0,\mathbf{f}) \cong S\).
\xlemm
\pf 
As in the proof of \thref{HypLocTate} we may assume that $\tilde w \in W_{\aff}$. Since the \(\mathcal{I}\)-orbit \(\Fl_w^\circ(\mathbf{a}_0,\mathbf{f})\subseteq \Fl_w(\mathbf{a}_0,\mathbf{f})\) is a \(\Gm\)-stable open subscheme containing the fixed point of $\mathcal{S}_w$, we have \(\mathcal{S}_w\cap \Fl_{w}(\mathbf{a}_0,\mathbf{f})  = \mathcal{S}_w\cap \Fl_{w}^\circ(\mathbf{a}_0,\mathbf{f})\).
Fix a reduced decomposition \(\dot{w} =s_1\ldots s_n\), and for each \(1\leq i\leq n\) consider the unique positive affine root \(\alpha_i\) corresponding to the wall separating the alcoves \(s_1\ldots s_{i-1}(\mathbf{a}_0)\) and \(s_1\ldots s_i(\mathbf{a}_0)\). 
By considering the Demazure resolution \(D(\dot{w})\to \Fl_w(\mathbf{a}_0,\mathbf{f})\), which is an isomorphism over \(\Fl_w^\circ(\mathbf{a}_0,\mathbf{f})\), one can show that \(\Fl_w^\circ(\mathbf{a}_0,\mathbf{f}) = L^+\Uu_{\alpha_1}\cdots L^+\Uu_{\alpha_i}\cdot w\) where the \(\Uu_{\alpha_i}\) are the corresponding affine root groups of \(G_{\mathbf{f}}\), cf. \cite[(5.11)]{AGLR:Local}.
Now, \cite[Corollary 5.6]{HainesNgo:Alcoves} shows that there exists a regular cocharacter repelling all the \(L^+\Uu_{\alpha_j}\). 
By the above, this cocharacter then also repels \(\Fl_w^\circ(\mathbf{a}_0,\mathbf{f})\), and we conclude that $\mathcal{S}_w \cap \Fl_w(\mathbf{a}_0,\mathbf{f}) \cong S$.
\xpf

For every regular cocharacter $\lambda \in X_*(T)$, we have the hyperbolic localization functor $(q_s^+)_! (p_s^+)^*u^! \colon \DM(\mathcal{I}\backslash \Fl_{\mathbf{f}}) \rightarrow \DM(\coprod_{W/W_{\mathbf{f}}} S)$ (the dependence on $\lambda$ is not reflected in the notation). The following conservativity result, used in the proof of \thref{theo-WT-BD} below, is similar to \cite[Propositions 6.4, 6.6]{AGLR:Local}.

\lemm \thlabel{lemm-hyp-loc-conservative}
Let $u^! \colon \DM(\mathcal{I} \backslash \Fl_{\mathbf{f}}) \rightarrow \DM(\Fl_{\mathbf{f}})$ be the forgetful functor.
\begin{enumerate}
    \item \label{item--lemm-hyp-loc-conservative1} Let $M \in \DM(\mathcal{I} \backslash \Fl_{\mathbf{f}})$ be a bounded motive. Then $u^!M$ is Tate (resp. compact, resp. anti-effective) if and only if $(q_s^+)_! (p_s^+)^* u^! M \in \DM(\coprod_{W/W_{\mathbf{f}}} S)$ is Tate (resp. compact, resp. anti-effective) for every regular cocharacter $\lambda$. 
    \item \label{item--lemm-hyp-loc-conservative2} Let $\phi \colon M_1 \rightarrow M_2$ be a map between bounded motives in $\DM(\mathcal{I} \backslash \Fl_{\mathbf{f}})$. Then $\phi$ is an isomorphism if and only if $(q_s^+)_! (p_s^+)^* u^! \phi$ is an isomorphism for every regular cocharacter $\lambda$.
\end{enumerate}
\xlemm

\pf
We prove both parts by noetherian induction on the support of the motives in $\DM(\mathcal{I} \backslash \Fl_{\mathbf{f}})$, starting with part (1). One direction is proved in \thref{HypLocTate}. For the other direction, let $\Fl_w^\circ(\mathbf{a}_0,\mathbf{f})$ be a maximal Schubert cell in the support of $M$. By \thref{TrivialIntersection} we may choose $\lambda \in X_*(T)$ such that \(\mathcal{S}_w\cap \Fl_w(\mathbf{a}_0,\mathbf{f})\cong S\). Then the hypothesis implies that the *-restriction of \(i_w^*u^!M\) to this copy of $S$ is Tate (resp. compact, resp. anti-effective), where \(i_w\colon \Fl_w^\circ(\mathbf{a}_0,\mathbf{f})\to \Fl_\mathbf{f}\) is the inclusion. Then we conclude by \thref{lemm--DMorbit} that \(i_w^* u^! M\) is Tate (resp. compact, resp. anti-effective). Now the motive $M'$ defined by the fiber sequence $M' \rightarrow u^! M \rightarrow i_{w*} i_w^* u^! M$ is $\mathcal{I}$-equivariant, and \((q^+_{s})_!(p^+_{s})^*M'\) is Tate (resp. compact, resp. anti-effective) for all regular cocharacters $\lambda$. Thus, we may conclude the proof of (1) by induction. The proof of part (2) is similar, using conservativity of $u^!$ and the fact that an isomorphism over \(\mathcal{S}_w\cap \Fl_w(\mathbf{a}_0,\mathbf{f})\cong S\) spreads out to an isomorphism over $\Fl_w^\circ(\mathbf{a}_0,\mathbf{f})$ by $\mathcal{I}$-equivariance and \thref{lemm--DMorbit}.
\xpf

We can now prove the Whitney--Tateness of \(\Gr_{\mathcal{G}_{\mathbf{f}}}(P)\).

\theo \thlabel{theo-WT-BD} Let $\mathbf{f}$ be a facet contained in the closure of $\mathbf{a}_0$ and let $P \in \{\emptyset, \zero, \ux, \ux \cup \zero \}$.
\begin{enumerate}
    \item \label{item--theo-WT-BD1} The stratification of \(\Gr_{\mathcal{G}_{\mathbf{f}}}(P)\) in \thref{defi-strat-Gr} is anti-effective universally Whitney--Tate.
    \item \label{item--theo-WT-BD2} The nearby cycles functor $\Upsilon_{P}$ in \thref{defi--nearby-hecke} preserves locally compact anti-effective Tate motives.
    \item \label{item--theo-WT-BD3} For $P_i \in \{\emptyset, \zero, \ux, \ux \cup \zero \}$ and $\mathcal{F}_i \in \DTM(\Gr_{\mathcal{G}_{\mathbf{f}}}(P_i))$, the K\"unneth map \refeq{Upsilon-box-map-multiple} is an isomorphism.
\end{enumerate}
\xtheo

\pf
First, we note that if $P = \emptyset$ then $ \Gr_{\mathcal{G}_{\mathbf{f}}}(\emptyset) = \Gm\sqcup \A^0$ is universally anti-effective Whitney--Tate by a computation involving relative purity. Moreover, (2) follows in this case since $\Upsilon_{\emptyset}$ is trivial on the constant family $C$ by \thref{lemm-A1-Upsilon}. Here we use that $\DTM(C_\eta)^{\xantic}$ is generated (under finite colimits, retracts and shifts) by $\Z_{C_\eta}$ and its twists.

The proofs of (1) and (2) for $P \in \{\zero, \ux, \ux \cup \zero\}$ are essentially the same; for simplicity we will 
only write down the case where $P = \ux$, so \(\Gr_{\mathcal{G}_{\mathbf{f}}}(\uP) = \Gr_{\mathcal{G}_{\mathbf{f}}}\).
By \cite[Proposition 3.7]{CassvdHScholbach:MotivicSatake} (and \cite[Corollary 3.10]{CassvdHScholbach:MotivicSatake} when $P =\ux \cup \zero$), $(\Gr_{\mathcal{G}_{\mathbf{f}}})_\eta$ and $(\Gr_{\mathcal{G}_{\mathbf{f}}})_s$ are anti-effective universally Whitney--Tate. 
In particular, we only need to understand what happens when we pass from $(\Gr_{\mathcal{G}_{\mathbf{f}}})_\eta$ to $(\Gr_{\mathcal{G}_{\mathbf{f}}})_s$. Thus, fix some \(\mu\in X_*(T)^+\) and consider the anti-effective Tate motive \(\Ff:=j_{\mu,*}(\Z)\), where \(j_\mu\colon \Gr^\mu\times C_\eta\to \Gr\times C_\eta \cong (\Gr_{\mathcal{G}_{\mathbf{f}}})_{\eta}\) is the associated immersion. We must show that $i^*j_*\Ff \in \DTM((\Gr_{\mathcal{G}_{\mathbf{f}}})_s)^{\anti}$, where $j$ and $i$ are defined below.

Note that the motive $\Ff$ is \((L^+G\times C_\eta)\)-equivariant, and in particular \(\mathbf{G}_{m,C_\eta}\)-equivariant for the action induced by any regular cocharacter \(\lambda \in X_*(T)\). We then have the following commutative diagram
\[\begin{tikzcd}
	(\Gr_{\mathcal{G}_{\mathbf{f}}})^0_{s} \arrow[d, "i^0"] & (\Gr_{\mathcal{G}_{\mathbf{f}}})^+_{s} \arrow[d, "i^+"] \arrow[l, "q^+_{s}"'] \arrow[r, "p^+_{s}"] & (\Gr_{\mathcal{G}_{\mathbf{f}}})_{s} \arrow[d, "i"]\\
	(\Gr_{\mathcal{G}_{\mathbf{f}}})^0 & (\Gr_{\mathcal{G}_{\mathbf{f}}})^+ \arrow[l, "q^+"'] \arrow[r, "p^+"] & \Gr_{\mathcal{G}_{\mathbf{f}}}\\
	(\Gr_{\mathcal{G}_{\mathbf{f}}})^0_{\eta} \arrow[u, "j^0"'] & (\Gr_{\mathcal{G}_{\mathbf{f}}})^+_{\eta} \arrow[u, "j^+"'] \arrow[l, "q^+_{\eta}"'] \arrow[r, "p^{+}_{\eta}"] & (\Gr_{\mathcal{G}_{\mathbf{f}}})_{\eta}, \arrow[u, "j"']
\end{tikzcd}\]
where the four small squares are cartesian by \cite[Corollary 1.16]{Richarz:Spaces}. Using base change and hyperbolic localization for the \(\GmC\)-equivariant motive \(j_*\Ff\), we compute 
\begin{align*}
q^+_{s!} p^{+*}_{s}( i^*j_*\Ff)  &\cong  q^+_{s!}  i^{+*} ( p^{+*} j_* \Ff) \cong i^{0*}q^+_! (p^{+*} j_*\Ff) \cong i^{0*}q^-_*(p^{-!}j_*\Ff) \\ 
 & \cong i^{0*}q^-_*(j^-_*p^{-!}_{\eta}\Ff) \cong i^{0*}j^0_*(q^-_{\eta *} p^{-!}_{\eta}\Ff) \cong i^{0*}j^0_*(q^+_{\eta!} p^{+*}_{\eta}\Ff).
\end{align*}

By \thref{HypLocTate}, the hyperbolic localization functor \(q^+_{\eta!}p^{+*}_{\eta}\) preserves anti-effective stratified Tate motives. Also, \(i^{0*}j^0_*\) preserves anti-effective stratified Tate motives by \thref{lemm-closure-of-fixed}\refit{lemm-closure-of-fixed3} (or \thref{lemm-Gm-on-Flf}\refit{lemm-Gm-on-Flf5} when $P = \zero$) and the corresponding result for \(\A^1=\Gm\sqcup \A^0\).
Thus, we see that \(q^+_{s!} p^{+*}_{s}( i^*j_*\Ff) \in \DTM(	(\Gr_{\mathcal{G}_{\mathbf{f}}})^0_{s})^{\anti}\).
As \(i^*j_*\Ff\) is $\mathcal{I}$-equivariant by smooth base change, we conclude by \thref{lemm-hyp-loc-conservative}\refit{lemm-hyp-loc-conservative1} that $i^*j_*\Ff\in \DTM(	\Gr_{\mathcal{G}_{\mathbf{f}}})_{s})^{\anti}$. Thus, \(\Gr_{\mathcal{G}_{\mathbf{f}}}\) is anti-effective Whitney--Tate stratified. 

The condition that this stratification is universal means that the base change map \cite[Eqn.~(2.4)]{CassvdHScholbach:MotivicSatake} is an isomorphism. Again, we only need to understand what happens when passing from $(\Gr_{\mathcal{G}_{\mathbf{f}}})_\eta$ to $(\Gr_{\mathcal{G}_{\mathbf{f}}})_s$. Since hyperbolic localization commutes with base change, then by \thref{lemm-hyp-loc-conservative}\refit{lemm-hyp-loc-conservative2} we again reduce to the case of \(\A^1=\Gm\sqcup \A^0\).

For (2) a similar argument applies, using \thref{lemm-hyp-loc-conservative}\refit{lemm-hyp-loc-conservative1} and that $\Psi_P$ commutes with hyperbolic localization by \thref{prop-Gm-Upsilon}, thus reducing us to the case of $\A^1$ where the result follows from \thref{lemm-A1-Upsilon}. Except for preservation of compactness, (2) also follows from (1) by \thref{Upsilon preserves Tate}.

For (3), let $X = \Gr_{\mathcal{G}_{\mathbf{f}}}(P_1) \times_C \cdots \times_C \Gr_{\mathcal{G}_{\mathbf{f}}}(P_n)$ and let $\mathcal{F}_i \in \DTM(\Gr_{\mathcal{G}_{\mathbf{f}}}(P_i))$.
Using the componentwise left actions of $L^+\mathcal{G}_{\mathbf{f}}^{\BD}$ on $X$, any $n$-tuple of regular cocharacters $(\lambda_1, \ldots, \lambda_n)$ induces a Zariski-locally linearizable $\GmC$-action on $X$
by \thref{lemm-Zar-lin}. For such a $\GmC$-action, the reduced closure of any connected component of $X_\eta^0$ is isomorphic to $C$ by \thref{lemm-Gm-on-Flf}\refit{lemm-Gm-on-Flf5} and \thref{lemm-closure-of-fixed}\refit{lemm-closure-of-fixed3}. Also, the $\calF_i$ and their box product are each $\Gm$-monodromic since strata in affine flag varieties are $\Gm$-stable.

The construction of the K\"unneth map \refeq{Upsilon-box-map-multiple} involves nearby cycles over $\A^n$, but this is only done in order to keep track of the monodromy. By inspecting the proof of \thref{Upsilon.boxtimes} and \thref{coro--nearby-pullback} one sees that the underlying map of motives in \refeq{Upsilon-box-map-multiple} can be constructed using a base change map over $\A^1$ (followed by compatibility of $\Upsilon$ with smooth pullback). In particular, to check that is an isomorphism we may also apply hyperbolic localization, since the latter commutes with $\Upsilon$ by \thref{prop-Gm-Upsilon} and with the formation of base change maps. This reduces us to the K\"unneth map for the trivial family over $\A^1$. Since the formation of K\"unneth maps commutes with colimits, we reduce to the case where $\calF_i = \Z_{\Gm}$ for all $i$,  in which case we have an isomorphism by \thref{lemm-A1-Upsilon}.
\xpf

\rema
Now that we know $\DTM(\Gr_{\mathcal{G}_{\mathbf{f}}})^{\antic}$ is well-defined, then since hyperbolic localization commutes with base change, it follows that \thref{HypLocTate} is also true for $\DTM(\Gr_{\mathcal{G}_{\mathbf{f}}})^{\antic}$.
\xrema

\subsection{The central functor}
\label{sect--centralfunctor}

\subsubsection{A functor to the center}
\theo\thlabel{CentralityIso}
For $\mathcal{A} \in \DTM(L^+G \backslash LG / L^+G)$ and $\mathcal{F} \in \DTM(L^+G_{\mathbf{f}} \backslash LG / L^+G_{\mathbf{f}})$ there are canonical isomorphisms
\begin{equation} \label{eqn--central} \mathsf{Z}_{\mathbf{f}}(\mathcal{A}) \star \mathcal{F} \cong \Upsilon_{\ux \cup \zero}(\mathcal{A} \boxtimes \mathcal{F} \boxtimes \Z_{C_\eta}) \cong  \mathcal{F} \star \mathsf{Z}_{\mathbf{f}}(\mathcal{A})\end{equation} in $\NilpQ \DTM(L^+G_{\mathbf{f}} \backslash LG / L^+G_{\mathbf{f}})$. The composites of these isomorphisms are compatible with the associativity constraint on $\Ho (\DTM(L^+G_{\mathbf{f}} \backslash LG / L^+G_{\mathbf{f}}))$ in the sense that the following diagram commutes for $\mathcal{A} \in \DTM(L^+G \backslash LG / L^+G)$ and $\mathcal{F}_1$, $\mathcal{F}_2 \in \DTM(L^+G_{\mathbf{f}} \backslash LG / L^+G_{\mathbf{f}})$.
$$\xymatrix{
 & \mathsf{Z}_{\mathbf{f}}(\mathcal{A}) \star (\mathcal{F}_1 \star \mathcal{F}_2) \ar[r]^{\eqref{eqn--central}} & (\mathcal{F}_1 \star \mathcal{F}_2) \star \mathsf{Z}_{\mathbf{f}}(\mathcal{A}) \ar[rd] &  \\ (\mathsf{Z}_{\mathbf{f}}(\mathcal{A}) \star \mathcal{F}_1) \star \mathcal{F}_2 \ar[ru] \ar[rd]^{\eqref{eqn--central}} &  &  & \mathcal{F}_1 \star (\mathcal{F}_2 \star \mathsf{Z}_{\mathbf{f}}(\mathcal{A}))  \\
 & (\mathcal{F}_1 \star \mathsf{Z}_{\mathbf{f}}(\mathcal{A})) \star \mathcal{F}_2  \ar[r] & \mathcal{F}_1 \star (\mathsf{Z}_{\mathbf{f}}(\mathcal{A}) \star \mathcal{F}_2)  \ar[ru]^{\eqref{eqn--central}} & 
}$$

\xtheo

\pf
The isomorphisms are constructed as in \cite[Proposition 3.2.1]{AcharRiche:Central}. Briefly, the first isomorphism comes from applying \thref{theo-WT-BD}\refit{theo-WT-BD3} and \thref{lemm-conv-nearby} in the case $P_1 = \ux$ and $P_2 = \zero$. For the second we take $P_1 = \zero$ and $P_2 = \ux$, noting that we have a canonical isomorphism $(\Gr_{\mathcal{G}_{\mathbf{f}}}^\BD)_{\eta} = \Gr \times \Fl_{\mathbf{f}} \times C_\eta \cong \Fl_{\mathbf{f}} \times \Gr \times C_\eta$ which swaps the factors. We also use the fact that $\Psi_{\zero}$ is trivial on the sheaves in \thref{nearby trivial family}, where \cite[Corollary 3.10]{CassvdHScholbach:MotivicSatake} ensures that the hypothesis on the K\"unneth map is satisfied. The commutativity of the diagram is proved exactly as in \cite[Theorem 3.2.3]{AcharRiche:Central}, using \thref{theo-WT-BD}\refit{theo-WT-BD3} for $n > 2$ factors along with \thref{rema--multiple-factors}.
\xpf

\rema 
If $\mathbf{f} = \mathbf{f}_0$ then $\mathsf{Z}_{\mathbf{f}_0} = \id$ by \thref{nearby trivial family} and \cite[Corollary 3.10]{CassvdHScholbach:MotivicSatake}. Thus, \thref{CentralityIso} endows $\Ho(\DTM(L^+G \backslash LG / L^+G))$ with the structure of a symmetric monoidal category. In \thref{prop--fusion-v-nearby} we compare this to the commutativity constraint on $\MTM(L^+G \backslash LG / L^+G)$ constructed in \cite[Proposition 5.51]{CassvdHScholbach:MotivicSatake} via the fusion interpretation.
\xrema

\rema \thlabel{rem--functor-to-center}
For general $\mathbf{f}$, the isomorphisms in \thref{CentralityIso} endow $\mathsf{Z}_{\mathbf{f}}$ with a factorization through the Drinfeld center 
$\mathcal{Z} (\Ho (\DTM(L^+G_{\mathbf{f}} \backslash LG / L^+G_{\mathbf{f}})))$, cf. \cite[\S 3.1]{AcharRiche:Central}: $$\mathsf{Z}_{\mathbf{f}} \colon \Ho(\DTM(L^+G \backslash LG / L^+G)) \rightarrow \mathcal{Z}(\Ho (\DTM(L^+G_{\mathbf{f}} \backslash LG / L^+G_{\mathbf{f}}))).$$ 
\xrema

\subsubsection{Monoidality}
\theo \thlabel{theo-moniodal-functor}
For $\mathcal{A}_1$, $\mathcal{A}_2 \in \DTM(L^+G \backslash LG / L^+G)$ there is a canonical isomorphism
\begin{equation} \label{eqn--monoidal} \mathsf{Z}_{\mathbf{f}}(\mathcal{A}_1) \star \mathsf{Z}_{\mathbf{f}}(\mathcal{A}_2) \cong \mathsf{Z}_{\mathbf{f}}(\mathcal{A}_1 \star \mathcal{A}_2)\end{equation} in $\NilpQ \DTM(L^+G_{\mathbf{f}} \backslash LG / L^+G_{\mathbf{f}})$. These isomorphisms are compatible with the centrality constraints \eqref{eqn--central} in the sense that the following diagram commutes for $\mathcal{A}_1$, $\mathcal{A}_2 \in \DTM(L^+G \backslash LG / L^+G)$ and $\mathcal{F} \in \DTM(L^+G_{\mathbf{f}} \backslash LG / L^+G_{\mathbf{f}})$ (we have omitted the associativity isomorphisms, but see \cite[Eqn.~(3.4.1)]{AcharRiche:Central} for the full diagram).
$$\xymatrix{
\mathsf{Z}_{\mathbf{f}}(\mathcal{A}_1 \star \mathcal{A}_2) \star \mathcal{F} \ar[rr]^{\eqref{eqn--central}} \ar[d]^{\eqref{eqn--monoidal}} & & \mathcal{F} \star \mathsf{Z}_{\mathbf{f}}(\mathcal{A}_1 \star \mathcal{A}_2) \\ \mathsf{Z}_{\mathbf{f}}(\mathcal{A}_1) \star \mathsf{Z}_{\mathbf{f}}(\mathcal{A}_2) \star \mathcal{F} \ar[r]^{\eqref{eqn--central}} & \mathsf{Z}_{\mathbf{f}}(\mathcal{A}_1) \star \mathcal{F} \star \mathsf{Z}_{\mathbf{f}}(\mathcal{A}_2) \ar[r]^{\eqref{eqn--central}} & \mathcal{F} \star \mathsf{Z}_{\mathbf{f}}(\mathcal{A}_1) \star \mathsf{Z}_{\mathbf{f}}(\mathcal{A}_2) \ar[u]^{\eqref{eqn--monoidal}}
}$$
\xtheo

\pf
The isomorphism \eqref{eqn--monoidal} comes from applying \thref{theo-WT-BD}\refit{theo-WT-BD3} and \thref{lemm-conv-nearby} in the case $P_1= P_2 = \ux$. The commutativity of the diagram is proved exactly as in \cite[Theorem 3.4.1]{AcharRiche:Central}.
\xpf

\rema
\thref{theo-moniodal-functor} endows $\mathsf{Z}_{\mathbf{f}}$ with the structure of a \emph{monoidal} functor to the Drinfeld center of $\Ho (\DTM(L^+G_{\mathbf{f}} \backslash LG / L^+G_{\mathbf{f}}))$.
\xrema

\subsubsection{Change of facet}
\theo \thlabel{thm--change-of-facet}
Let $\mathbf{c}$ be a facet contained in the closure of $\mathbf{f}$ and let $\pi \colon \Fl_{\mathbf{f}} \rightarrow \Fl_{\mathbf{c}}$ be the projection. 
\begin{enumerate}
    \item For $\mathcal{A} \in \DTM(L^+G \backslash LG / L^+G)$ there is a canonical isomorphism
      $$\pi_* \mathsf{Z}_{\mathbf{f}}(\mathcal{A}) \cong \mathsf{Z}_{\mathbf{c}}(\mathcal{A})$$
    of the underlying non-equivariant motives in $\NilpQ \DTM(\Fl_{\mathbf{c}})$.
    \item For $\mathcal{A}_1$, $\mathcal{A}_2 \in \DTM(L^+G \backslash LG / L^+G)$ there is a canonical isomorphism
    \begin{equation} \label{eqn--change-of-facet} \pi_* (\mathsf{Z}_{\mathbf{f}}(\mathcal{A}_1) \star \mathsf{Z}_{\mathbf{f}}(\mathcal{A}_2)) \cong \mathsf{Z}_{\mathbf{c}}(\mathcal{A}_1) \star \mathsf{Z}_{\mathbf{c}}(\mathcal{A}_2)
    \end{equation}
    of the underlying non-equivariant motives in $\NilpQ \DTM(\Fl_{\mathbf{c}})$.
    \item There is a commutative diagram of isomorphisms in $\Ho (\DTM(L^+G_{\mathbf{c}} \backslash LG / L^+G_{\mathbf{c}}))$ as follows:
        $$\xymatrix{
     \pi_*(\mathsf{Z}_{\mathbf{f}}(\mathcal{A}_1) \star \mathsf{Z}_{\mathbf{f}}(\mathcal{A}_2)) \ar[d]^{\eqref{eqn--change-of-facet}} \ar[r]^{\pi_* \eqref{eqn--central}} & \pi_*(\mathsf{Z}_{\mathbf{f}}(\mathcal{A}_2) \star \mathsf{Z}_{\mathbf{f}}(\mathcal{A}_1)) \ar[d]^{\eqref{eqn--change-of-facet}} \\  \mathsf{Z}_{\mathbf{c}}(\mathcal{A}_1) \star \mathsf{Z}_{\mathbf{c}}(\mathcal{A}_2) \ar[r]^{\eqref{eqn--central}}  &  \mathsf{Z}_{\mathbf{c}}(\mathcal{A}_2) \star \mathsf{Z}_{\mathbf{c}}(\mathcal{A}_1).
    }$$
\end{enumerate}
\xtheo

\rema
\thlabel{Z section}
If $\mathbf{f} = \mathbf{a}_0$ and $\mathbf{c} = \mathbf{f}_0$ then part (1) provides a canonical identification of $\pi_* \mathsf{Z}_{\mathbf{a}_0}$ with the forgetful functor $\DTM(L^+G \backslash LG / L^+G) \rightarrow \DTM(\Gr_G)$.
\xrema

\pf Part (1) follows from the fact that $\Upsilon$ commutes with pushforward along the proper map $\Gr_{\mathcal{G}_{\mathbf{f}}} \rightarrow \Gr_{\mathcal{G}_{\mathbf{c}}}$. The isomorphism \eqref{eqn--change-of-facet} in part (2) is constructed by combining part (1) and \eqref{eqn--monoidal}. Alternatively, \eqref{eqn--change-of-facet} may be constructed as follows. By combining part (1) with base change and the K\"unneth formula for the proper map $\pi \widetilde\times  \pi \colon LG \times^{L^+G_{\mathbf{f}}} \Fl_{\mathbf{f}} \rightarrow LG \times^{L^+G_{\mathbf{c}}} \Fl_{\mathbf{c}}$ we get an isomorphism of non-equivariant motives,
$$(\pi \widetilde\times  \pi)_{!} ( \mathsf{Z}_{\mathbf{f}}(\mathcal{A}_1) \widetilde \boxtimes \mathsf{Z}_{\mathbf{f}}(\mathcal{A}_2)) \cong ( \mathsf{Z}_{\mathbf{c}}(\mathcal{A}_1) \widetilde \boxtimes \mathsf{Z}_{\mathbf{c}}(\mathcal{A}_2)).$$ Applying $m_!$ also gives an isomorphism \eqref{eqn--change-of-facet}, which agrees with the previous construction by the same arguments as in \cite[Lemma 3.4.3]{AcharRiche:Central}. Finally, part (3) follows as in \cite[Lemma 3.4.4]{AcharRiche:Central} by compatibility of $\Upsilon$ with proper pushforward (e.g, along $\Gr_{\mathcal{G}_{\mathbf{f}}} \rightarrow \Gr_{\mathcal{G}_{\mathbf{c}}}$).
\xpf

\subsubsection{Mixed Tate motives and the motivic Satake equivalence}
\label{sect--MTM}
In this subsection we assume that $S$ satisfies the Beilinson--Soulé conjecture (this is unnecessary when working with reduced motives).
Recall from \refsect{mixed Tate motives} that $\DTM(S)$ carries a t-structure with heart $\MTM(S)$. By \cite[Lemma 2.13]{CassvdHScholbach:MotivicSatake} the $\Z(k)$ for $k \in \Z$ compactly generate $\MTM(S)$. By \cite[Lemma 2.15]{CassvdHScholbach:MotivicSatake} the strata $\Fl_w^{\circ}(\mathbf{f}, \mathbf{f})$ also satisfy the Beilinson--Soulé conjecture. We normalize such that $\Z[\dim \Fl_w^{\circ}(\mathbf{f}, \mathbf{f})] \in \MTM(\Fl_w^{\circ}(\mathbf{f}, \mathbf{f}))$. By \cite[Lemma 2.15]{CassvdHScholbach:MotivicSatake}, the t-structures on these strata glue to a t-structure on $\DTM(\Fl_{\mathbf{f}})$, as outlined in \refsect{mixed Tate motives}. By \cite[Lemma 2.26]{CassvdHScholbach:MotivicSatake}, there is also a t-structure on $\DTM(L^+G_{\mathbf{f}} \backslash \Fl_{\mathbf{f}})$ characterized by t-exactness of the forgetful functor to $\DTM(\Fl_{\mathbf{f}})$.

\defilemm
Let $\mathbf{f}$ be a facet and let $P \in \{\emptyset, \ux, \zero, \ux \cup \zero\}$.
\begin{enumerate}
\item There is a t-structure on $\DTM(\Gr_{\mathcal{G}_{\mathbf{f}}}(P))$ characterized by the property that the aisle $\DTM(\Gr_{\mathcal{G}_{\mathbf{f}}}(P))^{\leq 0}$ is generated by the objects $\iota_{!} \Z_{X}(k)[\dim_S X]$, for all the strata $\iota : X \r \Gr_{\mathcal{G}_{\mathbf{f}}}(P)$ in \thref{defi-strat-Gr} and $k \in \Z$. Hence the abelian subcategory $\MTM(\Gr_{\mathcal{G}_{\mathbf{f}}}(P)) \subset \DTM(\Gr_{\mathcal{G}_{\mathbf{f}}}(P))$ of mixed Tate motives is well-defined. The same is true for the Hecke stack $\Hck_{\mathcal{G}_{\mathbf{f}}}(P)$, where the t-structure is characterized by t-exactness of the forgetful functor to $\DTM(\Gr_{\mathcal{G}_{\mathbf{f}}}(P))$.
\item The forgetful functor $\MTM(\Hck_{\mathcal{G}_{\mathbf{f}}}(P)) \rightarrow \MTM(\Gr_{\mathcal{G}_{\mathbf{f}}}(P))$ is fully faithful, and the image is stable under subquotients.
\item The perverse truncation functors preserve anti-effective motives, so that the full abelian subcategory $\MTM(\Hck_{\mathcal{G}_{\mathbf{f}}}(P))^{\anti} \subset \MTM(\Hck_{\mathcal{G}_{\mathbf{f}}}(P))$ is well-defined.
\end{enumerate}
\xdefilemm

\pf We first observe that the stratification $\A^1 = \Gm \sqcup \A^0$ is Whitney--Tate and cellular relative to $S$ in the sense of \cite[Remark 2.10]{CassvdHScholbach:MotivicSatake}, so that $\MTM(C)$ is well-defined by \cite[Lemma 2.15]{CassvdHScholbach:MotivicSatake}. Next, we note that the Schubert cells $\Fl_w^\circ(\mathbf{f}, \mathbf{f})$ are smooth, and the finer stratifications by Iwahori-orbits are cellular and Whitney--Tate by \cite[Proposition 3.7]{CassvdHScholbach:MotivicSatake}. By \cite[Remark 2.10]{CassvdHScholbach:MotivicSatake}, this implies that the stratification of $\Gr_{\mathcal{G}_{\mathbf{f}}}(P)$ in \thref{defi-strat-Gr} is admissible relative to $C$ in the sense of \cite[Definition 2.9]{CassvdHScholbach:MotivicSatake}. The t-structures in part (1) 
are then constructed as in \cite[Lemma 2.15]{CassvdHScholbach:MotivicSatake} (for $\Gr_{\mathcal{G}_{\mathbf{f}}}(P)$) and in \cite[Lemma 2.26]{CassvdHScholbach:MotivicSatake} (for $\Hck_{\mathcal{G}_{\mathbf{f}}}(P)$). In the latter case we also need \cite[Proposition A.4.9]{RicharzScholbach:Intersection}, which implies that the kernel of $L^{+}\mathcal{G}_{\mathbf{f}}^{\BD} \rightarrow  L^{+,(0)}\mathcal{G}_{\mathbf{f}}^{\BD}$ is split pro-unipotent. Part (2) is proved exactly as in \cite[Lemma 4.30]{CassvdHScholbach:MotivicSatake}, except we replace the stratification of $\A^I$ by partial diagonals by the stratification $\A^1 = \Gm \sqcup \A^0$. Finally, the construction in (1) carries over to anti-effective motives since the stratification of $\Gr_{\mathcal{G}_{\mathbf{f}}}(P)$ is anti-effective, giving (3).
\xpf

\prop \thlabel{prop--satake-exact}
Let $\mathbf{f} = \mathbf{f}_0$, so that $\mathcal{G}_{\mathbf{f}_0} = G_C$. Let $P = \ux \cup \zero$, and let $j \colon \Hck_{G_C}(P)_{\eta} \rightarrow \Hck_{G_C}(P)$ and $i \colon \Hck_{G_C}(P)_{s} \rightarrow \Hck_{G_C}(P)$ be the inclusions.
\begin{enumerate}
\item \label{item--satake-exact1} The following nearby cycles functor is t-exact
$$\Upsilon_{P}((-) \boxtimes \Z_{C_\eta})  \colon \DTM_{L^+G \times L^+G}(\Gr \times \Gr) \rightarrow \DTM_{L^+G}(\Gr).$$
\item \label{item--satake-exact2} The following functors are also t-exact $$j_!((-)  \boxtimes \Z_{C_\eta}[1]), \: j_*((-)  \boxtimes \Z_{C_\eta}[1]) \colon \DTM_{L^+G \times L^+G}(\Gr \times \Gr) \rightarrow \DTM(\Hck_{G_C}(P)).$$
\item \label{item--satake-exact3}For $\mathcal{A} \in \MTM_{L^+G \times L^+G}(\Gr \times \Gr)$ we have a canonical isomorphism
$$
\Upsilon_{P}(\mathcal{A} \boxtimes \Z_{C_\eta}) \cong \pH^{-1} (i^*j_{!*} (\mathcal{A} \boxtimes \Z_{C_\eta}[1])).
$$
\end{enumerate}
\xprop

\pf By continuity we may restrict to bounded motives. To prove (1), we may compose with the t-exact and conservative constant term functor $\CT_B[\deg]$ by \cite[Proposition 5.10]{CassvdHScholbach:MotivicSatake}. This commutes with $\Upsilon_{P}$ by \thref{prop-Gm-Upsilon}, where we apply $\CT_{B \times B}[\deg_{B \times B}] \boxtimes \id_{\DTM(C_\eta)}$ over $C_\eta$. (In contrast to \cite{CassvdHScholbach:MotivicSatake}, we will not include the degree shift in the definition of \(\CT_B\).) This crucially uses that $\Gr_{G_C}(P)_s = \Gr$.
By \thref{lemm-closure-of-fixed}\refit{lemm-closure-of-fixed3}, the reduced locus of the closure of each connected component of $(\Gr_{G_C}(P))^0_\eta$ is isomorphic to $C$. This reduces (1) to the case where $G$ is trivial, in which case the result follows from \thref{nearby trivial family} (where we set $Y=S$).

To prove (2), we note that $i^*j_*$ and $i^!j_!$ commute with constant terms (both assertions follow from base change since $\CT_B = q^+_! p^{+*} = q^-_*p^{-!}$). The same is true for $i^!j_*$ and $i^*j_!$ since both of these are zero. Thus, we may similarly apply constant terms to reduce to the case where $G$ is trivial, in which case (2) follows from relative purity on $\A^1 = \Gm \sqcup \A^0$.

To prove (3),  we argue as in \cite[Lemma 9.1.9]{AcharRiche:Central}. 
Let $\mathcal{F} = \mathcal{A} \boxtimes \Z_{C_\eta}[1]$.
The canonical map $i^* j_* \calF \r \Upsilon_P(\calF)$ (\thref{trivial monodromy}) is a map between objects in $\DTM(\Gr)$, where the latter lies in cohomological degree $-1$ by part (1).
Applying $\pH^{-1}$ then gives a map $\pH^{-1} i^* j_* \calF \r \Upsilon_P(\calF)$ which we claim is an isomorphism.
Indeed, applying the t-exact and conservative constant term functor $\CT_B$ commutes with everything in sight.
Thus, to check this is an isomorphism, we may assume $G = T$, in which case the claim holds by \thref{nearby trivial family} and \thref{lemm-closure-of-fixed}\refit{lemm-closure-of-fixed3}.
By (2) and the exact triangle $j_! \mathcal{F} \rightarrow j_* \mathcal{F} \rightarrow i_*i^* j_* \mathcal{F}$, we may identify $\pH^{-1}(i^*j_* \mathcal{F})$ with the kernel of $j_! \mathcal{F} \rightarrow j_* \mathcal{F}$. This is the same as the kernel of $j_! \mathcal{F} \rightarrow j_{!*} \mathcal{F}$, and we conclude by using the exact triangle $j_! \mathcal{F} \rightarrow j_{!*} \mathcal{F} \rightarrow i_*i^* j_{!*} \mathcal{F}$.
\xpf

\prop \thlabel{prop--fusion-v-nearby}
Let $\mathcal{A}_1$, $\mathcal{A}_2 \in \MTM(L^+G \backslash LG / L^+G)$. Then the commutativity isomorphism $\pH^0 (\mathcal{A}_1 \star \mathcal{A}_2) \cong\pH^0( \mathcal{A}_2 \star \mathcal{A}_1)$ constructed via the fusion product in \cite[Proposition 5.51]{CassvdHScholbach:MotivicSatake}, omitting the sign modifications in \cite[Remark 5.48]{CassvdHScholbach:MotivicSatake}, agrees with the perverse truncation of  \eqref{eqn--central} for $\mathbf{f} = \mathbf{f}_0$.
\xprop

\pf
The unmodified commutativity isomorphism in \cite[Proposition 5.51]{CassvdHScholbach:MotivicSatake} is constructed as follows. Let $I=\{1,2\}$ and consider the Beilinson--Drinfeld Grassmannian $\Gr_{G_C, I} \rightarrow C^I$. Let $\Delta C \subset C^I$ be the diagonal and let $C^\circ \subset C^I$ be the complement. Then the fusion isomorphism implies that
$$\restr{\Gr_{G_C, I}}{\Delta C} = \Gr \times C, \quad \restr{\Gr_{G_C, I}}{C^\circ} = \Gr \times \Gr \times C^\circ.$$ Let $i_I \colon \restr{\Gr_{G_C, I}}{\Delta C} \rightarrow \Gr_{G_C, I}$ and $j_I \colon \restr{\Gr_{G_C, I}}{C^\circ} \rightarrow \Gr_{G_C, I}$ be the inclusions. There is a canonical isomorphism
$$i_I^* \circ j_{I, !*} (\pH^0 (\mathcal{A}_1 \boxtimes \mathcal{A}_2) \boxtimes \Z_{\C^\circ}[2] ) \cong \pH^0 (\mathcal{A}_1 * \mathcal{A}_2) \boxtimes \Z_{\Delta C}[2].$$ The commutativity isomorphism is constructed by observing that $i_I^* \circ j_{I, !*}$ is invariant under the isomorphism of $\Gr_{G_C, I}$ which swaps the coordinates of $C^I$. 

To relate this to nearby cycles, note that there is a canonical isomorphism
$$\Gr_{G_C, I} \cong C \times \Gr_{G_C}^{\BD}$$ which identifies $\restr{\Gr_{G_C, I}}{\Delta C}$ with $\{0\} \times \Gr_{G_C}^{\BD}$. This is constructed using the additive structure on $C$, cf. \cite[Lemma 3.3.6]{AcharRiche:Central}. Since the extra factor of $C$ may be absorbed into the base $S$, the unmodified commutativity isomorphism may alternatively be constructed using the canonical isomorphism  $$i^* j_{!*}(\pH^0 (\mathcal{A}_1 \boxtimes \mathcal{A}_2) \boxtimes \Z_{C_\eta}[1] )) \cong  \pH^0 (\mathcal{A}_1 * \mathcal{A}_2)[1].$$ By \thref{prop--satake-exact}\refit{satake-exact3}, this agrees with the perverse truncation of the commutativity isomorphism constructed using nearby cycles.
\xpf

\rema
The computation of constant terms $\CT_B[\deg] \colon \MTM(L^+G \backslash LG / L^+G)  \rightarrow \MTM(L^+T \backslash LT / L^+T)$ using nearby cycles in \cite[Proposition 5.36]{CassvdHScholbach:MotivicSatake} assembles to a monoidal structure on $\CT_B[\deg]$. By arguments similar to the proof of \thref{prop--fusion-v-nearby}, this monoidal structure coincides with the one constructed in \cite[Proposition 5.51]{CassvdHScholbach:MotivicSatake}, again omitting the sign modifications in \cite[Remark 5.48]{CassvdHScholbach:MotivicSatake}. Similar remarks apply to the fiber functor to $\MTM(C)$ obtained by pushforward, cf.~\cite[Proposition 3.3.8]{AcharRiche:Central}.
\xrema

\subsubsection{Centrality}
\theo \thlabel{theo-E2}
Let $\mathcal{A}$, $\mathcal{B} \in \MTM(L^+G \backslash LG / L^+G)$ be such that $\mathcal{A} \boxtimes \mathcal{B}$ is mixed Tate. Then the following diagram commutes in $\Ho (\DTM(L^+G_{\mathbf{f}} \backslash LG / L^+G_{\mathbf{f}}))$, where the subscripts indicate the facet to which an isomorphism is applied.
    \begin{equation} \label{eqn-central-iso} \xymatrix{
     \mathsf{Z}_{\mathbf{f}}(\mathcal{A}) \star \mathsf{Z}_{\mathbf{f}}(\mathcal{B}) \ar[d]^{\eqref{eqn--central}_{\mathbf{f}}} \ar[r]^-{\eqref{eqn--monoidal}_{\mathbf{f}}} & \mathsf{Z}_{\mathbf{f}}(\mathcal{A} \star \mathcal{B}) \ar[d]^{\mathsf{Z}_{\mathbf{f}} \circ \eqref{eqn--central}_{\mathbf{f}_0}} \\ \mathsf{Z}_{\mathbf{f}}(\mathcal{B}) \star \mathsf{Z}_{\mathbf{f}}(\mathcal{A}) \ar[r]^-{\eqref{eqn--monoidal}_{\mathbf{f}}}  & \mathsf{Z}_{\mathbf{f}}(\mathcal{B} \star \mathcal{A})
    } \end{equation}

\xtheo

\rema
To construct the isomorphisms $\eqref{eqn--central}_{\mathbf{f}}$ and $\mathsf{Z}_{\mathbf{f}} \circ \eqref{eqn--central}_{\mathbf{f}_0}$ we must choose which factor should be considered central for the purpose of commuting it with the other factor. By \thref{prop--fusion-v-nearby} the isomorphism $\mathsf{Z}_{\mathbf{f}} \circ \eqref{eqn--central}_{\mathbf{f}_0}$ is independent of this choice, and the proof will show that $\eqref{eqn--central}_{\mathbf{f}}$ is also independent of this choice, which amounts to choosing an order of the coordinates on $C^2$ for which to perform iterated nearby cycles.
\xrema

\rema
Note that $\mathcal{A} \boxtimes \mathcal{B}$ is always mixed Tate if one restricts to coefficients in a field. Even if $\mathcal{A} \boxtimes \mathcal{B}$ is not mixed Tate, it is reasonable to expect commutativity of the diagram analogous to \eqref{eqn-central-iso} where one takes the perverse truncation of all convolution products. The proof of this commutativity in \cite[Theorem 3.5.1]{AcharRiche:Central} uses t-exactness of $\Upsilon_{\ux \cup \zero}[-1]$, which is currently unavailable in the motivic setting.
\xrema

\pf
The proof is analogous to \cite[Theorem 3.5.1]{AcharRiche:Central}. We will construct all of the necessary isomorphisms; commutativity may then be checked as in \emph{loc.~cit.} 
We fix $I = \{1,2\}$ throughout the proof. Recall the global convolution Grassmannian $m_{i_0} \colon \widetilde \Gr_{\mathcal{G}_{\mathbf{f}}, I} \rightarrow \Gr_{\mathcal{G}_{\mathbf{f}}, I}$ from \thref{def--i_0BDconv}, where $i_0 = 2$. In \thref{lemm--bundle-prod} we also constructed $L_I^+\mathcal{G}_{\mathbf{f}}$-torsors
$$\Gr_{\mathcal{G}_{\mathbf{f}}} \times \Gr_{\mathcal{G}_{\mathbf{f}}} \xleftarrow{p} (\Gr_{\mathcal{G}_{\mathbf{f}}} \times C)^{(\infty)} \times_{C^2} (\Gr_{\mathcal{G}_{\mathbf{f}}} \times C) \xrightarrow{q} \widetilde \Gr_{\mathcal{G}_{\mathbf{f}}, I}.$$ Let us denote the nearby cycles functor on $\Gr_{\mathcal{G}_{\mathbf{f}}, I}$ (resp. $\widetilde \Gr_{\mathcal{G}_{\mathbf{f}}, I}$) by $\Upsilon_I$ (resp. $\Upsilon_{\tilde I}$). 

Let $\mathcal{A}' = \mathcal{A} \boxtimes \Z_{\C_\eta} \in \DTM((\Gr_{\mathcal{G}_{\mathbf{f}}})_\eta)$, and define $\mathcal{B}'$ likewise. While \cite[\S 4]{CassvdHScholbach:MotivicSatake} works over the base $C^2$, it applies equally well over $C_\eta^2 = \Gm \times \Gm$. Thus, by \cite[Lemma 4.16]{CassvdHScholbach:MotivicSatake} we may form the twisted product $\mathcal{A}' \widetilde \boxtimes \mathcal{B}' \in \DM((\widetilde \Gr_{\mathcal{G}_{\mathbf{f}}, I})_\eta)$. Associated to the usual order on $I$  have the iterated nearby cycles functor $\Upsilon_{\tilde{I}, 1} \circ \Upsilon_{\tilde{I}, 2}$. If $\Delta$ is the inclusion of diagonal over $C_\eta^2$ then we also have the nearby cycles functor $\Upsilon_{(\ux, \ux)} \circ \Delta^*$, since the fiber of $\widetilde \Gr_{\mathcal{G}_{\mathbf{f}}, I}$ over $\Delta$ is $\Gr_{\mathcal{G}_{\mathbf{f}}}(\ux, \ux)$.

By \thref{lemm-composition} and \thref{coro--nearby-pullback} there are maps
\begin{equation} \label{eqn--NearbyTop} \Upsilon_{\tilde{I}, 1} \circ \Upsilon_{\tilde{I}, 2} (\mathcal{A}' \widetilde \boxtimes \mathcal{B}') \leftarrow \Upsilon_{\tilde I}(\mathcal{A}' \widetilde \boxtimes \mathcal{B}') \rightarrow \Upsilon_{(\ux, \ux)}\Delta^*(\mathcal{A}' \widetilde \boxtimes \mathcal{B}').\end{equation}
We claim that these are isomorphisms. Both maps are compatible with smooth base change, and $p$ and $q$ admit sections Zariski-locally over $C^2_\eta$ and $C^2_s$ (over $C_\eta^2$ this is \cite[Lemma 4.8]{CassvdHScholbach:MotivicSatake}). Thus, it suffices to check that the analogous maps are isomorphisms for 
$\mathcal{A}' \boxtimes \mathcal{B}' \in \DM( \Gr_{\mathcal{G}_{\mathbf{f}}}^2)_\eta$.
For this we proceed as in \thref{theo-WT-BD}, using \thref{lemm-hyp-loc-conservative} and hyperbolic localization to reduce to the case of Tate motives supported on $\Gm^2 \subset \A^2$. By \thref{Upsilon.boxtimes} we have $\Upsilon_{\A^2/ \A^2}(\Z_{\Gm^2}) = \Z_S$, so that we indeed have isomorphisms by tracing through the construction of these maps.

We also claim that there are canonical isomorphisms
\begin{equation} \label{eqn--Conv} \Upsilon_{\tilde I}(\mathcal{A}' \widetilde \boxtimes \mathcal{B}') \cong \Upsilon_{\ux}(\mathcal{A}') \widetilde \boxtimes \Upsilon_{\ux}(\mathcal{B}') \cong \mathsf{Z}_{\mathbf{f}}(\mathcal{A}) \widetilde \boxtimes \mathsf{Z}_{\mathbf{f}}(\mathcal{B}).\end{equation} 
The second is a definition, and the first is proved similarly to \thref{theo-WT-BD}\refit{theo-WT-BD3}, using $p$ and $q$ to reduce to showing that the K\"unneth map $\Upsilon_{\ux}(\mathcal{A}') \boxtimes \Upsilon_{\ux}(\mathcal{B}') \rightarrow \Upsilon_{ \Gr_{\mathcal{G}_{\mathbf{f}}}^2}(\mathcal{A}' \boxtimes \mathcal{B}')$ in \thref{Upsilon.boxtimes} is an isomorphism. This latter fact can also be proved using hyperbolic localization, so this constructs \eqref{eqn--Conv}. 

Since the maps in \eqref{eqn--NearbyTop} are compatible with proper pushforward, we have a commutative diagram of isomorphisms
\begin{equation} \label{eqn--centrality} \xymatrix{ 
m_{s!} (\Upsilon_{\tilde{I}, 1} \circ \Upsilon_{\tilde{I}, 2} (\mathcal{A}' \widetilde \boxtimes \mathcal{B}') ) \ar[d]^-{\sim} & \ar[l]_-{\sim} m_{s!} \Upsilon_{\tilde I}(\mathcal{A}' \widetilde \boxtimes \mathcal{B}') \ar[d]^-{\sim} \ar[r]^-\sim & m_{s!} \Upsilon_{(\ux, \ux)}\Delta^*(\mathcal{A}' \widetilde \boxtimes \mathcal{B}') \ar[d]^-{\sim} \\
 \Upsilon_{I, 1} \circ \Upsilon_{I, 2} (m_{\eta !}(\mathcal{A}' \widetilde \boxtimes \mathcal{B}')) & \ar[l]_-{\sim} \Upsilon_{I}m_{\eta !}(\mathcal{A}' \widetilde \boxtimes \mathcal{B}') \ar[r]^-{\sim} & \Upsilon_{\ux}\Delta^*m_{\eta !}(\mathcal{A}' \widetilde \boxtimes \mathcal{B}').
}\end{equation}

We now make several observations about this diagram, similar to those in the proof of \cite[Theorem 3.5.1]{AcharRiche:Central}. First, we have a canonical isomorphism $\Delta^*(\mathcal{A}' \widetilde \boxtimes \mathcal{B}') \cong A \widetilde \boxtimes \mathcal{B} \boxtimes \Z_{\C_\eta}$, and the resulting automorphism of $\mathsf{Z}_{\mathbf{f}}(\mathcal{A}) \widetilde \boxtimes \mathsf{Z}_{\mathbf{f}}(\mathcal{B})$ obtained by combining \eqref{eqn--NearbyTop} and \eqref{eqn--Conv} is the identity. 

Second, if we take $S = \Gm$, $Y = (\Gr_{\mathcal{G}_{\mathbf{f}}})_\eta$, $X = \Gr_{\mathcal{G}_{\mathbf{f}}} \times S$, $N = \mathcal{A}'$, and $M = \mathcal{B}' \boxtimes \Z_{S}$, then the map \refeq{Upsilon vs outer box} is an isomorphism. Indeed, this can be checked after applying hyperbolic localization, which reduces us to the constant family $Y \times \A^1$, where the result follows from \thref{nearby trivial family} and \cite[Corollary 3.10]{CassvdHScholbach:MotivicSatake}. Combining \refeq{Upsilon vs outer box} with smooth pullback along $p$ and $q$ and pushforward along $m$ gives a canonical isomorphism

\begin{equation} \label{eqn--centrality-2} \Upsilon_{I, 1} \circ \Upsilon_{I, 2} (m_{\eta !}(\mathcal{A}' \widetilde \boxtimes \mathcal{B}')) \cong \Upsilon_{\ux \cup \zero}(\mathcal{A} \boxtimes \mathsf{Z}_{\mathbf{f}}(\mathcal{B}) \boxtimes \Z_{C_\eta}).
\end{equation} Recall that the right side of \eqref{eqn--centrality-2} is the key motive in the construction of $\eqref{eqn--central}_{\mathbf{f}}$. 

Third, by base change we have a canonical isomorphism
\begin{equation} \label{eqn--centrality-3} \Upsilon_{\ux}\Delta^*m_{\eta !}(\mathcal{A}' \widetilde \boxtimes \mathcal{B}') \cong \mathsf{Z}_{\mathbf{f}}(\mathcal{A} \star \mathcal{B}).\end{equation} 
Via the identification $m_{s!} \Upsilon_{(\ux, \ux)}\Delta^*(\mathcal{A}' \widetilde \boxtimes \mathcal{B}') \cong \mathsf{Z}_{\mathbf{f}}(\mathcal{A}) \star \mathsf{Z}_{\mathbf{f}}(\mathcal{B})$ coming from \eqref{eqn--Conv}, the right vertical map of \eqref{eqn--centrality} is identified with $\eqref{eqn--monoidal}_{\mathbf{f}}$.

Let $j \colon \restr{ \Gr_{\mathcal{G}_{\mathbf{f}}, I} }{C_\eta^2 \setminus \Delta_\eta} \rightarrow (\Gr_{\mathcal{G}_{\mathbf{f}}, I})_\eta$ be the inclusion over the locus with pairwise distinct nonzero coordinates. 
Note that $m$ is an isomorphism over $C_\eta^2 \setminus \Delta_\eta$. 
By \cite[Theorem 5.46]{CassvdHScholbach:MotivicSatake} and its proof (which shows that a global convolution product is mixed Tate as soon as the box product is mixed Tate), the motive $m_{\eta !}(\mathcal{A}' \widetilde \boxtimes \mathcal{B}')[2]$ lies in the Satake category $\Sat^{G, I}$. Then by \cite[Proposition 5.21]{CassvdHScholbach:MotivicSatake}, $$m_{\eta !}(\mathcal{A}' \widetilde \boxtimes \mathcal{B}')[2] \cong j_{!*} (\mathcal{A} \boxtimes \mathcal{B} \boxtimes \Z_{C_\eta^2 \setminus \Delta_\eta}[2]).$$ 
There is an isomorphism of $\Gr_{\mathcal{G}_{\mathbf{f}}, I}$ which swaps the two coordinates of $C^2$. Over $C_\eta^2 \setminus \Delta_\eta$ this isomorphism swaps the two factors of $\Gr$ in $ \restr{ \Gr_{\mathcal{G}_{\mathbf{f}}, I} }{C_\eta^2 \setminus \Delta_\eta} = \Gr \times \Gr \times C_\eta^2 \setminus \Delta_\eta$, and over $C^2_s$ it is the identity map. Thus, we get a canonical isomorphism $$ \Upsilon_{I}m_{\eta !}(\mathcal{A}' \widetilde \boxtimes \mathcal{B}') \cong \Upsilon_{I}m_{\eta !}(\mathcal{B}' \widetilde \boxtimes \mathcal{A}') .$$ Similarly we get an isomorphism
$$ \Upsilon_{\ux}\Delta^*m_{\eta !}(\mathcal{A}' \widetilde \boxtimes \mathcal{B}') \cong  \Upsilon_{\ux}\Delta^*m_{\eta !}(\mathcal{B}' \widetilde \boxtimes \mathcal{A}') $$ which identifies with ${\mathsf{Z}_{\mathbf{f}} \circ \eqref{eqn--central}_{\mathbf{f}_0}}$ under \eqref{eqn--centrality-3}. On the other hand, the swapping isomorphism sends $\Upsilon_{I, 1} \circ \Upsilon_{I, 2}$ to $\Upsilon_{I, 2} \circ \Upsilon_{I, 1}$. The resulting isomorphism $$\Upsilon_{I, 1} \circ \Upsilon_{I, 2} (m_{\eta !}(\mathcal{A}' \widetilde \boxtimes \mathcal{B}')) \cong \Upsilon_{I, 2} \circ \Upsilon_{I, 1} (m_{\eta !}(\mathcal{B}' \widetilde \boxtimes \mathcal{A}'))$$ is identified under \eqref{eqn--centrality-2} (and its analogue for $\Upsilon_{I, 2} \circ \Upsilon_{I, 1}$) with $\eqref{eqn--monoidal}_{\mathbf{f}}$. This finally brings all of the isomorphisms in \eqref{eqn-central-iso} into the picture, the commutativity of which now follows from that of \eqref{eqn--centrality} and all of the identifications we have made.
\xpf

\section{Wakimoto motives}
\label{sect--Wakimotives}
The classical Wakimoto sheaves, due to Mirkovi\'c, are Iwahori-equivariant sheaves on the affine flag variety that categorify the Bernstein elements in the Iwahori--Hecke algebra, cf.~\cite[§5]{AcharRiche:Central}. 
They are used in \cite{ArkhipovBezrukavnikov:Perverse} to construct the so-called Arkhipov--Bezrukavnikov equivalence, which relates Iwahori--Whittaker sheaves on \(\Fl\) to \(\hat{G}\)-equivariant coherent sheaves on the Springer resolution of \(\hat{G}\), where \(\hat{G}\) is the Langlands dual group of \(G\).
The treatment of these Wakimoto sheaves in \cite{ArkhipovBezrukavnikov:Perverse} crucially uses t-exactness of nearby cycles, which is currently not available in the motivic setting.
A similar problem arises when using diamonds, such as in \cite{ALWY:Gaitsgory}.
To solve this problem, the authors of loc.~cit.~refine the study of Wakimoto sheaves, and use them to show t-exactness of the central functor from sheaves on the affine Grassmannian to sheaves on the full affine flag variety.
In this section, we will define and study similar objects in the motivic setting, with a view towards a motivic Arkhipov--Bezrukavnikov equivalence, and use them as in \cite{ALWY:Gaitsgory} to deduce that \(\mathsf{Z}_{\mathbf{a}_0}\) is t-exact.
The general strategy will follow \cite[§4]{AcharRiche:Central}, with the necessary modifications to suit our purposes.

\subsection{Standard and costandard motives}

Before we define the motivic Wakimoto functors, let us recall the standard and costandard functors, along with some basic properties.
Let \(\iota\colon \Fl^\dagger\to \Fl\) be the stratification by the \(\mathcal{I}\)-orbits $\Fl^w$, for \(w\in W\). We also write \(\iota_w\colon \mathcal{I}\backslash\Fl^w\to \mathcal{I}\backslash\Fl\) for the induced map after modding out $\calI$. We denote by \(h_w\colon \mathcal{I}\backslash\Fl^w\to \mathcal{I}\backslash S\) the quotient of the structure map.

\defi\thlabel{Definition CoStandard}
For any \(w\in W\), the standard functor is defined as
\[\Delta_w\colon \DM(\mathcal{I}\backslash S) \to \DM(\mathcal{I}\backslash \Fl)\colon \Ff\mapsto \iota_{w!}h_w^*\Ff [l(w)],\]
Similarly, the costandard functor for \(w\) is defined as 
\[\nabla_w\colon \DM(\mathcal{I}\backslash S) \to \DM(\mathcal{I} \backslash \Fl)\colon \Ff\mapsto \iota_{w*}h_w^*\Ff [l(w)].\]
\xdefi

By \cite[Proposition 3.7]{CassvdHScholbach:MotivicSatake}, these functors send Tate motives on \(\mathcal{I}\backslash S\) to stratified Tate motives on \(\mathcal{I}\backslash\Fl\), and we will usually view \(\Delta_w\) and \(\nabla_w\) as functors \(\DTM_{\mathcal{I}}(S)\to \DTM_{\mathcal{I}}(\Fl)\).
The following proposition is the affine analogue of \cite[Lemma 4.21]{EberhardtKelly:Mixed}, and we follow its proof.

\prop\thlabel{StandardExact}
Both functors \(\iota_*,\iota_!\colon \DTM(\Fl^\dagger) \to \DTM(\Fl)\) are t-exact with respect to the t-structures in \refsect{mixed Tate motives}.
The same is true when considering \(\mathcal{I}\)-equivariant motives on \(\Fl\). Therefore, for any \(w\in W\), the standard and costandard functors \(\Delta_w\) and \(\nabla_w\) are t-exact.
\xprop

\pf
By homotopy invariance and $\Fl^w = \A^{l(w)}_S$, the t-exactness of $\iota_!$ is equivalent to the t-exactness of the $\Delta_w$ etc.

We will show the t-exactness of \(\Delta_w\), the case of \(\nabla_w\) being similar. To this end, we consider the (co)standard functors as functors $\DTM(S) \to \DTM(\Fl)$, i.e., we forget the equivariance.

Since \(\Delta_w\) is defined via !-pushforward, and lying in \(\DTM^{\leq 0}(\Fl)\) is a condition on the *-pullback along the stratification, it is clear that \(\Delta_w\) is right t-exact. 

To show the left t-exactness, let \(M\in \DTM^{\geq 0}(S)\).
We will show \(\Delta_w(M)\in \DTM^{\geq 0}(\Fl)\) by induction on \(l(w)\).
If \(l(w)=0\), then \(\iota_w\) is a closed immersion, pushforward along which is t-exact.
If \(l(w)>0\), we can find \(v\in W\) and a simple reflection \(s\) such that \(w=vs\) and \(l(w)=l(v)+1\).
Consider the codimension 1 facet \(\mathbf{f}_s\) in the closure of \(\mathbf{a}_0\) corresponding to \(s\), and the projection \(\pi\colon \Fl\to \Fl_{\mathbf{f}_s}\), which is an étale-locally trivial \(\P\)-fibration.
As in \cite[(5.1.2)]{RicharzScholbach:Intersection}, there is a distinguished triangle
\[\pi^*\pi_* \Delta_v(M) \to \Delta_v(M) \to \Delta_w(M).\]
Now for any \(x\in W\), we have the cartesian diagram 
\[\begin{tikzcd}
	\Fl^{xs} \cup \Fl^x \arrow[d, "p"] \arrow[r, "k"] & \Fl \arrow[d, "\pi"] \\
	\Fl_{x\cdot W_s}^\circ(\mathbf{a}_0,\mathbf{f}_s)  \arrow[r, "i"] & \Fl_{\mathbf{f}_s},
\end{tikzcd}\]
where \(p\) is a \(\P\)-fibration, and we write \(W_s:=W_{\mathbf{f}_s}\) for simplicity.
Moreover, \(\Fl_{x\cdot W_s}^\circ(\mathbf{a}_0,\mathbf{f}_s)\) is an affine space, so that we have an equivalence \(\DTM(\Fl_{x\cdot W_s}^\circ(\mathbf{a}_0,\mathbf{f}_s)) \cong \DTM(S)\), which is t-exact up to a shift.
Applying \(k^!\) to the distinguished triangle above and using base change, we get an exact triangle
\[p^*p_*k^! \Delta_v(M) \to k^!\Delta_v(M) \to k^! \Delta_w(M).\]
This reduces us to the claim that \(\cofib(p^*p_* \to \id)\) is left t-exact as an endofunctor on \(\DTM(\P)\) stratified as $\P = \A^1 \sqcup S$, which is an easy computation.
\xpf

\prop
\thlabel{Prop-Convolutions of standards and costandards}
Let \(v,w\in W\) and \(M,N\in \DTM_{\mathcal{I}}(S)\).
\begin{enumerate}
	\item\label{item--Convolutions of standards} If \(l(vw) = l(v) + l(w)\), there are natural isomorphisms \(\Delta_v(M)\star \Delta_w(N) \cong \Delta_{vw}(M\otimes N)\)
and \(\nabla_v(M)\star \nabla_w(N) \cong \nabla_{vw}(M\otimes N)\).
These are compatible with associativity as in \cite[Lemma 4.1.4]{AcharRiche:Central}.
	\item\label{item--Invertible for convolution} There exist isomorphisms \(\Delta_v(\Z)\star \nabla_{v^{-1}}(\Z) \cong \Delta_e(\Z)(-l(v)) \cong \nabla_{v^{-1}}(\Z) \star \Delta_v(\Z)\).
\end{enumerate}
\xprop
\pf
We follow the proof of \cite[Lemma 4.1.4]{AcharRiche:Central}. 
The isomorphism for standards in \refit{Convolutions of standards} can be proven verbatim as in loc.~cit. For costandards we additionally use \cite[Corollary 3.10]{CassvdHScholbach:MotivicSatake} to make sure the *-pushforward along an open embedding commutes with the (twisted) exterior product; note that, while the proof in loc.~cit.~only states a K\"unneth map for $*$-pushforward with $\Z$, it immediately follows for general Tate motives.

We review the argument for \refit{Invertible for convolution} in order to explain the appearance of the extra Tate twist.
Recall that forgetting the equivariance gives a fully faithful functor \(\MTM_{\mathcal{I}}(\Fl)\to \MTM(\Fl)\) \cite[Proposition 3.2.20]{RicharzScholbach:Intersection}. 
Since \(\Delta_e(\Z)(-l(v))\in \MTM_{\mathcal{I}}(\Fl)\) and lying in a certain degree can be checked after forgetting the equivariance, it is thus enough to prove \refit{Invertible for convolution} after forgetting the \(\mathcal{I}\)-action. 
Hence, we will again abuse notation and view \(\Delta_v\) and \(\Delta_{v^{-1}}\) as functors \(\DTM(S)\to \DTM(\Fl)\).
By \refit{Convolutions of standards}, we can assume that \(w\in \Omega\) or \(w=s\) is a simple reflection.
The case \(w\in \Omega\) is easy, as then \(l(w)=0\), so that \(\Fl^{\leq w}\cong S\cong \Fl^{\leq w^{-1}}\) are closed in \(\Fl\).

On the other hand, for a simple reflection \(s\), the convolution product \(\Delta_s(\Z)\star \nabla_s(\Z)\) is supported on \(\Fl^{\leq s}=\Fl^s\cup \Fl^e\), as \(s^2=e\).
Letting \(k\colon \Fl^{\leq s} \to \Fl\) denote the closed immersion, we get a localization triangle
\(k_*\Z\to \iota_{e*} \Z \to \Delta_s(\Z)\) in \(\DTM(\Fl)\).
Convolving then gives an exact triangle 
\[k_*\Z\star \nabla_s(\Z)\to \nabla_s(\Z) \to \Delta_s(\Z)\star \nabla_s(\Z).\]
The leftmost object in this triangle can be identified with \(k_*\Z[1]\): 
by \cite[Lemma 3.12 (i)]{RicharzScholbach:Motivic}, it is the *-pushforward of \(\Z[1]\) along the convolution map \(m'\colon \Fl^{\leq s} \widetilde{\times} \Fl^{s} \to \Fl^{\leq s}\).
As in the proof of \cite[Proposition 3.19]{RicharzScholbach:Motivic}, this map can be identified with \(\P \times \P \setminus \Delta(\P) \to \P\times \P \xrightarrow{\operatorname{pr}_2} \P\), which is a Zariski-locally trivial \(\A^1\)-bundle.
Hence the natural map \(\Z[1]\to m'_*m'^{*}\Z[1]\) is an isomorphism, which implies that \(k_*\Z\star \nabla_s(\Z) \cong k_*\Z[1]\).
We claim moreover that \(\Delta_s(\Z)\star \nabla_s(\Z)\) is supported at the origin.
Since \[\Hom_{\DM(\Fl^{\leq s})}(\Z[1],\iota_{s,*}\Z[1])\cong \Hom_{\DM(\A^1)}(\Z,\Z)\cong \Z,\]
this will imply
that the left arrow is isomorphic (up to an invertible scalar) to the natural inclusion \(k_*\Z[1]\hookrightarrow \nabla_s(\Z)\), the cone of which is \(\Delta_e(\Z)(-1)\).
Thus we get \(\Delta_s(\Z)\star \nabla_s(\Z)\cong \Delta_e(\Z)(-1)\) (non-canonically), and similarly one can show \(\nabla_s(\Z)\star \Delta_s(\Z)\cong \Delta_e(\Z)(-1)\).

We are left to prove the claim that \(\Delta_s(\Z)\star \nabla_s(\Z)\) is supported at the origin.
Let \(q\colon LG\to \Fl\) be the quotient map as usual, and consider the diagram
\[\begin{tikzcd}
	U \arrow[r, "a"] \arrow[d, "j_U"] & Y \arrow[rr, "b"] \arrow[d, "j_Y"] && \Fl^s \arrow[d, "j"]\\
	q^{-1}(\Fl^s) \overset{\mathcal{I}}{\times} \Fl^{s} \arrow[r, "f"] & q^{-1}(\Fl^s) \overset{\mathcal{I}}{\times} \Fl^{\leq s} \arrow[r, "g"] & q^{-1}(\Fl^{\leq s}) \overset{\mathcal{I}}{\times} \Fl^{\leq s} \arrow[r, "m"] & \Fl^{\leq s}
\end{tikzcd}\]
with cartesian squares,
where \(f\) and \(g\) are the natural open immersions, and \(m\) is the (restriction of the) multiplication map. 
Denoting by \(\overline{q}\colon q^{-1}(\Fl^s) \overset{\mathcal{I}}{\times} \Fl^{\leq s}\to \Fl^s\) the quotient map on the first factor, we see as in \cite[Lemma 4.1.4]{AcharRiche:Central} that
\[q^{-1}(\Fl^s) \overset{\mathcal{I}}{\times} \Fl^{\leq s} \xrightarrow{(\overline{q}, m\circ g)} \Fl^s\times \Fl^{\leq s}\cong \A^1\times \P\]
is an isomorphism, cf.~also \cite[(3.20)]{RicharzScholbach:Motivic}.
Now, in the diagram above, \(U\subset Y\) is an open subset with complement isomorphic to \(\A^1\). We denote the inclusion by \(c\colon Y\setminus U\subset Y\).
Note that \(\Delta_s(\Z)\star \nabla_s(\Z)=m_!g_!(\Z\widetilde{\boxtimes} \nabla_s(\Z)) = m_!g_!f_*(\Z\widetilde{\boxtimes} \Z)\), where we again use \cite[Corollary 3.10]{CassvdHScholbach:MotivicSatake} to ensure the commutation of \(f_*\) with the twisted exterior product.
The restriction of this motive to \(\Fl^s\subset \Fl^{\leq s}\) is given by \(j^*m_!g_!f_*(\Z\widetilde{\boxtimes} \Z) \cong b_!a_*j_U^*(\Z\widetilde{\boxtimes} \Z)\), which we want to show vanishes.
Let \(\Ff:=j_Y^*(\Z\widetilde{\boxtimes} \Z)\). 
By localization, we have an exact triangle 
\[b_!c_!c^!\Ff\cong c^!\Ff \to b_!\Ff\to b_!a_*a^*\Ff.\]
Note that $Y \cong \A^2$, the map $b \colon \A^2 \rightarrow \Fl^s \cong \A^1$ identifies with projection onto the second factor, and $c \colon \A^1 \rightarrow \A^2$ identifies with the diagonal inclusion. Hence $b \circ c \cong \id$, and also \(c^!\) and \(b_!\) induce the same equivalence $M \mapsto M(-1)[-2]$ on Tate motives. From this it follows that $c^!\Ff \to b_!\Ff$ is an isomorphism.
\xpf

In contrast to the affine Grassmannian the convolution product on \(\Fl\) is not t-exact, even when working with coefficients in a field. 
However, we can still show t-exactness properties when convolving a standard and a costandard motive, which will be used to show the Wakimoto motives lie in \(\MTM\).
Classically, this is shown by using t-exactness properties of affine pushforwards.
Since this is not available in the motivic setting, we give a different argument.

Recall that \(\mathcal{I}\backslash S\) is a placid prestack (in the sense of \cite[§A.2]{RicharzScholbach:Motivic}), so that \(\DM(\mathcal{I}\backslash S)\) is monoidal by composing the exterior product from \cite[Corollary A.15]{RicharzScholbach:Motivic} with pullback along \(\mathcal{I}\backslash S \to (\mathcal{I} \times \mathcal{I}) \backslash S\).
After forgetting the equivariance, this monoidal structure is exactly the usual tensor product on \(\DM(S)\), and in particular it preserves the subcategory of Tate motives.
Recall also from \cite[Definition and Lemma 4.11]{CassvdHScholbach:MotivicSatake} that convolution preserves Tate motives, so that we may view convolution with some \(\Ff\in \DTM_{\mathcal{I}}(\Fl)\) as an endofunctor of \(\DTM_{\mathcal{I}}(\Fl)\).

\lemm \thlabel{lemm--leftright}
Let \(v\in W\).
Convolving with \(\Delta_v(\Z)\) on the left or right is left t-exact, while convolving with \(\nabla_v(\Z)\) on the left or right is right t-exact.
\xlemm
\pf
We prove the lemma for right convolution with standard motives, the other cases being similar (see also \cite[Proposition 4.6]{Achar:ModularPerverseSheavesII}).
Choosing a reduced expression for \(v\), we can moreover assume by \thref{Prop-Convolutions of standards and costandards} that \(v=s\) is a simple reflection.

Let \(\Ff\in \DTM^{\geq 0}_{\mathcal{I}}(\Fl)\), which we may assume to be bounded.
Then \(\Ff\) is a finite extension of costandard objects \(\nabla_w(\Ff_w)\), for \(\Ff_w \in \DTM^{\geq 0}_{\mathcal{I}}(S)\), so we may assume \(\Ff = \nabla_w(\Ff_w)\) for some \(w\in W\).
We now distinguish two cases.
If \(l(ws)<l(w)\), then by \thref{Prop-Convolutions of standards and costandards} we have
\[\nabla_w(\Ff_w)\star \Delta_s(\Z) \cong \nabla_{ws}(\Ff_w) \star \nabla_s(\Z) \star \Delta_s(\Z) \cong \nabla_{ws}(\Ff_w(-1)),\]
which lies in \(\DTM^{\geq 0}_{\mathcal{I}}(\Fl)\) by \thref{StandardExact}.

On the other hand, \(\nabla_w(\Ff_w)\star \Delta_s(\Z)\) can always be viewed as the *-pushforward of \(h_w^*(\Ff_w[l(w)]) \widetilde{\boxtimes} \Delta_s(\Z)\) along the restricted convolution map
\[\Fl^w \widetilde{\times} \Fl^{\leq s} \to \Fl.\]
But under the assumption that \(l(w)<l(ws)\), the convolution map induces an isomorphism \(\Fl^w\widetilde{\times} \Fl^{\leq s} \to \Fl^w \cup \Fl^{ws}\).
As lying in degree \(\geq 0\) is a condition on the !-pullback, the desired left t-exactness can be checked by using base change, after which it follows from relative purity applied to the immersions of strata on the smooth scheme \(\Fl^w \cup \Fl^{ws}\).
\xpf

Combining this with \thref{StandardExact}, we get the following corollary.

\coro\thlabel{Wakimoto are perverse}
Let \(v,w\in W\). For any \(M,N\in \MTM_{\mathcal{I}}(S)\) such that \(M\otimes N\in \MTM_{\mathcal{I}}(S)\), the convolutions \(\Delta_v(M)\star \nabla_w(N)\) and \(\nabla_w(N)\star \Delta_v(M)\) lie in \(\MTM_{\mathcal{I}}(\Fl)\).
\xcoro

\pf
The convolution \(\Delta_v(M)\star \nabla_w(N)\) agrees with both \(\Delta_v(M \otimes N) \star \nabla_w(\Z)\) and \(\Delta_v(\Z) \star \nabla_w(M\otimes N)\). The other convolution is similar.
\xpf

\subsection{Wakimoto functors}

We now define the Wakimoto functors, the essential image of which will be called the Wakimoto motives.
Recall that we fixed some split maximal torus and Borel \(T\subseteq B\subseteq G\).
In particular, we have the dominant cocharacters \(X_*(T)^+\subseteq X_*(T)\), and the dominance order \(\trianglelefteq_B\) on \(X_*(T)\) given by
\[\lambda\trianglelefteq_B \mu \iff \mu-\lambda\in X_*(T)^+.\]

We will define a Wakimoto functor for each cocharacter in \(X_*(T)\), which agrees (up to a twist) with the costandard functors for dominant cocharacters.
Moreover, this collection of functors should be compatible with the addition on \(X_*(T)\) and the convolution product.
In order to implement this idea on the level of \(\infty\)-categories we use the following lemma.

\lemm\thlabel{monoidal functors of inf categories}
Let \(M\) be a discrete monoid, $\Cc$ a commutative algebra object in $\PrLomegaSt$ (cf.~\refsect{monoidal aspects} for notation), and $\Dd$ an algebra object in $\Mod_\Cc(\PrLSt)$.
Then there are equivalences
\[\Fun^{\t,\L}_{\Cc} (\Fun(M,\Cc), \Dd) \cong \Fun^{\t,\L} (\Fun(M,\Ani),\Dd)\cong \Fun^\t(M,\Dd).\]
Here, $\Ani$ denotes the \ii-category of anima (also called \(\infty\)-groupoids or spaces), \(\Fun^{\t,L}_{\Cc}\) denotes the \(\infty\)-category of monoidal, colimit-preserving, \(\Cc\)-linear functors, and \(\Fun(M,\Cc)\) is equipped with the Day convolution product.
\xlemm

\pf
For the first equivalence, we note that \(\Fun(M,\Cc)\) is a product of copies of \(\Cc\).
Since \(\Cc\) is compactly generated, it is dualizable in \(\PrLSt\).
Hence, tensoring with \(\Cc\) (in \(\PrLSt\)) preserves limits, and in particular products.

The second equivalence is given by precomposing with the (symmetric monoidal) Yoneda embedding, as in \cite[Example 3.2.7]{Robalo:Theorie}.
\xpf

Below, we will use that for any discrete commutative monoid $M$ (such as $X_*(T)^+$) the (derived) group completion (i.e., its image under the left adjoint of the forgetful functor $\Grp_{E_\infty}(\Ani) \r \Mon_{E_\infty}(\Ani)$) is actually discrete and therefore agrees with the usual (underived) Grothendieck completion. This is, say, a (very) special case of \cite[Proposition~6]{Nikolaus:Group}: for $M$ being discrete, the filtered colimit $M_m := \colim (M \stackrel m \r M \stackrel m \r \dots)$ is discrete for any $m \in M$, as is $M_S := (M_{\{m_1, \dots, m_n\}})_{m_{n+1}}$ for a finite subset $S = \{m_1 < \dots < m_{n+1}\}$ inside a fixed well-ordered set of generators of $M$. Hence $M_\infty := \colim_S M_S$ is also discrete, so the fundamental groups of $M_\infty$ (at all base points) are trivial, so $M_\infty$ is a model for the group completion by loc.~cit.

In order to apply \thref{monoidal functors of inf categories}, we want to extend
\begin{equation}\label{functor to extend}
X_*(T)^+\to \DTM_{\mathcal{I}}(\Fl) \colon \mu \mapsto \nabla_{t(\mu)}(\Z)(\langle 2\rho,\mu\rangle)
\end{equation}
to a monoidal functor \(X_*(T)\to \DTM_{\mathcal{I}}(\Fl)\), where \(2\rho\) is the sum of the positive roots of \(G\). (Here we use the \ii-categorical monoidal structure given by convolution, cf.~\refsect{convolution structures}).
Since \(\{\nabla_{t(\mu)}(\Z)(\langle 2\rho,\mu\rangle)\mid \mu\in X_*(T)^+\}\subseteq\MTM_{\mathcal{I}}(\Fl)\) by \thref{StandardExact}, which is the nerve of an ordinary category (as opposed to an \ii-category), \thref{Prop-Convolutions of standards and costandards} yields a monoidal structure on \eqref{functor to extend} (see also \cite[§4.2.3]{AcharRiche:Central}).
Now, $X_*(T)$ is the Grothendieck completion of $X_*(T)^+$ or, equivalently, the (derived) group completion. Moreover, \(\nabla_{t(\mu)}(\Z)(\langle 2\rho,\mu\rangle)\) is invertible with inverse \(\Delta_{t(-\mu)}(\Z)\), so we get a unique extension to a monoidal functor $X_*(T) \to \DTM_\calI(\Fl)$; the essential image of this functor still lies in \(\MTM_{\mathcal{I}}(\Fl)\) by \thref{Wakimoto are perverse}.

Note that \(\DTM_{\mathcal{I}}(S)\) is compactly generated by \cite[Theorem 7.2.1]{Krause:Localization}, as the full subcategory of \(\DM(\mathcal{I}\backslash S)\) closed under colimits and containing the Tate twists, which is compactly generated by \cite[Lemma 2.3.6]{RicharzScholbach:Intersection}.
Applying \thref{monoidal functors of inf categories} to the situation where \(M=X_*(T)\), \(\Cc = \DTM_{\mathcal{I}}(S)\) and \(\Dd = \DTM_{\mathcal{I}}(\Fl)\) then permits the following definition.

\defi
\thlabel{Definition Wakimoto}
The \emph{full Wakimoto functor} \[\JWak\colon \Rep_{\hat T}(\DTM_{\mathcal{I}}(S)) \cong \Fun(X_*(T),\DTM_{\mathcal{I}}(S)) \to \DTM_{\mathcal{I}}(\Fl)\] 
is the monoidal functor which, under the equivalence from \thref{monoidal functors of inf categories}, corresponds to the functor \(X_*(T)\to \DTM_{\mathcal{I}}(\Fl)\) constructed above.

Moreover, for any \(\mu\in X_*(T)\), we define \(\JWak_\mu\) by precomposing \(\JWak\) with the natural embedding \(\DTM_{\mathcal{I}}(S) \to \Rep_{\hat T}(\DTM_{\mathcal{I}}(S))\) corresponding to \(\mu\in X_*(T) = X^*(\hat{T})\).
We call \(\JWak_{\mu}\) the \emph{Wakimoto functor associated to \(\mu\)}.
In particular, we have a decomposition \(\JWak \cong \coprod_{\mu\in X_*(T)} \JWak_\mu\).
\xdefi

\rema \thlabel{Observations about Wakimoto functor}
We record  some basic properties of the \(\JWak_{\mu}\), whose analogues in the classical setting are well-known.
\begin{enumerate}
	\item Since the definition above involves the dominant cocharacters specifically, the Wakimoto functors really depend on the choice of Borel \(B\).
	\item For any \(\lambda\in X_*(T)\) and \(\mu\in X_*(T)^+\) such that \(\mu \trianglerighteq_B \lambda\), there is an isomorphism of functors
	\[\JWak_{\lambda} \cong \nabla_{t(\mu)}(\Z)(\langle 2\rho,\mu\rangle) \star \Delta_{t(\lambda-\mu)}\cong \Delta_{t(\lambda-\mu)}\star \nabla_{t(\mu)}(\Z)(\langle 2\rho,\mu\rangle).\] 
	\item By \thref{Wakimoto are perverse}, each \(\JWak_\mu\) is t-exact.
	Hence \(\JWak\) is t-exact as well.
	\item \label{Monoidality of Wakimoto} Monoidality of \(\JWak\) implies that for \(\lambda,\mu\in X_*(T)\) and \(M,N\in \DTM_{\mathcal{I}}(S)\), we have an isomorphism
	\[\JWak_{\lambda}(M) \star \JWak_{\mu}(N) \cong \JWak_{\lambda+\mu}(M\otimes N).\]
\end{enumerate} 
\xrema

\rema
Our definition of the Wakimoto functors differs slightly from \cite{AcharRiche:Central}, even in the mixed setting.
By the twist appearing in \thref{Prop-Convolutions of standards and costandards}\eqref{item--Invertible for convolution}, the functors \(\JWak_\mu\) also have to involve some additional twist.
But unlike \cite[§5.3.2]{AcharRiche:Central}, we cannot form half-twists, as we are working with motives.
Our particular choice to twist the Wakimoto functors for dominant, rather than anti-dominant, cocharacters, is motivated by the Bernstein elements in the generic Iwahori--Hecke algebra, and we refer to Section \ref{sect--generic IH algebra} for more details.
\xrema

\lemm\thlabel{semi-orthogonality of Wakimoto}
Let \(\lambda,\mu\in X_*(T)\), and \(M,N\in \DTM_{\mathcal{I}}(S)\).
If \[\Maps_{\DTM_{\mathcal{I}}(\Fl)}(\JWak_{\lambda}(M),\JWak_{\mu}(N)) \neq 0,\] then \(\lambda-\mu\) is a sum of positive roots.
\xlemm
\pf
By convolving with some sufficiently dominant cocharacter, we may assume \(\lambda,\mu\in X_*(T)^+\) are dominant.
Then the Wakimoto functors agree with costandard functors up to a twist.
Hence, by adjunction, the mapping spaces above vanish if \(\iota_{t(\mu)}^*\iota_{t(\lambda)*}\) is trivial.
By \cite[Lemma 4.1.2]{AcharRiche:Central}, this is the case unless \(\lambda-\mu\) is a sum of positive roots (for \(T\subseteq B\)).
\xpf

Recall that we introduced the Wakimoto functors in order to show t-exactness of \(\mathsf{Z}_{\mathbf{a}_0}\). 
As the full subcategory of \(\DTM_{\mathcal{I}}(\Fl)\) spanned by the essential image of \(\JWak\) is not stable and does not contain the essential image of \(\mathsf{Z}_{\mathbf{a}_0}\), we will also consider the following categories.

\defi
We let 
$$\Wak^+ \subseteq \Wak\subseteq \DTM_{\mathcal{I}}(\Fl)$$ be the full stable cocomplete subcategory generated by the images of \(\JWak_\mu\) for \(\mu\in X_*(T)^+\), resp.~\(\mu\in X_*(T)\).

We also let \(\Wak_{\operatorname{bd}}\subseteq \DTM_{\mathcal{I}}(\Fl)\) be the full subcategory generated under extensions by the images of \(\JWak_{\mu}\) for \(\mu\in X_*(T)\).
In particular, each object in \(\Wak\) is a colimit of objects in \(\Wak_{\operatorname{bd}}\).
\xdefi

By \thref{Observations about Wakimoto functor} \eqref{Monoidality of Wakimoto}, these are monoidal subcategories (with respect to the convolution product).
Note that any object in \(\Wak_{\operatorname{bd}}\) is bounded, but it is not clear whether a bounded motive in \(\Wak\) also lies in \(\Wak_{\operatorname{bd}}\).

For a motive \(\Ff\in \DTM_{\mathcal{I}}(\Fl)\), we let \(!\text{-}\Supp(\Ff):=\{w\in W\mid \iota_w^!\Ff\neq 0\}\), and we similarly define \(*\text{-}\Supp(\Ff)\).

\prop\thlabel{Wak+criterion}
An object \(\Ff\in \DTM_{\mathcal{I}}(\Fl)\) lies in \(\Wak^+\) if and only if \(!\text{-}\Supp(\Ff)\subset \{t(\mu)\mid \mu\in X_*(T)^+\}\).
\xprop
\pf
Since the Schubert cells \(\Fl^w\cong \A^{l(w)}\) are affine spaces, this follows from localization and homotopy invariance.
\xpf

In order to show the central motives lie in \(\Wak\), we need another criterion, which in turn uses the following lemma.

\lemm\thlabel{SubsetSupport}
Let \(\Ff\in \DTM_{\mathcal{I}}(\Fl)\) be bounded. Then there exists a finite subset \(A_{\Ff}\subset W\) such that
\[*\text{-}\Supp(\Delta_w(M)\star \Ff)\subset w\cdot A_{\Ff}, \qquad !\text{-}\Supp(\nabla_w(M)\star \Ff)\subset w\cdot A_{\Ff},\]
\[*\text{-}\Supp(\Ff\star \Delta_w(M))\subset A_{\Ff} \cdot w, \qquad !\text{-}\Supp(\Ff\star \nabla_w(M))\subset A_{\Ff}\cdot w,\]
for any \(M\in \DTM_{\mathcal{I}}(S)\) and \(w\in W\).
\xlemm
\pf
The proofs of \cite[Lemma 4.4.2, Proposition 4.4.4]{AcharRiche:Central} go through without change, using only geometric arguments and basic properties of the 6-functor-formalism.
\xpf

\prop
If a bounded object \(\Ff\in \DTM_{\mathcal{I}}(\Fl)\) satisfies \(\Delta_{t(\nu)}(\Z)\star \Ff\cong \Ff\star \Delta_{t(\nu)}(\Z)\) for all \(\nu\in -X_*(T)^+\), then \(\Ff\in \Wak_{\operatorname{bd}}\).
\xprop
\pf
Let \(A_{\Ff}\) be a finite subset of \(W\) obtained via \thref{SubsetSupport}.
Then we can find some \(\nu\in -X_*(T)^+\) such that
\[t(\nu)\cdot A_{\Ff}\subseteq \{t(\mu)\cdot w\mid \mu\in -X_*(T)^+_{\mathrm{reg}}, w\in W_0\}, \quad A_{\Ff}\cdot t(\nu)\subseteq \{w\cdot t(\mu)\mid \mu\in -X_*(T)^+_{\mathrm{reg}}, w\in W_0\},\]
where \(X_*(T)^+_{\mathrm{reg}}\) denotes the regular dominant cocharacters.
As every \(W_0\)-orbit in \(X_*(T)\) has a unique (anti-)dominant representative, we get
\[(t(\nu)\cdot A_{\Ff})\cap (A_{\Ff}\cdot t(\nu))\subseteq \{t(\mu)\mid \mu\in -X_*(T)^+_{\mathrm{reg}}\}.\]
Along with the assumption on \(\Ff\) and \thref{SubsetSupport}, we deduce that
\[*\text{-}\Supp(\JWak_\nu(\Z)\star \Ff)\subseteq \{t(\mu)\mid \mu\in -X_*(T)^+_{\mathrm{reg}}\}.\]
By localization and stratified Tateness of \(\Ff\), we see that \(\JWak_\nu(\Z)\star \Ff\) is an extension of standard objects associated to anti-dominant cocharacters.
Since \(\Ff\) is moreover bounded, only finitely many such standard objects are nontrivial.
As these are Wakimoto motives, and \(\Wak\) is closed under convolution, we see that \(\Ff \cong \JWak_{-\nu}\star \JWak_\nu\star \Ff\in \Wak_{\operatorname{bd}}\).
\xpf

Using the centrality isomorphism from \thref{CentralityIso}, the following corollary is immediate.

\coro\thlabel{ExistenceWakimotoFiltration}
For any \(\Ff\in \DTM_{L^+G}(\Gr)\), the object \(\mathsf{Z}_{\mathbf{a}_0}(\Ff)\in \DTM_{\mathcal{I}}(\Fl)\) lies in \(\Wak\).
\xcoro

Let us be more precise about what it means for a motive to lie in \(\Wak\).
The next lemma says that any object in \(\Wak\) admits a \emph{Wakimoto filtration}. 
We consider the partial order \(\preceq\) on \(X_*(T)\), where \(\lambda\preceq \mu\) if \(\mu-\lambda\) is a sum of positive coroots.

\defilemm\thlabel{WakimotoFiltration}
We view \(X_*(T)\) as an \(\infty\)-category induced by the partially ordered set \((X_*(T),\preceq)\).
Then there is a functor
\[\Fil\colon \Wak \to \Fun(X_*(T),\Wak)\colon \Ff\mapsto \left(\lambda\mapsto \Fil_{\lambda}(\Ff)\right),\]
along with a natural transformation \(\Fil_\lambda\to \id\) for each \(\lambda\), such that \(\Fil_\lambda(\Ff)\to \Ff\) is the final morphism out of an object in the full subcategory of \(\DTM_{\mathcal{I}}(\Fl)\) generated under extensions and colimits by the images of \(\JWak_{\mu}\) for \(\mu\preceq \lambda\).
Moreover, there is an isomorphism \(\Ff\cong \colim_{\lambda\in X_*(T)} \Fil_\lambda(\Ff)\) for each \(\Ff\in \Wak\).
\xdefilemm
\pf
Let \(\Ff\in \Wak_{\operatorname{bd}}\). 
Then \(\Ff\) lies in the subcategory of \(\DTM_{\mathcal{I}}(\Fl)\) generated by the essential images of finitely many Wakimoto functors, so that there exists \(\mu\in X_*(T)^+\) sufficiently dominant for which \(\Ff\star \JWak_{\mu}(\Z)\in \Wak^+\).
The desired filtration for \(\Ff\star \JWak_{\mu}(\Z)\) is obtained by sending non-dominant cocharacters to \(0\), and \(\lambda\in X_*(T)^+\) to \(\iota_{\preceq \lambda,*}\iota_{\preceq \lambda}^! \Ff\), where \(\iota_{\preceq \lambda}\colon \Fl^{\leq t(\lambda)}\to \Fl\) denotes the inclusion , and the transition maps for \(\lambda'\preceq \lambda\) arise via localization.
Convolving with the monoidal inverse of \(\JWak_{\mu}(\Z)\), as well as shifting the filtration by \(-\mu\), then gives the desired filtration for \(\Ff\).
This is clearly functorial for uniformly bounded \(\Ff\), and independent of the choice of \(\mu\).
Passing to the colimit then defines \(\Fil\colon \Wak \to \Fun(X_*(T),\Wak)\) in general, and it clearly satisfies the required properties.
\xpf

We note that the construction above is different, but equivalent to the construction of the Wakimoto filtration in \cite{AcharRiche:Central,ALWY:Gaitsgory}. 

\coro\thlabel{WakimotoGraded}
View \(X_*(T)\) as a discrete \(\infty\)-category.
Then the functor \(\Fil\) from \thref{WakimotoFiltration} induces a functor 
\[\Grad\colon \Wak \to \Fun(X_*(T),\DTM_{\mathcal{I}}(S)),\]
such that each \(\Ff\) is a (possibly infinite) extension of all \(\JWak_{\lambda}(\Grad_\lambda(\Ff))\)'s.
\xcoro
We call \(\Grad(\Ff)\) the \emph{associated graded of the Wakimoto filtration of \(\Ff\)}.
\pf
Fix \(\lambda\in X_*(T)\), and consider the cofiber \(\Ff'\) of the natural map \(\colim_{\mu\prec\lambda} \Fil_{\mu}(\Ff) \to \Fil_{\lambda}(\Ff)\), which is functorial.
After convolving with \(\JWak_{\nu}(\Z)\) for sufficiently dominant \(\nu\in X_*(T)^+\), this \(\Ff'\) lies in the image of the fully faithful embedding \(\iota_{t(\lambda) *}\colon \DTM_{\mathcal{I}}(\Fl^{t(\lambda)}) \to \DTM_{\mathcal{I}}(\Fl)\).
By homotopy invariance, we moreover have \(\DTM_{\mathcal{I}}(\Fl^{t(\lambda)})\cong \DTM_{\mathcal{I}}(S)\).
It follows that $\Ff'$ identifies canonically with $\JWak_{\lambda}(\Grad_\lambda(\Ff))$ for some $\Grad_\lambda(\Ff) \in \DTM_{\mathcal{I}}(S)$, which concludes the proof.
\xpf

The following property of this associated graded is immediate.

\coro\thlabel{WakimotoSplitting}
There is a natural isomorphism
\[\Grad\circ \JWak\cong \id\]
of endofunctors on \(\Fun(X_*(T),\DTM_{\mathcal{I}}(S))\).
In other words, \(\Grad\) defines a splitting of \(\JWak\).
\xcoro

Finally, we mention some monoidality properties of \(\Grad\) (although we do not construct a full \(\infty\)-categorical monoidal structure). 

\prop\thlabel{Monoidality of associated graded}
For \(\lambda,\mu,\nu\in X_*(T)\), there exists a canonical morphism
\[\alpha_{\lambda,\mu}\colon \Grad_\lambda(-) \otimes \Grad_\mu(-) \to \Grad_{\lambda+\mu}(-\star-)\]
of functors \(\Wak\times \Wak \to \DTM(S)\),
such that \(\bigoplus_{\lambda+\mu = \nu} \alpha_{\lambda,\mu}\) is an isomorphism.
\xprop
\pf
The proof is similar to \cite[Lemma 4.7.4 and Proposition 4.7.5]{AcharRiche:Central}.
However, it becomes simpler, as the failure of t-exactness of the convolution product is not an issue for us.
For \(\Ff,\Ff'\in \Wak\), there are natural morphisms
\[\Fil_\lambda(\Ff)\star \Fil_\mu(\Ff') \to \Ff\star \Ff' \to \cofib(\Fil_{\lambda+\mu}(\Ff\star \Ff') \to \Ff\star \Ff').\]
By \thref{semi-orthogonality of Wakimoto}, this composition must vanish, so that we get a canonical morphism
\[\Fil_\lambda(\Ff)\star \Fil_\mu(\Ff') \to \Fil_{\lambda+\mu}(\Ff\star \Ff').\]
By using \thref{semi-orthogonality of Wakimoto} again, the morphisms
\[\colim_{\lambda'\prec\lambda} \Fil_{\lambda'} (\Ff)\star \Fil_\mu(\Ff') \to \JWak_{\lambda+\mu}(\Grad_{\lambda+\mu}(\Ff\star \Ff'))\]
and \[\Fil_\lambda(\Ff) \star \colim_{\mu'\prec\mu} \Fil_\mu(\Ff') \to \JWak_{\lambda+\mu}(\Grad_{\lambda+\mu}(\Ff\star \Ff'))\]
vanish, which induces a canonical morphism
\[\alpha_{\lambda,\mu}\colon \Grad_\lambda(\Ff) \otimes \Grad_\mu(\Ff') \to \Grad_{\lambda+\mu}(\Ff\star \Ff').\]

To check that \(\bigoplus_{\lambda+\mu = \nu}\alpha_{\lambda,\mu}\) is an isomorphism, we may restrict to \(\Wak_{\operatorname{bd}}\) by continuity.
If \(\Ff= \JWak_{\lambda}(M)\) for some \(M\in \DTM(S)\), the claim follows from \thref{Observations about Wakimoto functor} \eqref{Monoidality of Wakimoto}.
For general bounded objects, this then follows by induction on the Wakimoto-filtration of \(\Ff\).
\xpf

\subsection{t-exactness of the central functor}

We are now ready to show \(\mathsf{Z}_{\mathbf{a}_0}\) is t-exact.
The idea is to use \thref{ExistenceWakimotoFiltration} to write a (bounded) central motive as an extension of Wakimoto motives. 
It then remains to see that for mixed Tate motives on \(\Gr\), the resulting Wakimoto motives appearing are also mixed Tate. 
This will be deduced using the geometric constant term functors on the affine flag variety as in \cite[§6.1]{HainesRicharz:TestFunctions} and \cite[§6.3]{AGLR:Local}.
In order to compare constant term functors for different affine flag varieties, we introduce a new notation, different from \cite[§5.1]{CassvdHScholbach:MotivicSatake}.
Namely, fix a cocharacter \(\lambda\in X_*(T)\), 
and recall the notation from \refsect{GmActions} and \thref{HypLocCentral}. In particular, the fixed points $\Fl^0$ depend on the chosen $\lambda$, and we have the associated Levi subgroup \(M\subseteq G\) and parabolic subgroup \(P\subseteq G\), defined as the fixed points, resp.~attractor, of the \(\Gm\)-action on \(G\) induced by \(\lambda\).
In order to define the constant term functors with values in equivariant motives, we need to know the \(L^+M_{\mathbf{f}}\)-action on \(\Fl\) restricts to actions on \(\Fl^{\pm}\) and \(\Fl^0\).
This follows from the following general statement.

\lemm
Let \(H/S\) be a group scheme acting on an ind-scheme \(X\), and \(\Gm\to H\) a homomorphism.
This yields a \(\Gm\)-action on \(X\) by restriction, and a \(\Gm\)-action on \(H\) by conjugation.
Then the action of \(H^0\subseteq H\) on \(X\) restricts to an action on \(X^0\), and similarly for attractors and repellers.
\xlemm
\pf
Let \(S'\) be any \(S\)-scheme, and \(h\in \Hom_S(S',H^0)\) and \(x\in \Hom_S(S',X^0)\) two \(S'\)-valued points.
Then for any \(g\in \Hom_S(S',\Gm)\), we have \(ghg^{-1}=h\) and \(gx=x\).
It follows that \(ghx = ghg^{-1}gx = hx\), which proves the statement about the fixed points. 
The case of attractors and repellers can be handled similarly.
\xpf

Now consider the maps
\[L^+M_{\mathbf{a}_0}\backslash \Fl^0 \xleftarrow{\overline{q}^\pm_{\mathbf{a}_0}} L^+M_{\mathbf{a}_0}\backslash \Fl^\pm \xrightarrow{\overline{p}^\pm_{\mathbf{a}_0}} L^+M_{\mathbf{a}_0}\backslash \Fl,\]
obtained by quotienting out the \(L^+M_{\mathbf{a}_0}\)-action from the usual hyperbolic localization diagram.

\defi
The geometric constant term functor associated to \(P\) is
\[\CT_P^{\mathbf{a}_0} := (\overline{q}_{\mathbf{a}_0}^+)_!(\overline{p}_{\mathbf{a}_0}^+)^* \colon \DM(\mathcal{I}\backslash \Fl) \to \DM(L^+M_{\mathbf{a}_0}\backslash \Fl^0).\]
\xdefi

Note that this really only depends on the parabolic \(P\), rather than the cocharacter \(\lambda\).
By \cite[Proposition 2.5]{CassvdHScholbach:MotivicSatake}, this functor is canonically equivalent to \((\overline{q}_{\mathbf{a}_0}^-)_*(\overline{p}_{\mathbf{a}_0}^-)^!\).
Moreover, these constant term functors preserve (anti-effective) stratified Tate motives by \thref{HypLocTate}. 
Recall also the usual constant term functor on the affine Grassmannian from \cite[Definition 5.3]{CassvdHScholbach:MotivicSatake}
\[\CT_P^{\mathbf{f}_0}:= (\overline{q}_{\mathbf{f}_0}^+)_!(\overline{p}_{\mathbf{f}_0}^+)^*
\colon \DM(L^+G\backslash \Gr_G) \to \DM(L^+M\backslash \Gr_{M}),\]
which similarly preserves stratified Tate motives.
(In contrast to loc.~cit., we have not included a shift in the definition of the constant term functor.) 

Since in this section, we are mostly interested in the Iwahori-level, rather than general parahoric, we will denote the central functor \(\mathsf{Z}_{\mathbf{a}_0}\) from \thref{defi--nearby-hecke} by \(\mathsf{Z}_G\).
This notation will make it easier to compare central functors for different groups.
For example, by \thref{prop-Gm-Upsilon}, there is a natural equivalence of functors
\begin{equation}\label{central functor commutes with constant terms}
	\CT_P^{\mathbf{a}_0} \circ \mathsf{Z}_{G} \cong \mathsf{Z}_{M} \circ \CT_P^{\mathbf{f}_0},
\end{equation}
where we implicitly use the inclusion \(\Gr_M\subseteq \Fl^0\) on reduced loci from \thref{HypLocCentral}.

Now, assume that \(\lambda\) is regular anti-dominant (for the fixed Borel \(B\)), in which case \(M=T\) is the maximal torus and \(P=B^-\) is the opposite Borel.
As before, cf.~\thref{lemm-Gm-on-Flf}, we can write the repeller for the induced \(\Gm\)-action on \(\Fl\) as \(\coprod_{w\in W} \mathcal{S}_w\), a coproduct of semi-infinite orbits, keeping the notation from Section \ref{subsubsec-WT}.
These semi-infinite orbits are exactly the \(LU\)-orbits in \(\Fl\), where \(U\subseteq B\) is the unipotent radical.
Since the fixed point locus of this \(\Gm\)-action consists of discrete points (up to taking reduced subschemes) indexed by \(W\), we will also write \(\CT_{B^-}^w\colon \DTM_{\mathcal{I}}(\Fl)\to \DTM_{L^+T}(S)\) for the value of \(\CT_{B^-}^{\mathbf{a}_0}\) at such a point \(w\in W\).

To check which objects appear in the Wakimoto filtrations of central motives, we will apply the constant term functors, along with the following results, compare \cite[§4.5.3]{AcharRiche:Central}.

\prop\thlabel{Constant term of convolution}
For \(\mu\in X_*(T)^+\) and \(w\in W\), there is an isomorphism
\begin{equation}\label{eq-CT of convolution}
	\CT_{B^-}^{t(\mu)w}(\nabla_{t(\mu)}(\Z)\star -) \cong \CT_{B^-}^w(-)[\langle 2\rho,\mu\rangle]
\end{equation}
of functors \(\DTM_{\mathcal{I}}(\Fl) \to \DTM_{L^+T}(S)\).
\xprop
\pf
The proof is analogous to \cite[Proposition 3.23]{ALWY:Gaitsgory}.
Let \(i_{w}\colon \mathcal{S}_w\to \Fl\) be the inclusion, \(f_w\colon \mathcal{S}_w\to S\) the structure map, and \(\Ff\in \DTM_{\mathcal{I}}(\Fl)\) any object.
Fix also a lift \(z^{\mu}\in N_G(T)(\Z\rpot{t})\subseteq LG(\Z)\) of \(t(\mu)\in W\), and denote its image in \(\Fl(\Z)\) the same way.
By hyperbolic localization, we have 
\[\CT_{B^-}^{t(\mu)w}(\nabla_{t(\mu)}(\Z)\star \Ff) \cong f_{t(\mu)w,*}i_{t(\mu)w}^!(m'_*(\Z[\langle 2\rho,\mu\rangle]\widetilde{\boxtimes} \Ff)),\]
where \(m'\colon \Fl^{t(\mu)}\widetilde{\times} \Fl \to \Fl\) is the restriction of the convolution map.
Now, we have \(\mathcal{S}_{t(\mu)}\supseteq \Fl^{t(\mu)}\): for this we can assume \(S\) is the spectrum of a field, in which case it can be shown as in \cite[Lemma 2.12]{ALWY:Gaitsgory}.
As \(\Fl^{t(\mu)} = L^+Uz^\mu\) (again by \cite[Lemma 2.12]{ALWY:Gaitsgory}), this implies the fiber product \(\mathcal{S}_{t(\mu)w} \times_{\Fl}  (\Fl^{t(\mu)}\widetilde{\times} \Fl)\) agrees with \(L^+U z^\mu L^+U \overset{L^+U}{\times} S_w\).

Now, consider the following diagram with the right square cartesian:
\[\begin{tikzcd}
	&L^+U z^\mu L^+U \times \mathcal{S}_w \arrow[dl, "a"'] \arrow[d, "b"] &\\
	\Fl^{t(\mu)} \times \mathcal{S}_w \arrow[d, "h_{t(\mu)}\times f_w"'] &L^+U z^\mu L^+U \overset{L^+U}{\times} \mathcal{S}_w \arrow[d, "m''"] \arrow[r, "\widetilde{i}"]& \Fl^{t(\mu)} \widetilde{\times} \Fl \arrow[d, "m'"]\\
	S&\mathcal{S}_{t(\mu)w} \arrow[r, "i_{t(\mu)w}"] \arrow[l, "f_{t(\mu)w}"']& \Fl.
\end{tikzcd}\]
By base change, it suffices to compute
\[f_{t(\mu)w,*}m''_*\widetilde{i}^!(\Z[\langle2\rho,\mu\rangle]\widetilde{\boxtimes} \Ff).\]
By definition, we have an isomorphism 
\[b^! \widetilde{i}^!(\Z[\langle2\rho,\mu\rangle]\widetilde{\boxtimes} \Ff) \cong a^!(\Z[\langle2\rho,\mu\rangle]\boxtimes i_w^!(\Ff)).\]
Replacing \(L^+U\) by a suitable split unipotent truncation, purity yields a similar isomorphism when replacing \(a^!,b^!\) by \(a^*,b^*\).
In that case, \(a\) and \(b\) are torsors under split unipotent groups, hence Zariski-locally trivial by \cite[Proposition A.6]{RicharzScholbach:Intersection}.
We can then compute
\[\CT_{B^-}^{t(\mu)w}(\nabla_{t(\mu)}(\Z)\star \Ff) \cong f_{t(\mu)w,*}m''_*b_*b^*\widetilde{i}^!(\Z[\langle2\rho,\mu\rangle]\widetilde{\boxtimes} \Ff)\]
\[\cong (h_{t(\mu)}\times f_w)_* a_*a^*(\Z[\langle2\rho,\mu\rangle]\boxtimes i_w^!(\Ff)) \cong \CT_{B^-}^w(\Ff)[\langle 2\rho,\mu\rangle],\]
where the last isomorphism uses homotopy invariance and \(\Fl^{t(\mu)}\cong \A^{l(t(\mu))}\).
\xpf

Since the maximal reductive quotient of \(\mathcal{I}\) is naturally isomorphic to the torus \(T\), \cite[Proposition 2.2.11]{RicharzScholbach:Intersection} gives us canonical equivalences
\[\DM(\mathcal{I}\backslash S) \cong \DM(T\backslash S) \cong \DM(L^+T\backslash S),\]
identifying the subcategories of Tate motives.

\coro\thlabel{Constant term of Wakimoto}
Let \(\mu\in X_*(T)\). 
Then \(\CT_{B^-}^{\mathbf{a}_0}\circ \JWak_\mu\) takes values in motives supported on the connected component of \((\Fl)^0\) corresponding to \(t(\mu)\), and there is an isomorphism
\[\CT_{B^-}^{t(\mu)}\circ \JWak_\mu\cong \id(\langle 2\rho,\mu\rangle)[\langle 2\rho,\mu\rangle]\] 
of functors \(\DTM_{\mathcal{I}}(S)\to \DTM_{\mathcal{I}}(S)\).
\xcoro
\pf

Let \(\nu\in X_*(T)^+\) be such that \(\nu \trianglerighteq_B \mu\), and \(\Ff\in \DTM_{\mathcal{I}}(S)\).
By \thref{Observations about Wakimoto functor} and \thref{Constant term of convolution}, there exist isomorphisms
\[\CT_{B^-}^w(\mathbf{J}_\mu^B(\Ff)) \cong \CT_{B^-}^{t(\nu-\mu)w}(\JWak_{\nu}(\Z)\star \Delta_e(\Ff))(-\langle 2\rho,\nu-\mu\rangle)[-\langle 2\rho,\nu-\mu\rangle]\]
\[\cong \CT_{B^-}^{t(-\mu)w}(\Delta_e(\Ff))(\langle 2\rho,\mu\rangle)[\langle 2\rho,\mu\rangle].\]
Since \(\Fl^e\subset \mathcal{S}_e\), the rightmost term vanishes unless \(w=t(\mu)\), and we clearly have \(\CT_{B^-}^e\circ \Delta_e\cong \id\).
\xpf

With all the preparations out of the way, we can now prove the main theorem of this section.
As usual, let \(\deg:=\langle 2\rho,-\rangle\colon\Gr_T\to \Z\) be the (locally constant) degree function.

\theo \thlabel{Grad-CT}
There is a natural isomorphism of functors
\[\Grad\circ \mathsf{Z}_{G} \cong \CT_{B^-}^{\mathbf{f}_0}(-\deg)[-\deg].\]
In particular, \(\mathsf{Z}_G=\mathsf{Z}_{\mathbf{a}_0}\) is t-exact.
\xtheo
\pf
Consider the commutative diagram
\begin{equation}\label{central functors and constant terms}\begin{tikzcd}[column sep=huge]
		\NilpQ \DTM_{L^+T}(\Fl^0) & \NilpQ \Wak \arrow[l, "\CT_{B^-}^{\mathbf{a}_0}{[}-\deg{]}"']\\
		\DTM_{L^+T}(\Gr_T) \arrow[u, "\mathsf{Z}_{T}"] & \DTM_{L^+G}(\Gr_G) \arrow[l, "\CT_{B^-}^{\mathbf{f}_0}{[}-\deg{]}"] \arrow[u, "\mathsf{Z}_{G}"']
\end{tikzcd}\end{equation}
arising from \thref{prop-Gm-Upsilon} and \thref{ExistenceWakimotoFiltration}.
By \thref{nearby trivial family}, \(\mathsf{Z}_{T}\) agrees with the pushforward along the open and closed immersion (on reduced loci) \(\Gr_T\subseteq \Fl^0\) (with trivial monodromy).
Moreover, the lower arrow is t-exact by \cite[Proposition 5.10]{CassvdHScholbach:MotivicSatake} (since \(\lambda\) is anti-dominant).

Recall from \thref{WakimotoGraded} that every \(\Ff\in \Wak\) is an extension of \(\JWak_{\mu}(\Grad_\mu(\Ff))\) for varying \(\mu\).
Since we know \(\CT_{B^-}^{\mathbf{a}_0}(\JWak_{\mu}(\Grad_\mu(\Ff)))\) is supported on \(t(\mu)\in \Fl^0\) with value \(\Grad_\mu(\Ff)(\langle 2\rho,\mu\rangle)[\langle 2\rho,\mu\rangle]\) by \thref{Constant term of Wakimoto}, the natural isomorphism \(\Grad\circ \mathsf{Z}_{G} \cong \CT_{B^-}^{\mathbf{f}_0}(-\deg)[-\deg]\) follows from \eqref{central functors and constant terms} (additionally using \thref{WakimotoSplitting}).

Again using the fact that \(\mathsf{Z}_G\) is an extension of \(\JWak_{\mu}(\Grad_\mu(\mathsf{Z}_G(-)))\)'s, each of which is t-exact by t-exactness of \(\CT_{B^-}^{\mathbf{f}_0}(-\deg)[-\deg]\), we conclude that \(\mathsf{Z}_G\) itself is t-exact.
\xpf

\subsection{The case of general parahorics}

Until now, this section was only concerned with the case where \(\mathbf{f}=\mathbf{a}_0\) is an alcove.
Following an idea of Achar, we now use these results to show that the central functor \(\mathsf{Z}_{\mathbf{f}}\) is in fact t-exact for general facets.
For the rest of this section, fix a facet \(\mathbf{f}\) in the closure of \(\mathbf{a}_0\).
We start with the analogues of \thref{Definition CoStandard} and \thref{StandardExact}, and use similar notations \(\iota_{\overline{w}}\) and \(h_{\overline{w}}\) for \(\overline{w}\in W/W_{\mathbf{f}}\). 

\defi
For \(\overline{w}\in W/W_{\mathbf{f}}\), the standard functor is
\[\Delta_{\overline{w}}^{\mathbf{f}}\colon \DM(\mathcal{I}\backslash S) \to \DM(\mathcal{I\backslash \Fl_{\mathbf{f}}})\colon \Ff\mapsto \iota_{\overline{w}!}h^*_{\overline{w}}\Ff[l(\overline{w})],\]
while the costandard functor is
\[\nabla_{\overline{w}}^{\mathbf{f}}\colon \DM(\mathcal{I}\backslash S) \to \DM(\mathcal{I\backslash \Fl_{\mathbf{f}}})\colon \Ff\mapsto \iota_{\overline{w}*}h^*_{\overline{w}}\Ff[l(\overline{w})].\]
\xdefi

These functors preserve Tate motives by \cite[Proposition 3.7]{CassvdHScholbach:MotivicSatake}. In the following proposition, which is an adaptation of  \cite[Proposition 4.7]{Achar:ModularPerverseSheavesII} to the motivic setting, we consider Tate motives on \(\Fl_{\mathbf{f}}\) with respect to the stratification by \(\mathcal{I}\)-orbits.

\prop\thlabel{Prop standards for general facets}
For any \(\overline{w}\in W/W_{\mathbf{f}}\), the functors \(\Delta_{\overline{w}}^{\mathbf{f}},\nabla_{\overline{w}}^{\mathbf{f}}\colon \DTM_{\mathcal{I}}(S) \to \DTM_{\mathcal{I}}(\Fl_{\mathbf{f}})\) are t-exact.
\xprop
\pf
We will consider the case of standard functors, as costandard functors can be handled similarly.
Note also that standard functors are automatically right t-exact, and that we may forget about the \(\mathcal{I}\)-equivariance.

Since \(\mathbf{f}\) lies in the closure of \(\mathbf{a}_0\), we can consider the projection \(\pi\colon \Fl\to \Fl_{\mathbf{f}}\), which is a Zariski-locally trivial fibration with typical fiber \(L^+G_{\mathbf{f}}/\mathcal{I}\) \cite[Proposition 4.3.13 (i)]{RicharzScholbach:Intersection} (the Zariski local triviality follows from the fact that \(LG\to \Fl_{\mathbf{f}}\) has Zariski local sections).
Let \(d\) be the relative dimension of the smooth \(S\)-scheme \(L^+G_{\mathbf{f}}/\mathcal{I}\).
Since \(\pi\) is \(\mathcal{I}\)-equivariant and smooth surjective, \(\pi^*[d]\) preserves stratified Tate motives and is conservative and t-exact by \cite[Proposition 4.3.9 (i)]{RicharzScholbach:Intersection}.

Now, let \(w_{\mathbf{f}}\in W_{\mathbf{f}}\) be the longest element, and consider the two maps \(\Fl\overset{\mathcal{I}}{\times}\Fl^{\leq w_{\mathbf{f}}}\to \Fl\), given by the projection on the first factor, and the convolution respectively.
These maps agree after composition with \(\pi\), so they induce a map \[\Fl\overset{\mathcal{I}}{\times}\Fl^{\leq w_{\mathbf{f}}} \to \Fl \times_{\Fl_{\mathbf{f}}} \Fl,\]
which is an isomorphism (compare \cite[Proof of Lemma 4.3]{Achar:ModularPerverseSheavesII}).
Let \(w\in W\) be the smallest length representative of \(\overline{w}\), so that \(\pi_!\circ\Delta_w \cong \Delta_{\overline{w}}^{\mathbf{f}}\).
Then base change and \thref{Prop-Convolutions of standards and costandards}\refit{Convolutions of standards} yield an isomorphism 
\[\pi^*[d](\Delta_{\overline{w}}^{\mathbf{f}})(-) \cong \Delta_w(\Z) \star \pi^*[d](\Delta_{\overline{e}}^{\mathbf{f}})(-) .\]
Now $\Delta_{\overline{e}}^{\mathbf{f}}(-)$ is t-exact, thus so is $\pi^*[d](\Delta_{\overline{e}}^{\mathbf{f}})(-)$, 
and hence \(\pi^*[d](\Delta_{\overline{w}}^{\mathbf{f}})(-)\) is left t-exact by \thref{lemm--leftright}.
This implies that \(\Delta_{\overline{w}}^{\mathbf{f}}\) is left t-exact by t-exactness and conservativity of \(\pi^*[d]\).
\xpf

Now, consider the convolution product
\[-\star-\colon \DM(\mathcal{I}\backslash \Fl) \times \DM(\mathcal{I}\backslash \Fl_{\mathbf{f}}) \to \DM(\mathcal{I}\backslash \Fl_{\mathbf{f}}),\]
which preserves stratified Tate motives by the same argument as in \cite[Theorem 3.17]{RicharzScholbach:Motivic}.
At least on homotopy categories, it clearly determines an action of the monoidal category \(\DTM_{\mathcal{I}}(\Fl)\) on \(\DTM_{\mathcal{I}}(\Fl_{\mathbf{f}})\).

\lemm\thlabel{leftright general facets}
Let \(v\in W\) be any element. Then left-convolving with \(\Delta_v(\Z)\) (resp.~\(\nabla_v(\Z)\)) is a left (resp.~right) t-exact endofunctor of \(\DTM_{\mathcal{I}}(\Fl_{\mathbf{f}})\). 
\xlemm
\pf
As usual, we only prove the case of standard functors.
Let \(\Ff\in \DTM_{\mathcal{I}}(\Fl_{\mathbf{f}})\) be any object.
The diagram
\[\begin{tikzcd}
	\Fl \widetilde{\times} \Fl \arrow[r, "m"] \arrow[d, "\id \widetilde{\times} \pi"] & \Fl\arrow[d, "\pi"]\\
	\Fl \widetilde{\times} \Fl_{\mathbf{f}} \arrow[r, "m"] & \Fl_{\mathbf{f}}
\end{tikzcd}\]
is cartesian. 
As in the proof of \thref{Prop standards for general facets}, $\pi$ is smooth and \(\pi^*[d]\) is t-exact and conservative, where $d$ is the relative dimension.  Then \(\pi^*(\Delta_v(\Z)\star \Ff)[d] \cong \Delta_v(\Z)\star \pi^*(\Ff)[d]\), so the proposition follows from \thref{lemm--leftright}.
\xpf

The idea of the following proof was explained to us by Louren\c{c}o.

\lemm \thlabel{lemm--minimalorbit}
There exists a Borel $B$ containing \(T\) such for any dominant \(\mu\) (with respect to this Borel), \(t(\mu)\) is minimal in its right \(W_{\mathbf{f}}\)-orbit.
\xlemm

\pf It is equivalent to require that the length functions on $W$ and $W/W_{\mathbf{f}}$ assign the same value to $t(\mu)$ and $\overline{t(\mu)}$, respectively.
If \(\mathbf{f}_0\) lies in the closure of \(\mathbf{f}\), then by \cite[Remark 4.2.16]{RicharzScholbach:Intersection} we may take $B$ to be the unique Borel whose associated Weyl chamber contains \(\mathbf{a}_0\). In general, it suffices to consider the case where $\mathbf{f}$ is a vertex. Let $\Sigma^+(\mathbf{f})$ be the set of affine roots \(\alpha\) which vanish at $\mathbf{f}$ and such that $\alpha(\mathbf{a}_0) > 0$. By \cite[(1.8), (1.10)]{Richarz:Schubert}, we need to find $B$ such that if $\mu$ is $B$-dominant, then $\alpha(\mathbf{a}_0 + t(\mu)) > 0$ for all $\alpha \in \Sigma^+(\mathbf{f})$ (note that the normalization of the Kottwitz map in op.~cit. differs from ours by a sign). 

Each $\alpha \in \Sigma^+(\mathbf{f})$ vanishes on some affine hyperplane in the standard apartment $\mathscr{A}$. The complement of this hyperplane is a disjoint union of connected components $\mathbf{H}_{\alpha}^+ \sqcup \mathbf{H}_{\alpha}^-$, where $\mathbf{H}_{\alpha}^+$ is determined by the requirement that it contains $\mathbf{a}_0$.
Then $\cap_{\alpha \in \Sigma^+(\mathbf{f})}\mathbf{H}_{\alpha}^+$ is an open affine cone with vertex at $\mathbf{f}$, and we need to choose $B$ so that its Weyl chamber $C$ satisfies $\mathbf{a}_0 +C \subset \cap_{\alpha \in \Sigma^+(\mathbf{f})}\mathbf{H}_{\alpha}^+$. 
 Let $\tau$ be the translation of $\mathscr{A}$ which sends $\mathbf{f}$ to $\mathbf{f}_0$. Then the boundary of the closure of $\tau(\mathbf{H}_{\alpha}^+)$ is a hyperplane equal to the vanishing locus of some non-affine root in the root system associated with $(G,T)$. Thus, for $B$ we may pick any Borel whose Weyl chamber is contained in $\cap_{\alpha \in \Sigma^+(\mathbf{f})}\tau(\mathbf{H}_{\alpha}^+)$.
\xpf

\rema
If $\mathbf{f}$ is a non-special vertex then the gradients of the affine roots vanishing at $\mathbf{f}$ do not give all the (non-affine) roots associated to $(G,T)$, even up to scaling. Hence there will be at least two possible Borels in \thref{lemm--minimalorbit} in this case. For example, if $\mathbf{f}$ is a right-angle vertex in the $B_2$ affine root system, then there are two possible Borels since each Weyl chamber has angle $\pi/4$.
\xrema

These results allow us to generalize (part of) \thref{Grad-CT} to general facets.

\theo\thlabel{t-exactness general facets}
The central functor \(\mathsf{Z}_{\mathbf{f}}\colon \DTM_{L^+G}(\Gr) \to \DTM_{L^+G_{\mathbf{f}}}(\Fl_{\mathbf{f}})\) is t-exact.
\xtheo
\pf
Since \(\mathbf{f}\) lies in the closure of \(\mathbf{a}_0\), we can consider the projection \(\pi\colon \Fl\to \Fl_{\mathbf{f}}\). Let us choose $B$ as in \thref{lemm--minimalorbit}; this is harmless since \(\mathsf{Z}_{\mathbf{f}}\) does not depend on a choice of Borel.
Then for $\mu$ dominant, \(\Fl^{t(\mu)}\) has the same dimension as its image in \(\Fl_{\mathbf{f}}\) by \thref{prop--summary}\refit{summary3}. 
This implies that \(\pi_*\nabla_{t(\mu)}\cong \nabla^{\mathbf{f}}_{\overline{t(\mu)}}\), which is t-exact by \thref{Prop standards for general facets}.
In particular, for any \(\lambda\in X_*(T)\) and \(\mu\in X_*(T)^+\) with \(\mu\trianglerighteq_B \lambda\), pushing forward the Wakimoto functors gives \[\pi_*\JWak_\lambda \cong \pi_*(\Delta_{t(\lambda-\mu)}(\Z)\star \nabla_{t(\mu)}(\langle 2\rho,\mu\rangle)) \cong \Delta_{t(\lambda-\mu)}(\Z)\star \nabla^{\mathbf{f}}_{\overline{t(\mu)}}(\langle 2\rho,\mu\rangle),\]
which is left t-exact by \thref{leftright general facets}.

Since we have \(\mathsf{Z}_{\mathbf{f}} \cong \pi_*\circ \mathsf{Z}_{\mathbf{a}_0}\) (after forgetting the equivariance), the t-exact Wakimoto filtration on \(\mathsf{Z}_{\mathbf{a}_0}\) from \thref{Grad-CT} shows that \(\mathsf{Z}_{\mathbf{f}}\) is left t-exact, at least when considering the stratification of \(\Fl_{\mathbf{f}}\) by \(\mathcal{I}\)-orbits.
By replacing \(B\) by its opposite Borel \(B^-\) (cf.~\thref{rmk-choice of borel}), an argument dual to the above shows that \(\mathsf{Z}_{\mathbf{f}}\) is right t-exact.
It remains to observe that by \cite[Proposition 3.2.22]{RicharzScholbach:Intersection}, if a motive \(\Ff\in \DTM_{L^+G_{\mathbf{f}}}(\Fl_{\mathbf{f}})\) is mixed Tate when viewed as a motive in \(\DTM_{\mathcal{I}}(\Fl_{\mathbf{f}})\), then we already have \(\Ff\in \MTM_{L^+G_{\mathbf{f}}}(\Fl_{\mathbf{f}})\).
\xpf

\section{Geometrization of generic Hecke algebras}\label{sect--generic IH algebra}
In this section, we introduce generic parahoric Hecke algebras as decategorifications, i.e., Grothendieck groups, of appropriate categories of Tate motives. We show how this specializes to classical Hecke algebras, and also interpolates between generic Iwahori--Hecke algebras \cite{Vigneras:Algebres,Vigneras:Pro-p} and generic spherical Hecke algebras as considered in \cite[§6.3]{CassvdHScholbach:MotivicSatake}.
We finally identify the centers of these generic parahoric Hecke algebras in many cases.
\iffalse
In \cite[§6.3]{CassvdHScholbach:MotivicSatake}, we defined the generic spherical Hecke algebra (generalizing \cite[Definition 6.2.2]{PepinSchmidt:Generic} in the case \(\GL_2\)) using generators and relations, and showed it agrees with the Grothendieck ring of a suitable category of anti-effective (mixed) Tate motives on the affine Grassmannian.
Similarly, one can define a generic Iwahori--Hecke algebras in terms of generators and relations \cite[Theorem 2.1]{Vigneras:Pro-p}, and again we can show that it agrees with the Grothendieck ring of a suitable category of motives on the full affine flag variety, cf.~\thref{comparison with Vigneras}.
The goal of this section is to generalize this to general parahoric levels, by defining generic Hecke algebras as Grothendieck rings of categories of motives, and show they indeed recover the classical parahoric algebras after specialization.
We then use the functoriality of these Grothendieck rings to construct central elements in these generic parahoric Hecke algebras.
\fi

Recall that for a monoidal triangulated category $C$, its Grothendieck group $K_0(C)$ is a ring.

\subsection{Generic parahoric Hecke algebras}

Let \(\mathbf{f}\) be a facet in the Bruhat-Tits building of \(G\), which we assume to lie in the closure of \(\mathbf{a}_0\).
Recall from \thref{monoidal structure locc} that the convolution product turns the \ii-category \(\DTM_{L^+G_{\mathbf{f}}}(\Fl_{\mathbf{f}})^{\antilc}\) of anti-effective locally compact equivariant stratified Tate motives on the partial affine flag variety of \(\mathbf{f}\) into a monoidal triangulated category.

\defi\thlabel{Defi generic parahoric}
The \emph{generic parahoric Hecke algebra} at \(G_{\mathbf{f}}\)-level is the Grothendieck ring
\[\calH^{\mathbf{f}}(\qq):= K_0(\DTM_{L^+G_{\mathbf{f}}}(\Fl_{\mathbf{f}})^{\antilc}).\]
\xdefi 

This is naturally a \(\Z[\qq]\)-algebra, where the variable \(\qq\) acts via the negative Tate twist \((-1)\).
As a \(\Z[\qq]\)-module, $\calH^{\bbf}$ is free with basis indexed by \(W_{\mathbf{f}}\backslash W/W_{\mathbf{f}}\) (Corollaries \ref{coro_K0 DTM} and \ref{Grothendieck-stratified-equivariant}).
(In particular, it is independent of the base scheme $S$ fixed in \thref{notation S}.)

\exam
If \(\mathbf{f}=\mathbf{f}_0\) is hyperspecial, we write $\calH^{\bbf_0}(\qq) =:\calH^{\sph}(\qq)$, and refer to it as the \emph{generic spherical Hecke algebra}.
According to \cite[Theorem 6.32 and Corollary 6.36]{CassvdHScholbach:MotivicSatake}, this recovers the definition in terms of generators and relations in \cite[Definition 6.34]{CassvdHScholbach:MotivicSatake}.
\xexam

On the other extreme, we have the case where \(\mathbf{f}=\mathbf{a}_0\) is an alcove, for which we also write \(\calH^{\mathbf{a}_0}(\qq) = \calH^{\mathcal{I}}(\qq)\).
We now relate this case of \thref{Defi generic parahoric} to the generic Hecke algebra defined by Vignéras \cite[Theorem 2.1]{Vigneras:Pro-p}, which was defined as the free $\Z[\qq]$-module with basis $\{T_w\}_{w \in W}$ endowed with the unique $\Z[\qq]$-algebra structure satisfying:
\begin{itemize}
  \item (the braid relations) $T_w T_{w'} = T_{ww'}$ if $w, w' \in W$ and $l(ww') = l(w) + l(w')$.
  \item (the quadratic relations) $T_s^2 = \qq + (\qq-1)T_s$ if $s \in S_{\aff}$.
\end{itemize}

Recall the standard sheaves \(\Delta_w(\Z)\) from \thref{Definition CoStandard}.

\prop\thlabel{comparison with Vigneras}
\thlabel{Geometrization of generic IHalgebra}
The generic Iwahori--Hecke algebra \(\calH^{\mathcal{I}}(\qq)\), with basis \((-1)^{l(w)}[\Delta_w(\Z)]\) satisfies the braid relations and the quadratic relations. 
In particular, it agrees with Vignéras' definition.
\xprop
\pf
The braid relations are satisfied by \thref{Prop-Convolutions of standards and costandards}.
For the quadratic relations, let $s \in S_{\aff}$. Since $s^2 = \id < s$, the reduced fibers of the convolution map  $m \colon \Fl^{\leq s} \tilde \times \Fl^{\leq s} \rightarrow \Fl^{\leq s}$ are all isomorphic to $\P$. This follows from the arguments in the proof of \cite[Lemma 9]{Faltings:Loops} since $m$ can be viewed as a Demazure map for the non-reduced word $ss$.  We must determine the fibers of the restricted convolution map $m' \colon \Fl^{s} \tilde \times \Fl^{s} \rightarrow \Fl^{\leq s}$. We have $(m')^{-1}(\{e\}) = m^{-1}(\{e\}) \setminus (\{e\}, \{e\}) \cong \A^1$ and $(m')^{-1}(\{s\}) = m^{-1}(\{s\}) \setminus (\{e\}, \{s\}) \cup (\{s\}, \{e\}) \cong \Gm$. Thus $\iota_e^*(\Delta_s(\Z) * \Delta_s(\Z)) \cong \Z(-1)[0]$ and $\iota_s^*(\Delta_s(\Z) * \Delta_s(\Z)) \cong \Z[1] \oplus \Z(-1)[0]$.
By a straightforward computation this implies that our map respects the quadratic relations.
\xpf

Recall that the \emph{classical Hecke algebra} $\calH^\bbf := \calH^\bbf(q)$ is defined as the set of bi-\(G_{\mathbf{f}}(\Fq\pot{t})\)-equivariant compactly supported functions on \(G(\Fq\rpot{t})\) with values in \(\Z\), equipped with the convolution product.

\prop\thlabel{generic parahoric}
There is a canonical ring isomorphism
\[\calH^{\mathbf{f}}(\qq) \otimes_{\Z[\qq], \qq\mapsto q} \Z \stackrel \cong \r \calH^{\mathbf{f}}.\]
\xprop
\pf
To construct a map $\calH^\bbf (\qq) \r \calH^\bbf$, we may choose $S = \Spec \Fq$.
We pass to motives with rational coefficients (this does not even affect the Grothendieck groups by \thref{Rationalization and K0}), which have étale descent.
To any motive \(M\in \DTM_{L^+G_{\mathbf{f}}}(\Fl_{\mathbf{f}})^{\antilc}\), we associate a compactly supported function on \(G_{\mathbf{f}}(\mathcal{O}_F)\backslash G(F) / G_{\mathbf{f}}(\mathcal{O}_F)\), by *-pulling back \(M\) to the basepoints of the strata in \(\Fl_{\mathbf{f}}\), and then applying the trace of Frobenius (e.g., cf.~\cite{Cisinski:SurveyCoho}).
Since the trace of Frobenius on \(\Z(-1)\) is \(q\), these functions indeed take values in \(\Z\). The resulting map is a ring morphism by a similar argument as in \cite[Lemma 5.6.1]{Zhu:Introduction}.
Clearly, \(\qq-q\) is contained in the kernel of this morphism.

The resulting map $\calH^{\mathbf{f}}(\qq) \otimes_{\Z[\qq], \qq\mapsto q} \Z \r \calH^{\mathbf{f}}$ sends the standard motive $\iota_{w!}(\Z[l(w)])$ (where \(\iota_w\colon \Fl_w^{\circ}(\mathbf{f},\mathbf{f}) \to \Fl_{\mathbf{f}}\) is the inclusion) to the characteristic function of the double coset of $w$, up to a factor of \(\pm 1\). That is, the map preserves bases (of free $\Z$-modules) and is therefore an isomorphism.
\xpf

\rema\thlabel{polynomials via cellularity}
One could also use the recent result of Haines on the cellularity of the fibers of convolution morphisms \cite[§9]{Haines:Pavings} to define generic Hecke algebras at parahoric level in terms of generators and relations.
However, maps such as \eqref{Z blah} below become much less clear in such an approach.
\xrema

\subsection{Construction of central elements}

Recall that the central functor \(\mathsf{Z}_{\mathbf{f}}\) from \thref{defi--nearby-hecke} preserves locally compact anti-effective Tate motives (\thref{theo-WT-BD}). By Theorems~\ref{theo-moniodal-functor}, \ref{theo-E2}, passing to Grothendieck groups, it induces a ring homomorphism
\begin{equation}\label{Z blah}
	\varphi^\bbf\colon\calH^{\sph}(\qq)\to Z(\calH^{\bbf}(\qq)),
\end{equation}
where $Z$ denotes the center of a ring.
By \thref{thm--change-of-facet}, \(\varphi^{\bbf}\) is injective for any \(\bbf\).
Using results on centers of classical Hecke algebras, we now show this is an isomorphism in case \(\bbf = \mathbf{a}_0\).

Recall that classically, \(\calH^{\mathcal{I}}\) contains a (non-unique) maximally commutative subalgebra generated by the Bernstein elements, which are in turn categorified by the Wakimoto sheaves.
Using the motivic Wakimoto functor from \thref{Definition Wakimoto}, we adapt this as follows.

\defi\thlabel{defi-Bernstein}
For \(\lambda\in X_*(T)\), we define the \emph{Bernstein element} 
\[\theta_\lambda := [\JWak_{\lambda}(\Z)(-\frac{\langle 2\rho,\lambda\rangle + l(t(\lambda))}{2})] \in \calH^{\mathcal{I}}(\qq).\]
We denote by \(\mathcal{A}(\qq)\subseteq \mathcal{H}^{\mathcal{I}}(\qq)\) the \(\Z[\qq]\)-subalgebra generated by \(\theta_\lambda\) for \(\lambda\in X_*(T)\).
\xdefi

Note that since \(\mathcal{A}(\qq)\) is defined via Wakimoto motives, it depends on the choice of Borel \(B\).
Recall from \thref{Observations about Wakimoto functor} that the Wakimoto functors for dominant cocharacters involved a twist, to ensure that the full Wakimoto functor is monoidal.
The drawback of this is that \(\JWak_{\mu}(\Z)\) is not anti-effective in general, which is why we reversed the twist in the definition above.
This implies that the \(\theta_\lambda\) are not invertible in \(\mathcal{H}_G^{\mathcal{I}}(\qq)\), only after inverting \(\qq\).

\rema
In order to compare our results to \cite{Vigneras:Algebres}, we set all parameters \(q_s\) to \(\qq\) in loc.~cit.
Then \(\theta_\lambda\) agrees with the element \((-1)^{l(t(\lambda))}E_{\lambda}\) from \cite[Théorème 1.5]{Vigneras:Algebres}, as long as one replaces the fixed Borel \(B\) in loc.~cit.~by its opposite \(B^-\).
Indeed, for dominant and anti-dominant cocharacters this follows from the definitions, and for general cocharacters it follows from the equality
\[\theta_\lambda \cdot \theta_\mu = \qq^{\frac{l(t(\lambda)) + l(t(\mu)) - l(t(\lambda+\mu))}{2}}\cdot\theta_{\lambda + \mu}\]
induced by \thref{Prop-Convolutions of standards and costandards}, and similar relations for the \(E_\lambda\) in \cite{Vigneras:Algebres}.
In particular, \(\mathcal{A}(\qq)\) agrees with the algebra \(A\cap H\) from \cite[(1.6.1)]{Vigneras:Algebres}.
This also shows that we indeed have \(\theta_\lambda\in \calH^{\mathcal{I}}(\qq)\), without having to invert \(\qq\).
\xrema

\prop\thlabel{generic Wakimoto}
The set \(\{\theta_\lambda\mid \lambda\in X_*(T)\}\) forms a \(\Z[\qq]\)-basis of \(\mathcal{A}(\qq)\), which in turn is a maximal commutative subalgebra of \(\mathcal{H}^{\mathcal{I}}(\qq)\).
\xprop

\pf
Recall from \thref{Prop-Convolutions of standards and costandards} that \(\Delta_s(\Z)\star \nabla_s(\Z) = \Delta_e(\Z)(-1)\), for any \(s\in S_{\aff}\).
Hence, writing \([\nabla_s(\Z)]\) as a linear combination of \([\Delta_e(\Z)]\) and \([\Delta_s(\Z)]\), the quadratic relations from \thref{comparison with Vigneras} yield \([\nabla_s(\Z)] = [\Delta_s(\Z)] + \qq -1\).
Applying \thref{Prop-Convolutions of standards and costandards} again, this allows us to determine the classes of all costandard objects, and we see that \([\nabla_w(\Z)]\) is a linear combination of \(\{T_v\mid v\leq w\}\), and that the coefficient of \(T_{w}\) appearing is \(\pm 1\).
Hence the linear independence of \(\{\theta_\lambda\mid \lambda\in X_*(T)\}\) follows from \thref{Geometrization of generic IHalgebra}.
The fact that \(\{\theta_\lambda\mid \lambda\in X_*(T)\}\) also span \(\mathcal{A}(\qq)\) follows essentially by definition, since  \(\{\theta_\lambda\mid \lambda\in X_*(T)\}\) is closed under multiplication, up to a power of \(\qq\).

Since \(\theta_\lambda\cdot \theta_\mu = \theta_\mu \cdot \theta_\lambda\), we see that \(\mathcal{A}(\qq)\) is a commutative algebra, and we want to show it is maximal with respect to this property.
Now, for any generic object, which is denoted by \((\qq)\), and any prime power \(q\), we denote by \((q)\) the specialization along \(\qq\mapsto q\).
Then we can consider the diagram
\[\begin{tikzcd}
	\mathcal{A}(\qq) \arrow[d] \arrow[r] & \mathcal{H}^{\mathcal{I}}(\qq) \arrow[d]\\
	\prod_q \mathcal{A}(q) \arrow[r] & \prod_q \mathcal{H}^{\mathcal{I}}(q),
\end{tikzcd}\]
where \(\mathcal{A}(q)\subseteq \mathcal{H}^{\mathcal{I}}(q)\) is the subalgebra generated by the classical Bernstein elements \(E_\lambda\) from \cite[Théorème 1.5]{Vigneras:Algebres}, and the vertical maps are the product of all specialization maps, obtained by mapping \(\qq\) to a prime power \(q\).
Then all maps above are injective, and the diagram is cartesian as \(\{\theta_\lambda\mid \lambda \in X_*(T)\}\) can be extended to a basis of \(\calH^{\mathcal{I}}(\qq)\) (similar to \cite[Théorème 1]{Vigneras:Algebres}).
On the other hand, all specialization maps \(\mathcal{A}(\qq) \to \mathcal{A}(q)\) and \(\mathcal{H}^{\mathcal{I}}(\qq)\to \mathcal{H}^{\mathcal{I}}(q)\) are surjective.
Hence, the fact that \(\mathcal{A}(\qq)\subseteq \mathcal{H}^{\mathcal{I}}(\qq)\) is maximal commutative will follow from the similar assertion for usual Hecke algebras with integral coefficients.
However, since the images of \(\theta_\lambda\) can be extended to a basis of \(\calH^{\mathcal{I}}(q)\), it even suffices to show it for Iwahori--Hecke algebras with complex coefficients.
This is well-known, and follows from the explicit description of the multiplication of the basis elements above, such as in \cite[Corollary 5.47]{Vigneras:Pro-p}.
\xpf

Having decategorified the Wakimoto motives, we can now identify the center of \(\calH^{\mathcal{I}}(\qq)\).

\theo\thlabel{generic Bernstein iso}
The map $\varphi^{\mathbf a_0}$ in \eqref{Z blah} is an isomorphism.
\xtheo

\pf
By (a suitable variation of) \cite[(1.6.5)]{Vigneras:Algebres}, the elements $\sum_{\lambda \in W_0(\mu)} \theta_\lambda$ for $\mu \in X_*(T)^+$ form a basis of $Z(\calH^\calI(\qq))$. On the other hand, by the motivic Satake equivalence \cite[Theorem 6.32]{CassvdHScholbach:MotivicSatake} and the generic Satake isomorphism \cite[Corollary 6.36] {CassvdHScholbach:MotivicSatake}, the classes of the rational IC-motives $[\IC_{\mu, \Q}] \in K_0(\DTM_{L^+G}(\Gr_G)^{\antilc})$ for $\mu \in X_*(T)^+$ form a basis of $\mathcal{H}^{\sph}(\qq)$. Now by \thref{ExistenceWakimotoFiltration}, $[\mathsf{Z}_{G} (\IC_{\mu, \Q})]$ is a $\Z[\qq^{\pm 1}]$-linear combination of the $\theta_\nu$ for $\nu$ in the convex hull of $W_0(\mu)$. It suffices to show this is actually a $\Z[\qq]$-linear combination, and that the coefficient of $\theta_{\nu}$ equals $1$ for all $\nu \in W_0(\mu)$.

By \thref{Grad-CT}, $\Grad_{\nu} (\mathsf{Z}_G(\IC_{\mu, \Q}))\cong \CT_{B^-}^{t(\nu)}(\IC_{\mu, \Q})(- \langle 2 \rho, \nu \rangle) [- \langle 2 \rho, \nu \rangle]$, where $\lambda$ is a regular antidominant cocharacter. Let $\mathcal{S}_\nu^+$ be the attractor associated to $\nu$ in $\Gr_G$ and let $n_{\mu, \nu}$ be the number of irreducible components of $\mathcal{S}_\nu^+ \cap \Gr_G^\mu$. Note that $n_{\mu, \nu} = 1$ if $\nu \in W_0(\mu)$, and $n_{\mu, \nu}$ is constant for $\nu$ in the same $W_0$-orbit. Now by \cite[Lemma 5.42]{CassvdHScholbach:MotivicSatake} we have $\CT_{B^-}^{t(\nu)}(\IC_{\mu, \Q}) [- \langle 2 \rho, \nu \rangle] \cong \Q^{\oplus n_{\mu, \nu}}(-\dim (\mathcal{S}_\nu^+ \cap \Gr^\mu))$. Using that $\dim (\mathcal{S}_\nu^+ \cap \Gr^\mu) = \langle \rho, \mu-v \rangle$ \cite[Theorem 3.2]{MirkovicVilonen:Geometric}, a straightforward computation shows
$$[\JWak_\nu(\Grad_{\nu} (\mathsf{Z}_G(\IC_{\mu, \Q})))] = n_{\mu, \nu} \qq^{\frac{l(\mu) - l(\nu)}{2}} \theta_\nu.$$
These coefficients clearly satisfy the desired properties.
\xpf

\bibliographystyle{alphaurl}
\bibliography{bib}

\end{document}